\documentclass[12pt,twoside,english]{extarticle}
\usepackage{lmodern}
\usepackage{helvet}
\usepackage[T1]{fontenc}
\usepackage[latin9]{inputenc}
\usepackage{geometry}
\usepackage{color}
\usepackage{babel}
\usepackage{float}
\usepackage{mathrsfs}
\usepackage{amsmath}
\usepackage{amsthm}
\usepackage{amssymb}
\usepackage{graphicx}
\usepackage{setspace}
\PassOptionsToPackage{normalem}{ulem}
\usepackage{ulem}
\usepackage[unicode=true,pdfusetitle,
 bookmarks=true,bookmarksnumbered=false,bookmarksopen=false,
 breaklinks=true,pdfborder={0 0 0},pdfborderstyle={},backref=false,colorlinks=true]
 {hyperref}

\makeatletter

\providecommand{\tabularnewline}{\\}

\numberwithin{equation}{section}
\numberwithin{table}{section}
\numberwithin{figure}{section}
\usepackage[natbibapa]{apacite}
\theoremstyle{plain}
\newtheorem{assumption}{\protect\assumptionname}
\theoremstyle{plain}
\newtheorem{prop}{\protect\propositionname}[section]
\theoremstyle{plain}
\newtheorem{lem}{\protect\lemmaname}[section]
\theoremstyle{plain}
\newtheorem{thm}{\protect\theoremname}[section]
\theoremstyle{definition}
\newtheorem{defn}{\protect\definitionname}[section]
\theoremstyle{plain}
\newtheorem{lyxalgorithm}{\protect\algorithmname}
\theoremstyle{plain}
\newtheorem{cor}{\protect\corollaryname}[section]
\theoremstyle{remark}
\newtheorem{rem}{\protect\remarkname}[section]

\usepackage{graphicx}
\usepackage{amsmath}
\usepackage{dcolumn}
\usepackage{listings}
\usepackage{color}
\usepackage{setspace}
\usepackage{ulem}  
\definecolor{hellgelb}{rgb}{1,1,0.8}
\definecolor{colKeys}{rgb}{0,0,1}
\definecolor{colIdentifier}{rgb}{0,0,0}
\definecolor{colComments}{rgb}{1,0,0}
\definecolor{colString}{rgb}{0,0.5,0}

\usepackage{lmodern}
\usepackage[T1]{fontenc}
\usepackage{geometry}
\usepackage{fancyhdr}
\usepackage{amsthm}
\usepackage{thmtools}
\usepackage{amsmath}
\usepackage{amssymb}
\usepackage{esint}
\usepackage{multirow}
\usepackage{mathrsfs} 
\usepackage{textgreek}
\usepackage{amsfonts}
\usepackage{amssymb}
\usepackage{mathrsfs}

\usepackage[USenglish]{isodate}
\usepackage{chngcntr}

\numberwithin{equation}{section}
\numberwithin{table}{section}
\numberwithin{assumption}{section}
\counterwithout{figure}{section}
\counterwithout{table}{section}

\makeatother

  \providecommand{\algorithmname}{Algorithm}
  \providecommand{\assumptionname}{Assumption}

  \providecommand{\corollaryname}{Corollary}
  
  \providecommand{\definitionname}{Definition}

  \providecommand{\lemmaname}{Lemma}

  \providecommand{\propositionname}{Proposition}
  \providecommand{\remarkname}{Remark}
  
  \providecommand{\theoremname}{Theorem}
 \providecommand{\corollaryname}{Corollary}
 \providecommand{\theoremname}{Theorem}
 
\usepackage{thmtools}

\newtheoremstyle{MyTheoremstyle}
  {\topsep} 
  {\topsep} 
  {} 
  {} 
  {\bfseries} 
  {.} 
  {.90em} 
  {} 
\theoremstyle{MyTheoremstyle} 
\theoremstyle{MyTheoremstyle} 
\theoremstyle{MyTheoremstyle} 
\theoremstyle{MyTheoremstyle} 
\theoremstyle{MyTheoremstyle}

\declaretheoremstyle[
    headfont=\bfseries,
    notefont=\normalfont,
    bodyfont=\itshape,
    headpunct=\newline,
    headformat={%
        \makebox{\NAME\ \NUMBER\ }{\NOTE}%
    },
]{theorem}

\newlength{\spacelength}
\declaretheoremstyle[
    headfont=\bfseries,
    notefont=\normalfont,
    bodyfont=\itshape,
    headpunct=\newline,
    headformat={%
        \makebox[0pt][l]{\NAME\ \NUMBER\ }\hskip-\spacelength{\NOTE}%
    },
]{theore}

 \usepackage{fancyhdr}
\newcommand{\simpleheading}{
\newgeometry{verbose,tmargin=2cm,bmargin=2cm,lmargin=2cm,rmargin=2cm,headheight=1cm,headsep=1cm,footskip=1cm}
\setlength{\headheight}{25pt}
\fancyhead{{\normalsize \textsc{}}}}

\usepackage{booktabs}
\usepackage{float}
\usepackage{graphicx}
\usepackage{epstopdf}
\usepackage{morefloats}
\usepackage[referable]{threeparttablex}
\usepackage{footnote}

\usepackage{setspace} 
\usepackage{graphicx,color}

\usepackage{listings}
\RequirePackage[T1]{fontenc}
\RequirePackage{ae,fancyvrb}
\DefineVerbatimEnvironment{Sinput}{Verbatim}{fontshape=sl}
\DefineVerbatimEnvironment{Soutput}{Verbatim}{}
\DefineVerbatimEnvironment{Scode}{Verbatim}{fontshape=sl}

\title{\bf Continuous Record Asymptotic Framework for Inference in Strucutral Change Models\footnote{%
We wish to thank Mark Podolskij for helpful suggestions. We also thank Christian Gouriéroux and Viktor Todorov for related discussions on some aspects of this project. We are grateful to Iván Fernández-Val, Hiro Kiado, Zhongjun Qu as well as seminar partecipants at the Boston University for valuable comments and suggestions. Seong Yeon Chang, Andres Sagner and Yohei Yamamoto have provided generous help with computer programming. We thank Yunjong Eo and James Morley for shraing their programs with us. Casini gratefully acknowledges partial financial support from the Bank of Italy.}}
\author{
\textsc{\textcolor{MyBlue}{Alessandro Casini}}\thanks{Corresponding author at: Department of Economics and Finance, University of Rome Tor Vergata, Via Columbia 2, Rome 00133, IT. 
Email: 
\texttt{\textcolor{MyBlue}{{alessandro.casini@uniroma2.it}}}.} 
\\
\small{University of Rome Tor Vergata}
\and
\textsc{\textcolor{MyBlue}{Pierre Perron}}\thanks{Department of Economics, Boston University, 270 Bay State Road, Boston, MA 02215, US. 
Email: 
\texttt{\textcolor{MyBlue}{\mbox{perron@bu.edu}}}.} 
\\
\small{Boston University}
}

\date{\small{\today}} 

 \makeatletter

\numberwithin{equation}{section}
\makeatother

\usepackage{babel}
\addto\captionsenglish{%
}

\usepackage{color}
\usepackage{xcolor}

\renewcommand*{\thesection}{\arabic{section}}

\definecolor{MyRed}{rgb}{0.8,0,0}
\definecolor{MyBlue}{rgb}{0,0,0.7}
\definecolor{Green}{rgb}{0,0.5,0}
\definecolor{hellgelb}{rgb}{1,1,0.8}
\definecolor{colKeys}{rgb}{0,0,1}
\definecolor{colIdentifier}{rgb}{0,0,0}
\definecolor{colComments}{rgb}{1,0,0}
\definecolor{colString}{rgb}{0,0.5,0}

\definecolor{MyLightRed}{rgb}{2.2,0.2,0.4} 
\definecolor{MyLightRed2}{rgb}{0.6,0.2,0.3} 
\definecolor{MyLightRed2temp}{rgb}{0.6,0.2,0.3}
\definecolor{MyLightRed3}{rgb}{0.8,0.1,0.1} 
\definecolor{MyRed}{rgb}{0.7,0.0,0}

\definecolor{MyLigthBlue13}{rgb}{0,0.2,0.7}
 \definecolor{MyLigthBlack}{rgb}{0.2,0.25,0.3} 

\hypersetup{%
  bookmarks=true
  pdftitle = {Title Here},
  pdfsubject = {},
  pdfkeywords = {Keywords},
  pdfauthor = {Alessandro Casini},}
\hypersetup{ 
  colorlinks = {true}, 
  citecolor = {MyBlue},
  linkcolor = {MyRed},
  filecolor= {Green},	
  urlcolor = {MyRed},
  hyperindex = {true},
  linktocpage = {true},
  linkbordercolor ={1 0 0},
 citebordercolor ={0 1 1},
 urlbordercolor ={0 1 1}
}

\usepackage{bookmark}

\usepackage[bottom]{footmisc}
\usepackage{multibib}
\newcites{ReferencesSupp}{References}

\makeatother

\providecommand{\algorithmname}{Algorithm}
\providecommand{\assumptionname}{Assumption}
\providecommand{\corollaryname}{Corollary}
\providecommand{\definitionname}{Definition}
\providecommand{\lemmaname}{Lemma}
\providecommand{\propositionname}{Proposition}
\providecommand{\remarkname}{Remark}
\providecommand{\theoremname}{Theorem}

\begin{document}
\setcounter{page}{0}
\title{\textbf{Continuous Record Asymptotics for Change-Point Models}\thanks{We wish to thank Mark Podolskij for helpful suggestions. We also thank
Christian Gouri\'{e}roux, Hashem Pesaran, Myung Hwan Seo and Viktor
Todorov for related discussions on some aspects of this project. We
are grateful to Iv\'{a}n Fern\'{a}ndez-Val, Hiro Kaido, Zhongjun
Qu as well as seminar participants at Boston University and participants
at 2018 NBER-NSF Time Series Conference and 11th Annual SoFiE Conference
for valuable comments and suggestions. Seong Yeon Chang, Andres Sagner
and Yohei Yamamoto have provided generous help with computer programming.
We thank Yunjong Eo and James Morley for sharing their programs. Casini
gratefully acknowledges partial financial support from the Bank of
Italy.}}
\maketitle
\begin{abstract}
{\footnotesize{}In the context of a linear regression model with a
single break point, we develop a continuous record asymptotic framework
to build inference methods for the break date. We have $T$ observations
with a sampling frequency $h$ over a fixed time horizon $\left[0,\,N\right],$
and let $T\rightarrow\infty$ with $h\downarrow0$ while keeping the
time span $N$ fixed. We consider the least-squares estimate of the
break date and establish consistency and convergence rate. We provide
a limit theory for shrinking magnitudes of shifts and locally increasing
variances. The asymptotic distribution corresponds to the location
of the extremum of a function of the quadratic variation of the regressors
and of a Gaussian centered martingale process over a certain time
interval. We can account for the asymmetric informational content
provided by the pre- and post-break regimes and show how the location
of the break and shift magnitude are key ingredients in shaping the
distribution. We consider a feasible version based on plug-in estimates,
which provides a very good approximation to the finite sample distribution.
We use the concept of Highest Density Region to construct confidence
sets. Overall, our method is reliable and delivers accurate coverage
probabilities and relatively short average length of the confidence
sets. Importantly, it does so irrespective of the size of the break.}{\footnotesize\par}
\end{abstract}

\onehalfspacing
\thispagestyle{empty}
\allowdisplaybreaks

\pagebreak{}

\section{Introduction}

In the context of a linear regression model with a single break point,
we develop a continuous record asymptotic framework and inference
methods for the break date. Our model is specified in continuous time
but estimated with discrete-time observations using a least-squares
method. We have $T$ observations with a sampling frequency $h$ over
a fixed time horizon $\left[0,\,N\right],$ where $N=Th$ denotes
the time span of the data. We consider a continuous record asymptotic
framework whereby $T$ increases by shrinking the time interval $h$
to zero while keeping time span $N$ fixed. We impose very mild conditions
on an underlying continuous-time model assumed to generate the data,
basically continuous Itô semimartingales. 

An extensive amount of research addressed change-point problems under
the classical large-$N$ asymptotics. Early contributions are \citet{hinkley:71},
\citet{bhattacharya:87}, and \citet{yao:87}, who adopted a Maximum
Likelihood (ML) approach, and for linear regression models, \citet{bai:97RES}
and \citet{bai/perron:98}.  See the reviews of \citet{csorgo/horvath:97},
\citet{aue/Horvath:13}, \citet{casini/perron_Oxford_Survey} and
references therein. In this literature, the resulting large-$N$
limit theory for the estimate of the break date depends on the exact
distributions of the regressors and disturbances. Therefore, a so-called
shrinkage asymptotic theory was adopted whereby the magnitude of
the shift, say $\delta_{T},$ converges to zero which leads to an
invariant limit distribution.

We study a general change-point problem under a continuous record
asymptotic framework and develop inference procedures based on the
derived asymptotic distribution. We establish consistency at rate-$T$
convergence for the least-squares estimate of the break date, assumed
to occur at time $N_{b}^{0}$. Given the fast rate of convergence,
we introduce a limit theory with shrinking magnitudes of shifts and
increasing variance of the residual process local to the change-point.
The asymptotic distribution corresponds to the location of the extremum
of a function of the (quadratic) variation of the regressors and of
a Gaussian centered martingale process over some time interval. It
is characterized by some notable aspects. With the time horizon $\left[0,\,N\right]$
fixed, we can account for the asymmetric informational content provided
by the pre- and post-break observations, i.e., the time span and the
position of the break date $N_{b}^{0}$ convey useful information
about the finite-sample distribution. In contrast, this is not achievable
under the large-$N$ shrinkage asymptotic framework because both pre-
and post-break segments expand proportionately as $T$ increases and,
given the mixing assumptions imposed, only the neighborhood around
the break date remains relevant. Further, the domain of the extremum
depends on the position of the break $N_{b}^{0}$ relative to $N$,
or total span, and thus the distribution is asymmetric, in general.
The degree of asymmetry increases as the true break point moves away
from mid-sample. This holds unless the magnitude of the break is large,
in which case the density is symmetric irrespective of the location
of the break. This accords with simulation evidence which documents
that the break point estimate is less precise and the coverage rates
of the confidence intervals less reliable when the break is not at
mid-sample {[}see, e.g., \citet{chang/perron:18}{]}. When the shift
magnitude is small, the  density displays three modes. As the shift
magnitude increases, this tri-modality vanishes. We show empirically
that all of these features are shared by the finite-sample distribution
of the least-squares estimator of the break date. Hence, the continuous
record asymptotics theory provides an accurate approximation to the
finite-sample distribution of the break date estimator. In contrast,
we show that the large-$N$ shrinkage asymptotic distribution of \citet{yao:87}
and \citet{bai:97RES} provides a poor approximation to the finite-sample
distribution and does not share any of those features as we discuss.

Our asymptotics can be seen as intermediate between the shrinkage
asymptotics and more recent approaches relying on weak identification
{[}see e.g., \citet{elliott/mueller:07} and \citet{elliott/mueller/watson:15}{]}.
On the one hand, using the usual shrinking condition of \citet{yao:87}
and \citet{bai:97RES} for which the break magnitude, say $\delta_{T},$
goes to zero at a rate slower than $O(T^{-1/2})$ leads to underestimation
of the uncertainty about the break date. On the other hand, the weak
identification condition of \citet{elliott/mueller:07} for which
$\delta_{T}$ goes to zero at a fast rate (i.e., $\delta_{T}=O(T^{-1/2})$,
so that the change-point cannot be consistently estimated) leads to
overstating the uncertainty. This has opposite consequences for the
confidence intervals of the break date. Confidence sets have poor
coverage probabilities when the break is small under \citeauthor{bai:97RES}'s
framework while they can be too wide under that of \citet{elliott/mueller:07}.
In this paper, the key is not to focus our asymptotic experiment on
shrinking condition on $\delta_{T}$ but to make assumptions on the
signal-to-noise ratio $\delta_{T}/\sigma_{t}$ instead, where $\sigma_{t}$
is the volatility of the errors. We require $\delta_{T}$ to go to
zero at a slower rate than that of \citet{elliott/mueller:07}---to
guarantee strong identification---and require $\sigma_{t}$ to increase
without bound when $t$ approaches the break date $T_{b}^{0}$. This
offers a new characterization of the uncertainty without compromising
strong identification and consistency of the model parameters needed
to conduct inference.

Despite the effort devoted to the construction of confidence intervals
for the break date {[}see e.g., \citet{bai/perron:98}, \citet{elliott/mueller:07}
and \citet{eo/morley:15}{]}, what is still missing is a method that,
for both large and small breaks, achieves both accurate coverage rates
and satisfactory average lengths. The most popular method is \citet{bai/perron:98}
which yields confidence intervals that are relatively short but have
good coverage only when the magnitude of the break is not small. However,
both small and large breaks are relevant for empirical work; breaks
that are statistically small can still be practically relevant.

Given the peculiar properties (e.g., multi-modality and asymmetry)
of the continuous record asymptotic distribution, we propose a non-standard
inference related to Bayesian analyses. We use the concept of Highest
Density Region to construct confidence sets for the break date. Our
method has good coverage and length across all break magnitudes. This
has important implications for empirical work because the user can
be confident that our confidence interval includes the true value
across all break sizes. For small breaks, the length of the confidence
intervals from any method can be quite large for some models. However,
our confidence interval is still informative because it reveals that
there is high uncertainty about the change-point. The same information
cannot be provided by existing methods either because they do not
have good coverage unless the break is not small {[}e.g., \citet{bai/perron:98}{]}
or because they have a large length even when the break is not small
{[}e.g., \citet{elliott/mueller:07}{]}.

We use the continuous record asymptotics to provide an alternative
approximation to the finite-sample distribution of the least-squares
estimator based on discrete-time data. This creates no contradiction
since asymptotic theory is intended as a thought experiment used to
obtain approximations to the distribution of estimators or test statistics.
The continuous record asymptotics has proven to be useful in other
discrete-time settings such as in the context of unit roots {[}cf.
\citet{phillips:87a} and \citet{perron:1991ecma}{]} and nonparametric
regression {[}cf. \citet{brown/low:1996}{]}. \citet{brown/low:1996}
showed the asymptotic equivalence between a nonparametric regression
problem and a white noise with drift problem. Our results show the
asymptotic equivalence between discrete-time and continuous-time regression
models with a change-point. 

Recent work in change-point analysis has focused on estimation when
the number of change-points is allowed to increase with the sample
size {[}e.g., \citet{fryzlewicz:14}{]} and when the change-point
is allowed to approach the start and end sample point. A growing literature
has also considered change-points in a high-dimensional setting {[}e.g.,
\citet{lee/seo/shin:16}, \citet{leonardi/buhlmann:16}, \citet{wang/lin/willett:19}
and \citet{wang/yu/rinaldo/willett:19}{]}. This work is mainly concerned
with consistent estimation of the change-point dates and development
of corresponding computational algorithms.  Our focus is on asymptotic
theory and inference within the classical change-point model with
a single break. Our results can also have useful implications for
the growing literature on inference in high-dimensional change-point
analysis and for the literature on threshold regression {[}see, e.g.,
\citet{hansen:00ecma} and \citet{hildago/lee/seo:19}{]}.

This paper relates to other work by the authors, namely \citeauthor{casini/perron_SC_BP_Lap}
(\citeyear{casini/perron_Lap_CR_Single_Inf}; \citeyear{casini/perron_SC_BP_Lap})\nocite{casini/perron_Lap_CR_Single_Inf}.
\citet{casini/perron_Lap_CR_Single_Inf} used the asymptotic results
developed in this paper and proposed a new Generalized Laplace estimator
of the break date under a continuous record asymptotic framework.
 \citet{casini/perron_SC_BP_Lap} analyzed the Generalized Laplace
method under classical asymptotics and focused on the theoretical
relationship between the asymptotic distribution of frequentist and
Bayesian estimators of the break point. Finally, \citet{chambers/taylor:19}
considered both deterministic one-time and continuous stochastic parameter
change in a continuous-time autoregressive model while \citet{casini_CR_Test_Inst_Forecast}
introduced continuous-time asymptotics to test for forecast failure.
Recently, \citet{casini:change-point-spectra} considered testing
and estimating change-points in a locally stationary process using
frequency-domain methods. 

The paper is organized as follows. Section \ref{Section: Model-and-Assumptions}
introduces the model and the estimation method. Section \ref{Section Consistency-and-Rate }
contains results about the consistency and rate of convergence for
fixed shifts. Section \ref{Section Asymptotic Distribution: Continuous Case}
develops the asymptotic theory. We compare our limit theory with the
finite-sample distribution in Section \ref{Section Approximation-to-the}.
Section \ref{Section Inference Methods} describes how to construct
the confidence sets, with simulation results  reported in Section
\ref{Section Small-Sample-Effectiveness-of}. Section \ref{Section Conclusions}
provides brief concluding remarks. The Supplement {[}\citet{casini/perron_CR_Single_Break_Supp}{]}
contains the proofs as well as additional material. 

\section{Model and Assumptions \label{Section: Model-and-Assumptions}}

We denote the transpose of a matrix $A$ by $A'$ and the $\left(i,\,j\right)$
elements of $A$ by $A^{\left(i,j\right)}$. We use $\left\Vert \cdot\right\Vert $
to denote the Euclidean norm of a linear space, i.e., $\left\Vert x\right\Vert =\left(\sum_{i=1}^{p}x_{i}^{2}\right)^{1/2}$
for $x\in\mathbb{R}^{p}.$ We use $\left\lfloor \cdot\right\rfloor $
to denote the largest smaller integer function. A sequence $\left\{ u_{kh}\right\} _{k=1}^{T}$
is $i.i.d.$ (resp., $i.n.d$) if the $u_{kh}$ are independent and
identically (resp., non-identically) distributed. We use $\overset{P}{\rightarrow},\,\Rightarrow,$
and $\overset{\mathcal{L}-\mathrm{s}}{\Rightarrow}$ to denote convergence
in probability, weak convergence and stable convergence in law, respectively.
For semimartingales $\left\{ S_{t}\right\} _{t\geq0}$ and $\left\{ R_{t}\right\} _{t\geq0}$,
we denote their covariation process by $\left[S,\,R\right]_{t}$ and
their predictable counterpart by $\left\langle S,\,R\right\rangle _{t}$.
The symbol ``$\triangleq$'' denotes definitional equivalence. 

Consider a change-point model with a single break point:
\begin{align}
Y_{t} & =D_{t}'\nu^{0}+Z_{t}'\delta_{Z,1}^{0}+e_{t},\quad(t=0,\,1,\ldots,\,T_{b}^{0})\label{Original SC Model}\\
Y_{t} & =D_{t}'\nu^{0}+Z_{t}'\delta_{Z,2}^{0}+e_{t},\quad(t=T_{b}^{0}+1,\ldots,\,T),\nonumber 
\end{align}
 where $Y_{t}$ is the dependent variable, $D_{t}$ and $Z_{t}$ are,
respectively, $q\times1$ and $p\times1$ vectors of regressors and
$e_{t}$ is an unobservable disturbance. The vector-valued parameters
$\nu^{0},\,\delta_{Z,1}^{0}$ and $\delta_{Z,2}^{0}$ are unknown
with $\delta_{Z,1}^{0}\neq\delta_{Z,2}^{0}$. Our main purpose is
to develop inference methods for the unknown change-point date $T_{b}^{0}$
when $T+1$ observations on $\left(Y_{t},\,D_{t},\,Z_{t}\right)$
are available. Before moving to the re-parametrization of the model,
we discuss the underlying continuous-time model assumed to generate
the data. The discrete-time variables are assumed to be generated
from the continuous-time processes $\{D_{s},\,Z_{s},\,e_{s}\}_{s\geq0}$
defined on a filtered probability space $(\Omega,\,\mathscr{F},\,(\mathscr{F}_{s})_{s\geq0},\,P).$
We observe realizations of $\{Y_{s},\,D_{s},\,Z_{s}\}$ at discrete
points of time. 

The sampling occurs at regularly spaced time intervals of length $h$
within a fixed time horizon $\left[0,\,N\right]$ where $N$ denotes
the span of the data. We observe $\left\{ _{h}Y_{kh},\,{}_{h}D_{kh},\,{}_{h}Z_{kh};\,k=0,\,1,\ldots,\,T=N/h\right\} $.
$_{h}D_{kh}\in\mathbb{R}^{q}$ and $_{h}Z_{kh}\in\mathbb{R}^{p}$
are random vector step functions which jump only at times $0,\,h,\ldots,\,Th$.
We shall allow $_{h}D_{kh}$ and $_{h}Z_{kh}$ to include both predictable
processes and locally-integrable semimartingles, though the case with
predictable regressors is more delicate and discussed in the supplement.
The discretized processes $_{h}D_{kh}$ and $_{h}Z_{kh}$ are assumed
to be adapted to the increasing and right-continuous filtration $\left\{ \mathscr{F}_{t}\right\} _{t\geq0}$.
For any process $X$ we denote its ``increments'' by $\Delta_{h}X_{k}=X_{kh}-X_{\left(k-1\right)h}$.
For $k=1,\ldots,T$, let $\Delta_{h}D_{k}\triangleq\mu_{D,k}h+\Delta_{h}M_{D,k}$
and $\Delta_{h}Z_{k}\triangleq\mu_{Z,k}h+\Delta_{h}M_{Z,k}$ where
the ``drifts'' $\mu_{D,t}\in\mathbb{R}^{q},\,\mu_{Z,t}\in\mathbb{R}^{p}$
are $\mathscr{F}_{t-h}$-measurable (exact assumptions will be given
below), and $M_{D,k}\in\mathbb{R}^{q},\,M_{Z,k}\in\mathbb{R}^{p}$
are continuous local martingales with finite conditional covariance
matrix $P$-a.s., $\mathbb{E}(\Delta_{h}M_{D,t}\Delta_{h}M_{D,t}'|\,\mathscr{F}_{t-h})=\Sigma_{D,t-h}\Delta t$
and $\mathbb{E}(\Delta_{h}M_{Z,t}\Delta_{h}M'_{Z,t}|\,\mathscr{F}_{t-h})=\Sigma_{Z,t-h}\Delta t$
($\Delta t$ and $h$ are used interchangeably). Let $\lambda_{0}\in\left(0,\,1\right)$
denote the fractional break date (i.e., $T_{b}^{0}=\left\lfloor T\lambda_{0}\right\rfloor $).
Via the Doob-Meyer Decomposition, model \eqref{Original SC Model}
can be expressed as 
\begin{align}
\Delta_{h}Y_{k} & \triangleq\begin{cases}
\left(\Delta_{h}D_{k}\right)'\nu^{0}+\left(\Delta_{h}Z_{k}\right)'\delta_{Z,1}^{0}+\Delta_{h}e_{k}^{*}, & \left(k=1,\ldots,\left\lfloor T\lambda_{0}\right\rfloor \right)\\
\left(\Delta_{h}D_{k}\right)'\nu^{0}+\left(\Delta_{h}Z_{k}\right)'\delta_{Z,2}^{0}+\Delta_{h}e_{k}^{*}, & \left(k=\left\lfloor T\lambda_{0}\right\rfloor +1,\ldots,\,T\right),
\end{cases}\label{Eq. Model 1, CT}
\end{align}
 where the error process $\left\{ \Delta_{h}e_{t}^{*},\,\mathscr{F}_{t}\right\} $
is a continuous local martingale difference sequence with conditional
variance $\mathbb{E}[\left(\Delta_{h}e_{t}^{*}\right)^{2}|\,\mathscr{F}_{t-h}]=\sigma_{e,t-h}^{2}\Delta t$
$P$-a.s. finite. The underlying continuous-time data-generating process
can thus be represented (up to $P$-null sets) in integral equation
form as
\begin{align}
D_{t}=D_{0}+\int_{0}^{t}\mu_{D,s}ds+\int_{0}^{t}\sigma_{D,s}dW_{D,s},\quad & Z_{t}=Z_{0}+\int_{0}^{t}\mu_{Z,s}ds+\int_{0}^{t}\sigma_{Z,s}dW_{Z,s},\label{Model Regressors Integral Form}
\end{align}
where $\sigma_{D,t}$ and $\sigma_{Z,t}$ are the instantaneous covariance
processes taking values in $\mathcal{M}_{q}^{\textrm{càdlàg}}$ and
$\mathcal{M}_{p}^{\textrm{càdlàg}}$ {[}the space of $p\times p$
positive definite real-valued matrices whose elements are càdlàg{]};
$W_{D}$ (resp., $W_{Z}$) is a $q$ (resp., $p$)-dimensional standard
Wiener process; $e^{*}=\left\{ e_{t}^{*}\right\} _{t\geq0}$ is a
continuous local martingale which is orthogonal (in a martingale sense)
to $\left\{ D_{t}\right\} _{t\geq0}$ and $\left\{ Z_{t}\right\} _{t\geq0}$;
and $D_{0}$ and $Z_{0}$ are $\mathscr{F}_{0}$-measurable random
vectors. In \eqref{Model Regressors Integral Form}, $\int_{0}^{t}\mu_{D,s}ds$
is a continuous adapted process with finite variation paths and $\int_{0}^{t}\sigma_{D,s}dW_{D,s}$
corresponds to a continuous local martingale. 
\begin{assumption}
\label{Assumption 1, CT}(i) $\mu_{D,t},\,\mu_{Z,t},\,\sigma_{D,t}$
and $\sigma_{Z,t}$ satisfy $P$-a.s., $\sup_{\omega\in\Omega,\,0<t\leq\tau_{T}}\left\Vert \mu_{D,t}\left(\omega\right)\right\Vert <\infty$,
$\sup_{\omega\in\Omega,}$ $_{0<t\leq\tau_{T}}\left\Vert \mu_{Z,t}\left(\omega\right)\right\Vert <\infty$,
$\sup_{\omega\in\Omega,\,0<t\leq\tau_{T}}\left\Vert \sigma_{D,t}\left(\omega\right)\right\Vert <\infty$
and $\sup_{\omega\in\Omega,\,0<t\leq\tau_{T}}\left\Vert \sigma_{Z,t}\left(\omega\right)\right\Vert <\infty$
for some localizing sequence $\left\{ \tau_{T}\right\} $ of stopping
times. Also, $\sigma_{D,\,s}$ and $\sigma_{Z,s}$ are càdlàg; (ii)
$\int_{0}^{t}\mu_{D,s}ds$ and \textup{$\int_{0}^{t}\mu_{Z,s}ds$
}belong to the class of continuous adapted finite variation processes;
(iii) $\int_{0}^{t}\sigma_{D,s}dW_{D,s}$ and $\int_{0}^{t}\sigma_{Z,s}dW_{Z,s}$
are continuous local martingales with $P$-a.s. finite positive definite
conditional variances (or spot covariances) defined by $\Sigma_{D,t}=\sigma_{D,t}\sigma'_{D,t}$
and $\Sigma_{Z,t},=\sigma_{Z,t}\sigma'_{Z,t}$, which for all $t<\infty$
satisfy $\int_{0}^{t}\Sigma_{D,s}^{\left(j,j\right)}ds<\infty$ $\left(j=1,\ldots,\,q\right)$
and $\int_{0}^{t}\Sigma_{Z,s}^{\left(j,j\right)}ds<\infty$ $\left(j=1,\ldots,\,p\right)$.
Furthermore, for every $j=1,\ldots,\,q,$ $r=1,\ldots,\,p$, and $k=1,\ldots,\,T$,
$h^{-1}\int_{\left(k-1\right)h}^{kh}\Sigma_{D,s}^{\left(j,j\right)}ds$
and $h^{-1}\int_{\left(k-1\right)h}^{kh}\Sigma_{Z,s}^{\left(r,r\right)}ds$
are bounded away from zero and infinity, uniformly in $k$ and $h$;
(iv) $e_{t}^{*}$ is such that $e_{t}^{*}\triangleq\int_{0}^{t}\sigma_{e,s}dW_{e,s}$
with $0<\sigma_{e,t}^{2}<\infty$, where $W_{e}$ is a one-dimensional
standard Wiener process. Furthermore, $\left\langle e,\,D\right\rangle _{t}=\left\langle e,\,Z\right\rangle _{t}=0$
identically for all $t\geq0$.
\end{assumption}
Part (i) restricts the processes to be locally bounded and part (ii)
requires the drifts to be adapted finite variation processes. These
are standard regularity conditions in the high-frequency statistics
literature {[}cf. \citet{barndorff/shephard:04}{]}. Part (iii) requires
the regressors to have finite integrated covariance. We rule out jump
processes; hence, our results are not expected to provide good approximations
for high-frequency data but for data sampled at lower frequencies.
In principle, ultra high frequency data would be used which, however,
are essentially available only for financial variables, which involve
a host of issues that we cannot handle (e.g., market-microstructure,
bid-ask spread, volatility jumps, non-continuous sampling, etc.).
\begin{assumption}
\label{Assumption 2}$D,\,Z,\,e$ and $\Sigma^{0}\triangleq\left\{ \Sigma_{\cdot,t},\,\sigma_{e,t}\right\} _{t\geq0}$
have $P$-a.s. continuous sample paths. 
\end{assumption}
An interesting issue is whether the theoretical results to be derived
for model \eqref{Eq. Model 1, CT} are applicable to classical structural
change models for which an increasing span of data is assumed. This
requires establishing a connection between the assumptions imposed
on the stochastic processes in both settings. Roughly, the classical
long-span setting uses approximation results valid for weakly dependent
data; e.g., ergodic and mixing procesess. Such assumptions are not
needed under our fixed-span asymptotics. Nonetheless, we can impose
restrictions on the probabilistic properties of the latent volatility
processes in our model and thereby guarantee that ergodic and mixing
properties are inherited by the corresponding observed processes.
This follows from Theorem 3.1 in \citet*{genon-catalot/jeantheau/laredo:00}
together with Proposition 4 in \citet{carrasco/chen:02}. For example,
these results imply that the observations $\left\{ Z_{kh}\right\} _{k\geq1}$
(with fixed $h$) can be viewed (under certain conditions) as a hidden
Markov model which inherits the ergodic and mixing properties of $\left\{ \sigma_{Z,t}\right\} _{t\geq0}$.
Hence, our model encompasses those considered in the structural change
literature that uses a long-span asymptotic setting. We shall extend
model \eqref{Eq. Model 1, CT} to allow for predictable processes
(e.g., a constant and/or lagged dependent variable) in the supplement.
\begin{assumption}
\label{Assumption 3 Break Date} $N_{b}^{0}=N\lambda_{0}$ for some
$\lambda_{0}\in\left(0,\,1\right)$.
\end{assumption}
It is useful to re-parametrize model \eqref{Eq. Model 1, CT}. Let
$y_{kh}=\Delta_{h}Y_{k},$ $x_{kh}=(\Delta_{h}D'_{k},\,\Delta_{h}Z'_{k})'$,
$z_{kh}=\Delta_{h}Z_{k}$, $e_{kh}=\Delta_{h}e_{k}^{*},$ $\beta^{0}=((\pi^{0}),\,(\delta_{Z,1}^{0})')'$
and $\delta_{Z}^{0}=\delta_{Z,2}^{0}-\delta_{Z,1}^{0}$. \eqref{Eq. Model 1, CT}
can be expressed as:
\begin{align}
y_{kh} & =x_{kh}'\beta^{0}+e_{kh},\qquad & (k=1,\ldots,\,T_{b}^{0})\label{Model (4), scalar format, CT}\\
y_{kh} & =x_{kh}'\beta^{0}+z_{kh}'\delta_{Z}^{0}+e_{kh},\qquad & (k=T_{b}^{0}+1,\ldots,\,T),\nonumber 
\end{align}
where the true parameter $\theta^{0}=((\beta^{0})',\,(\delta_{Z}^{0})')'$
takes value in a compact space $\Theta\subset\mathbb{R}^{\textrm{dim}\left(\theta\right)}$.
Also, define $z_{kh}=R'x_{kh}$, where $R$ is a $\left(q+p\right)\times p$
known matrix with full column rank. We consider a partial structural
change model for which $R=\left(0,\,I\right)'$ with $I$ an identity
matrix.

Finally, we write the model in matrix format which will be useful
for the derivations. Let $Y=(y_{h},\,\ldots,\,y_{Th})',\,X=(x_{h},\,\ldots,\,x_{Th})'$,
$e=(e_{h},\,\ldots,\,e_{Th})',$ $X_{1}=(x_{h},\,\ldots,\,x_{T_{b}h},\,0,\,\ldots,\,0)'$,
$X_{2}=(0,\,\ldots,\,0,\,x_{\left(T_{b}+1\right)h},\ldots,\,x_{Th})'$
and $X_{0}=(0,\,\ldots,\,0,\,x_{\left(T_{b}^{0}+1\right)h},\ldots,\,x_{Th})'$.
Note that the difference between $X_{0}$ and $X_{2}$ is that the
latter uses $T_{b}$ rather than $T_{b}^{0}$. Define $Z_{1}=X_{1}R,$
$Z_{2}=X_{2}R$ and $Z_{0}=XR$. \eqref{Model (4), scalar format, CT}
in matrix format is: $Y=X\beta^{0}+Z_{0}\delta_{Z}^{0}+e$. We consider
the least-squares estimator of $T_{b}$, i.e., the minimizer of $S_{T}\left(T_{b}\right)$,
the sum of squared residuals when regressing $Y$ on $X$ and $Z_{2}$
over all possible partitions, namely: $\widehat{T}_{b}^{\textrm{LS}}=\mathrm{argmin}_{p+q\leq T_{b}\leq T}S_{T}\left(T_{b}\right)$.
It is straightforward to show that $\widehat{T}_{b}^{\textrm{LS}}=\mathrm{\mathrm{argmin}}_{p+q\leq T_{b}\leq T}Q_{T}\left(T_{b}\right)$
where $Q_{T}\left(T_{b}\right)\triangleq\widehat{\delta}'_{T_{b}}\left(Z_{2}'MZ_{2}\right)\widehat{\delta}_{T_{b}}$,
$\widehat{\delta}_{T_{b}}$ is the least-squares estimator of $\delta_{Z}^{0}$
when regressing $Y$ on $X$ and $Z_{2}$, and $M=I-X\left(X'X\right)^{-1}X'$.
For brevity, we will write $\widehat{T}_{b}^{\textrm{}}$ for $\widehat{T}_{b}^{\textrm{LS}}$
with the understanding that $\widehat{T}_{b}$ is a sequence indexed
by $T$. Let $\widehat{\delta}=\widehat{\delta}{}_{\widehat{T}_{b}}$.
The estimate of the break fraction is then $\widehat{\lambda}_{b}=\widehat{T}_{b}/T$.
Both in practice and for theoretical analyses, a trimming parameter
$\pi\in\left(0,\,1/2\right)$ is  applied to restrict the minimization
over the interval $\left[T\pi,\,\left(1-\pi\right)T\right]$. 

\section{Consistency and Convergence Rate under Fixed Shifts\label{Section Consistency-and-Rate } }

We now establish the consistency and convergence rate of the least-squares
estimator under fixed shifts. Under the classical large-$N$ asymptotics,
related results have been established by \citet{bai:97RES} and \citet{bai/perron:98}.
Early important results for a mean-shift appeared in \citet{yao:87}
and \citet{bhattacharya:87} for an $\mathit{i.i.d.}$ series, \citet{bai:94a}
for linear processes and \citet{picard:85} for a Gaussian autoregressive
model. 
\begin{assumption}
\label{Assumption 4 Eigenvalue}There exists an $l_{0}$ such that
for all $l>l_{0},$ the matrices $\left(lh\right)^{-1}\sum_{k=1}^{l}x_{kh}x'_{kh},$
$\left(lh\right)^{-1}\sum_{k=T-l+1}^{T}x_{kh}x'_{kh},$  $\left(lh\right)^{-1}\sum_{k=T_{b}^{0}-l+1}^{T_{b}^{0}}x_{kh}x'_{kh},$
and $\left(lh\right)^{-1}\sum_{k=T_{b}^{0}+1}^{T_{b}^{0}+l}x_{kh}x'_{kh},$
have minimum eigenvalues bounded away from zero in probability.
\end{assumption}

\begin{assumption}
\label{Assumption 5 Identification}Let $Q_{0}\left(T_{b},\,\theta^{0}\right)\triangleq\mathbb{E}\left[Q_{T}\left(T_{b},\,\theta^{0}\right)-Q_{T}\left(T_{b}^{0},\,\theta^{0}\right)\right]$.
There exists a $T_{b}^{0}$ such that $Q_{0}\left(T_{b}^{0},\,\theta^{0}\right)>\sup_{\left(T_{b},\,\theta^{0}\right)\notin\mathbf{B}}Q_{0}\left(T_{b},\,\theta^{0}\right),$
for every open set $\mathbf{B}$ that contains $\left(T_{b}^{0},\,\theta^{0}\right)$. 
\end{assumption}
Assumption \ref{Assumption 4 Eigenvalue} is similar to A2 in \citet{bai/perron:98}
and requires enough variation around the break point and at the beginning
and end of the sample. The factor $h^{-1}$ normalizes the observations
so that the assumption is implied by a weak law of large numbers.
Assumption \ref{Assumption 5 Identification} is a standard uniqueness
identification condition. We then have the following results.
\begin{prop}
\label{Proposition (1) - Consistency}Under Assumption \ref{Assumption 1, CT}-\ref{Assumption 3 Break Date}
and \ref{Assumption 4 Eigenvalue}-\ref{Assumption 5 Identification},
$\widehat{\lambda}_{b}\overset{P}{\rightarrow}\lambda_{0}$.
\end{prop}

\begin{prop}
\label{Proposition 2, (Rate of Convergence)}Under Assumption \ref{Assumption 1, CT}-\ref{Assumption 3 Break Date}
and \ref{Assumption 4 Eigenvalue}-\ref{Assumption 5 Identification}
for any $\varepsilon>0$, there exists a $K>0$ such that for all
large $T$, $P(T\left|\widehat{\lambda}_{b}-\lambda_{0}\right|>K)<\varepsilon.$
\end{prop}
We have the same $T$-convergence rate as under large-$N$ asymptotics.
Let $\theta^{0}=((\beta^{0})',\,(\delta_{1}^{0})',$ $(\delta_{2}^{0})')'$.
The fast $T$-rate of convergence implies that the least-squares estimate
of $\theta^{0}$ is the same as when $\lambda_{0}$ is known. A natural
estimator for $\theta^{0}$ is $\mathrm{\mathrm{\mathrm{argmin}}_{\beta\in\mathbb{R}^{\mathit{p+q}},\delta\in\mathbb{R}^{\mathit{p}}}}||Y-X\beta-\widehat{Z}_{2}\delta||^{2}$,
where we use $\widehat{T}_{b}$ instead of $T_{b}$ in the construction
of $\widehat{Z}_{2}$. Then we have the following result, akin to
an extension of corresponding results in Section 3 of \citet{barndorff/shephard:04}.
As a matter of notation, let $\Sigma^{*}\triangleq\left\{ \mu_{\cdot,t},\,\Sigma_{\cdot,t},\,\sigma_{e,t}\right\} _{t\geq0}$
and denote expectation taken with respect to $\Sigma^{*}$ by $\mathbb{E}^{*}$.
\begin{prop}
\label{Proposition OLS Asymtptoc Distribu}Under Assumption \ref{Assumption 1, CT}-\ref{Assumption 3 Break Date}
and \ref{Assumption 4 Eigenvalue}-\ref{Assumption 5 Identification},
we have as $T\rightarrow\infty$ ($N$ fixed), conditionally on $\Sigma^{*},$
$(\sqrt{T/N}(\widehat{\beta}-\beta^{0}),\,\sqrt{T/N}(\widehat{\delta}-\delta_{Z}^{0}))'\overset{d}{\rightarrow}\mathscr{MN}\left(0,\,V\right)$
where $\mathscr{MN}$ denotes a mixed Gaussian distribution, with
\begin{align*}
V & \triangleq\overline{V}^{-1}\underset{T\rightarrow\infty}{\lim}T\begin{bmatrix}\sum_{k=1}^{T}\mathbb{E}^{*}\left(x_{kh}x'_{kh}e_{kh}^{2}\right) & \sum_{k=T_{b}^{0}}^{T}\mathbb{E}^{*}\left(x_{kh}z'_{kh}e_{kh}^{2}\right)\\
\sum_{k=T_{b}^{0}}^{T}\mathbb{E}^{*}\left(x_{kh}z'_{kh}e_{kh}^{2}\right) & \sum_{k=T_{b}^{0}}^{T}\mathbb{E}^{*}\left(z_{kh}z'_{kh}e_{kh}^{2}\right)
\end{bmatrix}\overline{V}^{-1},
\end{align*}
and 
\begin{align*}
\overline{V} & \triangleq\underset{T\rightarrow\infty}{\lim}\begin{bmatrix}\sum_{k=1}^{T}\mathbb{E}^{*}\left(x_{kh}x'_{kh}\right) & \sum_{k=T_{b}^{0}}^{T}\mathbb{E}^{*}\left(x_{kh}z'_{kh}\right)\\
\sum_{k=T_{b}^{0}}^{T}\mathbb{E}^{*}\left(x_{kh}x'_{kh}\right) & \sum_{k=T_{b}^{0}}^{T}\mathbb{E}^{*}\left(z_{kh}z'_{kh}\right)
\end{bmatrix}.
\end{align*}
\end{prop}
 $V$ can be random because   $\Sigma_{\cdot,t}$ and $\sigma_{e,t}$
can be stochastic. Under fixed shifts, Proposition \ref{Proposition (1) - Consistency}-\ref{Proposition OLS Asymtptoc Distribu}
shows the asymptotic equivalence of discrete and continuous-time regression
models with a change-point, a result corresponding to \citet{brown/low:1996}
for nonparametric regression. 

\section{\label{Section Asymptotic Distribution: Continuous Case}Asymptotic
Distribution under a Continuous Record}

We now present results about the limiting distribution of the least-squares
estimate of the break date under a continuous record framework. As
in the classical large-$N$ asymptotics, it depends on the exact distribution
of the data and the errors for fixed break sizes {[}c.f., \citet{hinkley:71}{]}.
This has forced researchers to consider a shrinkage asymptotic theory
where the size of the shift is made local to zero as $T$ increases,
an approach developed by \citet{picard:85} and \citet{yao:87}. We
continue with this avenue. Given the consistency result, we know that
there exists some $h^{*}$ such that for all $h<h^{*}$ with high
probability $\eta Th\leq\widehat{N}_{b}\leq\left(1-\eta\right)Th$,
for $\eta>0$ such that $\lambda_{0}\in\left(\eta,\,1-\eta\right)$.
By Proposition \ref{Proposition 2, (Rate of Convergence)}, $\widehat{N}_{b}-N_{b}^{0}=O_{p}\left(T^{-1}\right)$,
i.e., $\widehat{N}_{b}$ is in a shrinking neighborhood of $N_{b}^{0}$.
With a certain rescaling of the objective function one can first obtain
the shrinkage asymptotic distribution of \citet{bai:97RES}. However,
this is unsatisfactory for two reasons. First, as we show below {[}see
also Casini and Perron (\citeyear{casini/perron_SC_BP_Lap}; \citeyear{casini/perron_Lap_CR_Single_Inf}){]},
the shrinkage asymptotic distribution provides a poor approximation
to the finite-sample distribution of the least-squares estimator.
Second, the latter point also explains the poor coverage properties
of the confidence intervals derived from the shrinkage asymptotic
distribution when the magnitude of the break is not large. Some related
results were obtained by Jiang et al. \citeyearpar{jiang/wang/yu:16}
for a simple location model. Their approach is, however, quite restrictive
and no feasible inference procedure suggested. See the supplement
for a more complete discussion.

We begin with the following assumption which specifies that i) we
use a shrinking condition on $\delta_{Z}^{0}$; ii) we introduce a
locally increasing variance condition on the residual process. The
first is similarly used under classical large-$N$ asymptotics, while
the second is new and useful in our context in order to accurately
capture the relevant uncertainty in the change-point problem. We
do not impose restrictions only on $\delta_{Z}^{0}$ but also on the
ratio $\delta_{Z}^{0}/\sigma_{t}$ when $t$ is close to $T_{b}^{0}$.
We refer to $\delta_{Z}^{0}/\sigma_{t}$ as the signal-to-noise ratio.
Controlling this ratio rather than just $\delta_{Z}^{0}$ allows for
an alternative characterization of the uncertainty about the change-point
date in order to obtain an asymptotic distribution which provides
a better approximation of the finite-sample distribution of the estimator.
To emphasize that $\delta_{Z}^{0}$ depends on the sample-size we
denote it by $\delta_{h}$. 
\begin{assumption}
\label{Assumption 6 - Small Shifts}Let $\delta_{h}=\delta^{0}h^{1/4}$,
$\delta^{0}\in\mathbb{R}^{p}$ and assume that for all $t\in\left(N_{b}^{0}-\epsilon,\,N_{b}^{0}+\epsilon\right),$
with $\epsilon\downarrow0$ and $T^{1-\kappa}\epsilon\rightarrow B<\infty$,
$0<\kappa<1/2$, $\mathbb{E}[\left(\Delta_{h}e_{t}^{*}\right)^{2}|\,\mathscr{F}_{t-h}]=\sigma_{h,t-h}^{2}\Delta t$
$P$-a.s., where $\sigma_{h,t}\triangleq\sigma_{h}\sigma_{e,t}$,
$\sigma_{h}\triangleq\overline{\sigma}h^{-1/4}$ and $\overline{\sigma}\triangleq\int_{0}^{N}\sigma_{e,s}^{2}ds$.
\end{assumption}
Note that the localization parameter $\delta^{0}$ in the definition
of $\delta_{h}$ is different from the fixed parameter $\delta_{Z}^{0}$
since $h\rightarrow0$. The rate $1/4$ in the conditions $\delta_{h}=O(h^{1/4})$
and $\sigma_{h}=O(h^{-1/4})$ is for tractability. One can show that
consistency also holds for a rate faster than $1/4$, though slower
than $\kappa.$ However, for the derivation of the limiting distribution
one needs $\delta_{h}/\sigma_{h}=O(h^{1/2})$ and $O(\delta_{h})=O(\sigma_{h}^{-1})$
with $\kappa<1/2.$ The vector of scaled true parameters is $\theta{}_{h}\triangleq((\beta^{0})',\,\delta'_{h})'$.
Define 
\begin{align}
\Delta_{h}\widetilde{e}_{t} & \triangleq\begin{cases}
\Delta_{h}e_{t}^{*}, & t\notin\left(N_{b}^{0}-\epsilon,\,N_{b}^{0}+\epsilon\right)\\
h^{1/4}\Delta_{h}e_{t}^{*}, & t\in\left(N_{b}^{0}-\epsilon,\,N_{b}^{0}+\epsilon\right)
\end{cases}.\label{Eq. eps WN}
\end{align}
We shall refer to $\left\{ \Delta_{h}\widetilde{e}_{t},\,\mathscr{F}_{t}\right\} $
as the normalized residual process. Under this framework, the rate
of convergence of $\widehat{N}_{b}$ is now $T^{1-\kappa}$ with $0<\kappa<1/2$.
Due to the fast rate of convergence of the change-point estimator,
the objective function oscillates too rapidly as $h\downarrow0$.
By scaling up the volatility of the errors around the change-point,
we make the objective function behave as if it were a function of
a standard diffusion process. The neighborhood in which the errors
have relatively higher variance is shrinking at a rate $1/T^{1-\kappa}$,
the rate of convergence of $\widehat{N}_{b}.$ Hence, in a neighborhood
of $N_{b}^{0}$ in which we study the limiting behavior of the break
point estimator, the rescaled criterion function is regular enough
so that a feasible limit theory can be developed. The rate of convergence
$T^{1-\kappa}$ is still sufficiently fast to guarantee a $\sqrt{T}$-consistent
estimation of the slope parameters, as stated in the following proposition.
Let $\left\langle Z_{\Delta},\,Z_{\Delta}\right\rangle \left(v\right)$
be the predictable quadratic variation process of $Z_{\Delta}$. The
process $\mathscr{W}\left(v\right)$ is, conditionally on $\mathscr{F}$,
a two-sided centered Gaussian martingale with independent increments
and variances given in Section \ref{subsection Description Limiting Process}
of the supplement.
\begin{prop}
\label{Prop 3 Asym}Under Assumption \ref{Assumption 1, CT}-\ref{Assumption 3 Break Date},
\ref{Assumption 4 Eigenvalue}-\ref{Assumption 5 Identification}
and \ref{Assumption 6 - Small Shifts}, (i) $\widehat{\lambda}_{b}\overset{P}{\rightarrow}\lambda_{0}$;
(ii) for every $\varepsilon>0$ there exists a $K>0$ such that for
all large $T,$ $P(T^{1-\kappa}\left|\widehat{\lambda}_{b}-\lambda_{0}\right|>K||\delta^{0}||^{-2}\overline{\sigma}^{2})<\varepsilon$;
and (iii) for $\kappa\in(0,\,1/4],$ $(\sqrt{T/N}(\widehat{\beta}-\beta^{0}),\,\sqrt{T/N}(\widehat{\delta}-\delta_{h}))'\overset{d}{\rightarrow}\mathscr{MN}\left(0,\,V\right)$
as $T\rightarrow\infty,$ with $V$ given in Proposition \ref{Proposition OLS Asymtptoc Distribu}.
\end{prop}
We first present a general result which shows that under Assumption
\ref{Assumption 6 - Small Shifts} one can obtain a shrinkage asymptotic
distribution similar to \citet{bai:97RES}. The latter exploits the
consistency of $\widehat{\lambda}_{b}$ and the fact that mixing conditions
implies that the regimes before and after $\lambda_{0}$ are asymptotically
independent.  Let $Z_{\Delta}\triangleq(0,\ldots,\,0,\,z_{\left(T_{b}+1\right)h},\ldots,\,z_{T_{b}^{0}h},\,0,\ldots,\,0)$
if $T_{b}<T_{b}^{0}$ and $Z_{\Delta}\triangleq(0,\ldots,\,0,\,z_{\left(T_{b}^{0}+1\right)h},\ldots,\,z_{T_{b}h},$
$0,\ldots,\,0)$ if $T_{b}>T_{b}^{0}$.
\begin{prop}
\label{Prop. Slow Time Scale}Under Assumption \ref{Assumption 1, CT}-\ref{Assumption 3 Break Date},
\ref{Assumption 4 Eigenvalue}-\ref{Assumption 5 Identification}
and \ref{Assumption 6 - Small Shifts}, 
\begin{align}
T^{1-\kappa}(\widehat{\lambda}_{b}-\lambda_{0}) & \overset{\mathcal{L}\mathrm{-}\mathrm{s}}{\Rightarrow}\underset{v\in\left(-\infty,\,\infty\right)}{\mathrm{argmax}}2(\delta^{0})'\mathscr{W}\left(v\right).\label{CR Asym Dist =00003D Bai 97}
\end{align}
\end{prop}
The distribution in Proposition \ref{Prop. Slow Time Scale} is different
from \citet{bai:97RES}. One can show that his distribution can be
obtained under a continuous record if Assumption \ref{Assumption 6 - Small Shifts}
is modified as follows: $\delta_{h}=\delta^{0}h^{\kappa/2}$, $T^{1-\kappa}\epsilon\rightarrow B<\infty$,
$0<\kappa\leq1/2$ and $\sigma_{h}\triangleq\overline{\sigma}h^{-\kappa/2}$.
This would result in, 
\begin{align}
T^{1-\kappa} & \left(\widehat{\lambda}_{b}-\lambda_{0}\right)\overset{\mathcal{L}\mathrm{-}\mathrm{s}}{\Rightarrow}\label{CR Asym Dist =00003D Bai 97-1}\\
 & \underset{v\in\left(-\infty,\,\infty\right)}{\mathrm{argmax}}\left\{ -\left(\delta^{0}\right)'\left\langle Z_{\Delta},\,Z_{\Delta}\right\rangle \left(v\right)\delta^{0}+2\left(\delta^{0}\right)'\mathscr{W}\left(v\right)\right\} .\nonumber 
\end{align}
The difference between \eqref{CR Asym Dist =00003D Bai 97} and \eqref{CR Asym Dist =00003D Bai 97-1}
is the presence of the drift (or deterministic) part $-\left(\delta^{0}\right)'$
$\left\langle Z_{\Delta},\,Z_{\Delta}\right\rangle \left(v\right)\delta^{0}$.
Without relating the magnitude of the break to the local variance
condition, the order of the stochastic part dominates that of the
deterministic part and so the latter vanishes asymptotically. The
distributions in \eqref{CR Asym Dist =00003D Bai 97}-\eqref{CR Asym Dist =00003D Bai 97-1}
share the same issues as Bai's and so they do not add any particular
insight. We therefore move to discuss how to obtain a more useful
continuous record asymptotic distribution. 

Consider the set $\mathscr{\mathcal{D}}\left(C\right)\triangleq\left\{ N_{b}:\,N_{b}\in\left\{ N_{b}^{0}+Ch^{1-\kappa}\right\} ,\,\left|C\right|<\infty\right\} $,
on the original time scale. Let $\psi_{h}\triangleq h^{1-k}$. Here
we use the same device as in \citeauthor{foster/nelson:96} (\citeyear{nelson/foster:94};
\citeyear{foster/nelson:96})\nocite{nelson/foster:94}. Different
scaling factors applied to an objective function can lead to different
asymptotic distributions. We normalize $Q_{T}(T_{b})$ by $\psi_{h}$,
where $\psi_{h}$ corresponds to the rate of convergence in Proposition
\ref{Prop 3 Asym}. The rate of convergence implicitly describes the
order of the terms in the expansion of $Q_{T}\left(T_{b}\right)-Q_{T}\left(T_{b}^{0}\right)$.
 
\begin{lem}
\label{Lemma 1}Under Assumption \ref{Assumption 1, CT}-\ref{Assumption 3 Break Date},
\ref{Assumption 4 Eigenvalue}-\ref{Assumption 5 Identification}
and \ref{Assumption 6 - Small Shifts},
\begin{align}
( & Q_{T}\left(T_{b}\right)-Q_{T}\left(T_{b}^{0}\right))/\psi_{h}\label{Eq. (1), Lemma 1}\\
 & =-\delta'_{h}\left(Z_{\Delta}'Z_{\Delta}/\psi_{h}\right)\delta_{h}+2\delta_{h}'\left(Z'_{\Delta}e/\psi_{h}\right)\mathrm{sgn}\left(T_{b}^{0}-T_{b}\right)+o_{p}\left(h^{1/2}\right).\nonumber 
\end{align}
\end{lem}
For brevity, we use the notation $\pm$ in place of $\mathrm{sgn}\left(T_{b}^{0}-T_{b}\right)$,
henceforth. The conditional first moment of the centered criterion
function $Q_{T}\left(T_{b}\right)-Q_{T}\left(T_{b}^{0}\right)$ is
of order $O\left(h^{1-\kappa}\right)$, i.e., it ``oscillates''
rapidly as $h\downarrow0$. Hence, in order to approximate the behavior
of $\{\widehat{T}_{b}-T_{b}^{0}\}$, we proceed as in Section 3 in
\citet{nelson/foster:94} and rescale ``time''. For any $C>0$,
let $L_{C}\triangleq N_{b}^{0}-Ch^{1-\kappa}$ and $R_{C}\triangleq N_{b}^{0}+Ch^{1-\kappa}$,
where $L_{C}$ and $R_{C}$ are the left and right boundary points
of $\mathcal{D}\left(C\right)$, respectively. We then have $\left|R_{C}-L_{C}\right|=O\left(Ch^{1-\kappa}\right)$.
Now, take the vanishingly small interval $\left[L_{C},\,R_{C}\right]$
on the original time scale, and stretch it into a time interval $\left[T^{1-\kappa}L_{C},\,T^{1-\kappa}R_{C}\right]$
on a new ``fast time scale''. Changing time scale simply means that
we rescale the objective function in such a way that it is of higher
order as $h\downarrow0$, i.e., it fluctuates less. This leads to
an asymptotic distribution that accounts for higher uncertainty. Yet,
under our framework it is still possible to consistently estimate
the break fraction and the regression coefficients so that inference
is feasible.

Since the criterion function is scaled by $\psi_{h}^{-1}$, all scaled
processes are $O_{p}\left(1\right)$. Now, let $N_{b}\left(v\right)=N_{b}^{0}-vh^{1-\kappa},\,v\in\left[-C,\,C\right]$.
Using Lemma \ref{Lemma 1} and Assumption \ref{Assumption 6 - Small Shifts}
(see the appendix), 
\begin{align*}
\psi_{h}^{-1} & \left(Q_{T}\left(T_{b}\left(v\right)\right)-Q_{T}\left(T_{b}^{0}\right)\right)=\\
 & -\delta'_{h}\left(\sum_{k=T_{b}\left(v\right)+1}^{T_{b}^{0}}\frac{z_{kh}}{\sqrt{\psi_{h}}}\frac{z'_{kh}}{\sqrt{\psi_{h}}}\right)\delta_{h}\pm2\left(\delta^{0}\right)'\sum_{k=T_{b}\left(v\right)+1}^{T_{b}^{0}}\frac{z_{kh}}{\sqrt{\psi_{h}}}\frac{\widetilde{e}{}_{kh}}{\sqrt{\psi_{h}}}+o_{p}\left(h^{1/2}\right).
\end{align*}
where $\widetilde{e}_{kh}\triangleq h^{1/4}e_{kh}.$ In addition,
in view of \eqref{Model Regressors Integral Form}, we let $dZ_{\psi,s}=\psi_{h}^{-1/2}\sigma_{Z,s}dW_{Z,s}$
for $s\in\left[N_{b}^{0}-vh^{1-\kappa},\,N_{b}^{0}+vh^{1-\kappa}\right]$.
Applying the time scale change $s\rightarrow t\triangleq\psi_{h}^{-1}s$
to all processes including $\Sigma^{0}$, we have $dZ_{\psi,t}=\sigma_{Z,t}dW_{Z,t}$
with $t\in\mathcal{T}\left(C\right)$, where $\mathcal{T}\left(C\right)\triangleq\{t:\,t\in\left[N_{b}^{0}+v\left\Vert \delta^{0}\right\Vert ^{2}/\overline{\sigma}^{2}\right],\,\left|v\right|\leq C\}$.
Therefore, 
\begin{align*}
\psi_{h}^{-1} & \left(Q_{T}\left(T_{b}\left(v\right)\right)-Q_{T}\left(T_{b}^{0}\right)\right)\\
 & =-\delta'_{h}\left(\sum_{k=T_{b}\left(v\right)+1}^{T_{b}^{0}}z_{\psi,kh}z'_{\psi,kh}\right)\delta_{h}\pm2\left(\delta^{0}\right)'\sum_{k=T_{b}\left(v\right)+1}^{T_{b}^{0}}z_{\psi,kh}\widetilde{e}_{\psi,kh}+o_{p}\left(h^{1/2}\right),
\end{align*}
with $NT_{b}\left(v\right)/T=N_{b}\left(v\right)=N_{b}^{0}+v$, where
$z_{\psi,kh}\triangleq z_{kh}/\sqrt{\psi_{h}}$ and $\widetilde{e}_{\psi,kh}\triangleq\widetilde{e}_{kh}/\sqrt{\psi_{h}}$.
Because of the change of time scale, all processes in the last display
are scaled up to be $O_{p}\left(1\right)$ and thus behave as diffusion-like
processes. On this new ``fast time scale'', we have $T^{1-\kappa}R_{C}-T^{1-\kappa}L_{C}=O\left(1\right)$
and $Q_{T}\left(T_{b}\left(v\right)\right)-Q_{T}\left(T_{b}^{0}\right)$
is restored to be $O_{p}\left(1\right)$. Observe that changing the
time scale does not affect any statistic which depends on observations
from $k=1$ to $k=\left\lfloor L_{C}/h\right\rfloor $ or from $k=\left\lfloor R_{C}/h\right\rfloor $
to $k=T$ (since these involve a positive fraction of data). However,
it does affect quantities which include observations that fall in
$\left[T_{b}h,\,T_{b}^{0}h\right]$ (assuming $T_{b}<T_{b}^{0}$).
In particular, on the original time scale, the processes $\left\{ D_{t}\right\} ,\,\left\{ Z_{t}\right\} $
and $\left\{ e_{t}\right\} $ are well-defined and scaled to be $O_{p}\left(1\right)$
while $Q_{T}\left(T_{b}\right)-Q_{T}\left(T_{b}^{0}\right)$ (asymptotically)
oscillates more rapidly than a simple diffusion-type process. On the
new ``fast time scale'', $\left\{ D_{t}\right\} ,\,\left\{ Z_{t}\right\} $
and $\left\{ e_{t}\right\} $ are not affected since they have the
same order in $\left[T^{1-\kappa}L_{C},\,T^{1-\kappa}R_{C}\right]$
as $h\downarrow0$. That is, the first conditional moments are $O\left(h\right)$
while the corresponding moments for $Q_{T}\left(T_{b}\right)-Q_{T}\left(T_{b}^{0}\right)$
on $\mathcal{T}\left(C\right)$ are restored to be $O\left(h\right)$.
As the continuous-time limit is approached, the rescaled criterion
function $\left(Q_{T}\left(T_{b}\left(v\right)\right)-Q_{T}\left(T_{b}^{0}\right)\right)/h^{1/2}$\textit{\textcolor{red}{{}
}}operates on a ``fast time scale'' on $\mathcal{T}\left(C\right)$. 

Our analysis is local; we examine the limiting behavior of the centered
and rescaled criterion function process in a neighborhood $\mathcal{T}\left(C\right)$
of the the true break date $N_{b}^{0}$ defined on a new time scale.
We first obtain the weak convergence results for the statistic $\left(Q_{T}\left(T_{b}\left(v\right)\right)-Q_{T}\left(T_{b}^{0}\right)\right)/h^{1/2}$
and then apply a continuous mapping theorem for the argmax functional.
However, it is convenient to work with a re-parametrized objective
function. Proposition \ref{Prop 3 Asym} allows us to use
\begin{align*}
\overline{Q}_{T}\left(\theta^{*}\right) & =\left(Q_{T}\left(\theta_{h},\,T_{b}\left(v\right)\right)-Q_{T}\left(\theta^{0},\,T_{b}^{0}\right)\right)/h^{1/2},
\end{align*}
where $\theta^{*}\triangleq\left(\theta'_{h},\,v\right)'$ with $T_{b}\left(v\right)\triangleq T_{b}^{0}+\left\lfloor v/h\right\rfloor $
and $T_{b}\left(v\right)$ is the time index on the ``fast time scale''.
The normalizing factor $\psi_{h}h^{1/2}$ allows us to change the
time scale and obtain an alternative asymptotic distribution. When
$v$ varies, $T_{b}\left(v\right)$ potentially visits all integers
between $1$ and $T$. Thus, on the new time scale, we need to introduce
the trimming parameter $\pi\in\left(0,\,1\right)$ which determines
the region where $T_{b}\left(v\right)$ can vary. We have the normalizations
$T_{b}\left(v\right)=T\pi$ if $T_{b}\left(v\right)\leq T\pi$ and
$T_{b}\left(v\right)=T\left(1-\pi\right)$ if $T_{b}\left(v\right)\geq T\left(1-\pi\right)$.
On the old time scale, $N_{b}\left(u\right)=N_{b}^{0}+u$ with $v\rightarrow\psi_{h}^{-1}u$,
so that $N_{b}\left(u\right)$ is in a vanishing neighborhood of $N_{b}^{0}$.
On $\mathcal{T}\left(C\right)$, we index the process $Q_{T}\left(\theta_{h},\,T_{b}\left(v\right)\right)-Q_{T}\left(\theta^{0},\,T_{b}^{0}\right)$
by two time subscripts: one referring to the time $T_{b}$ on the
original time scale and one referring to the time elapsed since $T_{b}h$
on the ``fast time scale''. For simplicity, we omit the former;
since the limiting distribution of the least-squares estimator will
now depend on the trimming we use the notation $\widehat{T}_{b,\pi}=T\widehat{\lambda}_{b,\pi}$
where $\widehat{\lambda}_{b,\pi}$ is the least-squares estimator
of the fractional break date associated to the fast time scale (i.e.,
associated to the  factor $\psi_{h}h^{1/2}$). 

The optimization problem is not affected by the change of time scale.
In fact, by Proposition \ref{Prop 3 Asym}, $u=Th(\widehat{\lambda}_{b}-\lambda_{0})=KO_{p}\left(h^{1-\kappa}\right)$
on the old time scale; whereas on the new ``fast time scale'',
$v=Th(\widehat{\lambda}_{b,\pi}-\lambda_{0})=O_{p}\left(1\right)$.
The maximization problem is not changed because $v/h$ can take any
value in $\mathbb{R}$. The process $Q_{T}\left(\theta_{h},\,T_{b}\left(v\right)\right)-Q_{T}\left(\theta^{0},\,T_{b}^{0}\right)$
is thus analyzed on a fixed horizon since $v$ now varies over $[(N\pi-N_{b}^{0})/(||\delta^{0}||^{-2}\overline{\sigma}^{2}),\,(N\left(1-\pi\right)-N_{b}^{0})/(\left\Vert \delta^{0}\right\Vert ^{-2}\overline{\sigma}^{2})]$.
Define  the modification to the set $\mathcal{D}\left(C\right)$
applicable to the new time scale by 
\begin{align*}
\mathcal{D}^{*}\left(C\right) & =\biggl\{\left(\beta^{0},\,\delta_{h},\,v\right):\,\left\Vert \theta^{0}\right\Vert \leq C;\,T_{b}\left(v\right)=T_{b}^{0}+vN^{-1}\left\Vert \delta^{0}\right\Vert ^{-2}\overline{\sigma}^{2};\\
 & \quad\frac{\left(N\pi-N_{b}^{0}\right)}{\left\Vert \delta^{0}\right\Vert ^{-2}\overline{\sigma}^{2}}\leq v\leq\frac{N\left(1-\pi\right)-N_{b}^{0}}{\left\Vert \delta^{0}\right\Vert ^{-2}\overline{\sigma}^{2}}\biggr\}.
\end{align*}
Let $\mathbb{D}\left(\mathcal{D}^{*}\left(C\right),\,\mathbb{R}\right)$
denote the space of all \textit{c\`{a}dl\`{a}g} functions from $\mathcal{D}^{*}\left(C\right)$
into $\mathbb{R}.$ Endow this space with the Skorokhod topology.
 Under a continuous record, we can apply limit theorems for statistics
involving (co)variation between regressors and errors. This enables
us to deduce the limiting process for $\overline{Q}_{T}\left(\theta^{*}\right)$,
 relying upon the work of Jacod \citeyearpar{jacod:94,jacod:97}
 and \citet{jacod/protter:98}.

To guide intuition, note that under the new re-parametrization, the
limit law of $\overline{Q}_{T}\left(\theta^{*}\right)$ is, according
to Lemma \ref{Lemma 1}, the same as the limit law of 
\begin{align*}
-h^{-1/2} & \delta_{h}'\left(Z_{\Delta}'Z_{\Delta}\right)\delta_{h}\pm2h^{-1/2}\delta_{h}'\left(Z'_{\Delta}e\right)\\
 & \overset{d}{\equiv}-\left(\delta^{0}\right)'\left(Z_{\Delta}'Z_{\Delta}\right)\delta^{0}\pm2h^{-1/2}\left(\delta^{0}\right)'h^{1/4}\left(Z'_{\Delta}h^{-1/4}\widetilde{e}\right),
\end{align*}
where $\overset{d}{\equiv}$ denotes (first order) equivalence in
law, and since (approximately) $e_{kh}\sim i.n.d.\,\mathscr{N}(0,$
$\sigma_{h,k-1}^{2}h),\,\sigma_{h,k}=\sigma_{h}\sigma_{e,k}$ then
$\widetilde{e}_{kh}\sim i.n.d.\,\mathscr{N}\left(0,\,\sigma_{e,k-1}^{2}h\right)$.
Hence, the limit law of $\overline{Q}_{T}\left(\theta^{*}\right)$
is, to first-order, equivalent to the law of 
\begin{align}
-\left(\delta^{0}\right)'\left(Z_{\Delta}'Z_{\Delta}\right)\delta^{0}\pm2\left(\delta^{0}\right)'\left(h^{-1/2}Z'_{\Delta}\widetilde{e}\right) & .\label{Eq Intuition Limit Law}
\end{align}
 We apply a law of large numbers to the first term and a stable convergence
in law under the Skorokhod topology to the second. Assumption \ref{Assumption 6 - Small Shifts}
combined with the normalizing factor $h^{-1/2}$ in $\overline{Q}_{T}\left(\theta^{*}\right)$
account for the discrepancy between the deterministic and stochastic
component in \eqref{Eq Intuition Limit Law}. 

Having outlined the main steps in the arguments used to derive the
continuous records limit distribution of the break date estimate,
we now state the main result of this section. The limiting process
is realized on a extension of the original probability space and we
relegate this description to Section \ref{subsection Description Limiting Process}
in the supplement.  
\begin{thm}
\label{Theorem 1}Under Assumption \ref{Assumption 1, CT}-\ref{Assumption 3 Break Date},
\ref{Assumption 4 Eigenvalue}-\ref{Assumption 5 Identification}
and \ref{Assumption 6 - Small Shifts},
\begin{align}
N & \left(\widehat{\lambda}_{b,\pi}-\lambda_{0}\right)\overset{\mathcal{L}\mathrm{-}\mathrm{s}}{\Rightarrow}\underset{v\in\mathcal{A}}{\mathrm{argmax}}\left\{ -\left(\delta^{0}\right)'\left\langle Z_{\Delta},\,Z_{\Delta}\right\rangle \left(v\right)\delta^{0}+2\left(\delta^{0}\right)'\mathscr{W}\left(v\right)\right\} ,\label{CR Asymptotic Distribution}
\end{align}
where 
\begin{align*}
\mathcal{A} & \triangleq\left[\frac{N\pi-N_{b}^{0}}{\left\Vert \delta^{0}\right\Vert ^{-2}\overline{\sigma}^{2}},\,\frac{N\left(1-\pi\right)-N_{b}^{0}}{\left\Vert \delta^{0}\right\Vert ^{-2}\overline{\sigma}^{2}}\right].
\end{align*}
\end{thm}
Note the differences between the results in Theorem \ref{Theorem 1}
and in Proposition \ref{Prop. Slow Time Scale}. First, on the fast
time scale, $\widehat{\lambda}_{b,\pi}$ behaves as an inconsistent
estimator for $\lambda_{0}$ for $N$ fixed, but it is consistent
as $N\rightarrow\infty$. On the original time scale $\widehat{\lambda}_{b}$
is not only consistent for $\lambda_{0}$ but it also enjoys a similar
asymptotic distribution as in \citet{bai:97RES}. Second, the asymptotic
distribution of $\widehat{\lambda}_{b,\pi}$ depends on the span of
the data and consequently on the trimming $\pi.$ Proposition \ref{Prop. Slow Time Scale},
in contrast, suggests that the span, the trimming and the location
of the break are irrelevant for the limiting behavior of the estimator.
This intuitively follows from the fact that under the original time
scale the break date estimator is consistent. We will show that indeed
the span of the data and the location of the break influence the finite-sample
properties of the least-squares estimator, and that Theorem \ref{Theorem 1}
provides a more useful approximation. An important implication of
Theorem \ref{Theorem 1} is that the precision of the estimator depends
more on the span $N$ than to the number of observations $T$. 

Unlike Bai's distribution, the distribution in Theorem \ref{Theorem 1}
involves the location of the maximum of a function of the (quadratic)
variation of the regressors and of a two-sided centered Gaussian martingale
process over the interval $[(N\pi-N_{b}^{0})/(\left\Vert \delta^{0}\right\Vert ^{-2}\overline{\sigma}^{2}),\,(N\left(1-\pi\right)-N_{b}^{0})/(||\delta^{0}||^{-2}\overline{\sigma}^{2})]$.
Notably, this domain depends on the true value of $N_{b}^{0}$ and
therefore the limit distribution is asymmetric, in general. The degree
of asymmetry increases as the true break point moves away from mid-sample.
This holds even when the distributions of the errors and regressors
are the same in the pre- and post-break regimes. The presence of the
trimming confirms that the span of the (trimmed) data affects the
 limit distribution. It is well-known that the least-squares estimator
of the break date can be sensitive to trimming {[}see \citet{bai/perron:03}
for some recommendations on the trimming choice{]}. Our asymptotic
theory accommodates this property of the least-squares estimator while
others do not.

Additional relevant remarks follow; more details are provided in the
supplement. The magnitude of the break plays a key role in determining
the density of the asymptotic distribution. More precisely, the density
displays interesting properties which change when the signal-to-noise
ratio as well as other parameters of the model change. Moreover, the
distribution in Theorem \ref{Theorem 1} is able to reproduce important
features of the small-sample results obtained via simulations {[}e.g.,
\citet{bai/perron:06}{]}. First, the second moments of the regressors
impact the asymptotic mean as well as the second-order behavior of
the break point estimator (e.g., the persistence of the regressors
influences the finite-sample performance of the estimator). Second,
the continuous record setting manages to preserve information about
the time span $N$ of the data, a clear advantage since the location
of the true break point matters for the small-sample distribution
of the estimator. It has been shown via simulations that in small-samples
the break point estimator tends to be imprecise if the break size
is small, and some bias arises if the break point is not at mid-sample.
In our framework, the (trimmed) time horizon $\left[N\pi,\,N\left(1-\pi\right)\right]$
is fixed and thus we can distinguish between the statistical content
of the segments $\left[N\pi,\,N_{b}^{0}\right]$ and $\left[N_{b}^{0},\,N\left(1-\pi\right)\right]$.
In contrast, this is not feasible under the classical shrinkage large-$N$
asymptotics because both the pre- and post-break segments increase
 proportionately and mixing conditions are imposed so that the only
relevant information is a neighborhood around the true break date.
Details on how to simulate the limiting distribution in Theorem \ref{Theorem 1}
are given in Section \ref{subsec:Simulation-of-the LD in Theorem 4.1}
of the supplement.

We further characterize the asymptotic distribution by exploiting
the ($\mathscr{F}$-conditionally) Gaussian property of the limit
process. The analysis also holds unconditionally if we assume that
the volatility processes are non-stochastic. Thus, as in the classical
setting, we begin with a second-order stationarity assumption within
each regime. The following assumption guarantees that the results
below remain valid without the need to condition on $\mathscr{F}.$
\begin{assumption}
\label{Assumtpion - Regimes}The process $\Sigma^{0}$ is (possibly
time-varying) deterministic; $\left\{ z_{kh},\,e_{kh}\right\} $ is
second-order stationary within each regime. For $k=1,\ldots,\,T_{b}^{0}$,
$\mathbb{E}(z_{kh}z'_{kh}|\,\mathscr{F}_{\left(k-1\right)h})=\Sigma_{Z,1}h$,
$\mathbb{E}(\widetilde{e}_{kh}^{2}|$ $\mathscr{F}_{\left(k-1\right)h})=\sigma_{e,1}^{2}h$
and $\mathbb{E}(z_{kh}z'_{kh}\widetilde{e}_{kh}^{2}|\,\mathscr{F}_{\left(k-1\right)h})=\Omega_{\mathscr{W},1}h^{2}$
while for $k=T_{b}^{0}+1,\ldots,\,T$, $\mathbb{E}(z_{kh}z'_{kh}|$
$\mathscr{F}_{\left(k-1\right)h})=\Sigma_{Z,2}h$, $\mathbb{E}(\widetilde{e}_{kh}^{2}|\,\mathscr{F}_{\left(k-1\right)h})=\sigma_{e,2}^{2}h$
and $\mathbb{E}(z_{kh}z'_{kh}\widetilde{e}_{kh}^{2}|\,\mathscr{F}_{\left(k-1\right)h})=\Omega_{\mathscr{W},2}h^{2}$.
\end{assumption}
Let $W_{i}^{*},$ $i=1,\,2,$ be two independent standard Wiener processes
defined on $[0,\,\infty),$ starting at the origin when $s=0.$ Let
\begin{align*}
\mathscr{V}\left(s\right) & =\begin{cases}
-\frac{\left|s\right|}{2}+W_{1}^{*}\left(s\right), & \textrm{if }s<0\\
-\frac{\left(\delta^{0}\right)'\Sigma_{Z,2}\delta^{0}}{\left(\delta^{0}\right)'\Sigma_{Z,1}\delta^{0}}\frac{\left|s\right|}{2}+\left(\frac{\left(\delta^{0}\right)'\Omega_{\mathscr{W},2}\left(\delta^{0}\right)}{\left(\delta^{0}\right)'\Omega_{\mathscr{W},1}\left(\delta^{0}\right)}\right)^{1/2}W_{2}^{*}\left(s\right), & \textrm{if }s\geq0.
\end{cases}
\end{align*}

\begin{thm}
\label{Theorem 2, Asymptotic Distribution immediate Stationary Regimes}Under
Assumption \ref{Assumption 1, CT}-\ref{Assumption 3 Break Date},
\ref{Assumption 4 Eigenvalue}-\ref{Assumption 5 Identification}
and \ref{Assumption 6 - Small Shifts}-\ref{Assumtpion - Regimes},
\begin{align}
\frac{\left(\left(\delta^{0}\right)'\left\langle Z,\,Z\right\rangle _{1}\delta^{0}\right)^{2}}{\left(\delta^{0}\right)'\Omega_{\mathscr{W},1}\delta^{0}}N\left(\widehat{\lambda}_{b,\pi}-\lambda_{0}\right) & \Rightarrow\underset{s\in\mathcal{A}^{*}}{\mathrm{argmax}}\mathscr{V}\left(s\right),\label{Equation (2) Asymptotic Distribution}
\end{align}
 where 
\begin{align*}
\mathcal{A}^{*} & \triangleq\left[\frac{N\pi-N_{b}^{0}}{\left\Vert \delta^{0}\right\Vert ^{-2}\overline{\sigma}^{2}}\frac{\left(\left(\delta^{0}\right)'\left\langle Z,\,Z\right\rangle _{1}\delta^{0}\right)^{2}}{\left(\delta^{0}\right)'\Omega_{\mathscr{W},1}\left(\delta^{0}\right)},\,\frac{N\left(1-\pi\right)-N_{b}^{0}}{\left\Vert \delta^{0}\right\Vert ^{-2}\overline{\sigma}^{2}}\frac{\left(\left(\delta^{0}\right)'\left\langle Z,\,Z\right\rangle _{1}\delta^{0}\right)^{2}}{\left(\delta^{0}\right)'\Omega_{\mathscr{W},1}\delta^{0}}\right].
\end{align*}
\end{thm}
Unlike the asymptotic distribution derived under classical large-$N$
asymptotics, the probability density  in \eqref{Equation (2) Asymptotic Distribution}
is not available in closed form. Furthermore, the limiting distribution
depends on unknown quantities. In the next section we explain how
one can derive a feasible counterpart. This will be useful to characterize
the main features of interest that will guide us in devising methods
to construct confidence sets for $T_{b}^{0}$.

\section{Feasible Approximations to the Finite-Sample Distributions\label{Section Approximation-to-the}}

In Section \ref{subsec:Evaluation-of-the} we propose a feasible version
of our limit theory and compare it with the finite-sample distribution.
In Section \ref{Subsec:Comparison-with Bai and EM} we discuss some
differences between our approach and others. Let 
\begin{align*}
\rho=\frac{\left(\left(\delta^{0}\right)'\left\langle Z,\,Z\right\rangle _{1}\delta^{0}\right)^{2}}{\left(\left(\delta^{0}\right)'\Omega_{\mathscr{W},1}\delta^{0}\right)},\quad\xi_{1}=\frac{\left(\delta^{0}\right)'\left\langle Z,\,Z\right\rangle _{2}\delta^{0}}{\left(\delta^{0}\right)'\left\langle Z,\,Z\right\rangle _{1}\delta^{0}},\quad & \xi_{2}=\frac{\left(\delta^{0}\right)'\Omega_{\mathscr{W},2}\delta^{0}}{\left(\delta^{0}\right)'\Omega_{\mathscr{W},1}\delta^{0}}.
\end{align*}

\subsection{\label{subsec:Evaluation-of-the}A Feasible Version of the Limit
Distribution}

\begin{figure}[H]
\includegraphics[width=18cm,height=7.5cm]{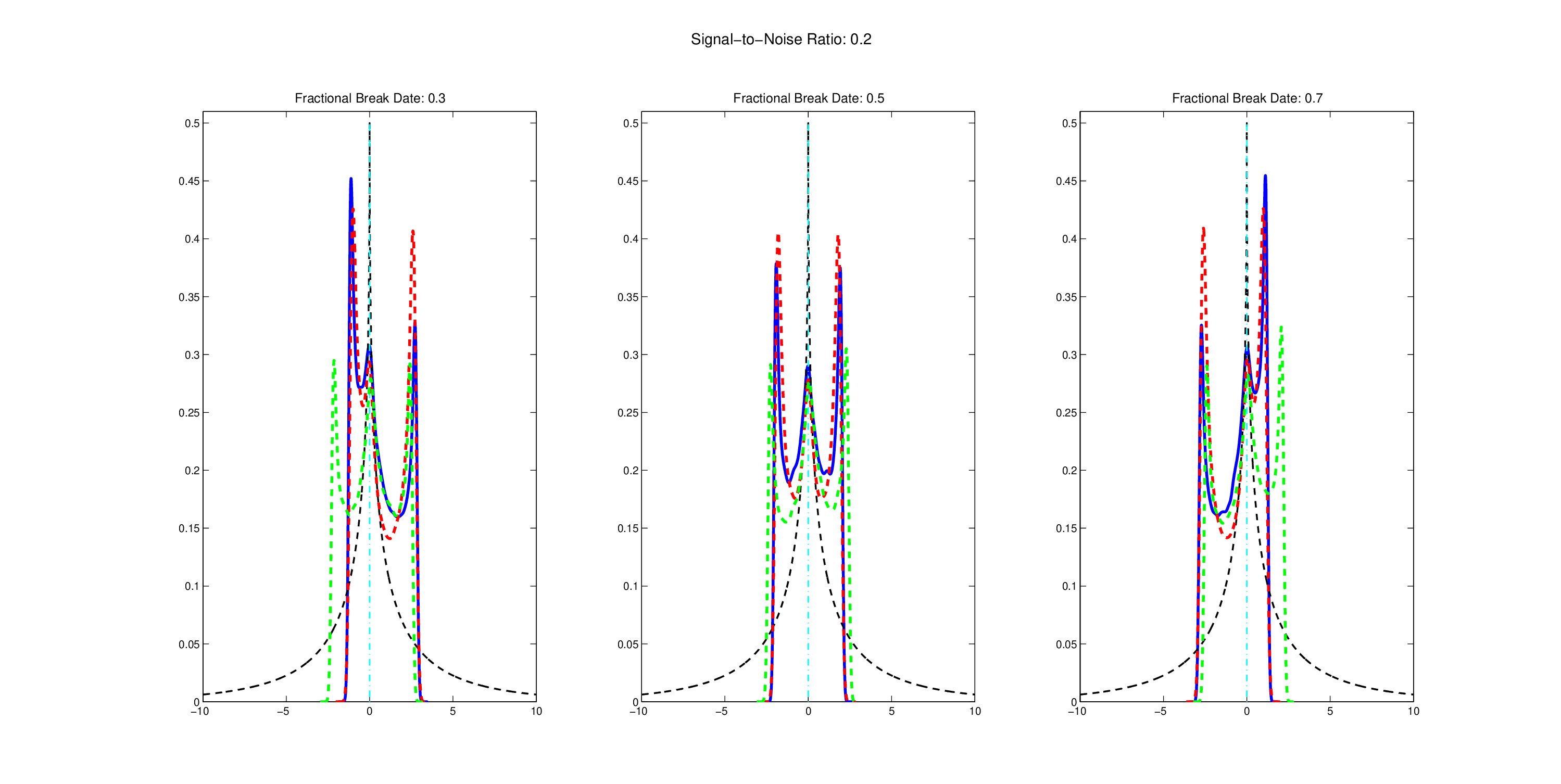}

{\footnotesize{}\caption{{\footnotesize{}\label{Fig19}}{\scriptsize{}The probability density
of $\rho\left(\widehat{T}_{b,\pi}-T_{b}^{0}\right)$ for model (5.1)
with break magnitude $\delta_{Z}^{0}=0.2$ and true break fraction
$\lambda_{0}=0.3,\,0.5$ and $0.7$ (the left, middle and right panel,
respectively). The signal-to-noise ratio is $\delta_{Z}^{0}/\sigma_{e}=\delta_{Z}^{0}$
since $\sigma_{e}^{2}=1$. The blue solid (green broken) line is the
density of the infeasible (reps. feasible) asymptotic distribution
derived under a continuous record, the black broken line is the density
of the asymptotic distribution from Bai (1997) and the red broken
line is the density of the finite-sample distribution.}}
}{\footnotesize\par}
\end{figure}

\begin{figure}[H]
\includegraphics[width=18cm,height=7.5cm]{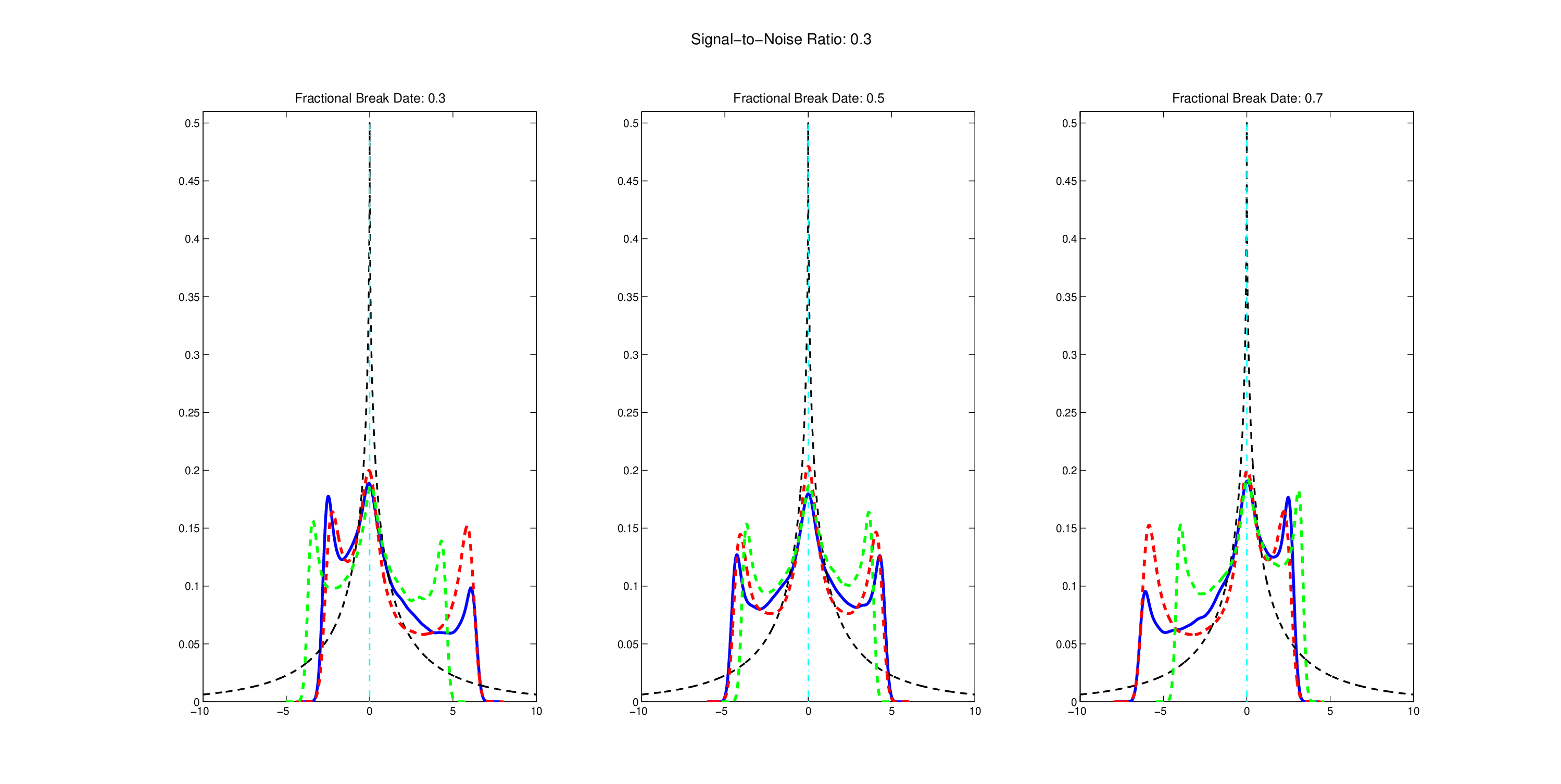}

{\footnotesize{}\caption{{\scriptsize{}\label{Fig20}The probability density of $\rho\left(\widehat{T}_{b,\pi}-T_{b}^{0}\right)$
for model (5.1) with break magnitude $\delta_{Z}^{0}=0.3$ and true
break fraction $\lambda_{0}=0.3,\,0.5$ and $0.7$ (the left, middle
and right panel, respectively). The signal-to-noise ratio is $\delta_{Z}^{0}/\sigma_{e}=\delta_{Z}^{0}$
since $\sigma_{e}^{2}=1$. The blue solid (green broken) line is the
density of the infeasible (reps. feasible) asymptotic distribution
derived under a continuous record, the black broken line is the density
of the asymptotic distribution from Bai (1997) and the red broken
line is the density of the finite-sample distribution.}}
}{\footnotesize\par}
\end{figure}

\begin{figure}[H]
\includegraphics[width=18cm,height=7.5cm]{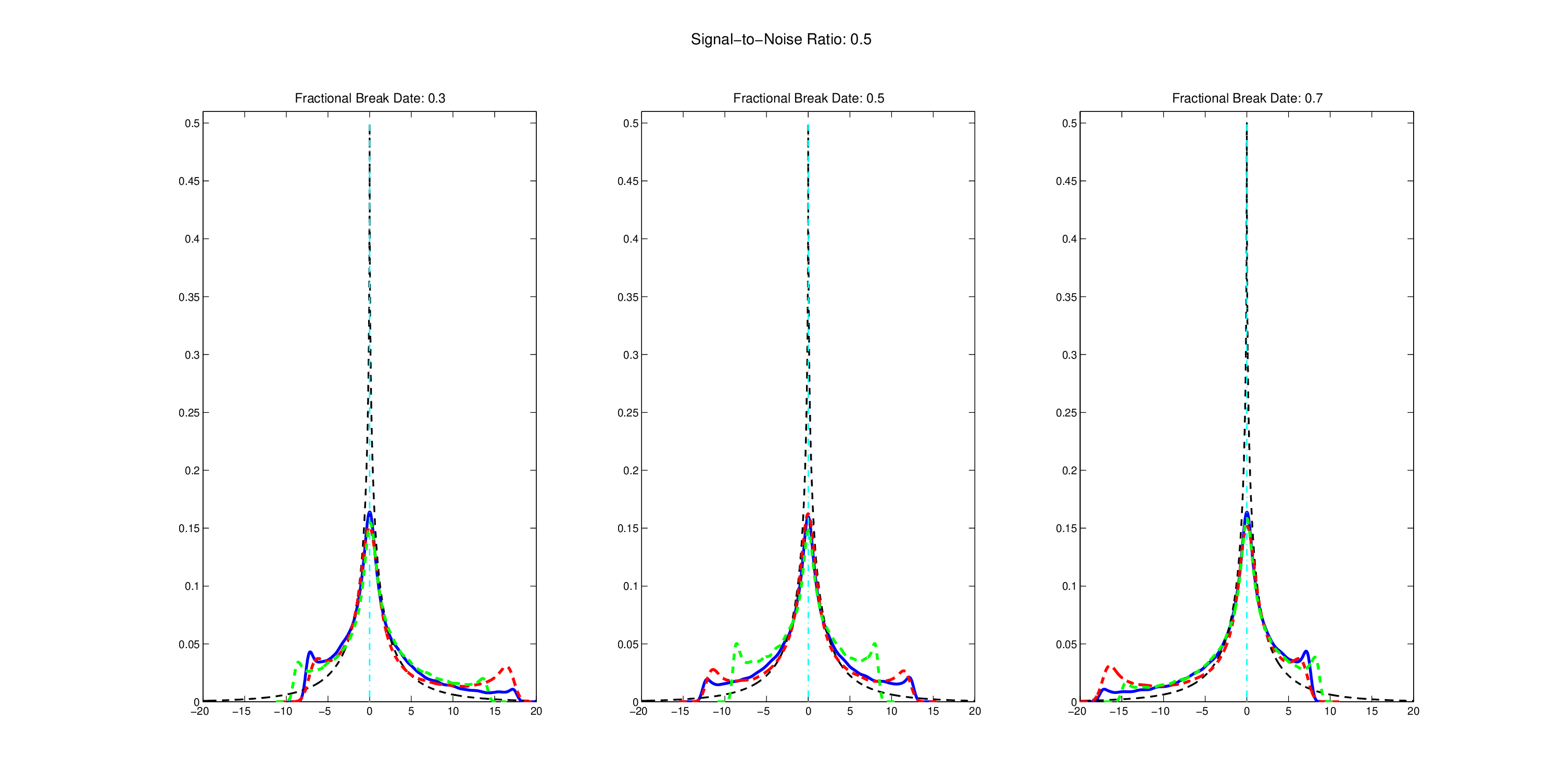}

{\footnotesize{}\caption{{\footnotesize{}\label{Fig21}}{\scriptsize{}The probability density
of $\rho\left(\widehat{T}_{b,\pi}-T_{b}^{0}\right)$ for model (5.1)
with break magnitude $\delta_{Z}^{0}=0.5$ and true break fraction
$\lambda_{0}=0.3,\,0.5$ and $0.7$ (the left, middle and right panel,
respectively). The signal-to-noise ratio is $\delta_{Z}^{0}/\sigma_{e}=\delta_{Z}^{0}$
since $\sigma_{e}^{2}=1$. The blue solid (green broken) line is the
density of the infeasible (reps. feasible) asymptotic distribution
derived under a continuous record, the black broken line is the density
of the asymptotic distribution from Bai (1997) and the red broken
line is the density of the finite-sample distribution.}}
}{\footnotesize\par}
\end{figure}

\begin{figure}[H]
\includegraphics[width=18cm,height=7.5cm]{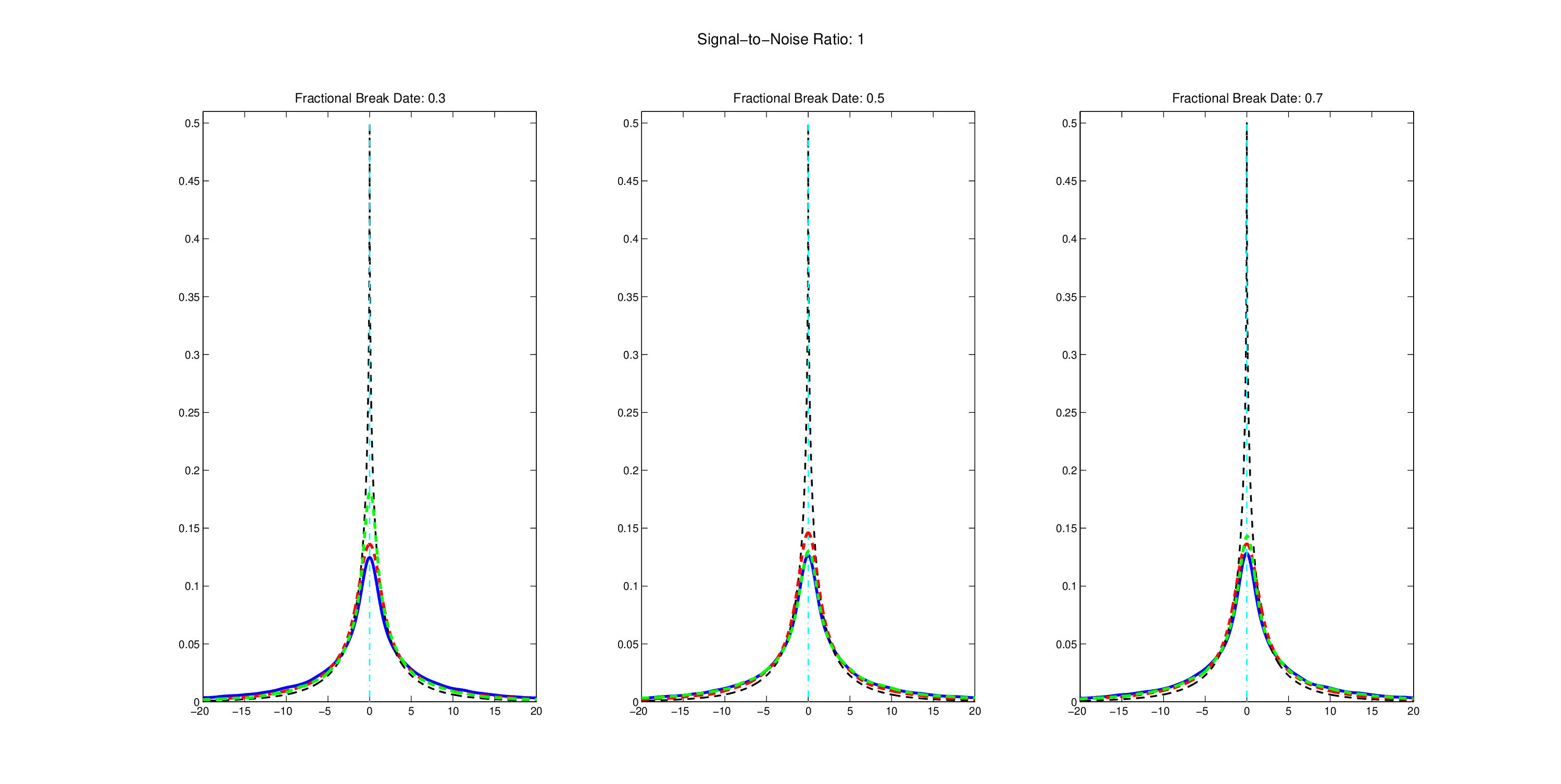}

{\footnotesize{}\caption{{\footnotesize{}\label{Fig23}}{\scriptsize{}The probability density
of $\rho\left(\widehat{T}_{b,\pi}-T_{b}^{0}\right)$ for model (5.1)
with break magnitude $\delta_{Z}^{0}=1$ and true break fraction $\lambda_{0}=0.3,\,0.5$
and $0.7$ (the left, middle and right panel, respectively). The signal-to-noise
ratio is $\delta_{Z}^{0}/\sigma_{e}=\delta_{Z}^{0}$ since $\sigma_{e}^{2}=1$.
The blue solid (green broken) line is the density of the infeasible
(reps. feasible) asymptotic distribution derived under a continuous
record, the black broken line is the density of the asymptotic distribution
from Bai (1997) and the red broken line is the density of the finite-sample
distribution.}}
}{\footnotesize\par}
\end{figure}

In order to use the continuous record asymptotic distribution in practice
one needs consistent estimates of the unknown quantities. In this
section, we compare the finite-sample distribution of the least-squares
estimator of the change-point date with a feasible version of the
continuous record asymptotic distribution obtained with plug-in estimates.
We obtain the finite-sample distribution of $\rho(\widehat{T}_{b,\pi}-T_{b}^{0})$
based on 100,000 simulations from the following model:
\begin{align}
Y_{t}=D'_{t}\nu^{0}+Z'_{t}\beta^{0}+Z'_{t}\delta_{Z}^{0}\boldsymbol{1}_{\left\{ t>T_{b}^{0}\right\} }+e_{t}, & \qquad t=1,\ldots,\,T,\label{Casin i-  Model for FS vs Asy Dist}
\end{align}
 where $Z_{t}=0.5Z_{t-1}+u_{t}$ with $u_{t}\sim i.i.d.\,\mathscr{N}\left(0,\,1\right)$
independent of $e_{t}\sim i.i.d.\,\mathscr{N}\left(0,\,\sigma_{e}^{2}\right)$,
$\sigma_{e}^{2}=1$, $\nu^{0}=1$, $Z_{0}=0$, $D_{t}=1$ for all
$t$, and $T=100.$ We set $\pi=0.05$, $T_{b}^{0}=\left\lfloor T\lambda_{0}\right\rfloor $
with $\lambda_{0}=0.3,\,0.5,\,0.7$ and consider different break sizes
$\delta_{Z}^{0}=0.2,\,0.3,\,0.5,\,1$. The infeasible continuous record
asymptotic distribution is computed assuming knowledge of the data
generating process (DGP) as well as of the model parameters, i.e.,
using Theorem \ref{Theorem 2, Asymptotic Distribution immediate Stationary Regimes}
where we set $N_{b}^{0}$ equal to its true value, $||\delta^{0}||^{-2}\overline{\sigma}^{2}=||\delta_{Z}^{0}||^{-2}\sigma_{e}^{2},$
and $\xi_{1},\,\xi_{2}$ and $\rho$ equal to their true values, respectively,
with $\delta_{Z}^{0}$ in place of $\delta^{0}$. Note that the scaling
$h^{\kappa/2}$ and $h^{-1/4}$ in the definition of $\delta_{h}$
and $\sigma_{h}$ respectively, cancel using the fact that they appear
in both numerator and denominator and applying a change in variables.
The feasible counterparts are constructed with plug-in estimates of
$\xi_{1},\,\xi_{2},\,\rho$ and $(N_{b}^{0}\left\Vert \delta^{0}\right\Vert ^{2}/\overline{\sigma}^{2})\rho$.
In practice we need to use a normalization for $N$. A common choice
is $N=1$. Then $\widehat{\lambda}_{b}^{\mathrm{}}=\widehat{T}_{b}/T$
 is a natural estimate of $\lambda_{0}$, using the consistency result
of $\widehat{\lambda}_{b}^{\mathrm{}}$  that holds in the setting
of Theorem \ref{Theorem 1} which can also be rationalized for large
$N$ under the conditions of Theorem \ref{Theorem 1}. In practice
this means that we approximate the distribution of the estimator $\widehat{\lambda}_{b,\pi}$
where $\pi$ is chosen by the researcher and we plug-in the estimator
$\widehat{\lambda}_{b}$ which can be based on any  trimming because
of the consistency property. Here, we set $\widehat{\lambda}_{b}$
equal to the least-squares estimator based on a trimming $0.15$,
which is also used for the other plug-in estimates. The estimates
of $\xi_{1}$ and $\xi_{2}$ are given, respectively, by 
\begin{align*}
\widehat{\xi}_{1}=\frac{\widehat{\delta}'\left(T-\widehat{T}_{b}^{\mathrm{}}\right)^{-1}\sum_{k=\widehat{T}_{b}^{\mathrm{}}+1}^{T}z_{kh}z'_{kh}\widehat{\delta}}{\widehat{\delta}'\left(\widehat{T}_{b}^{\mathrm{}}\right)^{-1}\sum_{k=1}^{\widehat{T}_{b}^{\mathrm{}}}z_{kh}z'_{kh}\widehat{\delta}} & ,\qquad\widehat{\xi}_{2}=\frac{\widehat{\delta}'\left(T-\widehat{T}_{b}^{\mathrm{}}\right)^{-1}\sum_{k=\widehat{T}_{b}^{\mathrm{}}+1}^{T}\widehat{e}_{kh}^{2}z_{kh}z'_{kh}\widehat{\delta}}{\widehat{\delta}'\left(\widehat{T}_{b}^{\mathrm{}}\right)^{-1}\sum_{k=1}^{\widehat{T}_{b}^{\mathrm{}}}\widehat{e}_{kh}^{2}z_{kh}z'_{kh}\widehat{\delta}},
\end{align*}
where $\widehat{\delta}$ is the least-squares estimator of $\delta_{h}$
and $\widehat{e}_{kh}$ are the least-squares residuals. Note that
in $\widehat{\xi}_{1}$ and $\widehat{\xi}_{2}$, the estimate $\widehat{\delta}$
appears in both numerator and denominator so that the scaling $h^{\kappa/2}$
in the definition of $\delta_{h}$ cancels. Use is made of the fact
that $\left\langle Z,\,Z\right\rangle _{1}$ is consistently estimated
by $\sum_{k=1}^{\widehat{T}_{b}^{\mathrm{}}}z_{kh}z'_{kh}/\widehat{\lambda}_{b}^{\mathrm{}}$
while $\Omega_{\mathscr{W},1}$ is consistently estimated by $T\sum_{k=1}^{\widehat{T}_{b}^{\mathrm{}}}\widehat{e}_{kh}^{2}z_{kh}z'_{kh}/\widehat{\lambda}_{b}^{\mathrm{}}$.
The method to estimate $\lambda_{0}\left\Vert \delta^{0}\right\Vert ^{2}\overline{\sigma}^{-2}\rho$
is less immediate because it involves manipulating the scaling of
each of the three estimates. Let $\vartheta=\left\Vert \delta^{0}\right\Vert ^{2}\overline{\sigma}^{-2}\rho$.
We use the following estimates for $\vartheta$ and $\rho$, respectively,
\begin{align*}
\widehat{\vartheta}= & \widehat{\rho}\left\Vert \widehat{\delta}\right\Vert ^{2}\left(T^{-1}\sum_{k=1}^{T}\widehat{e}_{kh}^{2}\right)^{-1},\qquad\widehat{\rho}=\frac{\left(\widehat{\delta}'\left(\widehat{T}_{b}^{\mathrm{}}\right)^{-1}\sum_{k=1}^{\widehat{T}_{b}^{\mathrm{}}}z_{kh}z'_{kh}\widehat{\delta}\right)^{2}}{\widehat{\delta}'\left(\widehat{T}_{b}^{\mathrm{}}\right)^{-1}\sum_{k=1}^{\widehat{T}_{b}^{\mathrm{}}}\widehat{e}_{kh}^{2}z_{kh}z'_{kh}\widehat{\delta}},
\end{align*}
Whereas we have $\widehat{\xi}_{i}\overset{p}{\rightarrow}\xi_{i}$
$\left(i=1,\,2\right)$, the corresponding approximations for $\widehat{\rho}$
and $\widehat{\vartheta}$ are given by $\widehat{\rho}/h^{\kappa}\overset{p}{\rightarrow}\rho$
and $\widehat{\vartheta}/h^{2\kappa}\overset{p}{\rightarrow}\vartheta$.
 However, before letting $T\rightarrow\infty$ we can apply a change
in variable using the fact that $\widehat{\lambda}_{b}-\lambda_{b}^{0}=O\left(h^{1-\kappa}\right)$
which result in the extra factor $h^{2\kappa}$ canceling. 
\begin{prop}
\label{Proposition CI}Under the conditions of Theorem \ref{Theorem 2, Asymptotic Distribution immediate Stationary Regimes},
\eqref{Equation (2) Asymptotic Distribution} holds when using $\widehat{\xi}_{1},\,\widehat{\xi}_{2},\,\widehat{\rho}$
and $\widehat{\vartheta}$ in place of $\xi_{1},\,\xi_{2},\,\rho$
and $\vartheta$, respectively. 
\end{prop}
The proposition implies that the limiting distribution can be simulated
using plug-in estimates. This allows feasible inference about the
break date. The results are presented in Figure \ref{Fig19}-\ref{Fig23}
which also plot the asymptotic distribution from \citet{bai:97RES}
and the infeasible distribution from Theorem \ref{Theorem 2, Asymptotic Distribution immediate Stationary Regimes}.
Here by signal-to-noise ratio we mean $\delta_{Z}^{0}/\sigma_{e}$
which, given $\sigma_{e}^{2}=1$, equals the break size $\delta_{Z}^{0}$.

Several interesting observations appear at the outset. The density
of the large-$N$ shrinkage asymptotic distribution does not depend
on the location of the break, and thus it is always unimodal and symmetric
about the origin. None of these features are shared by the density
derived under a continuous record. When the true break is at mid-sample
$\left(\lambda_{0}=0.5\right)$, the density function is symmetric
and centered at zero. However, when the signal-to-noise ratio is low,
the density features three modes. This tri-modality vanishes as the
signal-to-noise ratio increases. When $\delta_{Z}^{0}$ is low and
the break is not at mid-sample the density is asymmetric; for values
of $\lambda_{0}$ less (larger) than 0.5, the density is right (left)
skewed. When the signal is low and $\lambda_{0}$ is less (larger)
than 0.5, the density has highest mode at some value near $\widehat{\lambda}_{b}$
being close to the starting (end) sample point than centered at $\lambda_{0}.$
However, as in the case of $\lambda_{0}=0.5,$ when the signal-to-noise
ratio increases the highest mode is centered at a value which corresponds
to $\widehat{\lambda}_{b}$ being close to $\lambda_{0}$. Asymmetry
and multi-modality of the finite-sample distribution of the break
point estimator were also found by \citet{perron/zhu:05} and \citet{deng/perron:06}
in models with a trend. 

The interpretation of these features are straightforward. For example,
asymmetry reflects the fact that the span of the data and the actual
location of the break play a crucial role on the behavior of the estimator.
If the break occurs early in the sample there is a tendency to overestimate
the break date and vice-versa if the break occurs late in the sample.
The marked changes in the shape of the density as we raise $\delta_{Z}^{0}$
confirms that the magnitude of the shift matters a great deal as well.
The tri-modality of the density when the shift size is small reflects
the uncertainty in the data as to whether a structural change is present
at all; i.e., the least-squares estimator finds it easier to locate
the break at either the beginning or the end of the sample. Unlike
the shrinkage asymptotic distribution, the density of the feasible
version of the continuous record distribution provides a remarkably
good approximation to the infeasible one and thus also to the finite-sample
distribution. The extended working paper \citet{casini/perron_CR_Single_Break_Extended}
shows that the quality of the approximation is good for a variety
of models. 

\subsection{\label{Subsec:Comparison-with Bai and EM}Comparison with Other Approaches}

The figures reported above have shown that there is a high degree
of uncertainty when the break magnitude is not large. The classical
shrinkage asymptotics of \citet{bai:97RES} with $\delta_{T}$ required
to convergence to zero at a rate slower than $O(T^{-1/2})$ clearly
underestimates that degree of uncertainty and, as the figures show,
provides a poor approximation to the finite-sample behavior of the
least-squares estimator. In Section \ref{Section Small-Sample-Effectiveness-of}
we show that this issue is responsible for the poor coverage probabilities
of the confidence intervals introduced in \citet{bai:97RES} when
the break magnitude is small. On the other hand, \citet{elliott/mueller:07}
and \citet{elliott/mueller/watson:15} require $\delta_{T}$ to go
to zero at the fast rate $O(T^{-1/2})$ leading to weak identification.
The latter implies that the relevant quantities in the model become
inconsistent. This can be problematic for inference and indeed, their
inference often suffers from the opposite problem in that confidence
intervals for $\widehat{T}_{b}$ can be too large {[}Casini and Perron
(\citeyear{casini/perron_Oxford_Survey}; \citeyear{casini/perron_Lap_CR_Single_Inf})
and \citet{chang/perron:18}{]}. 

We impose conditions on the signal-to-noise ratio $\delta/\sigma$
rather than just on $\delta.$ Consider a simple location model with
a change $\delta$ in the mean and independent errors. What describes
the uncertainty about the break in this model is the ratio $\delta/\sigma$
where $\sigma$ is the volatility of the errors. We let $\delta$
go to zero at a not too fast rate while letting $\sigma$ increase
to infinity in a neighborhood of $T_{b}^{0}$. That is $\left(\delta_{T}/\sigma_{t}\right)\rightarrow0$
at rate $O(T^{-1/2})$ in a neighborhood of $T_{b}^{0}$. Interestingly,
this is the same rate Elliott and M{\"u}ller used for $\delta_{T}\rightarrow0.$
Away from $T_{b}^{0}$, we require $\left(\delta_{T}/\sigma_{t}\right)\rightarrow0$
at slower rate---similar to \citet{yao:87} and \citet{bai:97RES}.
The difference now is that we do not lose identification and all the
parameters in the model remain consistent.  Under continuous-time,
the variance of the processes is proportional to the sampling interval.
This allows us to trade-off the rate of convergence at which $\widehat{\lambda}_{b}$
approaches $\lambda_{0}$ with the variance of the errors in a neighborhood
of $T_{b}^{0}$ by letting $\sigma_{t}$ become large when $t$ is
close to $T_{b}^{0}$ {[}i.e., a change of time scale as in Foster
and Nelson (\citeyear{nelson/foster:94}, \citeyear{foster/nelson:96}){]}.
This offers a new characterization of higher uncertainty without losing
identification. 

\section{Highest Density Region-based Confidence Sets \label{Section Inference Methods}}

\begin{figure}[H]
\includegraphics[width=18cm,height=20cm]{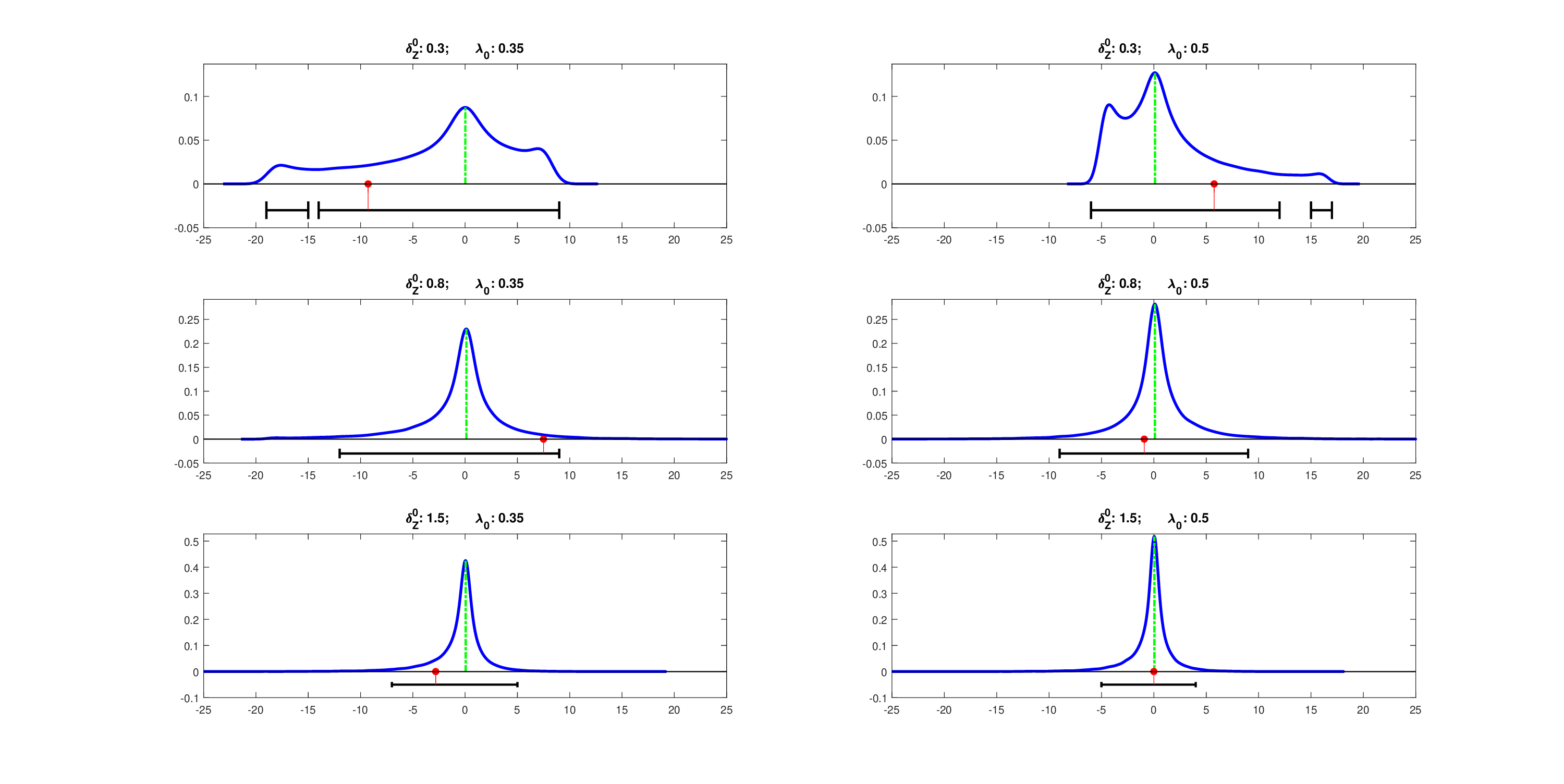}

{\footnotesize{}\caption{{\footnotesize{}\label{Fig30}}{\scriptsize{}Highest Density Regions
(HDRs) of the feasible probability density of $\rho\left(\widehat{T}_{b,\pi}-T_{b}^{0}\right)$
as described in Section \ref{Section Inference Methods}. The significance
level is $\alpha=0.05$, the true break point is $\lambda_{0}=0.3$
and $0.5$ (the left and right panels, respectively) and the break
magnitude is $\delta_{Z}^{0}=0.3,\,0.8$ and 1.5 (the top, middle
and bottom panels, respectively). The horizontal axis is the support
of $\rho\left(\widehat{T}_{b,\pi}-T_{b}^{0}\right)$. The red dot
is the true value of the break point. The union of the black lines
below the horizontal axis is the 95\% HDR confidence region.}}
}{\footnotesize\par}
\end{figure}

The features of the limit and finite-sample distributions suggest
that standard methods to construct confidence intervals may be inappropriate;
e.g., two-sided intervals around the estimated break date based on
the standard deviations of the estimate. Our suggested approach is
rather non-standard and relates to Bayesian methods. In our context,
the Highest Density Region (HDR) seems the most appropriate in light
of the asymmetry and, especially, the multi-modality of the distribution
for small break sizes. All that is needed to implement the procedure
is an estimate of the density function, using plug-in estimates as
explained in Section \ref{Section Approximation-to-the}. Choose
some significance level $0<\alpha<1$ and let $\widehat{P}_{T_{b}}$
denote the empirical counterpart of the probability distribution of
$\rho N(\widehat{\lambda}_{b,\pi}-\lambda_{b}^{0})$ as defined in
Theorem \ref{Theorem 2, Asymptotic Distribution immediate Stationary Regimes}.
Further, let $\widehat{p}_{T_{b}}$ denote the density function defined
by the Radon-Nikodym equation $\widehat{p}_{T_{b}}^{\mathrm{}}=d\widehat{P}_{T_{b}}^{\mathrm{}}/d\lambda_{\mathrm{L}},$
where $\lambda_{\mathrm{L}}$ denotes the Lebesgue measure. 
\begin{defn}
\textbf{Highest Density Region:} Assume that the density function
$f_{Y}\left(y\right)$ of some random variable $Y$ defined on a probability
space $(\Omega_{Y},\,\mathscr{F}_{Y},\,\mathbb{P}_{Y})$ and taking
values on the measurable space $\left(\mathcal{Y},\,\mathscr{Y}\right)$
is continuous and bounded. Then the $\left(1-\alpha\right)100\%$
Highest Density Region is a subset $\mathbf{S}(\kappa_{\alpha})$
of $\mathcal{Y}$ defined as $\mathbf{S}(\kappa_{\alpha})=\{y:\,f_{Y}\left(y\right)>\kappa_{\alpha}\}$
where $\kappa_{\alpha}$ the largest constant that satisfies $\mathbb{P}_{Y}(Y\in\mathbf{S}(\kappa_{\alpha}))\geq1-\alpha$. 
\end{defn}
The concept of HDR and of its estimation has an established literature
in statistics. The definition reported here is from \citet{hyndman:96};
see also \citet{samworth/wand:10} and \citeauthor{mason/polonik:09}
(\citeyear{mason/polonik:08}; \citeyear{mason/polonik:09}) for more
recent developments.\nocite{mason/polonik:08}
\begin{defn}
\textbf{Confidence Sets for $T_{b}^{0}$ under a Continuous Record:}\label{Def. Confidence Sets for Break Date}
Under Assumption \ref{Assumption 1, CT}-\ref{Assumption 3 Break Date},
\ref{Assumption 4 Eigenvalue}-\ref{Assumption 5 Identification}
and \ref{Assumption 6 - Small Shifts}-\ref{Assumtpion - Regimes},
a $\left(1-\alpha\right)100\%$ confidence set for $T_{b}^{0}$ is
a subset of $\left\{ 1,\ldots,\,T\right\} $ given by $C\left(\textrm{cv}_{\alpha}\right)=\left\{ T_{b}\in\left\{ 1,\ldots,\,T\right\} :\,T_{b}\in\mathbf{S}\left(\textrm{cv}_{\alpha}\right)\right\} ,$
where $\mathbf{S}\left(\textrm{cv}_{\alpha}\right)=\left\{ T_{b}:\,\widehat{p}_{T_{b}}>\textrm{cv}_{\alpha}\right\} $
and $\textrm{cv}_{\alpha}$ satisfies $\sup_{\textrm{cv}_{\alpha}\in\mathbb{R}_{+}}\widehat{P}_{T_{b}}\left(T_{b}\in\mathbf{S}\left(\textrm{cv}_{\alpha}\right)\right)\geq1-\alpha$.
\end{defn}
The confidence set $C(\textrm{cv}_{\alpha})$ has a frequentist interpretation
even though the concept of HDR is often encountered in Bayesian analyses
since it associates naturally to the derived posterior distribution,
especially when the latter is multi-modal. A feature of the confidence
set $C(\textrm{cv}_{\alpha})$ under our context is that, at least
when the size of the shift is small, it consists of the union of several
disjoint intervals. The appeal of using HDR is that one can directly
deal with such features. As the break size increases and the distribution
becomes unimodal, the HDR becomes equivalent to the standard way of
constructing confidence sets. In practice, one can proceed as follows.
\begin{lyxalgorithm}
\textbf{$\mathrm{\mathbf{Confidence\,sets\,for\,}}T_{b}^{0}\mathbf{:}$}\label{algorithmn Confidence-sets-for}
1) Estimate by least-squares the break point and the regression coefficients
from model \eqref{Model (4), scalar format, CT}; 2) Replace quantities
appearing in \eqref{Equation (2) Asymptotic Distribution} by consistent
estimators as explained in Section \ref{Section Approximation-to-the};
3) Simulate the limiting distribution $\widehat{P}_{T_{b}}^{\mathrm{}}$
from Theorem \ref{Theorem 2, Asymptotic Distribution immediate Stationary Regimes};
4) Compute the HDR of the empirical distribution $\widehat{P}_{T_{b}}^{\mathrm{}}$
and include the point $T_{b}$ in the level $1-\alpha$ confidence
set $C\left(\textrm{cv}_{\alpha}\right)$ if $T_{b}$ satisfies the
conditions in Definition \ref{Def. Confidence Sets for Break Date}.
\textcolor{white}{\uline{d}}
\end{lyxalgorithm}
This procedure will not deliver contiguous confidence sets when the
size of the break is small. Indeed, we find that in such cases, the
overall confidence set for $T_{b}^{0}$ consists in general of the
union of disjoint intervals if $\widehat{T}_{b}$ is not near the
tails of the sample. One is located around the estimate of the break
date, while the others are in the pre- and post-break regimes. To
provide an illustration, we consider a simple example involving a
single draw from a simulation experiment. Figure \ref{Fig30} reports
the HDR of the feasible limiting distribution of $\rho(\widehat{T}_{b,\pi}-T_{b}^{0})$
for a random draw from the model in \eqref{Casin i-  Model for FS vs Asy Dist}
with parameters $\nu^{0}=1$, $\beta^{0}=0$, unit variance and autoregressive
coefficient 0.6 for $Z_{t}$ and $\sigma_{e}^{2}=1.2$. We set $\lambda_{0}=0.35,\,0.5$
and $\delta_{Z}^{0}=0.3,\,0.8,\,1.5$. We use a trimming 0.15 for
the plug-in estimator $\widehat{T}_{b}$ and $\pi=0.05$ for $\widehat{T}_{b,\pi}$.
As explained in Section \ref{subsec:Evaluation-of-the}, we could
use any other trimming in place of 0.15. The results remain unchanged.
We set $T=100$ and the significance level is $\alpha=0.05$. Note
that the origin is at the estimated break date. The point on the horizontal
axis corresponds to the true break date. The black intervals on the
horizontal axis correspond to regions of high density. The resulting
confidence set is their union. Once a confidence region for $\rho(\widehat{T}_{b,\pi}-T_{b}^{0})$
is computed, it is straightforward to derive a 95\% confidence set
for $T_{b}^{0}.$ The top panel (left plot) reports results for the
case $\delta_{Z}^{0}=0.3$ and $\lambda_{0}=0.35$ and shows that
the HDR is composed of two disjoint intervals. The estimated break
date is $\widehat{T}_{b}=70$ and the implied 95\% confidence set
for $T_{b}^{0}$ is given by $C(\textrm{cv}_{0.05})=\left\{ 1,\ldots,12\right\} \cup\left\{ 18,\ldots100\right\} $.
This includes  $T_{b}^{0}$ and the overall length is 95 observations.
Table \ref{Table 1 HDR} reports for various methods whether $T_{b}^{0}$
is covered or not and the length of the confidence sets for this example.
The length of \citeauthor{bai:97RES}\textquoteright s (1997) confidence
interval is 55 but does not include $T_{b}^{0}$. \citeauthor{elliott/mueller:07}\textquoteright s
(2007) confidence set, denoted by $\widehat{U}_{T}.\textrm{eq}$ in
Table \ref{Table 1 HDR}, also does not include the true break date
at the 90\% confidence level, but does so at the 95\% and its length
is 95. Our method covers $T_{b}^{0}$ and has a relatively short length
across different $\delta_{Z}^{0}$.

\section{Small-Sample Properties of the HDR Confidence Sets\label{Section Small-Sample-Effectiveness-of}}

We now assess via simulations the finite-sample performance of the
method proposed to construct confidence sets for the break date. We
also make comparisons with alternative methods in the literature:
\citeauthor{bai:97RES}'s \citeyearpar{bai:97RES} approach based
on the large-$N$ shrinkage asymptotics; \citeauthor{elliott/mueller:07}\textquoteright s
\citeyearpar{elliott/mueller:07}, hereafter EM, method on\textcolor{red}{{}
}inverting \citeauthor{nyblom:89}'s \citeyearpar{nyblom:89} statistic;\textcolor{red}{{}
}the Inverted Likelihood Ratio (ILR) approach of \citet{eo/morley:15}.
We omit the technical details of these methods and refer to the original
sources or \citet{chang/perron:18} for a review and comparisons.
We consider two DGPs: M1 is $y_{t}=\beta^{0}+\delta_{Z}^{0}\boldsymbol{1}_{\left\{ t>T_{b}^{0}\right\} }+e_{t}$
with $\beta^{0}=1$ and $e_{t}\sim i.i.d.\,\mathscr{N}\left(0,\,1\right)$;
M2 is $y_{t}=\delta_{Z}^{0}\left(1-\nu^{0}\right)\boldsymbol{1}_{\left\{ t>T_{b}^{0}\right\} }+\nu^{0}y_{t-1}+e_{t}$
with $\nu^{0}=0.8$ and $e_{t}\sim i.i.d.\,\mathscr{N}\left(0,\,0.04\right)$.
Our companion paper \citet{casini/perron_CR_Single_Break_Extended}
includes extensive simulation results. We set the significance level
at $\alpha=0.05$, and the break occurs at date $\left\lfloor T\lambda_{0}\right\rfloor $,
where $\lambda_{0}=0.2,\,0.35,\,0.5$ and $T=200$ for M1 and $T=100$
for M2. The results are presented in Table \ref{Table M1}-\ref{Table M7}.
The last row in each table includes the rejection probability of a
5\%-level sup-Wald test using the asymptotic critical value in \citet{andrews:93},
which provides a measure of the magnitude of the break relative to
the noise. For models with predictable processes we use the two-step
procedure described in Section \ref{subsection Asymptotic-Results- Pre}.

Overall, the simulation results confirm previous findings about the
performance of existing methods. \citeauthor{bai:97RES}\textquoteright s
\citeyearpar{bai:97RES} method has a coverage rate below the nominal
level when the size of the break is small. Overall, our HDR method
and that of EM show accurate empirical coverage rates for all DGP
considered. However, EM\textquoteright s method almost always displays
confidence sets which are larger than those from the other approaches.
Over all DGPs considered, the average length of the HDR confidence
sets are 40\% to 70\% shorter than those obtained with EM\textquoteright s
approach when the size of the shift is moderate to high. The results
for M2, a change in mean with a lagged dependent variable and strong
correlation, are quite revealing. EM\textquoteright s method yields
confidence intervals that are very wide, increasing with the size
of the break and for large breaks covering nearly the entire sample.
This does not occur with the other methods. For instance, when $\lambda_{0}=0.5$
and $\delta_{Z}^{0}=2$, the average length from the HDR method is
8.34 compared to 93.71 with EM\textquoteright s. This concurs with
the results in \citet{chang/perron:18}.

In summary, the small-sample simulation results suggest that our continuous
record HDR-based inference provides accurate coverage probabilities
close to the nominal level and average lengths of the confidence sets
shorter relative to existing methods. It is also valid and reliable
under a wider range of DGPs including long-memory processes. Specifically
noteworthy is the fact that it performs well for all break sizes,
whether small or large.

\section{\label{Section Conclusions}Conclusions}

We examined a change-point model under a continuous record asymptotics.
With the time horizon $\left[0,\,N\right]$ fixed, we can account
for the asymmetric informational content provided by the pre- and
post-break samples. We derived a feasible counterpart of the continuous
record asymptotic distribution of the change-point estimator using
consistent plug-in estimates and showed that it provides accurate
approximations to the finite-sample distributions.  We used our limit
theory to construct confidence sets for the change-point date based
on the concept of Highest Density Region. Overall, it delivers accurate
coverage probabilities and relatively short average lengths of the
confidence sets. Importantly, it does so irrespective of the magnitude
of the break, whether large or small, a notoriously difficult problem
in the literature.

\newpage{}

\bibliographystyle{elsarticle-harv}
\bibliography{References_JoE}
\addcontentsline{toc}{section}{References}

\clearpage{}

\newpage{}

\clearpage 
\pagenumbering{arabic}
\renewcommand*{\thepage}{A-\arabic{page}}
\appendix

\clearpage\pagebreak{}

\begin{singlespace} 
\noindent 
\small

\begin{singlespace}
\noindent 
\end{singlespace}

\noindent 
\small

\setcounter{page}{1} \renewcommand{\thepage}{F-\arabic{page}}

\begin{singlespace}
\noindent 

\end{singlespace}

\begin{singlespace}
\noindent 

\end{singlespace}

\begin{singlespace}
\noindent 

\end{singlespace}

\begin{singlespace}
\noindent 

\end{singlespace}

\begin{singlespace}
\noindent 

\end{singlespace}

\pagebreak{}

\end{singlespace}

\begin{singlespace} 

\setcounter{page}{1} \renewcommand{\thepage}{T-\arabic{page}}

{\small{}}
\begin{table}[H]
{\small{}\caption{\label{Table 1 HDR}Coverage rate and length of the confidence set
for the example of Section \ref{Section Inference Methods}}
}{\small\par}
\begin{centering}
{\scriptsize{}}%
\begin{tabular}{ccccccc}
\hline 
 & \multicolumn{2}{c}{{\scriptsize{}$\delta_{Z}^{0}=0.3$}} & \multicolumn{2}{c}{{\scriptsize{}$\delta_{Z}^{0}=0.8$}} & \multicolumn{2}{c}{{\scriptsize{}$\delta_{Z}^{0}=1.5$}}\tabularnewline
 & {\scriptsize{}$\textrm{Cov.}$} & {\scriptsize{}$\textrm{Lgth.}$} & {\scriptsize{}$\textrm{Cov.}$} & {\scriptsize{}$\textrm{Lgth.}$} & {\scriptsize{}$\textrm{Cov.}$} & {\scriptsize{}$\textrm{Lgth.}$}\tabularnewline
\hline 
{\scriptsize{}$\lambda_{0}=0.35$} &  &  &  &  &  & \tabularnewline
{\scriptsize{}HDR} & {\scriptsize{}1} & {\scriptsize{}94} & {\scriptsize{}1} & {\scriptsize{}27} & {\scriptsize{}1} & {\scriptsize{}10}\tabularnewline
{\scriptsize{}Bai (1997)} & {\scriptsize{}0} & {\scriptsize{}55} & {\scriptsize{}0} & {\scriptsize{}13} & {\scriptsize{}1} & {\scriptsize{}8}\tabularnewline
{\scriptsize{}$\widehat{U}_{T}.\textrm{neq}$} & {\scriptsize{}1} & {\scriptsize{}95} & {\scriptsize{}1} & {\scriptsize{}37} & {\scriptsize{}1} & {\scriptsize{}24}\tabularnewline
 &  &  &  &  &  & \tabularnewline
{\scriptsize{}$\lambda_{0}=0.5$} &  &  &  &  &  & \tabularnewline
{\scriptsize{}HDR} & {\scriptsize{}1} & {\scriptsize{}82} & {\scriptsize{}1} & {\scriptsize{}14} & {\scriptsize{}1} & {\scriptsize{}4}\tabularnewline
{\scriptsize{}Bai (1997)} & {\scriptsize{}1} & {\scriptsize{}67} & {\scriptsize{}1} & {\scriptsize{}18} & {\scriptsize{}1} & {\scriptsize{}5}\tabularnewline
{\scriptsize{}$\widehat{U}_{T}.\textrm{neq}$} & {\scriptsize{}1} & {\scriptsize{}95} & {\scriptsize{}1} & {\scriptsize{}35} & {\scriptsize{}1} & {\scriptsize{}14}\tabularnewline
\hline 
\end{tabular}{\scriptsize\par}
\par\end{centering}
{\small{}~~~~~~~~~~~~~~~~~~~~~~~~~~~~}%
\begin{minipage}[t]{0.6\columnwidth}%
{\scriptsize{}Coverage rate and length of the confidence sets corresponding
to the example from Section }\ref{Section Inference Methods}{\scriptsize{}.
See also Figure \ref{Fig30}. The significance level is $\alpha=0.05$.
Cov. and Lgth. refer to the coverage rate and average size of the
confidence sets (i.e. average number of dates in the confidence sets),
respectively. Cov=1 if the confidence set includes $T_{b}^{0}$ and
Cov=0 otherwise. The sample size is $T=100.$}%
\end{minipage}{\small\par}
\end{table}
{\small\par}

{\small{}\pagebreak}{\small\par}

{\small{}}
\begin{table}[H]
{\small{}\caption{\label{Table M1}Small-sample coverage rate and length of the confidence
set for model M1}
}{\small\par}
\begin{centering}
{\scriptsize{}}%
\begin{tabular}{cccccccccc}
\hline 
 &  & \multicolumn{2}{c}{{\scriptsize{}$\delta_{Z}^{0}=0.3$}} & \multicolumn{2}{c}{{\scriptsize{}$\delta_{Z}^{0}=0.6$}} & \multicolumn{2}{c}{{\scriptsize{}$\delta_{Z}^{0}=1$}} & \multicolumn{2}{c}{{\scriptsize{}$\delta_{Z}^{0}=2$}}\tabularnewline
 &  & {\scriptsize{}$\textrm{Cov.}$} & {\scriptsize{}$\textrm{Lgth.}$} & {\scriptsize{}$\textrm{Cov.}$} & {\scriptsize{}$\textrm{Lgth.}$} & {\scriptsize{}$\textrm{Cov.}$} & {\scriptsize{}$\textrm{Lgth.}$} & {\scriptsize{}$\textrm{Cov.}$} & {\scriptsize{}$\textrm{Lgth.}$}\tabularnewline
\hline 
{\scriptsize{}$\lambda_{0}=0.5$} & {\scriptsize{}HDR} & {\scriptsize{}0.938} & {\scriptsize{}131.35} & {\scriptsize{}0.941} & {\scriptsize{}69.05} & {\scriptsize{}0.943} & {\scriptsize{}24.02} & {\scriptsize{}0.962} & {\scriptsize{}6.89}\tabularnewline
 & {\scriptsize{}Bai (1997)} & {\scriptsize{}0.842} & {\scriptsize{}114.24} & {\scriptsize{}0.855} & {\scriptsize{}51.58} & {\scriptsize{}0.911} & {\scriptsize{}19.75} & {\scriptsize{}0.964} & {\scriptsize{}5.70}\tabularnewline
 & {\scriptsize{}$\widehat{U}_{T}.\textrm{eq}$} & {\scriptsize{}0.946} & {\scriptsize{}146.23} & {\scriptsize{}0.943} & {\scriptsize{}76.13} & {\scriptsize{}0.948} & {\scriptsize{}33.45} & {\scriptsize{}0.930} & {\scriptsize{}14.59}\tabularnewline
 & {\scriptsize{}ILR} & {\scriptsize{}0.954} & {\scriptsize{}147.25} & {\scriptsize{}0.956} & {\scriptsize{}78.17} & {\scriptsize{}0.965} & {\scriptsize{}23.87} & {\scriptsize{}0.973} & {\scriptsize{}5.25}\tabularnewline
{\scriptsize{}$\lambda_{0}=0.35$} & {\scriptsize{}HDR} & {\scriptsize{}0.939} & {\scriptsize{}129.02} & {\scriptsize{}0.934} & {\scriptsize{}63.70} & {\scriptsize{}0.939} & {\scriptsize{}24.23} & {\scriptsize{}0.951} & {\scriptsize{}5.78}\tabularnewline
 & {\scriptsize{}Bai (1997)} & {\scriptsize{}0.855} & {\scriptsize{}111.45} & {\scriptsize{}0.855} & {\scriptsize{}49.52} & {\scriptsize{}0.914} & {\scriptsize{}19.39} & {\scriptsize{}0.956} & {\scriptsize{}5.62}\tabularnewline
 & {\scriptsize{}$\widehat{U}_{T}.\textrm{eq}$} & {\scriptsize{}0.933} & {\scriptsize{}148.74} & {\scriptsize{}0.933} & {\scriptsize{}75.94} & {\scriptsize{}0.933} & {\scriptsize{}33.08} & \multicolumn{1}{c}{{\scriptsize{}0.933}} & {\scriptsize{}14.43}\tabularnewline
 & {\scriptsize{}ILR} & {\scriptsize{}0.946} & {\scriptsize{}149.81} & {\scriptsize{}0.960} & {\scriptsize{}77.54} & {\scriptsize{}0.964} & {\scriptsize{}25.63} & {\scriptsize{}0.982} & {\scriptsize{}5.42}\tabularnewline
{\scriptsize{}$\lambda_{0}=0.2$} & {\scriptsize{}HDR} & {\scriptsize{}0.941} & {\scriptsize{}127.29} & {\scriptsize{}0.940} & {\scriptsize{}62.13} & {\scriptsize{}0.942} & {\scriptsize{}22.06} & {\scriptsize{}0.946} & {\scriptsize{}5.73}\tabularnewline
 & {\scriptsize{}Bai (1997)} & {\scriptsize{}0.863} & {\scriptsize{}110.12} & {\scriptsize{}0.911} & {\scriptsize{}53.14} & {\scriptsize{}0.931} & {\scriptsize{}20.20} & {\scriptsize{}0.967} & {\scriptsize{}5.67}\tabularnewline
 & {\scriptsize{}$\widehat{U}_{T}.\textrm{eq}$} & {\scriptsize{}0.950} & {\scriptsize{}158.98} & {\scriptsize{}0.951} & {\scriptsize{}97.12} & {\scriptsize{}0.950} & {\scriptsize{}35.26} & {\scriptsize{}0.950} & {\scriptsize{}13.99}\tabularnewline
 & {\scriptsize{}ILR} & {\scriptsize{}0.956} & {\scriptsize{}162.32} & {\scriptsize{}0.956} & {\scriptsize{}96.45} & {\scriptsize{}0.965} & {\scriptsize{}33.31} & {\scriptsize{}0.976} & {\scriptsize{}5.96}\tabularnewline
\hline 
\end{tabular}{\scriptsize\par}
\par\end{centering}
{\small{}}%
\noindent\begin{minipage}[t]{1\columnwidth}%
{\scriptsize{}The model is $y_{t}=\beta^{0}+\delta_{Z}^{0}\mathbf{1}_{\left\{ t>\left\lfloor T\lambda_{0}\right\rfloor \right\} }+e_{t},\,e_{t}\sim i.i.d.\,\mathscr{N}\left(0,\,1\right),\,T=200$.
Cov. and Lgth. refer to the coverage probability and the average length
of the confidence set (i.e., the average number of dates in the confidence
set). sup-W refers to the rejection probability of the sup-Wald test
using a 5\% size with the asymptotic critical value. The number of
simulations is 5,000.}%
\end{minipage}{\small\par}
\end{table}
{\small\par}

\pagebreak{}

\begin{table}[H]
\caption{\label{Table M7}Small-sample coverage rate and length of the confidence
sets for model M2}

\begin{centering}
{\scriptsize{}}%
\begin{tabular}{cccccccccc}
\hline 
 &  & \multicolumn{2}{c}{{\scriptsize{}$\delta_{Z}^{0}=1$}} & \multicolumn{2}{c}{{\scriptsize{}$\delta_{Z}^{0}=1.5$}} & \multicolumn{2}{c}{{\scriptsize{}$\delta_{Z}^{0}=2$}} & \multicolumn{2}{c}{{\scriptsize{}$\delta_{Z}^{0}=3$}}\tabularnewline
 &  & {\scriptsize{}$\textrm{Cov.}$} & {\scriptsize{}$\textrm{Lgth.}$} & {\scriptsize{}$\textrm{Cov.}$} & {\scriptsize{}$\textrm{Lgth.}$} & {\scriptsize{}$\textrm{Cov.}$} & {\scriptsize{}$\textrm{Lgth.}$} & {\scriptsize{}$\textrm{Cov.}$} & {\scriptsize{}$\textrm{Lgth.}$}\tabularnewline
\hline 
{\scriptsize{}$\lambda_{0}=0.5$} & {\scriptsize{}HDR} & {\scriptsize{}0.916} & {\scriptsize{}30.68} & {\scriptsize{}0.944} & {\scriptsize{}14.77} & {\scriptsize{}0.969} & {\scriptsize{}8.34} & {\scriptsize{}0.995} & {\scriptsize{}4.55}\tabularnewline
 & {\scriptsize{}Bai (1997)} & {\scriptsize{}0.793} & {\scriptsize{}12.87} & {\scriptsize{}0.877} & {\scriptsize{}7.11} & {\scriptsize{}0.929} & {\scriptsize{}4.78} & {\scriptsize{}0.973} & {\scriptsize{}2.957}\tabularnewline
 & {\scriptsize{}$\widehat{U}_{T}.\textrm{eq}$} & {\scriptsize{}0.951} & {\scriptsize{}91.64} & {\scriptsize{}0.955} & {\scriptsize{}93.94} & {\scriptsize{}0.959} & {\scriptsize{}93.71} & {\scriptsize{}0.961} & {\scriptsize{}90.34}\tabularnewline
 & {\scriptsize{}ILR} & {\scriptsize{}0.951} & {\scriptsize{}46.31} & {\scriptsize{}0.967} & {\scriptsize{}34.19} & {\scriptsize{}0.977} & {\scriptsize{}26.48} & {\scriptsize{}0.991} & {\scriptsize{}16.49}\tabularnewline
{\scriptsize{}$\lambda_{0}=0.35$} & {\scriptsize{}HDR} & {\scriptsize{}0.925} & {\scriptsize{}33.02} & {\scriptsize{}0.933} & {\scriptsize{}16.67} & {\scriptsize{}0.971} & {\scriptsize{}9.40} & {\scriptsize{}0.994} & {\scriptsize{}4.33}\tabularnewline
 & {\scriptsize{}Bai (1997)} & {\scriptsize{}0.804} & {\scriptsize{}13.00} & {\scriptsize{}0.876} & {\scriptsize{}7.11} & {\scriptsize{}0.923} & {\scriptsize{}4.94} & {\scriptsize{}0.974} & {\scriptsize{}2.93}\tabularnewline
 & {\scriptsize{}$\widehat{U}_{T}.\textrm{eq}$} & {\scriptsize{}0.952} & {\scriptsize{}91.22} & {\scriptsize{}0.945} & {\scriptsize{}92.61} & \multicolumn{1}{c}{{\scriptsize{}0.957}} & {\scriptsize{}92.48} & {\scriptsize{}0.964} & {\scriptsize{}93.08}\tabularnewline
 & {\scriptsize{}ILR} & {\scriptsize{}0.949} & {\scriptsize{}47.54} & {\scriptsize{}0.967} & {\scriptsize{}34.18} & {\scriptsize{}0.982} & {\scriptsize{}25.84} & {\scriptsize{}0.984} & {\scriptsize{}16.76}\tabularnewline
{\scriptsize{}$\lambda_{0}=0.2$} & {\scriptsize{}HDR} & {\scriptsize{}0.937} & {\scriptsize{}34.66} & {\scriptsize{}0.953} & {\scriptsize{}19.24} & {\scriptsize{}0.954} & {\scriptsize{}11.42} & {\scriptsize{}0.994} & {\scriptsize{}5.36}\tabularnewline
 & {\scriptsize{}Bai (1997)} & {\scriptsize{}0.832} & {\scriptsize{}13.64} & {\scriptsize{}0.885} & {\scriptsize{}7.19} & {\scriptsize{}0.931} & {\scriptsize{}4.92} & {\scriptsize{}0.971} & {\scriptsize{}2.91}\tabularnewline
 & {\scriptsize{}$\widehat{U}_{T}.\textrm{eq}$} & {\scriptsize{}0.944} & {\scriptsize{}89.64} & {\scriptsize{}0.951} & {\scriptsize{}89.58} & {\scriptsize{}0.956} & {\scriptsize{}88.22} & {\scriptsize{}0.961} & {\scriptsize{}85.95}\tabularnewline
 & {\scriptsize{}ILR} & {\scriptsize{}0.946} & {\scriptsize{}49.13} & {\scriptsize{}0.970} & {\scriptsize{}33.54} & {\scriptsize{}0.980} & {\scriptsize{}24.48} & {\scriptsize{}0.989} & {\scriptsize{}12.51}\tabularnewline
\hline 
\end{tabular}{\scriptsize\par}
\par\end{centering}
\noindent\begin{minipage}[t]{1\columnwidth}%
{\scriptsize{}The model is $y_{t}=\delta_{Z}^{0}\left(1-\nu^{0}\right)\mathbf{1}_{\left\{ t>\left\lfloor T\lambda_{0}\right\rfloor \right\} }+\nu^{0}y_{t-1}+e_{t},\,e_{t}\sim i.i.d.\,\mathscr{N}\left(0,\,0.04\right),\,\nu^{0}=0.8,\,T=100$.
The notes of Table \ref{Table M1} apply.}%
\end{minipage}
\end{table}

\newpage{}

\section{Supplemental Materials}

\setcounter{page}{1} \renewcommand{\thepage}{A-\arabic{page}}The
supplement for online publication {[}cf. \citet{casini/perron_CR_Single_Break_Supp}{]}
includes the followings: (i) it describes how to simulate the continuous
record limiting distribution; (ii) it describes the limiting process
in Theorem \ref{Theorem 1}; (iii) it extends the benchmark model
in Section \ref{Section: Model-and-Assumptions} to include predictable
processes; (iv) it includes all proofs of the results in the paper;
(v) it presents additional small-sample evaluations of the HDR confidence
sets.

\newpage{}

\end{singlespace} 


\pagebreak{}

\section*{}
\addcontentsline{toc}{part}{Supplemental Material}
\begin{center}
\Large{Supplemental Material to} 
\end{center}

\begin{center}
\title{\textbf{\Large{Continuous Record Asymptotics  for Change-Point Models}}} 
\maketitle
\end{center}
\medskip{} 
\medskip{} 
\medskip{} 
\thispagestyle{empty}

\begin{center}
$\qquad$ \textsc{\textcolor{MyBlue}{Alessandro Casini}} $\qquad$ \textsc{\textcolor{MyBlue}{Pierre Perron}}\\
\small{University of Rome Tor Vergata} $\quad$ \small{Boston University} 
\\
\medskip{}
\medskip{} 
\medskip{} 
\medskip{} 
\date{\small{\today} \\
\footnotesize{First Vesrion: \printdate{28.10.2015}}}
\medskip{} 
\medskip{} 
\medskip{} 
\end{center}
\begin{abstract}
{\footnotesize{}This supplemental material is structured as follows.
Section \ref{subsec:Simulation-of-the LD in Theorem 4.1} describes
how to simulate the continuous record limiting distribution}.{\footnotesize{}
Section \ref{subsection Description Limiting Process} describes the
limiting process in Theorem \ref{Theorem 1} of the main text.} {\footnotesize{}Section
\ref{SubSection: The-Extended-Model} extends the benchmark model
in Section \ref{Section: Model-and-Assumptions} of the main text
to include predictable processes. Section \ref{Section Mathematical-Proofs Supp}
includes all proofs of the results in the paper. Section \ref{sec:Additional-Simulations-Results}
presents additional small-sample evaluations of the HDR confidence
sets.}{\footnotesize\par}
\end{abstract}
\setcounter{page}{0}
\setcounter{section}{0}
\renewcommand*{\theHsection}{\the\value{section}}

\newpage{}

\begin{singlespace} 
\noindent 
\small

\allowdisplaybreaks


\renewcommand{\thepage}{S-\arabic{page}}   
\renewcommand{\thesection}{S.\Alph{section}}   
\renewcommand{\theequation}{S.\arabic{equation}}



\simpleheading

\section{\label{subsec:Simulation-of-the LD in Theorem 4.1}Simulation of
the Limiting Distribution in Theorem \ref{Theorem 1}}

We discuss how to simulate the limiting distribution in Theorem \ref{Theorem 1}
which is slightly different from simulating the limiting distribution
in Theorem \ref{Theorem 2, Asymptotic Distribution immediate Stationary Regimes}.
However, the idea is similar in that we replace unknown quantities
by consistent estimates. First, we replace $N_{b}^{0}$ by $\widehat{N}_{b}^{\mathrm{}}=\widehat{T}_{b}^{\mathrm{}}/T$.
The ratio $||\delta^{0}||^{2}/\overline{\sigma}^{2}$ is consistently
estimated by $||\widehat{\delta}||^{2}/(T^{-1}\sum_{k=1}^{T}\widehat{e}_{kh}^{2})$
because under the ``fast time scale'' $h^{1/2}\sum_{k=1}^{T}\widehat{e}_{kh}^{2}\overset{p}{\rightarrow}\overline{\sigma}^{2}$
(cf. Assumption \ref{Assumption 6 - Small Shifts}). Now consider
the term $\{-(\delta^{0})'\left\langle Z_{\Delta},\,Z_{\Delta}\right\rangle \left(v\right)\delta^{0}+2(\delta^{0})'\mathscr{W}\left(v\right)\}$.
For $v\leq0$, this can be consistently estimated by 
\begin{align}
-T^{1/2}\left[\left(\widehat{\delta}\right)'\left(\sum{}_{k=\widehat{T}_{b}^{\mathrm{}}+\left\lfloor v/h\right\rfloor }^{\widehat{T}_{b}^{\mathrm{}}}z_{kh}z'_{kh}\right)\widehat{\delta}-2\widehat{\delta}'\widehat{\mathscr{W}}_{h}\left(v\right)\right] & ,\label{eq. RHS Limit Process Th. 4.1}
\end{align}
 where $\widehat{\mathscr{W}_{h}}$ is a simple-size dependent sequence
of Gaussian processes whose marginal distribution is characterized
by $h^{1/2}T\sum_{k=\widehat{T}_{b}+\left\lfloor v/h\right\rfloor }^{\widehat{T}_{b}}\widehat{e}_{kh}^{2}z_{kh}z'_{kh}$
which is a consistent estimate of $\int_{v}^{0}\Omega_{Ze,s}ds$.
Thus, in the limit $\widehat{\mathscr{W}_{h}}\left(v\right)$ has
the same marginal distribution as $\mathscr{W}\left(v\right)$, and
it follows that the limiting distribution from Theorem \ref{Theorem 1}
can be simulated. The proposed method is valid under a continuous-record
asymptotic (i.e., under Assumption \ref{Assumption 6 - Small Shifts}
and the adoption of the ``fast time scale''). It can also be shown
to be valid under a fixed-shifts framework.

\section{Description of the Limiting Process in Theorem \ref{Theorem 1}\label{subsection Description Limiting Process}}

\begin{onehalfspace}
We describe the probability setup underlying the limit process of
Theorem \ref{Theorem 1}. Note that $Z'_{\Delta}e/h^{1/2}=h^{-1/2}\sum_{k=T_{b}+1}^{T_{b}^{0}}z_{kh}e_{kh}$
if $T_{b}\leq T_{b}^{0}.$ Consider an additional measurable space
$(\Omega^{*},\,\mathscr{F}^{*})$ and a transition probability $H\left(\omega,\,d\omega^{*}\right)$
from $\left(\Omega,\,\mathscr{F}\right)$ into $\left(\Omega^{*},\,\mathscr{F}^{*}\right)$.
Next, we can define the products $\widetilde{\Omega}=\Omega\times\Omega^{*}$,
$\widetilde{\mathscr{F}}=\mathscr{F\otimes\mathscr{F}}^{*}$, $\widetilde{P}\left(d\omega,\,d\omega^{*}\right)=P\left(d\omega\right)H\left(\omega,\,d\omega^{*}\right)$.
This defines an extension $(\widetilde{\Omega},\,\widetilde{\mathscr{F}},\,\widetilde{P})$
of the original space $(\Omega,\,\mathscr{F},\,\{\mathscr{F}_{t}\}_{t\geq0},\,P)$.
We also consider another filtration $\{\widetilde{\mathscr{F}}_{t}\}_{t\geq0}$
which takes the following product form $\widetilde{\mathscr{F}}_{t}=\cap_{s>t}\mathscr{F}_{s}\otimes\mathscr{F}{}_{s}^{*}$
where $\{\mathscr{F}{}_{t}^{*}\}_{t\geq0}$ is a filtration on $(\Omega^{*},\,\mathscr{F}^{*})$.
For the transition probability $H$, we consider the simple form $H\left(\omega,\,d\omega^{*}\right)=P^{*}\left(d\omega^{*}\right)$
for some probability measure $P^{*}$ on $(\Omega^{*},\,\mathscr{F}^{*})$.
This constitutes a ``very good'' product filtered extension. Next,
assume that $(\Omega^{*},\,\mathscr{F}^{*},\,(\mathscr{F}{}_{t}^{*})_{t\geq0},\,P^{*})$
supports $p$-dimensional $\{\mathscr{F}{}_{t}^{*}\}$-standard independent
Wiener processes $W^{i*}\left(v\right)$ $\left(i=1,\,2\right)$.
Finally, we postulate the process $\Omega_{Ze,t}$ with entries $\Sigma_{Z}^{\left(i,\,j\right)}\sigma_{e}^{2}$
to admit a progressively measurable $p\times p$ matrix-valued process
(i.e., a symmetric ``square-root'' process) $\sigma_{Ze}$, satisfying
$\Omega_{Ze}=\sigma_{Ze}\sigma'_{Ze}$, with the property that $||\sigma_{Ze}||^{2}\leq K||\Omega_{Ze}||$
for some $K<\infty$. Define the process $\mathscr{W}\left(v\right)=\mathscr{W}_{1}\left(v\right)$
if $v\leq0$, and $\mathscr{W}\left(v\right)=\mathscr{W}_{2}\left(v\right)$
if $v>0$, where $\mathscr{W}_{1}\left(v\right)=\int_{N_{b}^{0}+v}^{N_{b}^{0}}\sigma_{Ze,s}dW_{s}^{1*}$
and $\mathscr{W}_{2}\left(v\right)=\int_{N_{b}^{0}}^{N_{b}^{0}+v}\sigma_{Ze,s}dW_{s}^{2*}$
with components $\mathscr{W}^{\left(j\right)}\left(v\right)=\sum_{r=1}^{p}\int_{N_{b}^{0}+v}^{N_{b}^{0}}\sigma_{Ze,s}^{\left(jr\right)}dW_{s}^{1*\left(r\right)}$
if $v\leq0$ and $\mathscr{W}^{\left(j\right)}\left(v\right)=\sum_{r=1}^{p}\int_{N_{b}^{0}}^{N_{b}^{0}+v}\sigma_{Ze,s}^{\left(jr\right)}dW_{s}^{2*\left(r\right)}$
if $v>0$. The process $\mathscr{W}\left(v\right)$ is well defined
on the product extension $(\widetilde{\Omega},\,\widetilde{\mathscr{F}},\,\left\{ \widetilde{\mathscr{F}}_{t}\right\} _{t\geq0},\,\widetilde{P})$,
and furthermore, conditionally on $\mathscr{F}$, is a two-sided centered
continuous Gaussian process with independent increments and (conditional)
covariance 
\begin{align}
\widetilde{\mathbb{E}}\left(\mathscr{W}^{\left(u\right)}\left(v\right)\mathscr{W}^{\left(j\right)}\left(v\right)\right) & =\Omega_{\mathscr{W}}^{\left(u,j\right)}\left(v\right)=\begin{cases}
\Omega_{\mathscr{W},1}^{\left(u,j\right)}\left(v\right), & \textrm{ if }v\leq0\\
\Omega_{\mathscr{W},2}^{\left(u,j\right)}\left(v\right), & \textrm{ if }v>0
\end{cases},\label{Eq. Cov Matrix of W process, CT}
\end{align}
where $\Omega_{\mathscr{W},1}^{\left(u,j\right)}\left(v\right)=\int_{N_{b}^{0}+v}^{N_{b}^{0}}\Omega_{Ze,s}^{\left(u,j\right)}ds$
and $\Omega_{\mathscr{W},2}^{\left(u,j\right)}\left(v\right)=\int_{N_{b}^{0}}^{N_{b}^{0}+v}\Omega_{Ze,s}^{\left(u,j\right)}ds$.
Therefore, $\mathscr{W}\left(v\right)$ is conditionally on $\mathscr{F}$,
a continuous martingale with \textquotedblleft deterministic\textquotedblright{}
quadratic covariation process $\Omega_{\mathscr{W}}$. The continuity
of $\Omega_{\mathscr{W}}$ signifies that $\mathscr{W}\left(v\right)$
is not only conditionally Gaussian but also a.s. continuous. Precise
treatment of this result can be found in Section II.7 of \citeReferencesSupp{jacod/shiryaev:03}.
\end{onehalfspace}

\section{The Extended Model with Predictable Processes\label{SubSection: The-Extended-Model}}

\subsection{The Extended Model}

The assumptions on $D_{t}$ and $Z_{t}$ specify that they are continuous
semimartingale of the form \eqref{Model Regressors Integral Form}.
This precludes predictable processes, which are often of interest
in applications; e.g., a constant and/or a lagged dependent variable.
Technically, these require a separate treatment since the coefficients
associated with predictable processes are not identified under a fixed-span
asymptotic setting. Let
\begin{align*}
\tau_{1,k}=\mu_{1,h}h+\alpha_{1,h}Y_{\left(k-1\right)h},\qquad & \left(k\leq\left\lfloor T\lambda_{0}\right\rfloor \right),
\end{align*}
 and
\begin{align*}
\tau_{2,k}=\mu_{2,h}h+\alpha_{2,h}Y_{\left(k-1\right)h},\qquad & \left(k>\left\lfloor T\lambda_{0}\right\rfloor +1\right).
\end{align*}
We consider the following extended model:

\begin{flalign}
\Delta_{h}Y_{k} & \triangleq\begin{cases}
\tau_{1,k}+\left(\Delta_{h}D_{k}\right)'\nu^{0}+\left(\Delta_{h}Z_{k}\right)'\delta_{Z,1}^{0}+\Delta_{h}e_{k}^{*},\qquad & \left(k=1,\ldots,T_{b}^{0}\right)\\
\tau_{2,k}+\left(\Delta_{h}D_{k}\right)'\nu^{0}+\left(\Delta_{h}Z_{k}\right)'\delta_{Z,2}^{0}+\Delta_{h}e_{k}^{*},\qquad & \left(k=T_{b}^{0}+1,\ldots,\,T\right)
\end{cases}\label{Model Extended, Sect 1}
\end{flalign}
for some given initial value $Y_{0}$. We specify the parameters associated
with the constant and the lagged dependent variable as being of higher
order in $h$, or lower in $T$, as $h\downarrow0$ so that some fixed
true parameter values can be identified, i.e., $\mu_{1,h}\triangleq\mu_{1}^{0}h^{-1/2}$,
$\mu_{2,h}\triangleq\mu_{2}^{0}h^{-1/2}$, $\mu_{\delta,h}\triangleq\mu_{2,h}-\mu_{1,h}$,
$\alpha_{1,h}\triangleq\alpha_{1}^{0}h^{-1/2}$, $\alpha_{2,h}\triangleq\alpha_{2}^{0}h^{-1/2}$
and $\alpha_{\delta,h}\triangleq\alpha_{2,h}-\alpha_{1,h}$. Our framework
is then similar to the small-diffusion setting studied previously
{[}cf. \citeReferencesSupp{ibragimov/has:80}, \citeReferencesSupp{galtchouk/konev:01},
\citeReferencesSupp{laredo:00} and \citeReferencesSupp{sorensen/uchida:03}{]}.
With $\mu_{\cdot,h}$ and $\alpha_{\cdot,h}$ independent of $h$
and fixed, respectively, at the true values $\mu_{\cdot}^{0}$ and
$\alpha_{\cdot}^{0}$, the continuous-time model is then equivalent
to
\begin{align}
Y_{t} & =Y_{0}+\int_{0}^{t}\left(\mu_{1}^{0}+\mu_{\delta}^{0}\mathbf{1}_{\left\{ s>N_{b}^{0}\right\} }\right)ds+\int_{0}^{t}\left(\alpha_{1}^{0}+\alpha_{\delta}^{0}\mathbf{1}_{\left\{ s>N_{b}^{0}\right\} }\right)Y_{s}ds\label{Full Model}\\
 & \quad+D'_{t}\nu^{0}+\int_{0}^{t}\left(\delta_{Z,1}^{0}+\delta_{Z}^{0}\mathbf{1}_{\left\{ s>N_{b}^{0}\right\} }\right)'dZ_{s}+e_{t}^{*},\nonumber 
\end{align}
 for $t\in\left[0,\,N\right],$ where $Y_{t}=\sum_{k=1}^{\left\lfloor t/h\right\rfloor }\Delta_{h}Y_{k}$,
$D_{t}=\sum_{k=1}^{\left\lfloor t/h\right\rfloor }\Delta_{h}D_{k}$,
$Z_{t}=\sum_{k=1}^{\left\lfloor t/h\right\rfloor }\Delta_{h}Z_{k}$
and $e_{t}^{*}=\sum_{k=1}^{\left\lfloor t/h\right\rfloor }\Delta_{h}e_{k}^{*}$.
The results to be discussed below go through in this extended framework.
However, some additional technical details are needed. Hence, we treat
both cases with and without predictable components separately. Note
that the model and results can be trivially extended to allow for
more general forms of predictable processes (e.g., more lagged values),
at the expense of additional technical details of no substance. 
\begin{singlespace}

\subsection{Asymptotic\label{subsection Asymptotic-Results- Pre} Results for
the Model with Predictable Processes}
\end{singlespace}

\begin{onehalfspace}
In this section, we present asymptotic results allowing for predictable
processes that include a constant and a lagged dependent variable
among the regressors. Recall model \eqref{Model Extended, Sect 1}.
Let $\beta^{0}=(\mu_{1}^{0},\,\alpha_{1}^{0},\,(\nu^{0})',\,(\delta_{Z,1}^{0})')'$,
$\delta_{Z}^{0}=(\mu_{\delta}^{0},\,\alpha_{\delta}^{0},\,(\delta_{Z,2}^{0}-\delta_{Z,1}^{0})')'$,
$((\beta{}^{0})',\,((\delta_{Z}^{0})'))'\in\Theta_{0}$, and $x_{kh}=((\mu_{1,h}/\mu_{1}^{0})h,$
$(\alpha_{1,h}/\alpha_{1}^{0})Y_{\left(k-1\right)h}h,\,\Delta_{h}D'_{k},\,\Delta_{h}Z'_{k}$).
In matrix format, the model is $Y=X\beta^{0}+Z_{0}\delta_{Z}^{0}+e,$
where now $X$ is $T\times\left(p+q+2\right)$ and $Z_{0}=X\overline{R}$,
$\overline{R}\triangleq[(I_{2},\,0_{2\times p})',\,(0'_{\left(p+q\right)\times2},\,R)]'$,
with $R$ as defined in Section \ref{Section: Model-and-Assumptions}.
Natural estimates of $\beta^{0}$ and $\delta_{Z}^{0}$ minimize the
following criterion function,
\begin{align}
h^{-1}\sum_{k=1}^{T} & \left(\Delta_{h}Y{}_{k}-\beta'\int_{\left(k-1\right)h}^{kh}X_{s}ds-\delta'\int_{\left(k-1\right)h}^{kh}Z_{s}ds\right)^{2}\nonumber \\
 & =h^{-1}\sum_{k=1}^{T}\biggl(\Delta_{h}Y{}_{k}-\mu_{1}^{h}h-\alpha_{1}^{h}\int_{\left(k-1\right)h}^{kh}Y_{s}ds-\pi'\Delta_{h}D{}_{k}\label{QLS Full Criterion}\\
 & \quad-\delta_{Z,1}'\Delta_{h}Z{}_{k}\mathbf{1}\left\{ k\leq T_{b}\right\} -\delta_{Z,2}'\Delta_{h}Z{}_{k}\mathbf{1}\left\{ k>T_{b}\right\} \biggl)^{2}.\nonumber 
\end{align}
Hence, we define our LS estimator as the minimizer of the following
approximation to \eqref{QLS Full Criterion}:
\begin{align*}
h^{-1}\sum_{k=1}^{T}\biggl( & \Delta_{h}Y{}_{k}-\mu_{1}^{h}h-\alpha_{1}^{h}Y_{\left(k-1\right)h}h-\nu'\Delta_{h}D_{k}\\
 & \quad-\delta_{Z,1}'\Delta_{h}Z{}_{k}\mathbf{1}\left\{ k\leq T_{b}\right\} -\delta_{Z,2}'\Delta_{h}Z{}_{k}\mathbf{1}\left\{ k>T_{b}\right\} \biggl)^{2}.
\end{align*}
 Such approximations are common {[}cf. \citeReferencesSupp{christopeit:86},
\citeReferencesSupp{lai/wei:83} and \citeReferencesSupp{melnikov/novikov:88}
and the more recent work of \citeReferencesSupp{galtchouk/konev:01}{]}.
Define $\Delta_{h}\widetilde{Y}_{k}\triangleq h^{1/2}\Delta_{h}Y_{k}$
and $\Delta_{h}\widetilde{V}_{k}=h^{1/2}\Delta_{h}V_{k}(\nu^{0},\,\delta_{Z,1}^{0},\,\delta_{Z,2}^{0})$,
where 
\begin{align*}
\Delta_{h}V_{k}\left(\nu^{0},\,\delta_{Z,1}^{0},\,\delta_{Z,2}^{0}\right) & \triangleq\begin{cases}
\left(\nu^{0}\right)'\Delta_{h}D{}_{k}+\left(\delta_{Z,1}^{0}\right)'\Delta_{h}Z{}_{k}+\Delta_{h}e_{k}^{*}, & \textrm{if }k\leq T_{b}^{0}\\
\left(\nu^{0}\right)'\Delta_{h}D{}_{k}+\left(\delta_{Z,2}^{0}\right)'\Delta_{h}Z{}_{k}+\Delta_{h}e_{k}^{*}, & \textrm{if }k>T_{b}^{0}
\end{cases}.
\end{align*}
The small-dispersion format of our model is then
\begin{flalign}
 &  & \Delta_{h}\widetilde{Y}_{k} & =\left(\mu_{1}^{0}h+\alpha_{1}^{0}\widetilde{Y}_{\left(k-1\right)h}h\right)\mathbf{1}\left\{ k\leq T_{b}^{0}\right\} \label{Small-Dispersion Model}\\
 &  &  & \quad+\left(\mu_{2}^{0}h+\alpha_{2}^{0}\widetilde{Y}_{\left(k-1\right)h}h\right)\mathbf{1}\left\{ k>T_{b}^{0}\right\} +\Delta_{h}\widetilde{V}_{k}\left(\nu^{0},\,\delta_{Z,1}^{0},\,\delta_{Z,2}^{0}\right).\nonumber 
\end{flalign}
 This re-parametrization emphasizes that asymptotically our model
describes small disturbances to the approximate dynamical system
\begin{align}
d\widetilde{Y}_{t}^{0}/dt= & \left(\mu_{1}^{0}+\alpha_{1}^{0}\widetilde{Y}_{t}^{0}\right)\mathbf{1}\left\{ t\leq N_{b}^{0}\right\} +\left(\mu_{2}^{0}+\alpha_{2}^{0}\widetilde{Y}_{t}^{0}\right)\mathbf{1}\left\{ t>N_{b}^{0}\right\} .\label{Eq. Dynamical System}
\end{align}
The process $\left\{ \widetilde{Y}_{t}^{0}\right\} _{t\geq0}$ is
the solution to the underlying ordinary differential equation. The
LS estimate of the break point is then defined as $\widehat{T}_{b}\triangleq\arg\max_{T_{b}}Q_{T}\left(T_{b}\right),$
where
\begin{align*}
Q_{T}\left(T_{b}\right) & \triangleq Q_{T}\left(\widehat{\beta}\left(T_{b}\right),\,\widehat{\delta}\left(T_{b}\right),\,T_{b}\right)=\widehat{\delta}'\left(Z_{2}'MZ_{2}\right)\widehat{\delta},
\end{align*}
and the LS estimates of the regression parameters are
\begin{align*}
\widehat{\theta}^{\mathrm{}} & \triangleq\arg\min_{\theta\in\Theta_{0}}h\left(S_{T}\left(\beta,\,\delta,\,\widehat{T}_{b}\right)-S_{T}\left(\beta^{0},\,\delta_{Z}^{0},\,T_{b}^{0}\right)\right),
\end{align*}
where $S_{T}$ is the sum of square residuals. With the exception
of our small-dispersion assumption and consequent more lengthy derivations,
our analysis remains the same as in the model without predictable
processes. Hence, the asymptotic distribution of the break point estimator
is derived under the same setting as in Section \ref{Section Asymptotic Distribution: Continuous Case}.
We show that the limiting distribution is qualitatively equivalent
to that in Theorem \ref{Theorem 1}. 
\end{onehalfspace}
\begin{assumption}
\label{Assumption 1, Pre}Assumption \ref{Assumption 3 Break Date}
and \ref{Assumption 5 Identification} hold. Assumption \ref{Assumption 1, CT},
\ref{Assumption 2} and \ref{Assumption 4 Eigenvalue} now apply to
the last $p$ (resp. $q$) elements of the process $\left\{ Z_{t}\right\} _{t\geq0}$
(resp. $\left\{ D_{t}\right\} _{t\geq0}$). 
\end{assumption}
\begin{prop}
\label{Prop 1-2 Asym - Pre}Consider model \eqref{Model Extended, Sect 1}.
Under Assumption \ref{Assumption 1, Pre}: (i) $\widehat{\lambda}_{b}\overset{P}{\rightarrow}\lambda_{0}$;
(ii) for every $\varepsilon>0$ there exists a $K>0$ such that for
all large $T,$ $P(T|\widehat{\lambda}_{b}-\lambda_{0}|>K||\delta^{0}||^{-2}\overline{\sigma}^{2})<\varepsilon$.
\end{prop}
\begin{assumption}
\label{Assumption 5 - Small Shifts - Pre}Let $\delta_{h}=h^{1/4}\delta^{0}$
and for $i=1,\,2$ $\mu_{i}^{h}=h^{1/4}\mu_{i}^{0}$ and $\alpha_{i}^{h}=h^{1/4}\alpha_{i}^{0}$,
and assume that for all $t\in(N_{b}^{0}-\epsilon,\,N_{b}^{0}+\epsilon),$
with $\epsilon\downarrow0$ and $T^{1-\kappa}\epsilon\rightarrow B<\infty$,
$0<\kappa<1/2$, $\mathbb{E}[\left(\Delta_{h}e_{t}^{*}\right)^{2}|\,\mathscr{F}_{t-h}]=\sigma_{h,t}^{2}\Delta t$
$P$-a.s, where $\sigma_{h,t}\triangleq\sigma_{h}\sigma_{e,t}$ with
$\sigma_{h}\triangleq h^{-1/4}\overline{\sigma}$. 
\end{assumption}
Furthermore, define the normalized residual $\Delta_{h}\widetilde{e}_{t}$
as in \eqref{Eq. eps WN}. We shall derive a stable convergence in
distribution for $\overline{Q}_{T}\left(\cdot,\,\cdot\right)$ as
defined in Section \ref{Section Asymptotic Distribution: Continuous Case}.
The description of the limiting process is similar to the one presented
in the previous section. However, here we shall condition on the $\sigma$-field
$\mathscr{G}$ generated by all latent processes appearing in the
model. In view of its properties, the $\sigma$-field $\mathscr{F}$
admits a regular version of the $\mathscr{G}$-conditional probability,
denoted $H\left(\omega,\,d\omega^{*}\right)$. The limit process is
then realized on the extension $(\widetilde{\Omega},\,\widetilde{\mathscr{F}},\,\{\widetilde{\mathscr{F}}_{t}\}_{t\geq0},\,\widetilde{P})$
of the original filtered probability space as explained in Section
\ref{subsection Description Limiting Process}. We again introduce
a two-sided Gaussian process $\mathscr{W}_{Ze}\left(\cdot\right)$
with a different dimension in order to accommodate for the presence
of the predictable regressors in the first two columns of both $X$
and $Z$. That is, $\mathscr{W}_{Ze}\left(\cdot\right)$ is a $p$-dimensional
process which is $\mathscr{G}$-conditionally Gaussian and has $P$-a.s.
continuous sample paths. We then have the following theorem.
\begin{thm}
\label{Theorem - Asy Dist - Pre}Consider model \eqref{Small-Dispersion Model}.
Under Assumption \ref{Assumption 1, Pre}-\ref{Assumption 5 - Small Shifts - Pre}:
(i) $\widehat{\lambda}_{b}\overset{P}{\rightarrow}\lambda_{0}$; (ii)
for every $\varepsilon>0$ there exists a $K>0$ such that for all
large $T,$ $P(T^{1-\kappa}|\widehat{\lambda}_{b}-\lambda_{0}|>K||\delta^{0}||^{-2}\overline{\sigma}^{2})<\varepsilon$;
(iii) 
\begin{align}
N\left(\widehat{\lambda}_{b,\pi}-\lambda_{0}\right) & \overset{\mathcal{L}-\mathrm{s}}{\Rightarrow}\underset{v\in\left[\frac{N\pi-N_{b}^{0}}{\left\Vert \delta^{0}\right\Vert ^{-2}\overline{\sigma}^{2}},\,\frac{N\left(1-\pi\right)-N_{b}^{0}}{\left\Vert \delta^{0}\right\Vert ^{-2}\overline{\sigma}^{2}}\right]}{\mathrm{argmax}}\left\{ -\left(\delta^{0}\right)'\Lambda\left(v\right)\delta^{0}+2\left(\delta^{0}\right)'\mathscr{W}\left(v\right)\right\} ,\label{CR Asymptotic Distribution, Pre}
\end{align}
where $\Lambda\left(v\right)$ is a process given by 
\begin{align*}
\Lambda\left(v\right) & \triangleq\begin{cases}
\Lambda_{1}\left(v\right), & \mathrm{if}\,v\leq0\\
\Lambda_{2}\left(v\right), & \mathrm{if}\,v>0
\end{cases},\qquad\mathrm{with}
\end{align*}
\[
\quad\Lambda_{1}\left(v\right)\triangleq\begin{bmatrix}\int_{N_{b}^{0}+v}^{N_{b}^{0}}ds & \int_{N_{b}^{0}+v}^{N_{b}^{0}}\widetilde{Y}_{s}ds & 0_{1\times p}\\
\int_{N_{b}^{0}+v}^{N_{b}^{0}}\widetilde{Y}_{s}ds & \int_{N_{b}^{0}+v}^{N_{b}^{0}}\widetilde{Y}_{s}^{2}ds & 0_{1\times p}\\
0_{p\times1} & 0_{p\times1} & \left\langle Z,\,Z\right\rangle _{1}\left(v\right)
\end{bmatrix},
\]
and $\Lambda_{2}\left(v\right)$ is defined analogously, where $\left\langle Z,\,Z\right\rangle _{1}\left(v\right)$
is the $p\times p$ predictable quadratic covariation process of the
pair $(Z_{\Delta}^{\left(u\right)},\,Z_{\Delta}^{\left(j\right)}),\,3\leq u,j\leq p$
and $v\leq0$. The process $\mathscr{W}\left(v\right)$ is, conditionally
on $\mathscr{G},$ a two-sided centered Gaussian martingale with independent
increments. 
\end{thm}
When $v\leq0$, the limit process $\mathscr{W}\left(v\right)$ is
defined as follows,  
\begin{align*}
\mathscr{W}^{\left(j\right)}\left(v\right) & =\begin{cases}
\int_{N_{b}^{0}+v}^{N_{b}^{0}}dW_{e,s}, & j=1,\\
\int_{N_{b}^{0}+v}^{N_{b}^{0}}\widetilde{Y}_{s}dW_{e,s}, & j=2,\\
\mathscr{W}_{Ze}^{\left(j-2\right)}\left(v\right), & j=3,\ldots,\,p+2,
\end{cases}
\end{align*}
where $\mathscr{W}_{Ze}^{\left(i\right)}\left(v\right)\triangleq\sum_{r=1}^{p}\int_{N_{b}^{0}+v}^{N_{b}^{0}}\sigma_{Ze,s}^{\left(i,r\right)}dW_{s}^{1*\left(r\right)}$
($i=1,\ldots,\,p$) and analogously when $v>0$. That is, $\mathscr{W}_{Ze}\left(v\right)$
corresponds to the process $\mathscr{W}\left(v\right)$ used for the
benchmark model (and so are $W_{s}^{1*}$, $W_{s}^{2*}$ and $\Omega_{Ze,s}$
below). Its conditional covariance is given by 
\begin{align}
\widetilde{\mathbb{E}}\left(\mathscr{W}^{\left(u\right)}\left(v\right)\mathscr{W}^{\left(j\right)}\left(v\right)\right) & =\Omega_{\mathscr{W}}^{\left(u,j\right)}\left(v\right)=\begin{cases}
\Omega_{\mathscr{W},1}^{\left(u,j\right)}\left(v\right), & \textrm{ if }v\leq0\\
\Omega_{\mathscr{W},2}^{\left(u,j\right)}\left(v\right), & \textrm{ if }v>0
\end{cases},\label{Eq. Covariance Matrix of Limit process, Pre}
\end{align}
where $\Omega_{\mathscr{W},1}^{\left(u,j\right)}\left(v\right)=\int_{N_{b}^{0}+v}^{N_{b}^{0}}\sigma_{e,s}^{2}ds$,
if $u,j=1$; $\Omega_{\mathscr{W},1}^{\left(u,j\right)}\left(v\right)=\int_{N_{b}^{0}+v}^{N_{b}^{0}}\widetilde{Y}_{s}^{2}\sigma_{e,s}^{2}ds$,
if $u,j=2$; $\Omega_{\mathscr{W},1}^{\left(u,j\right)}\left(v\right)=\int_{N_{b}^{0}+v}^{N_{b}^{0}}\widetilde{Y}_{s}^{2}\sigma_{e,s}^{2}ds$,
if $1\leq u,j\leq2,\,u\neq j$; $\Omega_{\mathscr{W},1}^{\left(u,j\right)}\left(v\right)=0$,
if $u=1,\,2,\,j=3,\ldots,\,p$; $\Omega_{\mathscr{W},1}^{\left(u,j\right)}\left(v\right)=\int_{N_{b}^{0}+v}^{N_{b}^{0}}\Omega_{Ze,s}^{\left(u-2,j-2\right)}$
$ds$ if $3\leq u,j\leq p+2$; and similarly for $\Omega_{\mathscr{W},2}^{\left(u,j\right)}\left(v\right)$.
The asymptotic distribution is qualitatively the same as in Theorem
\ref{Theorem 1}. When the volatility processes are deterministic,
we have convergence in law under the Skorhokod topology to the same
limit process $\mathscr{W}\left(\cdot\right)$ with a Gaussian unconditional
law. The case with stationary regimes is described as follows.
\begin{assumption}
\label{Assumtpion - Regimes - Pre}$\Sigma^{*}=\{\mu_{\cdot,t},\,\Sigma_{\cdot,t},\,\sigma_{e,t}\}_{t\geq0}$
is deterministic and the regimes are stationary. 
\end{assumption}
Let $W_{i}^{*},$ $i=1,\,2,$ be two independent standard Wiener processes
defined on $[0,\,\infty),$ starting at the origin when $s=0.$ Let
\begin{align*}
\mathscr{V}\left(s\right) & =\begin{cases}
-\frac{\left|s\right|}{2}+W_{1}^{*}\left(s\right), & \textrm{if }s<0\\
-\frac{\left(\delta^{0}\right)'\Lambda_{2}\delta^{0}}{\left(\delta^{0}\right)'\Lambda_{1}\delta^{0}}\frac{\left|s\right|}{2}+\left(\frac{\left(\delta^{0}\right)'\Omega_{\mathscr{W},2}\delta^{0}}{\left(\delta^{0}\right)'\Omega_{\mathscr{W},1}\delta^{0}}\right)^{1/2}W_{2}^{*}\left(s\right), & \textrm{if }s\geq0.
\end{cases}
\end{align*}

\begin{cor}
\label{Corollary - Asymptotic Distribution immediate Stationary Regimes, Pre}Under
Assumption \ref{Assumption 1, Pre}-\ref{Assumtpion - Regimes - Pre},
\begin{align}
\frac{\left(\left(\delta^{0}\right)'\Lambda_{1}\delta^{0}\right)^{2}}{\left(\delta^{0}\right)'\Omega_{\mathscr{W},1}\delta^{0}} & N\left(\widehat{\lambda}_{b,\pi}-\lambda_{0}\right)\Rightarrow\underset{s\in\mathcal{A}_{1}}{\mathrm{argmax}}\mathscr{V}\left(s\right),\label{Equation (2) Asymptotic Distribution, Casini (2015a)-1}
\end{align}
where 
\begin{align*}
\mathcal{A}_{1} & =\left[\frac{N\pi-N_{b}^{0}}{\left\Vert \delta^{0}\right\Vert ^{-2}\overline{\sigma}^{2}}\frac{\left(\left(\delta^{0}\right)'\Lambda_{1}\delta^{0}\right)^{2}}{\left(\delta^{0}\right)'\Omega_{\mathscr{W},1}\delta^{0}},\,\frac{N\left(1-\pi\right)-N_{b}^{0}}{\left\Vert \delta^{0}\right\Vert ^{-2}\overline{\sigma}^{2}}\frac{\left(\left(\delta^{0}\right)'\Lambda_{1}\delta^{0}\right)^{2}}{\left(\delta^{0}\right)'\Omega_{\mathscr{W},1}\delta^{0}}\right].
\end{align*}
\end{cor}

In the next two corollaries, we assume stationary errors across regimes.
Corollary \ref{Corollary Mean} considers the basic case of a change
in the mean of a sequence of $\mathit{i.i.d.}$ random variables.
Let 
\begin{align*}
\mathscr{V}_{\mathrm{sta}}\left(s\right) & =\begin{cases}
-\frac{\left|s\right|}{2}+W_{1}^{*}\left(s\right), & \textrm{if }s<0\\
-\frac{\left(\delta^{0}\right)'\Lambda_{2}\delta^{0}}{\left(\delta^{0}\right)'\Lambda_{1}\delta^{0}}\frac{\left|s\right|}{2}+\left(\frac{\left(\delta^{0}\right)'\Lambda_{2}\delta^{0}}{\left(\delta^{0}\right)'\Lambda_{1}\delta^{0}}\right)^{1/2}W_{2}^{*}\left(s\right), & \textrm{if }s\geq0
\end{cases},
\end{align*}
\[
\quad\mathscr{V}_{\mu,\mathrm{sta}}\left(s\right)=\begin{cases}
-\frac{\left|s\right|}{2}+W_{1}^{*}\left(s\right), & \textrm{if }s<0\\
-\frac{\left|s\right|}{2}+W_{2}^{*}\left(s\right), & \textrm{if }s\geq0
\end{cases}.
\]

\begin{cor}
Under Assumption \ref{Assumption 1, Pre}-\ref{Assumtpion - Regimes - Pre}
and assuming that the second moments of the residual process are stationary
across regimes, $\sigma_{e,s}=\overline{\sigma}$ for all $0\leq s\leq N$,
\begin{align*}
\frac{\left(\delta^{0}\right)'\Lambda_{1}\delta^{0}}{\overline{\sigma}^{2}}N\left(\widehat{\lambda}_{b,\pi}-\lambda_{0}\right) & \Rightarrow\underset{s\in\mathcal{A}_{2}}{\mathrm{argmax}}\mathscr{V}_{\mathrm{sta}}\left(s\right),
\end{align*}
 where 
\begin{align*}
\mathcal{A}_{2} & =\left[\frac{N\pi-N_{b}^{0}}{\left\Vert \delta^{0}\right\Vert ^{-2}\overline{\sigma}^{2}}\frac{\left(\delta^{0}\right)'\Lambda_{1}\delta^{0}}{\overline{\sigma}^{2}},\,\frac{N\left(1-\pi\right)-N_{b}^{0}}{\left\Vert \delta^{0}\right\Vert ^{-2}\overline{\sigma}^{2}}\frac{\left(\delta^{0}\right)'\Lambda_{1}\delta^{0}}{\overline{\sigma}^{2}}\right].
\end{align*}
\end{cor}

\begin{cor}
\label{Corollary Mean}Under Assumption \ref{Assumption 1, Pre}-\ref{Assumtpion - Regimes - Pre},
with $\nu^{0}=0,$ $\delta_{Z,i}^{0}=0$, and $\alpha_{i}^{0}=0$
for $i=1,\,2$: 
\begin{align*}
\left(\delta^{0}/\overline{\sigma}\right)^{2}N\left(\widehat{\lambda}_{b,\pi}-\lambda_{0}\right) & \Rightarrow\underset{s\in\left[\left(N\pi-N_{b}^{0}\right)\left(\delta^{0}/\overline{\sigma}\right)^{2},\,\left(N\left(1-\pi\right)-N_{b}^{0}\right)\left(\delta^{0}/\overline{\sigma}\right)^{2}\right]}{\mathrm{argmax}}\mathscr{V}_{\mu,\mathrm{sta}}\left(s\right).
\end{align*}
\end{cor}
\begin{rem}
The last corollary reports the result for the simple case of a shift
in the mean of an $i.i.d.$ process. This case was recently considered
by \citeReferencesSupp{jiang/wang/yu:16} under a continuous-time
setting in their Theorem 4.2-(b) which is similar to our Corollary
\ref{Corollary Mean}. Our limit theory differs in many respects,
besides being obviously more general. \citeReferencesSupp{jiang/wang/yu:16}
only develop an infeasible distribution theory for the break date
estimator whereas we also derive a feasible version. This is because
we introduce an assumption about the drift in order to ``keep''
it in the asymptotics. The limiting distribution is also derived under
a different asymptotic experiment (cf. Assumption \ref{Assumption 5 - Small Shifts - Pre}
above and the change of time scale as discussed in Section \ref{Section Asymptotic Distribution: Continuous Case}).
A direct consequence is that the estimate of the break fraction is
shown to be consistent as $h\downarrow0$ whereas \citeReferencesSupp{jiang/wang/yu:16}
do not have such a result.
\end{rem}
The results are similar to those in the benchmark model. However,
the estimation of the regression parameters is more complicated because
of the identification issues about the parameters associated with
predictable processes. Nonetheless, our model specification allows
us to construct feasible estimators. Given the small-dispersion specification
in \eqref{Small-Dispersion Model}, we propose a two-step estimator.
In fact, \eqref{Eq. Dynamical System} essentially implies that asymptotically
the evolution of the dependent variable is governed by a deterministic
drift function given by $\mu_{1}^{0}+\alpha_{1}^{0}\widetilde{Y}_{t}^{0}$
(resp., $\mu_{2}^{0}+\alpha_{2}^{0}\widetilde{Y}_{t}^{0}$) if $t\leq N_{b}^{0}$
(resp., $t>N_{b}^{0}$). Thus, in a first step we construct least-squares
estimates of $\mu_{i}^{0}$ and $\alpha_{i}^{0}$ $\left(i=1,\,2\right)$.
Next, we subtract the estimate of the deterministic drift from the
dependent variable so as to generate a residual component that will
be used (after rescaling) as a new dependent variable in the second
step where we construct the least-squares estimates of the parameters
associated with the stochastic semimartigale regressors.
\begin{prop}
\label{Proposition Consistency QLS, Pre}Under Assumption \ref{Assumption 1, Pre}-\ref{Assumption 5 - Small Shifts - Pre},
as $h\downarrow0$, $\widehat{\theta}^{\mathrm{}}\overset{P}{\rightarrow}\theta^{0}$.
\end{prop}
The consistency of the estimate $\widehat{\theta}^{\mathrm{}}$ is
all that is needed to carry out our inference procedures about the
break point $T_{b}^{0}$ presented in Section \ref{Section Inference Methods}.
The relevance of the result is that even though the drifts cannot
in general be consistently estimated, we can, under our setting, estimate
the parameters entering the limiting distribution; i.e., $\mu_{i}^{0}$
and $\alpha_{i}^{0}$. 

\section{\label{Section Mathematical-Proofs Supp}Mathematical Proofs}

\subsection{Additional Notations}

For a matrix $A$, the orthogonal projection matrices $P_{A},\,M_{A}$
are defined as $P_{A}=A\left(A'A\right)^{-1}A'$ and $M_{A}=I-P_{A}$,
respectively. For a matrix $A$, we use the vector-induced norm, i.e.,
$\left\Vert A\right\Vert =\sup_{x\neq0}\left\Vert Ax\right\Vert /\left\Vert x\right\Vert .$
Also, for a projection matrix $P$, $\left\Vert PA\right\Vert \leq\left\Vert A\right\Vert .$
We denote the $d$-dimensional identity matrix by $I_{d}.$ When the
context is clear we omit the subscript notation in the projection
matrices. We denote the $\left(i,\,j\right)$-th element of the outer
product matrix $A'A$ as $\left(A'A\right)_{i,j}$ and the $i\times j$
upper-left (resp., lower-right) sub-block of $A'A$ as $\left[A'A\right]_{\left\{ i\times j,\cdot\right\} }$
(resp., $\left[A'A\right]_{\left\{ \cdot,i\times j\right\} }$). For
a random variable $\xi$ and a number $r\geq1,$ we write $\left\Vert \xi\right\Vert _{r}=(\mathbb{E}\left\Vert \xi\right\Vert ^{r})^{1/r}.$
$B$ and $C$ are generic constants that may vary from line to line;
we may sometime write $C_{r}$ to emphasize the dependence of $C$
on a number $r.$ For two scalars $a$ and $b$ the symbol $a\wedge b$
means the infimum of $\left\{ a,\,b\right\} $. The symbol ``$\overset{\textrm{u.c.p.}}{\Rightarrow}$''
signifies uniform locally in time convergence under the Skorokhod
topology and recall that it implies convergence in probability. The
symbol ``$\overset{d}{\equiv}$'' signifies equivalence in distribution.
We also use the same notations as detailed in Section \ref{Section: Model-and-Assumptions}. 

\subsection{Preliminary Lemmas \label{subsection: Preliminary Lemmas}}

Lemma \ref{Lemma A1, Casini (2015a)} is Lemma A.1 in \citeReferencesSupp{bai:97RES}.
Let $X_{\Delta}$ be defined as in the display equation after \eqref{eq. A.2.4, Casini (2015a)-2}.
\begin{lem}
\label{Lemma A1, Casini (2015a)} The following inequalities hold
$P$-a.s.:
\begin{align}
\left(Z_{0}'MZ_{0}\right) & -\left(Z_{0}'MZ_{2}\right)\left(Z'_{2}MZ_{2}\right)^{-1}\left(Z'_{2}MZ_{0}\right)\label{Equation (36)}\\
 & \geq R'\left(X'_{\Delta}X_{\Delta}\right)\left(X'_{2}X_{2}\right)^{-1}\left(X'_{0}X_{0}\right)R,\quad T_{b}<T_{b}^{0}\nonumber \\
\left(Z_{0}'MZ_{0}\right) & -\left(Z_{0}'MZ_{2}\right)\left(Z'_{2}MZ_{2}\right)^{-1}\left(Z'_{2}MZ_{0}\right)\label{Eq. 37}\\
 & \geq R'\left(X'_{\Delta}X_{\Delta}\right)\left(X'X-X'_{2}X_{2}\right)^{-1}\left(X'X-X'_{0}X_{0}\right)R,\quad T_{b}\geq T_{b}^{0}.\nonumber 
\end{align}
\end{lem}
The following lemma presents the uniform approximation to the instantaneous
covariation between continuous semimartingales. This will be useful
in the proof of the convergence rate of our estimator. Below, the
time window in which we study certain estimates is shrinking at a
rate no faster than $h^{1-\epsilon}$ for some $0<\epsilon<1/2.$
\begin{lem}
\label{Lemma, Spot Uniform Approx}Let $X_{t}$ (resp., $\widetilde{X}_{t}$)
be a $q$ (resp., $p$)-dimensional It\^o continuous semimartingale
defined on $\left[0,\,N\right]$. Let $\Sigma_{t}$ denote the time
$t$ instantaneous covariation between $X_{t}$ and $\widetilde{X}_{t}$.
Choose a fixed number $\epsilon>0$ and $\varpi$ satisfying $1/2-\epsilon\geq\varpi\geq\epsilon>0$.
Further, let $B_{T}\triangleq\left\lfloor N/h-T^{\varpi}\right\rfloor .$
Define the moving average of $\Sigma_{t}$ as $\overline{\Sigma}_{kh}\triangleq\left(T^{\varpi}h\right)^{-1}\int_{kh}^{kh+T^{\varpi}h}\Sigma_{s}ds,$
and let $\widehat{\Sigma}_{kh}\triangleq\left(T^{\varpi}h\right)^{-1}\sum_{i=1}^{\left\lfloor T^{\varpi}\right\rfloor }\Delta_{h}X_{k+i}\Delta_{h}\widetilde{X}'_{k+i}.$Then,
$\sup_{1\leq k\leq B_{T}}||\widehat{\Sigma}_{kh}-\overline{\Sigma}_{kh}||=o_{p}\left(1\right).$
Furthermore, for each $k$ and some $K>0$ with $N-K>kh>K$, $\sup_{T^{\epsilon}\leq T^{\varpi}\leq T^{1-\epsilon}}$
$||\widehat{\Sigma}_{kh}-\overline{\Sigma}_{kh}||=o_{p}\left(1\right)$. 
\end{lem}

\noindent \textit{Proof.} By a polarization argument, we can assume
that $X_{t}$ and $\widetilde{X}_{t}$ are univariate without loss
of generality, and by standard localization arguments, we can assume
that the drift and diffusion coefficients of $X_{t}$ and $\widetilde{X}_{t}$
are bounded. Then, by It\^o Lemma, 
\begin{align*}
\widehat{\Sigma}_{kh}-\overline{\Sigma}_{kh} & \triangleq\frac{1}{T^{\varpi}h}\sum_{i=1}^{\left\lfloor T^{\varpi}\right\rfloor }\int_{\left(k+i-1\right)h}^{\left(k+i\right)h}\left(X_{s}-X_{\left(k+i-1\right)h}\right)d\widetilde{X}{}_{s}\\
 & \quad+\frac{1}{T^{\varpi}h}\sum_{i=1}^{\left\lfloor T^{\varpi}\right\rfloor }\int_{\left(k+i-1\right)h}^{\left(k+i\right)h}\left(\widetilde{X}_{s}-\widetilde{X}_{\left(k+i-1\right)h}\right)dX{}_{s}.
\end{align*}
For any $l\geq1,$ $||\widehat{\Sigma}_{kh}-\overline{\Sigma}_{kh}||_{l}\leq K_{l}T^{-\varpi/2},$
which follows from standard estimates for continuous It\^o semimartignales.
By a maximal inequality, 
\[
\left\Vert \sup_{1\leq k\leq B_{T}}\left|\widehat{\Sigma}_{kh}-\overline{\Sigma}_{kh}\right|\right\Vert _{l}\leq K_{l}T^{1/l}T^{-\varpi/2},
\]
 which goes to zero choosing $l>2/\varpi$. This proves the first
claim. For the second, note that for $l\geq1,$ 
\begin{align*}
\left\Vert \sup_{T^{\epsilon}\leq T^{\varpi}\leq T^{1-\epsilon}}\left|\widehat{\Sigma}_{kh}-\overline{\Sigma}_{kh}\right|\right\Vert _{l} & =\left\Vert \sup_{1\leq T^{\varpi-\epsilon}\leq T^{1-2\epsilon}}\left|\widehat{\Sigma}_{kh}-\overline{\Sigma}_{kh}\right|\right\Vert _{l}\\
 & \leq K_{l}T^{\left(1-2\epsilon\right)/l}T^{-\epsilon/2},
\end{align*}
and choose $l>\left(2-4\epsilon\right)/\epsilon$ to verify the claim.
$\square$

\subsection{Preliminary Results\label{Subsection: Preliminaries Results for Asymptotic Theory}}

As it is customary in related contexts, we use a standard localization
argument as explained in Section 1.d in \citeReferencesSupp{jacod/shiryaev:03}
and thus we can replace Assumption \ref{Assumption 1, CT}-\ref{Assumption 2}
with the following stronger assumption. 
\begin{assumption}
\label{Assumption Localization}Let Assumption \ref{Assumption 1, CT}-\ref{Assumption 2}
hold. The process $\left\{ Y_{t},\,D_{t},\,Z_{t}\right\} _{t\geq0}$
takes value in some compact set, $\left\{ \sigma_{\cdot,t}\right\} _{t\geq0}$
is bounded càdlàg and the process $\left\{ \mu_{\cdot,t}\right\} $
is bounded càdlàg or càglàd.
\end{assumption}
The localization technique basically translates all the local conditions
into global ones. We next introduce concepts and results which will
be useful in some of the proofs below. 

\subsubsection{Approximate Variation, LLNs and CLTs}

We review some basic definitions about approximate covariation and
more general high-frequency statistics. Given a continuous-time semimartingales
$X=\left(X^{i}\right)_{1\leq i\leq d}\in\mathbb{R}^{d}$ with zero
initial value over the time horizon $\left[0,\,N\right],$ with $P$-a.s.
continuous paths, the covariation of $X$ over $\left[0,\,t\right]$
is denoted $\left[X,\,X\right]_{t}.$ The $\left(i,\,j\right)$-element
of the \textit{quadratic covariation process} $\left[X,\,X\right]_{t}$
is defined as\footnote{The reader may refer to \citeReferencesSupp{jacod/protter:12} or
\citeReferencesSupp{jacod/shiryaev:03} for a complete introduction
to the material of this section.} 
\begin{align*}
\left[X^{i},\,X^{j}\right]_{t} & =\underset{T\rightarrow\infty}{\textrm{plim}}\sum_{k=1}^{T}\left(X_{kh}^{i}-X_{\left(k-1\right)h}^{i}\right)\left(X_{kh}^{j}-X_{\left(k-1\right)h}^{j}\right),
\end{align*}
where plim denotes the probability limit of the sum. $\left[X,\,X\right]_{t}$
takes values in the cone of all positive semidefinite symmetric $d\times d$
matrices and is continuous in $t$, adapted and of locally finite
variation. Associated with this, we can define the $\left(i,j\right)$-element
of the approximate covariation\textbf{ }matrix as 
\begin{align*}
\sum_{k\geq1}\left(_{h}X_{kh}^{i}-{}_{h}X_{\left(k-1\right)h}^{i}\right)\left(_{h}X_{kh}^{j}-{}_{h}X_{\left(k-1\right)h}^{j}\right) & ,
\end{align*}
which consistently estimates the increments of the quadratic covariation
$\left[X^{i},\,X^{j}\right].$ It is an ex-post estimator of the covariability
between the components of $X$ over the time interval $\left[0,\,t\right]$.
More precisely, as $h\downarrow0$: 
\begin{align*}
\sum_{k\geq1}^{\left\lfloor t/h\right\rfloor }\left(X_{kh}^{i}-X_{\left(k-1\right)h}^{i}\right)\left(X_{kh}^{j}-X_{\left(k-1\right)h}^{j}\right) & \overset{P}{\rightarrow}\int_{0}^{t}\Sigma_{XX,s}^{\left(i,j\right)}ds,
\end{align*}
where $\Sigma_{XX,s}^{\left(i,j\right)}$ is referred to as the \textit{spot
(not integrated)} volatility. 

After this brief review, we turn to the statement of the asymptotic
results for some statistics to be encountered in the proofs below.
We simply refer to \citeReferencesSupp{jacod/protter:12}. More specifically,
Lemma \ref{Lemma LLNs 1}-\ref{Lemma LLNs 2} follow from their Theorem
3.3.1-(b), while Lemma \ref{=00005BJ=000026P,-Theorem-5.4.2=00005D}
follows from their Theorem 5.4.2.
\begin{lem}
\label{Lemma LLNs 1}Under Assumption \ref{Assumption Localization},
we have as $h\downarrow0,\,T\rightarrow\infty$ with $N$ fixed and
for any $1\leq i,\,j\leq p$, 

(i) $\left|\left(Z'_{2}e\right)_{i,1}\right|\overset{P}{\rightarrow}0$
where $\left(Z'_{2}e\right)_{i,1}=\sum_{k=T_{b}+1}^{T}z_{kh}^{\left(i\right)}e_{kh}$; 

(ii) $\left|\left(Z'_{0}e\right)_{i,1}\right|\overset{P}{\rightarrow}0$
where $\left(Z'_{0}e\right)_{i,1}=\sum_{k=T_{b}^{0}+1}^{T}z_{kh}^{\left(i\right)}e_{kh}$;

(iii) $\left|\left(Z'_{2}Z_{2}\right)_{i,j}-\int_{\left(T_{b}+1\right)h}^{N}\Sigma_{ZZ,s}^{\left(i,j\right)}ds\right|\overset{P}{\rightarrow}0$
where $\left(Z'_{2}Z_{2}\right)_{i,j}=\sum_{k=T_{b}+1}^{T}z_{kh}^{\left(i\right)}z_{kh}^{\left(j\right)}$; 

(iv) $\left|\left(Z'_{0}Z_{0}\right)_{i,j}-\int_{\left(T_{b}^{0}+1\right)h}^{N}\Sigma_{ZZ,s}^{\left(i,j\right)}ds\overset{P}{\rightarrow}\right|0$
where $\left(Z'_{0}Z_{0}\right)_{i,j}=\sum_{k=T_{b}^{0}+1}^{T}z_{kh}^{\left(i\right)}z_{kh}^{\left(j\right)}.$ 

For the following estimates involving $X,$ we have, for any $1\leq r\leq p$
and $1\leq l\leq q+p,$ 

(v) $\left|\left(Xe\right)_{l,1}\right|\overset{P}{\rightarrow}0$
where $\left(Xe\right)_{l,1}=\sum_{k=1}^{T}x_{kh}^{\left(l\right)}e_{kh}$; 

(vi) $\left|\left(Z'_{2}X\right)_{r,l}-\int_{\left(T_{b}+1\right)h}^{N}\Sigma_{ZX,s}^{\left(r,l\right)}ds\right|\overset{P}{\rightarrow}0$
where $\left(Z'_{2}X\right)_{r,l}=\sum_{k=T_{b}+1}^{T}z_{kh}^{\left(r\right)}x_{kh}^{\left(l\right)}$; 

(vii) $\left|\left(Z'_{0}X\right)_{r,l}-\int_{\left(T_{b}^{0}+1\right)h}^{N}\Sigma_{ZX,s}^{\left(r,l\right)}ds\right|\overset{P}{\rightarrow}0$
where $\left(Z'_{0}X\right)_{r,l}=\sum_{k=T_{b}^{0}+1}^{T}z_{kh}^{\left(r\right)}x_{kh}^{\left(l\right)}$. 

Further, for $1\leq u,\,d\leq q+p,$ 

(viii) $\left|\left(X'X\right)_{u,d}-\int_{0}^{N}\Sigma_{XX,s}^{\left(u,d\right)}ds\right|\overset{P}{\rightarrow}0$
where $\left(X'X\right)_{u,d}=\sum_{k=1}^{T}x_{kh}^{\left(u\right)}x_{kh}^{\left(d\right)}.$ 
\end{lem}

\begin{lem}
\label{Lemma LLNs 2} Under Assumption \ref{Assumption Localization},
we have as $h\downarrow0,\,T\rightarrow\infty$ with $N$ fixed, $\left|N_{b}^{0}-N_{b}\right|>\gamma>0$
and for any $1\leq i,\,j\leq p$, 

(i) with $\left(Z'_{\Delta}Z_{\Delta}\right)_{i,j}=\sum_{k=T_{b}^{0}+1}^{T_{b}}z_{kh}^{\left(i\right)}z_{kh}^{\left(j\right)}$,
we have 
\[
\begin{cases}
|\left(Z'_{\Delta}Z_{\Delta}\right)_{i,j}-\int_{\left(T_{b}+1\right)h}^{T_{b}^{0}h}\Sigma_{ZZ,s}^{\left(i,j\right)}ds|\overset{P}{\rightarrow}0, & \textrm{if }T_{b}<T_{b}^{0}\\
|\left(Z'_{\Delta}Z_{\Delta}\right)_{i,j}-\int_{T_{b}^{0}h}^{\left(T_{b}+1\right)h}\Sigma_{ZZ,s}^{\left(i,j\right)}ds|\overset{P}{\rightarrow}0, & \textrm{if }T_{b}>T_{b}^{0}
\end{cases};
\]

and for $1\leq r\leq p+q$,

(ii) with $\left(Z'_{\Delta}X_{\Delta}\right)_{i,r}=\sum_{k=T_{b}^{0}+1}^{T_{b}}z_{kh}^{\left(i\right)}x_{kh}^{\left(r\right)}$,
we have 
\[
\begin{cases}
|\left(Z'_{\Delta}X_{\Delta}\right)_{i,r}-\int_{\left(T_{b}+1\right)h}^{T_{b}^{0}h}\Sigma_{ZX,s}^{\left(i,r\right)}ds|\overset{P}{\rightarrow}0, & \textrm{if }T_{b}<T_{b}^{0}\\
|\left(Z'_{\Delta}X_{\Delta}\right)_{i,r}-\int_{T_{b}^{0}h}^{\left(T_{b}+1\right)h}\Sigma_{ZX,s}^{\left(i,r\right)}ds|\overset{P}{\rightarrow}0, & \textrm{if }T_{b}>T_{b}^{0}
\end{cases}.
\]
\end{lem}
Next, we turn to the central limit theorems, they all feature a
limiting process defined on an extension of the original probability
space $\left(\Omega,\,\mathscr{F},\,P\right).$ In order to avoid
non-useful repetitions, we present a general framework valid for all
statistics considered in the paper. The first step is to carry out
an extension of the original probability space $\left(\Omega,\,\mathscr{F},\,P\right).$
We accomplish this in the usual way. We first fix the original probability
space $(\Omega,\,\mathscr{F},\,\{\mathscr{F}_{t}\}_{t\geq0},\,P)$.
Consider an additional measurable space $\left(\Omega^{*},\,\mathscr{F}^{*}\right)$
and a transition probability $Q\left(\omega,\,d\omega^{*}\right)$
from $\left(\Omega,\,\mathscr{F}\right)$ into $\left(\Omega^{*},\,\mathscr{F}^{*}\right)$.
Next, we can define the products $\widetilde{\Omega}=\Omega\times\Omega^{*}$,
$\widetilde{\mathscr{F}}=\mathscr{F\otimes\mathscr{F}}^{*}$ and $\widetilde{P}\left(d\omega,\,d\omega^{*}\right)=P\left(d\omega\right)Q\left(\omega,\,d\omega^{*}\right)$.
This defines the extension $(\widetilde{\Omega},\,\widetilde{\mathscr{F}},\,\widetilde{P})$
of the original space $(\Omega,\,\mathscr{F},\,\{\mathscr{F}_{t}\}_{t\geq0},\,P)$.
Any variable or process defined on either $\Omega$ or $\Omega^{*}$
is extended in the usual to $\widetilde{\Omega}$ as follows: for
example, let $Y_{t}$ be defined on $\Omega$. Then we say that $Y_{t}$
is extended in the usual way to $\widetilde{\Omega}$ by writing $Y_{t}\left(\omega,\,\omega^{*}\right)=Y_{t}\left(\omega\right)$.
Further, we identify $\mathscr{F}_{t}$ with $\mathscr{F}_{t}\otimes\{\emptyset,\,\Omega^{*}\},$
so that we have a filtered space $(\widetilde{\Omega},\,\widetilde{\mathscr{F}},\,\{\mathscr{F}_{t}\}_{t\geq0},\,\widetilde{P}).$
Finally, as for the filtration, we can consider another filtration
$\{\widetilde{\mathscr{F}}_{t}\}_{t\geq0}$ taking the product form
$\widetilde{\mathscr{F}}_{t}=\cap_{s>t}\mathscr{F}_{s}\otimes\mathscr{F}{}_{s}^{*},$
where $\{\mathscr{F}{}_{t}^{*}\}_{t\geq0}$ is a filtration on $\left(\Omega^{*},\,\mathscr{F}^{*}\right)$.
As for the transition probability $Q$ we can consider the simple
form $Q\left(\omega,\,d\omega^{*}\right)=P^{*}\left(d\omega^{*}\right)$
for some probability measure on $\left(\Omega^{*},\,\mathscr{F}^{*}\right)$.
This defines the way a product filtered extension $(\widetilde{\Omega},\,\widetilde{\mathscr{F}},\,\{\widetilde{\mathscr{F}}_{t}\}_{t\geq0},\,\widetilde{P})$
of the original filtered space $(\Omega,\,\mathscr{F},\,\{\mathscr{F}_{t}\}_{t\geq0},\,P)$
is constructed in this paper. Assume that the auxiliary probability
space $(\Omega^{*},\,\mathscr{F}^{*},\,\{\mathscr{F}{}_{t}^{*}\}_{t\geq0},\,P^{*})$
supports a $p^{2}$-dimensional standard Wiener process $W_{s}^{\dagger}$
which is adapted to $\{\widetilde{\mathscr{F}}_{t}\}$. We need some
additional ingredients in order to describe the limiting process.
We choose a progressively measurable ``square-root'' process $\sigma_{Z}^{*}$
of the $\mathcal{M}_{p^{2}\times p^{2}}^{+}$-valued process $\widehat{\Sigma}_{Z,s}$,
whose elements are given by $\widehat{\Sigma}_{Z,s}^{\left(ij,kl\right)}=\Sigma_{Z,s}^{\left(ik\right)}\Sigma_{Z,s}^{\left(jl\right)}$.
Due to the symmetry of $\Sigma_{Z,s},$ the matrix with entries $(\sigma_{Z,s}^{*,\left(ij,kl\right)}+\sigma_{Z,s}^{*,\left(ji,kl\right)})/\sqrt{2}$
is a square-root of the matrix with entries $\widehat{\Sigma}_{Z,s}^{\left(ij,kl\right)}+\widehat{\Sigma}_{Z,s}^{\left(il,jk\right)}$.
Then the process $\mathscr{U}_{t}$ with components $\mathscr{U}_{t}^{\left(r,j\right)}=2^{-1/2}\sum_{k,l=1}^{p}\int_{0}^{t}(\sigma_{Z,s}^{\left(rj,kl\right)}+\sigma_{Z,s}^{\left(jr,kl\right)})dW_{s}^{\dagger\left(kl\right)}$
is, conditionally on $\mathscr{F}$, a continuous Gaussian process
with independent increments and (conditional) covariance $\widetilde{\mathbb{E}}(\mathscr{U}^{\left(r,j\right)}\left(v\right)\mathscr{U}^{\left(k,l\right)}\left(v\right)|\,\mathscr{F})=\int_{T_{b}^{0}h+v}^{T_{b}^{0}h}(\Sigma_{Z,s}^{\left(rk\right)}\Sigma_{Z,s}^{\left(jl\right)}+\Sigma_{Z,s}^{\left(rl\right)}\Sigma_{Z,s}^{\left(jk\right)})ds,$
where $v\leq0.$ The CLT of interest is as follows.
\begin{lem}
\label{=00005BJ=000026P,-Theorem-5.4.2=00005D}Let $Z$ be a continuous
Itô semimartingale satisfying Assumption \ref{Assumption Localization}.
Then, $\left(Nh\right)^{-1/2}$ $(Z_{2}'Z_{2}-(\left[Z,\,Z\right]_{Th}-\left[Z,\,Z\right]_{\left(T_{b}+1\right)h}))\overset{\mathcal{L}-\mathrm{s}}{\Rightarrow}\mathscr{U}.$
\end{lem}

\subsection{Proofs of the Results in Sections \ref{Section Consistency-and-Rate }
and \ref{Section Asymptotic Distribution: Continuous Case}}

\subsubsection{\label{Subsection Additional-Notation Spot Continuous case}Additional
Notation}

In some of the proofs we face a setting in which $N_{b}$ is allowed
to vary within a shrinking neighborhood of $N_{b}^{0}.$ Some estimates
only depend on observations in this window. For example, assume $T_{b}<T_{b}^{0}$
and consider $\sum_{k=T_{b}+1}^{T_{b}^{0}}x_{kh}x'_{kh}$. When $N_{b}$
is allowed to vary within a shrinking neighborhood of $N_{b}^{0},$
this sum approximates a local window of asymptotically shrinking size.
Introduce a sequence of integers $\left\{ l_{T}\right\} $ that satisfies
$l_{T}\rightarrow\infty$ and $l_{T}h\rightarrow0.$ Below we shall
establish a $T^{1-\kappa}$-rate of convergence of $\widehat{\lambda}_{b}$
toward $\lambda_{0},$ considering the case where $N_{b}-N_{b}^{0}=T^{-\left(1-\kappa\right)}$
for some $\kappa\in\left(0,\,1/2\right)$. Hence, define 
\begin{align}
\widehat{\Sigma}_{X}\left(T_{b},\,T_{b}^{0}\right) & \triangleq\sum_{k=T_{b}+1}^{T_{b}^{0}}x_{kh}x'_{kh}=\sum_{k=T_{b}^{0}+1-l_{T}}^{T_{b}^{0}}x_{kh}x'_{kh},\label{Spot one side XX}
\end{align}
where now $l_{T}=\left\lfloor T^{\kappa}\right\rfloor \rightarrow\infty$
and $l_{T}h=h^{1-\kappa}\rightarrow0.$ Note that $1/h^{1-\kappa}$
is the rate of convergence and the interpretation for $\widehat{\Sigma}_{X}\left(T_{b},\,T_{b}^{0}\right)$
is that it involves asymptotically an infinite number of observations
 falling in the shrinking (at rate $h^{1-\kappa}$) block $(\left(T_{b}-1\right)h,\,T_{b}^{0}h]$.
Other statistics involving the regressors and errors are defined similarly:
\begin{align}
\widehat{\Sigma}_{Xe}\left(T_{b},\,T_{b}^{0}\right) & \triangleq\sum_{k=T_{b}+1}^{T_{b}^{0}}x_{kh}e{}_{kh}=\sum_{k=T_{b}^{0}+1-l_{T}}^{T_{b}^{0}}x_{kh}e{}_{kh},\label{Spot one side Xe}
\end{align}
 and 
\begin{align}
\widehat{\Sigma}_{Ze}\left(T_{b},\,T_{b}^{0}\right) & \triangleq\sum_{k=T_{b}^{0}+1-l_{T}}^{T_{b}^{0}}z_{kh}e{}_{kh}.\label{Eq. spot one side Ze}
\end{align}
Further, we let $\overline{\Sigma}_{Xe}\left(T_{b},\,T_{b}^{0}\right)\triangleq h^{-\left(1-\kappa\right)}\int_{N_{b}}^{N_{b}^{0}}\Sigma_{Xe,s}ds$
and analogously when $Z$ replaces $X$. We also define 
\begin{align}
\widehat{\Sigma}_{h,X}\left(T_{b},\,T_{b}^{0}\right) & \triangleq h^{-\left(1-\kappa\right)}\sum_{k=T_{b}^{0}+1-l_{T}}^{T_{b}^{0}}x_{kh}x'_{kh}.\label{Eq. spot one side Ze-1 H}
\end{align}
The proofs of Section \ref{Section Asymptotic Distribution: Continuous Case}
are first given for the case where $\mu_{\cdot,t}$ from equation
\eqref{Model Regressors Integral Form} are identically zero. In the
last step, this is relaxed. Furthermore, throughout the proofs we
proceed conditionally on the processes $\mu_{\cdot,t}$ and $\Sigma_{t}^{0}$
(defined in Assumption \ref{Assumption 2}) so that they are treated
as if they were deterministic. This is a natural strategy since the
processes $\mu_{\cdot,t}$ are of higher order in $h$ and they do
not play any role for the asymptotic results {[}cf. \citeReferencesSupp{barndorff/shephard:04}{]}. 

\subsubsection{Proof of Proposition \ref{Proposition (1) - Consistency}}

\noindent \textit{Proof.} The concentrated sample objective function
evaluated at $\widehat{T}_{b}$ is $Q_{T}(\widehat{T}_{b})=\widehat{\delta}'_{T_{b}}\left(Z_{2}'MZ_{2}\right)\widehat{\delta}_{T_{b}}.$
We have 
\begin{align*}
\widehat{\delta}_{T_{b}} & =\left(Z'_{2}MZ_{2}\right)^{-1}\left(Z'_{2}MY\right)=\left(Z'_{2}MZ_{2}\right)^{-1}\left(Z'_{2}MZ_{0}\right)\delta_{Z}^{0}+\left(Z'_{2}MZ_{2}\right)^{-1}Z_{2}Me,
\end{align*}
 and $\widehat{\delta}_{T_{b}^{0}}=\left(Z'_{0}MZ_{0}\right)^{-1}\left(Z'_{0}MY\right)=\delta_{Z}^{0}+\left(Z'_{0}MZ_{0}\right)^{-1}\left(Z'_{0}Me\right)$.
Therefore,
\begin{align}
Q_{T}\left(T_{b}\right)-Q_{T}\left(T_{b}^{0}\right) & =\widehat{\delta}'_{T_{b}}\left(Z_{2}'MZ_{2}\right)\widehat{\delta}_{T_{b}}-\widehat{\delta}'_{T_{b}^{0}}\left(Z_{0}'MZ_{0}\right)\widehat{\delta}_{T_{b}^{0}}\label{eq. A.2.0, Casini (2015a)-2}\\
 & =\left(\delta_{Z}^{0}\right)'\left\{ \left(Z'_{0}MZ_{2}\right)\left(Z'_{2}MZ_{2}\right)^{-1}\left(Z'_{2}MZ_{0}\right)-Z_{0}'MZ_{0}\right\} \delta_{Z}^{0}\label{eq. A.2.1, Casini (2015a)-2}\\
 & \quad+g_{e}\left(T_{b}\right),\label{eq. A.2.2-2}
\end{align}
 where
\begin{align}
g_{e}\left(T_{b}\right) & =2\left(\delta_{Z}^{0}\right)'\left(Z'_{0}MZ_{2}\right)\left(Z'_{2}MZ_{2}\right)^{-1}Z_{2}Me-2\left(\delta_{Z}^{0}\right)'\left(Z'_{0}Me\right)\label{eq. A.2.3, Casini (2015a)-2}\\
 & \quad+e'MZ_{2}\left(Z'_{2}MZ_{2}\right)^{-1}Z_{2}Me-e'MZ_{0}\left(Z'_{0}MZ_{0}\right)^{-1}Z'_{0}Me.\label{eq. A.2.4, Casini (2015a)-2}
\end{align}
Denote 
\begin{flalign*}
 &  & X_{\Delta} & \triangleq X_{2}-X_{0}=\left(0,\,\ldots,\,0,\,x_{\left(T_{b}+1\right)h},\ldots,\,x_{T_{b}^{0}h},\,0,\ldots,\,\right)', & \textrm{for }T_{b}<T_{b}^{0}\\
 &  & X_{\Delta} & \triangleq-\left(X_{2}-X_{0}\right)=\left(0,\,\ldots,\,0,\,x_{\left(T_{b}^{0}+1\right)h},\ldots,\,x_{T_{b}h},\,0,\ldots,\,\right)', & \textrm{for }T_{b}>T_{b}^{0}\\
 &  & X_{\Delta} & \triangleq0, & \textrm{for }T_{b}=T_{b}^{0}.
\end{flalign*}
Observe that when $T_{b}^{0}\neq T_{b}$ we have $X_{2}=X_{0}+X_{\Delta}\textrm{sign}\left(T_{b}^{0}-T_{b}\right)$.
When the sign is immaterial, we simply write $X_{2}=X_{0}+X_{\Delta}$.
Next, let $Z_{\Delta}=X_{\Delta}R$, and define
\begin{align}
r\left(T_{b}\right) & \triangleq\frac{\left(\delta_{Z}^{0}\right)'\left\{ \left(Z_{0}'MZ_{0}\right)-\left(Z'_{0}MZ_{2}\right)\left(Z'_{2}MZ_{2}\right)^{-1}\left(Z'_{2}MZ_{0}\right)\right\} \delta_{Z}^{0}}{\left|T_{b}-T_{b}^{0}\right|}.\label{eq. A.2.5, Casini (2015a)-2}
\end{align}
We arbitrarily define $r\left(T_{b}\right)=\left(\delta_{Z}^{0}\right)'\delta_{Z}^{0}$
when $T_{b}=T_{b}^{0}$. We write \eqref{eq. A.2.0, Casini (2015a)-2}
as
\begin{align}
Q_{T}\left(T_{b}\right)-Q_{T}\left(T_{0}\right) & =-\left|T_{b}-T_{b}^{0}\right|r\left(T_{b}\right)+g_{e}\left(T_{b}\right),\qquad\textrm{for all }T_{b}.\label{eq. A.2.6, Casini (2015a)-2}
\end{align}
By definition, $\widehat{T}_{b}$ is an extremum estimator and thus
satisfies $g_{e}(\widehat{T}_{b})\geq\left|\widehat{T}_{b}-T_{b}^{0}\right|r(\widehat{T}_{b})$.
Therefore, 
\begin{align}
P\left(\left|\widehat{\lambda}_{b}-\lambda_{0}\right|>K\right) & =P\left(\left|\widehat{T}_{b}-T_{b}^{0}\right|>TK\right)\nonumber \\
 & \leq P\left(\sup_{\left|T_{b}-T_{b}^{0}\right|>TK}\left|g_{e}\left(T_{b}\right)\right|\geq\inf_{\left|T_{b}-T_{b}^{0}\right|>TK}\left|T_{b}-T_{b}^{0}\right|r\left(T_{b}\right)\right)\nonumber \\
 & \leq P\left(\sup_{p\leq T_{b}\leq T-p}\left|g_{e}\left(T_{b}\right)\right|\geq TK\inf_{\left|T_{b}-T_{b}^{0}\right|>TK}r\left(T_{b}\right)\right)\label{eq. last line of A.2.6.b-2}\\
 & =P\left(r_{T}^{-1}\sup_{p\leq T_{b}\leq T-p}\left|g_{e}\left(T_{b}\right)\right|\geq K\right),\nonumber 
\end{align}
where we recall that $p\leq T_{b}\leq T-p$ is needed for identification,
and $r_{T}\triangleq T\inf_{\left|T_{b}-T_{b}^{0}\right|>TK}r\left(T_{b}\right).$
Lemma \ref{Lemma r_CT, Prop 1-2} below shows that $r_{T}$ is positive
and bounded away from zero. Thus, it is sufficient to verify that
the stochastic component is negligible as $h\downarrow0$, i.e., 
\begin{align}
\sup_{p\leq T_{b}\leq T-p}\left|g_{e}\left(T_{b}\right)\right| & =o_{p}\left(1\right).\label{eq. A.2.7, Casini (2015a)-2}
\end{align}
The first term of $g_{e}\left(T_{b}\right)$ is
\begin{align}
2\left(\delta_{Z}^{0}\right)'\left(Z'_{0}MZ_{2}\right)\left(Z'_{2}MZ_{2}\right)^{-1/2}\left(Z'_{2}MZ_{2}\right)^{-1/2}Z_{2}Me & .\label{eq. A.2.8, Casini (2015a)-2}
\end{align}
Lemma \ref{=00005BJ=000026P,-Theorem-5.4.2=00005D} implies that for
any $1\leq j\leq p,$ $\left(Z_{2}e\right)_{j,1}/\sqrt{h}=O_{p}\left(1\right)$
and for any $1\leq i\leq q+p,$ $\left(Xe\right)_{i,1}/\sqrt{h}=O_{p}\left(1\right)$.
These hold because they both involve a positive fraction of the data.
Furthermore, from Lemma \ref{Lemma LLNs 1}, we also have that $Z'_{2}MZ_{2}$
and $Z'_{0}MZ_{2}$ are $O_{p}\left(1\right).$ Therefore, the supremum
of $\left(Z'_{0}MZ_{2}\right)\left(Z'_{2}MZ_{2}\right)^{-1/2}$ over
all $T_{b}$ is 
\begin{align*}
\sup_{T_{b}}\left(Z'_{0}MZ_{2}\right)\left(Z'_{2}MZ_{2}\right)^{-1}\left(Z'_{2}MZ_{0}\right)\leq Z'_{0}MZ_{0} & =O_{p}\left(1\right),
\end{align*}
 by Lemma \ref{Lemma LLNs 1}. By Assumption \eqref{Assumption 1, CT}-(iii),
$\left(Z'_{2}MZ_{2}\right)^{-1/2}Z'_{2}Me$ is $O_{p}\left(1\right)O_{p}(\sqrt{h})$
uniformly, which implies that \eqref{eq. A.2.8, Casini (2015a)-2}
is $O_{p}(\sqrt{h})$ uniformly over $p\leq T_{b}\leq T-p$. As for
the second term of \eqref{eq. A.2.3, Casini (2015a)-2}, $Z'_{0}Me=O_{p}(\sqrt{h}).$
The first term in \eqref{eq. A.2.4, Casini (2015a)-2} is uniformly
$o_{p}\left(1\right)$ and the same holds for the last term. Therefore,
combining these results, $\sup_{T_{b}}\left|g_{e}\left(T_{b}\right)\right|=O_{p}(\sqrt{h})$
uniformly when $|\widehat{\lambda}_{b}-\lambda_{0}|>K$. Therefore
for some $B>0$, these arguments combined with Lemma \ref{Lemma r_CT, Prop 1-2}
below result in $P(r_{B}^{-1}\sup_{p\leq T_{b}\leq T-p}\left|g_{e}\left(T_{b}\right)\right|\geq K)\leq\varepsilon,$
from which it follows that the right-hand side of \eqref{eq. last line of A.2.6.b-2}
is weakly smaller than $\varepsilon.$ This concludes the proof since
$\varepsilon>0$ was arbitrary. $\square$
\begin{lem}
\label{Lemma r_CT, Prop 1-2}For $B>0$, let $r_{B}=\inf_{\left|T_{b}-T_{b}^{0}\right|>TB}Tr\left(T_{b}\right).$
There exists a $\kappa>0$ such that for every $\varepsilon>0,$
there exists a $B<\infty$ such that $P\left(r_{B}\geq\kappa\right)\leq1-\varepsilon$,
i.e., $r_{B}$ is positive and bounded away from zero with high probability. 
\end{lem}
\noindent \textit{Proof.} Assume $T_{b}\leq T_{b}^{0}$ and observe
that $r_{T}\geq r_{B}$ for an appropriately chosen $B.$ From the
first inequality result in Lemma \ref{Lemma A1, Casini (2015a)},
\[
r\left(T_{b}\right)\geq\left(\delta_{Z}^{0}\right)'R'\left(X'_{\Delta}X_{\Delta}/\left(T_{b}^{0}-T_{b}\right)\right)\left(X'_{2}X_{2}\right)^{-1}\left(X'_{0}X_{0}\right)R\delta_{Z}^{0}.
\]
 When multiplied by $T,$ we have
\begin{align*}
Tr\left(T_{b}\right) & \geq T\left(\delta_{Z}^{0}\right)'R'\frac{X'_{\Delta}X_{\Delta}}{T_{b}^{0}-T_{b}}\left(X'_{2}X_{2}\right)^{-1}\left(X'_{0}X_{0}\right)R\delta_{Z}^{0}\\
 & =\left(\delta_{Z}^{0}\right)'R'\frac{X'_{\Delta}X_{\Delta}}{N_{b}^{0}-N_{b}}\left(X'_{2}X_{2}\right)^{-1}\left(X'_{0}X_{0}\right)R\delta_{Z}^{0}.
\end{align*}
Note that $0<K<B<h\left(T_{b}^{0}-T_{b}\right)<N$. Then, 
\[
Tr\left(T_{b}\right)\geq\left(\delta_{Z}^{0}\right)'R'\left(X'_{\Delta}X_{\Delta}/N\right)\left(X'_{2}X_{2}\right)^{-1}\left(X'_{0}X_{0}\right)R\delta_{Z}^{0},
\]
 and by standard estimates for It\^o semimartingales, $X'_{\Delta}X_{\Delta}=O_{p}\left(1\right)$
(i.e., use the Burkh\"{o}lder-Davis-Gundy  inequality, recalling
that $|\widehat{N}_{b}-N_{b}^{0}|>BN$). Hence, we conclude that $Tr\left(T_{b}\right)\geq\left(\delta_{Z}^{0}\right)'R'O_{p}\left(1/N\right)O_{p}\left(1\right)R\delta_{Z}^{0}$
$\geq\kappa>0$, where $\kappa$ is some positive constant. The last
inequality follows whenever $X'_{\Delta}X_{\Delta}$ is positive definite
since $R'X'_{\Delta}X_{\Delta}\left(X'_{2}X_{2}\right)^{-1}\left(X'_{0}X_{0}\right)R$
can be rewritten as $R'[\left(X'_{0}X_{0}\right)^{-1}+\left(X'_{\Delta}X_{\Delta}\right)^{-1}]R$.
According to Lemma \ref{Lemma LLNs 1}, $X'_{2}X_{2}$ is $O_{p}\left(1\right)$.
The same argument applies to $X'_{0}X_{0}$, which together with the
the fact that $R$ has full common rank in turn implies that we can
choose a $B>0$ such that $r_{B}=\inf_{\left|T_{b}-T_{b}^{0}\right|>TB}Tr\left(T_{b}\right)$
satisfies $P\left(r_{B}\geq\kappa\right)\leq1-\varepsilon.$ The case
with $T_{b}>T_{b}^{0}$ is similar and omitted. $\square$

\subsubsection{Proof of Proposition \ref{Proposition 2, (Rate of Convergence)}}

\noindent \textit{Proof.} Given the consistency result, one can restrict
attention to the local behavior of the objective function for those
values of $T_{b}$ in $\mathbf{B}_{T}\triangleq\left\{ T_{b}:\,T\eta\leq T_{b}\leq T\left(1-\eta\right)\right\} ,$
where $\eta>0$ satisfies $\eta\leq\lambda_{0}\leq1-\eta.$ By Proposition
\ref{Proposition (1) - Consistency}, the estimator $\widehat{T}_{b}$
will visit the set $\mathbf{B}_{T}$ with large probability as $T\rightarrow\infty.$
That is, for any $\varepsilon>0,$ $P\left(\widehat{T}_{b}\notin\mathbf{B}_{T}\right)<\varepsilon$
for sufficiently large $T.$ We show that for large $T,$ $\widehat{T}_{b}$
eventually falls in the set $\mathbf{B}_{K,T}\triangleq\left\{ T_{b}:\,\left|N_{b}-N_{b}^{0}\right|\leq KT^{-1}\right\} ,$
for some $K>0.$ For any $K>0,$ define the intersection of $\mathbf{B}_{T}$
and the complement of $\mathbf{B}_{K,T}$ by 
\[
\mathbf{D}_{K,T}\triangleq\Bigl\{ T_{b}:\,N\eta\leq N_{b}\leq N\left(1-\eta\right),\left.\left|N_{b}-N_{b}^{0}\right|>KT^{-1}\right\} .
\]
 Notice that
\begin{align*}
 & \left\{ \left|\widehat{\lambda}_{b}-\lambda_{0}\right|>KT^{-1}\right\} =\\
 & \quad\left\{ \left|\widehat{\lambda}_{b}-\lambda_{0}\right|>KT^{-1}\cap\widehat{\lambda}_{b}\in\left(\eta,\,1-\eta\right)\right\} \\
 & \quad\cup\left\{ \left|\widehat{\lambda}_{b}-\lambda_{0}\right|>KT^{-1}\cap\widehat{\lambda}_{b}\notin\left(\eta,\,1-\eta\right)\right\} \\
 & \quad\subseteq\left\{ \left|\widehat{\lambda}_{b}-\lambda_{0}\right|>K\left(T^{-1}\right)\cap\widehat{\lambda}_{b}\in\left(\eta,\,1-\eta\right)\right\} \cup\left\{ \widehat{\lambda}_{b}\notin\left(\eta,\,1-\eta\right)\right\} ,
\end{align*}
 and so
\begin{align*}
P\left(\left|\widehat{\lambda}_{b}-\lambda_{0}\right|>KT^{-1}\right) & \leq P\left(\widehat{\lambda}_{b}\notin\left(\eta,\,1-\eta\right)\right)\\
 & \quad+P\left(\left|\widehat{T}_{b}-T_{b}^{0}\right|>K\cap\widehat{\lambda}_{b}\in\left(\eta,\,1-\eta\right)\right),
\end{align*}
 and for large $T$, 
\begin{align*}
P\left(\left|\widehat{\lambda}_{b}-\lambda_{0}\right|>KT^{-1}\right) & \leq\varepsilon+P\left(\left|\widehat{\lambda}_{b}-\lambda_{0}\right|>KT^{-1}\cap\widehat{\lambda}_{b}\in\left(\eta,\,1-\eta\right)\right)\\
 & \leq\varepsilon+P\left(\sup_{T_{b}\in\mathbf{D}_{K,T}}Q_{T}\left(T_{b}\right)\geq Q_{T}\left(T_{b}^{0}\right)\right).
\end{align*}
Therefore it is enough to show that the second term above is negligible
as $h\downarrow0$. Suppose $T_{b}<T_{b}^{0}$. Since $\widehat{T}_{b}=\arg\max Q_{T}\left(T_{b}\right),$
it is enough to show that $P(\sup_{T_{b}\in\mathbf{D}_{K,T}}Q_{T}\left(T_{b}\right)\geq Q_{T}\left(T_{b}^{0}\right))<\varepsilon$.
Note that this implies $\left|T_{b}-T_{b}^{0}\right|>KN^{-1}.$ Therefore,
we have to deal with a setting where the time span in $\mathbf{D}_{K,T}$
between $N_{b}$ and $N_{b}^{0}$ is actually shrinking. The difficulty
arises from the quantities depending on the difference $\left|N_{b}-N_{b}^{0}\right|$.
We can rewrite $Q_{T}\left(T_{b}\right)\geq Q_{T}\left(T_{b}^{0}\right)$
as $g_{e}\left(T_{b}\right)/\left|T_{b}-T_{b}^{0}\right|\geq r\left(T_{b}\right),$
with $g_{e}\left(T_{b}\right)$ and $r\left(T_{b}\right)$ as defined
above. Thus, we need to show,
\begin{align*}
P & \left(\sup_{T_{b}\in\mathbf{D}_{K,T}}h^{-1}\frac{g_{e}\left(T_{b}\right)}{\left|T_{b}-T_{b}^{0}\right|}\geq\inf_{T_{b}\in\mathbf{D}_{K,T}}h^{-1}r\left(T_{b}\right)\right)<\varepsilon.
\end{align*}
By Lemma \ref{Lemma A1, Casini (2015a)},
\begin{align*}
\inf_{T_{b}\in\mathbf{D}_{K,T}}r\left(T_{b}\right) & \geq\inf_{T_{b}\in\mathbf{D}_{K,T}}\left(\delta_{Z}^{0}\right)'R'\frac{X'_{\Delta}X_{\Delta}}{\left|T_{b}-T_{b}^{0}\right|}\left(X'_{2}X_{2}\right)^{-1}\left(X'_{0}X_{0}\right)R\delta_{Z}^{0}.
\end{align*}
The asymptotic results used so far rely on statistics involving integrated
covariation between continuous semimartingales. However, since $\left|T_{b}-T_{b}^{0}\right|>K/N,$
the context is different and the same results do not apply because
the time horizon is decreasing as the sample size increases for quantities
depending on $\left|N_{b}-N_{b}^{0}\right|.$ Thus, we shall apply
asymptotic results for the local approximation of the covariation
between processes. Moreover, when $\left|T_{b}-T_{b}^{0}\right|>K/N$,
there are at least $K$ terms in this sum with asymptotically vanishing
moments. That is, for any $1\leq i,\,j\leq q+p$, we have $\mathbb{E}[x_{kh}^{\left(i\right)}x_{kh}^{\left(j\right)}|\,\mathscr{F}_{\left(k-1\right)h}]=\Sigma_{X,\left(k-1\right)h}^{\left(i,j\right)}h,$
and note that $x_{kh}/\sqrt{h}$ is i.n.d. with finite variance and
thus by Assumption \ref{Assumption 4 Eigenvalue} we can always choose
a $K$ large enough such that $\left(h\left|T_{b}-T_{b}^{0}\right|\right)^{-1}X'_{\Delta}X_{\Delta}=\left(h\left|T_{b}-T_{b}^{0}\right|\right)^{-1}\sum_{k=T_{b}+1}^{T_{b}^{0}}x_{kh}x'_{kh}=A>0$
for all $T_{b}\in\mathbf{D}_{K,T}$. This shows that $\inf_{T_{b}\in\mathbf{D}_{K,T}}h^{-1}r\left(T_{b}\right)$
is bounded away from zero. Note that for the other terms in $r\left(T_{b}\right)$
we can use the same arguments since they do not depend on $\left|N_{b}-N_{b}^{0}\right|$.
Hence,
\begin{align}
P\left(\sup_{T_{b}\in\mathbf{D}_{K,T}}h^{-1}\left(T_{b}^{0}-T_{b}\right)^{-1}g_{e}\left(T_{b}\right)\geq B/N\right) & <\varepsilon,\label{eq. (17),}
\end{align}
 for some $B>0.$ Consider the terms of $g_{e}\left(T_{b}\right)$
in \eqref{eq. A.2.3, Casini (2015a)-2}. When $T_{b}\in\mathbf{D}_{K,T}$,
$Z_{2}$ involves at least a positive fraction $N\eta$ of the data.
From Lemma \ref{Lemma LLNs 1}, as $h\downarrow0$, it follows that
\begin{align*}
h^{-1}\left(T_{b}^{0}-T_{b}\right)^{-1} & e'MZ_{2}\left(Z'_{2}MZ_{2}\right)^{-1}Z_{2}Me\\
 & =\left(T_{b}^{0}-T_{b}\right)^{-1}h^{-1}O_{p}\left(h^{1/2}\right)O_{p}\left(1\right)O_{p}\left(h^{1/2}\right)=\frac{O_{p}\left(1\right)}{T_{b}^{0}-T_{b}},
\end{align*}
 uniformly in $T_{b}$. Choose $K$ large enough so that the probability
that the right-hand size is larger than $B/N$ is less than $\varepsilon/4.$
A similar argument holds for the second term in \eqref{eq. A.2.4, Casini (2015a)-2}.
Next consider the first term of $g_{e}\left(T_{b}\right)$ in \eqref{eq. A.2.3, Casini (2015a)-2}.
Using $Z_{2}=Z_{0}\pm Z_{\Delta}$ we can deduce that
\begin{align}
\left(\delta_{Z}^{0}\right)' & \left(Z'_{0}MZ_{2}\right)\left(Z'_{2}MZ_{2}\right)^{-1}Z_{2}Me\nonumber \\
 & =\left(\delta_{Z}^{0}\right)'\left(\left(Z'_{2}\pm Z'_{\Delta}\right)MZ_{2}\right)\left(Z'_{2}MZ_{2}\right)^{-1}Z_{2}Me\nonumber \\
 & =\left(\delta_{Z}^{0}\right)'Z'_{0}Me\pm\left(\delta_{Z}^{0}\right)'Z'_{\Delta}Me\label{eq. (18), Casini (2015a)-1}\\
 & \quad\pm\left(\delta_{Z}^{0}\right)'\left(Z'_{\Delta}MZ_{2}\right)\left(Z'_{2}MZ_{2}\right)^{-1}Z_{2}Me,\nonumber 
\end{align}
from which it follows that
\begin{align}
 & \left|2\left(\delta_{Z}^{0}\right)'\left(Z'_{0}MZ_{2}\right)\left(Z'_{2}MZ_{2}\right)^{-1}Z_{2}Me-2\left(\delta_{Z}^{0}\right)'\left(Z'_{0}Me\right)\right|\nonumber \\
 & =\left|\left(\delta_{Z}^{0}\right)'Z'_{\Delta}Me\right|+\left|\left(\delta_{Z}^{0}\right)'\left(Z'_{\Delta}MZ_{2}\right)\left(Z'_{2}MZ_{2}\right)^{-1}\left(Z_{2}Me\right)\right|.\label{eq. 21, Casini (2015a)-1}
\end{align}
First, we can apply Lemma \ref{Lemma LLNs 1} {[}(vi) and (viii){]},
and Lemma \ref{Lemma LLNs 2} {[}(i)-(ii){]}, together with Assumption
\ref{Assumption 1, CT}-(iii), to terms that do not involve $\left|N_{b}-N_{b}^{0}\right|$,
i.e., 
\begin{align*}
h^{-1}\left(\delta_{Z}^{0}\right)'\left(Z'_{\Delta}MZ_{2}\right) & =h^{-1}\left(\delta_{Z}^{0}\right)'\left(Z'_{\Delta}Z_{2}\right)-h^{-1}\left(\delta_{Z}^{0}\right)'\left(Z'_{\Delta}X_{\Delta}\left(X'X\right)^{-1}X'Z_{2}\right)\\
 & =\frac{\left(\delta_{Z}^{0}\right)'\left(Z'_{\Delta}Z_{\Delta}\right)}{h}-\left(\delta_{Z}^{0}\right)'\left(\frac{Z'_{\Delta}X_{\Delta}}{h}\left(X'X\right)^{-1}X'Z_{2}\right).
\end{align*}
Consider $Z'_{\Delta}Z_{\Delta}$. By the same reasoning as above,
whenever $T_{b}\in\mathbf{D}_{K,T}$, $\left(Z'_{\Delta}Z_{\Delta}\right)/h\left(T_{b}^{0}-T_{b}\right)=O_{p}\left(1\right)$
for $K$ large enough. The term $Z'_{\Delta}X_{\Delta}/h\left(T_{b}^{0}-T_{b}\right)$
is also $O_{p}\left(1\right)$ uniformly. Thus, it follows from Lemma
\ref{=00005BJ=000026P,-Theorem-5.4.2=00005D} that the second term
of \eqref{eq. 21, Casini (2015a)-1} is $O_{p}(h^{1/2})$. Next, note
that $Z'_{\Delta}Me=Z'_{\Delta}e-Z'_{\Delta}X\left(X'X\right)^{-1}X'e.$
We can write
\begin{align*}
\frac{Z'_{\Delta}Me}{\left(T_{b}^{0}-T_{b}\right)h} & =\frac{1}{\left(T_{b}^{0}-T_{b}\right)h}\sum_{k=T_{b}+1}^{T_{b}^{0}}z_{kh}e_{kh}\\
 & \quad-\frac{1}{\left(T_{b}^{0}-T_{b}\right)h}\left(\sum_{k=T_{b}+1}^{T_{b}^{0}}z_{kh}x'_{kh}\right)\left(X'X\right)^{-1}\left(X'e\right).
\end{align*}
Note that the sequence $\{h^{-1/2}z_{kh}h^{-1/2}x_{kh}\}$ is i.n.d.
with finite mean identically in $k$. There are at least $K$ terms
in this sum, so $(\sum_{k=T_{b}+1}^{T_{b}^{0}}z_{kh}x'_{kh})/\left(T_{b}^{0}-T_{b}\right)h$
is $O_{p}\left(1\right)$ for a large enough $K$ in view of Assumption
\ref{Assumption 4 Eigenvalue}. Then,
\begin{align}
\frac{1}{\left(T_{b}^{0}-T_{b}\right)h}\left(\sum_{k=T_{b}+1}^{T_{b}^{0}}z_{kh}x'_{kh}\right)\left(X'X\right)^{-1}\left(X'e\right) & =O_{p}\left(1\right)O_{p}\left(1\right)O_{p}\left(h^{1/2}\right),\label{Eq. h^-1 ZX XX Xe}
\end{align}
when $K$ is large. Thus, 
\begin{align}
\frac{1}{\left(T_{b}^{0}-T_{b}\right)h}g_{e}\left(T_{b}\right) & =\frac{1}{\left(T_{b}^{0}-T_{b}\right)h}\left(\delta^{0}\right)'2Z'_{\Delta}e+\frac{O_{p}\left(1\right)}{T_{b}^{0}-T_{b}}+O_{p}\left(h^{1/2}\right).\label{eq. A.12}
\end{align}
We can now prove \eqref{eq. (17),} using \eqref{eq. A.12}. To this
end, we sneed a $K>0$, such that 
\begin{align}
P & \left(\sup_{T_{b}\in\mathbf{D}_{K,T}}\left\Vert \left(\delta_{Z}^{0}\right)'\frac{2}{h}\frac{1}{T_{b}^{0}-T_{b}}\sum_{k=T_{b}+1}^{T_{b}^{0}}z_{kh}e_{kh}\right\Vert >\frac{B}{4N}\right)\label{Eq. Haek Prop 2}\\
 & \leq P\left(\sup_{T_{b}\leq T_{b}^{0}-KN^{-1}}\left\Vert \frac{1}{h}\frac{1}{T_{b}^{0}-T_{b}}\sum_{k=T_{b}+1}^{T_{b}^{0}}z_{kh}e_{kh}\right\Vert >\frac{B}{8N\left\Vert \delta_{Z}^{0}\right\Vert }\right)<\varepsilon.\nonumber 
\end{align}
 Note that $\left|T_{b}-T_{b}^{0}\right|$ is bounded away from zero
in $\mathbf{D}_{K,T}$. Observe that $(z_{kh}/\sqrt{h})(e_{kh}/\sqrt{h})$
are independent in $k$ and have zero mean and finite second moments.
Hence, by the Hájek-Réiny inequality {[}see Lemma A.6 in \citeReferencesSupp{bai/perron:98}{]},
\begin{align*}
P & \left(\sup_{T_{b}\leq T_{b}^{0}-KN^{-1}}\left\Vert \frac{1}{T_{b}^{0}-T_{b}}\sum_{k=T_{b}+1}^{T_{b}^{0}}\frac{z_{kh}}{\sqrt{h}}\frac{e_{kh}}{\sqrt{h}}\right\Vert >\frac{B}{8\left\Vert \delta_{Z}^{0}\right\Vert N}\right)\\
 & \leq A\frac{64\left\Vert \delta_{Z}^{0}\right\Vert ^{2}N^{2}}{B^{2}}\frac{1}{KN^{-1}},
\end{align*}
where $A>0$. We can choose $K$ large enough such that the right-hand
side is less than $\varepsilon/4.$ Combining the above arguments,
we deduce the claim in \eqref{eq. (17),} which then concludes the
proof of Proposition \ref{Proposition 2, (Rate of Convergence)}.
$\square$

\subsubsection{Proof of Proposition \ref{Proposition OLS Asymtptoc Distribu}}

We focus on the case with $T_{b}\leq T_{0}.$ The arguments for the
other case are similar and omitted. From Proposition \ref{Proposition (1) - Consistency}
the distance $|\widehat{\lambda}_{b}-\lambda_{0}|$ can be made arbitrary
small. Proposition \ref{Proposition 2, (Rate of Convergence)} gives
the associated rate of convergence: $T(\widehat{\lambda}_{b}-\lambda_{0})=O_{p}\left(1\right).$
Given the consistency result for $\widehat{\lambda}_{b}$, we can
apply a restricted search. In particular, by Proposition \ref{Proposition 2, (Rate of Convergence)},
for large $T>\overline{T},$ we know that $\{T_{b}\notin\mathbf{D}_{K,T}\},$
or equivalently $|T_{b}-T_{b}^{0}|\leq K$, with high probability
for some  $K.$ Essentially, what we shall show is that from the results
of Proposition \ref{Proposition (1) - Consistency}-\ref{Proposition 2, (Rate of Convergence)}
the error in replacing $T_{b}^{0}$ with $\widehat{T}_{b}$ is stochastically
small and thus it does not affect the estimation of the parameters
$\beta^{0},\,\delta_{Z,1}^{0}$ and $\delta_{Z,2}^{0}$. Toward this
end, we first find a lower bound on the convergence rate for $\widehat{\lambda}_{b}$
that guarantees its estimation to be asymptotically independent from
that of the regression parameters. This result will also be used in
later proofs. We shall see that the rate of convergence established
in Proposition \ref{Proposition 2, (Rate of Convergence)} is strictly
faster than the lower bound. Below, we use $\widehat{T}_{b}$ in order
to construct $Z_{2}$ and define $\widehat{Z}_{0}\triangleq Z_{2}.$ 
\begin{lem}
\label{Lemma gamma OLS}Fix $\gamma\in\left(0,\,1/2\right)$ and some
constant $A>0.$ For all large $T>\overline{T},$ if $|\widehat{N}_{b}-N_{b}^{0}|\leq AO_{p}\left(h^{1-\gamma}\right)$,
then $X'(Z_{0}-\widehat{Z}_{0})=O_{p}(h^{1-\gamma})$ and $Z'_{0}(Z_{0}-\widehat{Z}_{0})=O_{p}(h^{1-\gamma})$.
\end{lem}
\noindent \textit{Proof.} Note that the setting of Proposition \ref{Proposition 2, (Rate of Convergence)}
satisfies the conditions of this lemma because $\widehat{N}_{b}-N_{b}^{0}=O_{p}\left(h\right)\leq AO_{p}\left(h^{1-\gamma}\right)$
as $h\downarrow0$. By assumption, there exists some constant $C>0$
such that $P(h^{\gamma}|\widehat{T}_{b}-T_{b}^{0}|>C)<\varepsilon.$
We have to show that although we only know $|\widehat{T}_{b}-T_{b}^{0}|\leq Ch^{-\gamma},$
the error when replacing $T_{b}^{0}$ by $\widehat{T}_{b}$ in the
construction of $Z_{2}$ goes to zero fast enough. This is achieved
because $|\widehat{N}_{b}-N_{b}^{0}|\rightarrow0$ at least at rate
$h^{1-\gamma}$, which is faster than the standard convergence rate
for regression parameters (i.e., $\sqrt{T}$-rate). Without loss of
generality we take $C=1$. We have
\begin{align*}
h^{-1/2}X'\left(Z_{0}-\widehat{Z}_{0}\right) & =h^{1/2-\gamma}\frac{1}{h^{1-\gamma}}\sum_{T_{b}^{0}-\left\lfloor T^{\gamma}\right\rfloor }^{T_{b}^{0}}x_{kh}z_{kh}.
\end{align*}
 Notice that, as $h\downarrow0$, the number of terms in the sum on
the right-hand side, for all $T>\overline{T}$, increases to infinity
at rate $1/h^{\gamma}$. Since $\widehat{N}_{b}$ approaches $N_{b}^{0}$
at rate $T^{-\left(1-\gamma\right)}$, the quantity $X'(Z_{0}-\widehat{Z}_{0})/h^{1-\gamma}$
is a consistent estimate of the so-called instantaneous or spot covariation
between $X$ and $Z$ at time $N_{b}^{0}$. Theorem 9.3.2 part (i)
in \citeReferencesSupp{jacod/protter:12} can be applied since the
``window'' is decreasing at rate $h^{1-\gamma}$ and the same factor
$h^{1-\gamma}$ is in the denominator. Thus, we have as $h\downarrow0,$
\begin{align}
X'_{\Delta}Z_{\Delta}/h^{1-\gamma} & \overset{P}{\rightarrow}\Sigma_{XX,N_{b}^{0}},\label{Eq. 1, Proof Corrollary 1 - X'X/h-2-1}
\end{align}
which implies that $h^{-1/2}X'\left(Z_{0}-\widehat{Z}_{0}\right)=O_{p}(h^{1/2-\gamma}).$
This shows that the order of the error in replacing $Z_{0}$ by $Z_{2}=\widehat{Z}_{0}$
goes to zero at a fast enough rate. That is, by definition we can
write 
\[
Y=X\beta^{0}+\widehat{Z}_{0}\delta_{Z}^{0}+\left(Z_{0}-\widehat{Z}_{0}\right)\delta_{Z}^{0}+e,
\]
 from which it follows that 
\[
X'\widehat{Z}_{0}=X'Z_{0}+o_{p}\left(1\right),\,X'\left(Z_{0}-\widehat{Z}_{0}\right)\delta_{Z}^{0}=o_{p}\left(1\right),
\]
and $Z_{0}'(Z_{0}-\widehat{Z}_{0})\delta_{Z}^{0}=o_{p}\left(1\right).$
To see this, consider for example
\begin{align*}
X'\left(\widehat{Z}_{0}-Z_{0}\right) & =\sum_{T_{b}^{0}-\left\lfloor T^{\gamma}\right\rfloor }^{T_{b}^{0}}x_{kh}z_{kh}=\frac{h^{1-\gamma}}{h^{1-\gamma}}\sum_{T_{b}^{0}-\left\lfloor T^{\gamma}\right\rfloor }^{T_{b}^{0}}x_{kh}z_{kh}=h^{1-\gamma}O_{p}\left(1\right),
\end{align*}
which clearly implies that $X'\widehat{Z}_{0}=X'Z_{0}+o_{p}\left(1\right)$.
The other case can be proved similarly. This concludes the proof of
the Lemma. $\square$

\textit{\noindent }Using Lemma \ref{Lemma gamma OLS}, the proof
of the proposition becomes simple.

\noindent \textit{Proof} \textit{of Proposition \ref{Proposition OLS Asymtptoc Distribu}.}
By standard arguments, 
\begin{align*}
\sqrt{T}\begin{bmatrix}\widehat{\beta}-\beta^{0}\\
\widehat{\delta}-\delta_{Z}^{0}
\end{bmatrix} & =\begin{bmatrix}X'X & X'\widehat{Z}_{0}\\
\widehat{Z}_{0}'X & \widehat{Z}_{0}'\widehat{Z}_{0}
\end{bmatrix}^{-1}\sqrt{T}\begin{bmatrix}X'e+X'\left(Z_{0}-\widehat{Z}_{0}\right)\delta_{Z}^{0}\\
\widehat{Z}_{0}'e+\widehat{Z}_{0}'\left(Z_{0}-\widehat{Z}_{0}\right)\delta_{Z}^{0}
\end{bmatrix},
\end{align*}
from which it follows that 
\begin{align*}
\begin{bmatrix}X'X & X'\widehat{Z}_{0}\\
\widehat{Z}_{0}'X & \widehat{Z}_{0}'\widehat{Z}_{0}
\end{bmatrix}^{-1}\frac{1}{h^{1/2}}X'\left(Z_{0}-\widehat{Z}_{0}\right)\delta_{Z}^{0} & =O_{p}\left(1\right)o_{p}\left(1\right)=o_{p}\left(1\right),
\end{align*}
 and a similar reasoning applies to $\widehat{Z}_{0}'(Z_{0}-\widehat{Z}_{0})\delta_{Z}^{0}.$
All other terms involving $\widehat{Z}_{0}$ can be treated in analogous
fashion. In particular, the $O_{p}\left(1\right)$ result above follows
from Lemma \ref{Lemma LLNs 1}-\ref{Lemma LLNs 2}. The rest of the
arguments (including mixed normality) follows from \citeReferencesSupp{barndorff/shephard:04}
and are omitted. $\square$

\subsubsection{Proof of Proposition \ref{Prop 3 Asym}}

\textit{\noindent Proof of part (i) of Proposition \ref{Prop 3 Asym}.}
Below $C$ is a generic positive constant which may change from line
to line. Let $\widetilde{e}$ denote the vector of normalized residuals
$\widetilde{e}_{t}$ defined by \eqref{Eq. eps WN}. Recall that $\widehat{T}_{b}=\arg\max_{T_{b}}Q_{T}\left(T_{b}\right)$,
$Q_{T}(\widehat{T}_{b})=\widehat{\delta}'_{T_{b}}\left(Z_{2}'MZ_{2}\right)\widehat{\delta}_{T_{b}}$,
and the decomposition 
\begin{align}
Q_{T}\left(T_{b}\right)-Q_{T}\left(T_{b}^{0}\right) & =\widehat{\delta}'_{T_{b}}\left(Z_{2}'MZ_{2}\right)\widehat{\delta}_{T_{b}}-\widehat{\delta}'_{T_{b}^{0}}\left(Z_{0}'MZ_{0}\right)\widehat{\delta}_{T_{b}^{0}}\label{eq. A.2.0, Prop 3}\\
 & =\delta_{h}'\left\{ \left(Z'_{0}MZ_{2}\right)\left(Z'_{2}MZ_{2}\right)^{-1}\left(Z'_{2}MZ_{0}\right)-Z_{0}'MZ_{0}\right\} \delta_{h}\label{eq. A.2.1, Prop 3}\\
 & \quad+g_{e}\left(T_{b}\right),\label{eq. A.2.2 Prop 3}
\end{align}
 where
\begin{align}
g_{e}\left(T_{b}\right) & =2\delta_{h}'\left(Z'_{0}MZ_{2}\right)\left(Z'_{2}MZ_{2}\right)^{-1}Z_{2}Me-2\delta_{h}'\left(Z'_{0}Me\right)\label{eq. A.2.3, Prop 3}\\
 & \quad+e'MZ_{2}\left(Z'_{2}MZ_{2}\right)^{-1}Z_{2}Me-e'MZ_{0}\left(Z'_{0}MZ_{0}\right)^{-1}Z'_{0}Me.\label{eq. A.2.4, Prop 3}
\end{align}
Since $g_{e}(\widehat{T}_{b})\geq|\widehat{T}_{b}-T_{b}^{0}|r(\widehat{T}_{b})$,
we have
\begin{align}
P & \left(\left|\widehat{\lambda}_{b}-\lambda_{0}\right|>K\right)\nonumber \\
 & =P\left(\left|\widehat{T}_{b}-T_{b}^{0}\right|>TK\right)\nonumber \\
 & \leq P\left(\sup_{\left|T_{b}-T_{b}^{0}\right|>TK}h^{-1/2}\left|g_{e}\left(T_{b}\right)\right|\geq\inf_{\left|T_{b}-T_{b}^{0}\right|>TK}h^{-1/2}\left|T_{b}-T_{b}^{0}\right|r\left(T_{b}\right)\right)\nonumber \\
 & \leq P\left(\sup_{p\leq T_{b}\leq T-p}h^{-1/2}\left|g_{e}\left(T_{b}\right)\right|\geq TK\inf_{\left|T_{b}-T_{b}^{0}\right|>TK}h^{-1/2}r\left(T_{b}\right)\right)\nonumber \\
 & =P\left(r_{T}^{-1}\sup_{p\leq T_{b}\leq T-p}h^{-1/2}\left|g_{e}\left(T_{b}\right)\right|\geq K\right),\label{eq. last line of A.2.6. Prop 3}
\end{align}
where $r_{T}=T\inf_{\left|T_{b}-T_{b}^{0}\right|>TK}h^{-1/2}r\left(T_{b}\right)$,
which is positive and bounded away from zero by Lemma \ref{Lemma r_CT, Asy, Prop 3}.
Thus, it is sufficient to verify that
\begin{align}
\sup_{p\leq T_{b}\leq T-p}h^{-1/2}\left|g_{e}\left(T_{b}\right)\right| & =o_{p}\left(1\right).\label{eq. A.2.7, Casini (2015a)-2-2}
\end{align}
 Consider the first term of $g_{e}\left(T_{b}\right)$:
\begin{align}
2\delta_{h}'\left(Z'_{0}MZ_{2}\right) & \left(Z'_{2}MZ_{2}\right)^{-1/2}\left(Z'_{2}MZ_{2}\right)^{-1/2}Z_{2}Me\label{eq. A.2.8, Prop 3}\\
 & \leq2h^{1/4}\left(\delta^{0}\right)'\left(Z'_{0}MZ_{2}\right)\left(Z'_{2}MZ_{2}\right)^{-1/2}\left(Z'_{2}MZ_{2}\right)^{-1/2}Z_{2}Me.\nonumber 
\end{align}
For any $1\leq j\leq p,$ $\left(Z_{2}\widetilde{e}\right)_{j,1}/\sqrt{h}=O_{p}\left(1\right)$
by Theorem \ref{=00005BJ=000026P,-Theorem-5.4.2=00005D}, and similarly,
for any $1\leq i\leq q+p,$ $\left(X\widetilde{e}\right)_{i}/\sqrt{h}=O_{p}\left(1\right)$.
Furthermore, from Lemma \ref{Lemma LLNs 1} we also have that $Z'_{2}MZ_{2}$
and $Z'_{0}MZ_{2}$ are $O_{p}\left(1\right).$ Therefore, the supremum
of $\left(Z'_{0}MZ_{2}\right)\left(Z'_{2}MZ_{2}\right)^{-1/2}$ over
all $T_{b}$ is such that
\begin{align*}
\sup_{T_{b}}\left(Z'_{0}MZ_{2}\right)\left(Z'_{2}MZ_{2}\right)^{-1}\left(Z'_{2}MZ_{0}\right) & \leq Z'_{0}MZ_{0}=O_{p}\left(1\right),
\end{align*}
by Lemma \ref{Lemma LLNs 1}. By Assumption \ref{Assumption 1, CT}-(iii),
$\left(Z'_{2}MZ_{2}\right)^{-1/2}Z_{2}M\widetilde{e}$ is $O_{p}\left(1\right)O_{p}(\sqrt{h})$
uniformly, which implies that \eqref{eq. A.2.8, Prop 3} is $O_{p}(\sqrt{h})$
uniformly over $p\leq T_{b}\leq T-p$. In view of Assumption \ref{Assumption 6 - Small Shifts}
{[}recall \eqref{Eq. eps WN}{]}, we need to study the behavior of
$\left(X'e\right)_{j,1}$ for $1\leq j\leq p+q.$ Note first that
$|\widehat{\lambda}_{b}-\lambda_{0}|>K$ or $N>|\widehat{N}_{b}-N_{b}^{0}|>KN$.
Then, by It\^o formula, proceeding as in the proof of Lemma \ref{Lemma, Spot Uniform Approx},
we have a standard result for the local volatility of a continuous
It\^o semimartingale; namely that for some $A>0$ (recall the condition
$T^{1-\kappa}\epsilon\rightarrow B>0$), 
\begin{align*}
\left\Vert \mathbb{E}\left(\frac{1}{\epsilon}\sum_{T_{b}^{0}-\left\lfloor T^{\kappa}\right\rfloor }^{T_{b}^{0}}x_{kh}\widetilde{e}_{kh}-\frac{1}{\epsilon}\int_{N_{b}^{0}-\epsilon}^{N_{b}^{0}}\Sigma_{Xe,s}ds|\,\mathscr{F}_{\left(T_{b}^{0}-1\right)h}\right)\right\Vert  & \leq Ah^{1/2}.
\end{align*}
From Assumption \ref{Assumption 1, CT}-(iv), since $\Sigma_{Xe,t}=0$
for all $t\geq0$, we have
\begin{align}
X'e & =\sum_{k=1}^{T_{b}^{0}-\left\lfloor T^{\kappa}\right\rfloor }x_{kh}\widetilde{e}_{kh}+h^{-1/4}\sum_{k=T_{b}^{0}-\left\lfloor T^{\kappa}\right\rfloor +1}^{T_{b}^{0}+\left\lfloor T^{\kappa}\right\rfloor }x_{kh}\widetilde{e}_{kh}+\sum_{k=T_{b}^{0}+\left\lfloor T^{\kappa}\right\rfloor +1}^{T}x_{kh}\widetilde{e}_{kh}\nonumber \\
 & =O_{p}\left(h^{1/2}\right)+h^{-1/4}O_{p}\left(h^{1-\kappa+1/2}\right)+O_{p}\left(h^{1/2}\right)=O_{p}\left(h^{1/2}\right).\label{Eq. Xej spot vol eps}
\end{align}
The same bound applies to $Z_{2}'e$ and $Z_{0}'e$. Thus, \eqref{eq. A.2.8, Prop 3}
is such that 
\begin{align*}
2h^{-1/2}h^{1/4}\left(\delta^{0}\right)' & \left(Z'_{0}MZ_{2}\right)\left(Z'_{2}MZ_{2}\right)^{-1/2}\left(Z'_{2}MZ_{2}\right)^{-1/2}Z_{2}Me\\
 & =2h^{-1/2}h^{1/4}\left\Vert \delta^{0}\right\Vert O_{p}\left(1\right)O_{p}\left(h^{1/2}\right)=O_{p}\left(1\right)O_{p}\left(h^{1/4}\right).
\end{align*}
As for the second term of \eqref{eq. A.2.3, Prop 3},
\begin{align*}
h^{-1/2}\delta_{h}'\left(Z'_{0}Me\right) & =2h^{-1/4}\left(\delta^{0}\right)'\left(Z'_{0}Me\right)=Ch^{-1/4}O_{p}\left(h^{1/2}\right)=CO_{p}\left(h^{1/4}\right),
\end{align*}
 using \eqref{Eq. Xej spot vol eps}. Again using \eqref{Eq. Xej spot vol eps},
the first term in \eqref{eq. A.2.4, Prop 3} is such that, uniformly
in $T_{b}$, 
\begin{align}
h^{-1/2} & e'MZ_{2}\left(Z'_{2}MZ_{2}\right)^{-1}Z_{2}Me\label{Eq. AAA5}\\
 & =h^{-1/2}BO_{p}\left(h^{1/2}\right)O_{p}\left(1\right)O_{p}\left(h^{1/2}\right)=O_{p}\left(h^{1/2}\right).\nonumber 
\end{align}
Similarly, the last term in \eqref{eq. A.2.4, Prop 3} is $O_{p}(\sqrt{h}).$
Therefore, combining these results, we have $h^{-1/2}\sup_{T_{b}}|g_{e}\left(T_{b}\right)|$
$=BO_{p}(h^{1/4})$, from which it follows that the right-hand side
of \eqref{eq. last line of A.2.6. Prop 3} is weakly smaller than
$\varepsilon.$ 
\begin{lem}
\label{Lemma r_CT, Asy, Prop 3}For $B>0$, let $r_{B,h}=\inf_{\left|T_{b}-T_{b}^{0}\right|>TB}Th^{-1/2}r\left(T_{b}\right).$
There exists an $A>0$ such that for every $\varepsilon>0,$ there
exists a $B<\infty$ such that $P\left(r_{B,h}\geq A\right)\leq1-\varepsilon$. 
\end{lem}

\noindent \textit{Proof.} Assume $N_{b}\leq N_{b}^{0},$ and observe
that $r_{T}\geq r_{B,h}$ for an appropriately chosen $B.$ From the
first inequality result in Lemma \ref{Lemma A1, Casini (2015a)},
 
\begin{align*}
Th^{-1/2}r\left(T_{b}\right) & \geq Th^{-1/2}h^{1/2}\left(\delta^{0}\right)'R'\frac{X'_{\Delta}X_{\Delta}}{T_{b}^{0}-T_{b}}\left(X'_{2}X_{2}\right)^{-1}\left(X'_{0}X_{0}\right)R\delta^{0}\\
 & =\left(\delta^{0}\right)'R'\left(X'_{\Delta}X_{\Delta}/\left(N_{b}^{0}-N_{b}\right)\right)\left(X'_{2}X_{2}\right)^{-1}\left(X'_{0}X_{0}\right)R\delta^{0}.
\end{align*}
Note that $B<h\left(T_{b}^{0}-T_{b}\right)<N$. Then 
\begin{align*}
Th^{-1/2}r\left(T_{b}\right) & \geq\left(\delta^{0}\right)'R'\left(X'_{\Delta}X_{\Delta}/N\right)\left(X'_{2}X_{2}\right)^{-1}\left(X'_{0}X_{0}\right)R\delta^{0}>A,
\end{align*}
 by the same argument as in Lemma \ref{Lemma r_CT, Prop 1-2}. Following
the same reasoning as in the proof of Lemma \ref{Lemma r_CT, Prop 1-2},
we can choose a $B>0$ such that $r_{B,h}=\inf_{\left|T_{b}-T_{b}^{0}\right|>TB}Th^{-1/2}r\left(T_{b}\right)$
satisfies $P(r_{B,h}\geq A)\leq1-\varepsilon.$ $\square$\textit{ }

\textit{\noindent Proof of part (ii) of Proposition \ref{Prop 3 Asym}.}
Suppose $T_{b}<T_{b}^{0}$. Let
\begin{align*}
\mathbf{D}_{K,T} & =\left\{ T_{b}:\,N\eta\leq N_{b}\leq N\left(1-\eta\right),\,\left|N_{b}-N_{b}^{0}\right|>K\left(T^{1-\kappa}\right)^{-1}\right\} .
\end{align*}
It is enough to show that $P(\sup_{T_{b}\in\mathbf{D}_{K,T}}Q_{T}\left(T_{b}\right)\geq Q_{T}(T_{b}^{0}))<\varepsilon.$
The difficulty is again to control the estimates that depend on $\left|N_{b}-N_{b}^{0}\right|$.
We shall show that
\begin{align*}
P & \left(\sup_{T_{b}\in\mathbf{D}_{K,T}}h^{-3/2}\frac{g_{e}\left(T_{b},\,\delta_{h}\right)}{\left|T_{b}-T_{b}^{0}\right|}\geq\inf_{T_{b}\in\mathbf{D}_{K,T}}h^{-3/2}r\left(T_{b}\right)\right)<\varepsilon.
\end{align*}
By Lemma \ref{Lemma A1, Casini (2015a)},
\begin{align*}
\inf_{T_{b}\in\mathbf{D}_{K,T}}r\left(T_{b}\right) & \geq\inf_{T_{b}\in\mathbf{D}_{K,T}}\delta'_{h}R'\frac{X'_{\Delta}X_{\Delta}}{T_{b}^{0}-T_{b}}\left(X'_{2}X_{2}\right)^{-1}\left(X'_{0}X_{0}\right)R\delta_{h},
\end{align*}
 and, since $\left|T_{b}-T_{b}^{0}\right|>KT^{\kappa},$ it is important
to consider $X'_{\Delta}X_{\Delta}=\sum_{k=T_{b}+1}^{T_{b}^{0}}x_{kh}x'_{kh}.$
We shall apply asymptotic results for the local approximation of the
covariation between processes. Consider 
\begin{align*}
\frac{X'_{\Delta}X_{\Delta}}{h\left(T_{b}^{0}-T_{b}\right)}=\frac{1}{h\left(T_{b}^{0}-T_{b}\right)}\sum_{k=T_{b}+1}^{T_{b}^{0}}x_{kh}x'_{kh} & .
\end{align*}
 By Theorem 9.3.2-(i) in \citeReferencesSupp{jacod/protter:12}, as
$h\downarrow0$
\begin{align}
\frac{1}{h\left(T_{b}^{0}-T_{b}\right)}\sum_{k=T_{b}+1}^{T_{b}^{0}}x_{kh}x'_{kh} & \overset{P}{\rightarrow}\Sigma_{XX,N_{b}^{0}},\label{Eq (1), Prop. 2 Casini (2015a) - Spot Estimator-2}
\end{align}
 since $\left|N_{b}-N_{b}^{0}\right|$ shrinks at a rate no faster
than $Kh^{1-\kappa}$ and $1/Kh^{1-\kappa}\rightarrow\infty$. By
Lemma \ref{Lemma, Spot Uniform Approx} this approximation is uniform,
establishing that
\begin{align*}
h^{-3/2} & \inf_{T_{b}\in\mathbf{D}_{K,T}}\left(\delta_{h}\right)'R'\frac{X'_{\Delta}X_{\Delta}}{T_{b}^{0}-T_{b}}\left(X'_{2}X_{2}\right)^{-1}\left(X'_{0}X_{0}\right)R\delta_{h}\\
 & =\inf_{T_{b}\in\mathbf{D}_{K,T}}\left(\delta^{0}\right)'R'\frac{X'_{\Delta}X_{\Delta}}{h\left(T_{b}^{0}-T_{b}\right)}\left(X'_{2}X_{2}\right)^{-1}\left(X'_{0}X_{0}\right)R\delta^{0},
\end{align*}
is bounded away from zero. Thus, it is sufficient to show that
\begin{align}
P\left(\sup_{T_{b}\in\mathbf{D}_{K,T}}h^{-3/2}\frac{g_{e}\left(T_{b},\,\delta_{h}\right)}{\left|T_{b}-T_{b}^{0}\right|}\geq B\right) & <\varepsilon,\label{eq. (17), Casini (2015a)-1-2}
\end{align}
 for some $B>0.$ Consider the terms of $g_{e}\left(T_{b}\right)$
in \eqref{eq. A.2.4, Prop 3}. Using $Z_{2}=Z_{0}\pm Z_{\Delta}$,
we deduce for the first term that
\begin{align}
\delta_{h}' & \left(Z'_{0}MZ_{2}\right)\left(Z'_{2}MZ_{2}\right)^{-1}Z_{2}Me\nonumber \\
 & =\delta_{h}'\left(\left(Z'_{2}\pm Z_{\Delta}\right)MZ_{2}\right)\left(Z'_{2}MZ_{2}\right)^{-1}Z_{2}Me\nonumber \\
 & =\delta_{h}'Z'_{0}Me\pm\delta_{h}'Z'_{\Delta}Me\pm\delta_{h}'\left(Z'_{\Delta}MZ_{2}\right)\left(Z'_{2}MZ_{2}\right)^{-1}Z_{2}Me.\label{eq. (18), Casini, Prop. 3}
\end{align}
First, we can apply Lemma \ref{Lemma LLNs 1} {[}(vi)-(viii){]}, together
with Assumption \ref{Assumption 1, CT}-(iii), to terms not involving
$\left|N_{b}-N_{b}^{0}\right|$. The third term is such that 
\begin{align}
K^{-1}h^{-\left(1-\kappa\right)}\left(Z'_{\Delta}MZ_{2}\right) & =\frac{Z'_{\Delta}Z_{2}}{Kh^{1-\kappa}}-\frac{Z'_{\Delta}X_{\Delta}}{Kh^{1-\kappa}}\left(X'X\right)^{-1}X'Z_{2}.\label{Eq. AAA1}
\end{align}
Consider $Z'_{\Delta}Z_{\Delta}$ (the argument for $Z'_{\Delta}X_{\Delta}$
is analogous). By Lemma \ref{Lemma, Spot Uniform Approx}, $Z'_{\Delta}Z_{\Delta}/Kh^{1-\kappa}$
uniformly approximates the moving average of $\Sigma_{ZZ,t}$ over
$(N_{b}^{0}-KT^{\kappa}h,\,N_{b}^{0}]$. Hence, as $h\downarrow0$,
\begin{align}
Z'_{\Delta}Z_{\Delta}/Kh^{1-\kappa} & =BO_{p}\left(1\right),\label{eq Zdelta X Pro 3}
\end{align}
 for some $B>0,$ uniformly in $T_{b}$. The second term in \eqref{Eq. AAA1}
is thus also $O_{p}\left(1\right)$ uniformly, using Lemma \ref{Lemma LLNs 1}.
Then, using \eqref{Eq. Xej spot vol eps} and \eqref{Eq. AAA1} into
the third term of \eqref{eq. (18), Casini, Prop. 3}, we have 
\begin{align}
\frac{1}{K} & h^{-\left(1-\kappa\right)-1/2}\left(\delta_{h}\right)'\left(Z'_{\Delta}MZ_{2}\right)\left(Z'_{2}MZ_{2}\right)^{-1}Z_{2}Me\label{AAA2, Prop}\\
 & \leq\frac{1}{K}h^{-1/4}\left(\delta^{0}\right)'\left(\frac{Z'_{\Delta}MZ_{2}}{h^{1-\kappa}}\right)\left(Z'_{2}MZ_{2}\right)^{-1}Z_{2}Me\nonumber \\
 & \leq h^{-1/4}\frac{Z'_{\Delta}MZ_{2}}{Kh^{1-\kappa}}O_{p}\left(1\right)O_{p}\left(h^{1/2}\right)\leq O_{p}\left(h^{1/4}\right),\nonumber 
\end{align}
 where $\left(Z'_{2}MZ_{2}\right)^{-1}=O_{p}\left(1\right)$. So the
right-and side of \eqref{AAA2, Prop} is less than $\varepsilon/4$
in probability. Therefore, for the second term of \eqref{eq. (18), Casini, Prop. 3},
\begin{align}
K^{-1} & h^{-\left(1-\kappa\right)-1/2}\delta_{h}'Z'_{\Delta}Me\nonumber \\
 & =\frac{h^{-1/2}}{Kh^{1-\kappa}}\delta'_{h}\sum_{k=T_{b}+1}^{T_{b}^{0}}z_{kh}e_{kh}-\frac{h^{-1/2}}{h^{1-\kappa}}\delta'_{h}\left(\sum_{k=T_{b}+1}^{T_{b}^{0}}z_{kh}x'_{kh}\right)\left(X'X\right)^{-1}\left(X'e\right)\nonumber \\
 & \leq\frac{h^{-1/2}}{Kh^{1-\kappa}}\delta'_{h}\sum_{k=T_{b}+1}^{T_{b}^{0}}z_{kh}e_{kh}-B\frac{1}{K}\frac{h^{-1/4}}{h^{1-\kappa}}\left(\delta^{0}\right)\left(\sum_{k=T_{b}+1}^{T_{b}^{0}}z_{kh}x'_{kh}\right)\left(X'X\right)^{-1}\left(X'e\right)\nonumber \\
 & \leq\frac{h^{-1/2}}{Kh^{1-\kappa}}\delta'_{h}\sum_{k=T_{b}+1}^{T_{b}^{0}}z_{kh}e_{kh}-h^{-1/4}O_{p}\left(1\right)O_{p}\left(h^{1/2}\right).\label{Eq. AAA3, Prop}
\end{align}
 Thus, using \eqref{eq. (18), Casini, Prop. 3}, \eqref{eq. A.2.3, Prop 3}
is such that
\begin{align*}
2\delta_{h}'Z'_{0}M & e\pm2\delta_{h}'Z'_{\Delta}Me\pm2\delta_{h}'\left(Z'_{\Delta}MZ_{2}\right)\left(Z'_{2}MZ_{2}\right)^{-1}Z_{2}Me-2\delta_{h}'\left(Z'_{0}Me\right)\\
 & =2\delta_{h}'Z'_{\Delta}Me\pm2\delta_{h}'\left(Z'_{\Delta}MZ_{2}\right)\left(Z'_{2}MZ_{2}\right)^{-1}Z_{2}Me\\
 & \leq\frac{h^{-1/2}}{Kh^{1-\kappa}}\left(\delta^{0}\right)'\sum_{k=T_{b}+1}^{T_{b}^{0}}z_{kh}\widetilde{e}_{kh}-h^{-1/4}O_{p}\left(1\right)O_{p}\left(h^{1/2}\right)+O_{p}\left(h^{-1/4}\right),
\end{align*}
in view of \eqref{AAA2, Prop} and \eqref{Eq. AAA3, Prop}. Next,
consider \eqref{eq. A.2.4, Prop 3}. We can use the decomposition
$Z_{2}=Z_{0}\pm Z_{\Delta}$ and show that all terms involving the
matrix $Z_{\Delta}$ are negligible. To see this, consider the first
term when multiplied by $K^{-1}h^{-\left(3/2-\kappa\right)}$, then
\begin{align}
K^{-1} & h^{-\left(3/2-\kappa\right)}e'MZ_{2}\left(Z'_{2}MZ_{2}\right)^{-1}Z_{2}Me\label{Eq. AAA6, Prop}\\
 & =K^{-1}h^{-\left(3/2-\kappa\right)}e'MZ_{0}\left(Z'_{2}MZ_{2}\right)^{-1}Z_{2}Me\nonumber \\
 & \quad\pm K^{-1}h^{-\left(3/2-\kappa\right)}e'MZ_{\Delta}\left(Z'_{2}MZ_{2}\right)^{-1}Z_{2}Me.\nonumber 
\end{align}
By the same argument as in \eqref{Eq. Xej spot vol eps}, $Z'_{2}Me=O_{p}(h^{1/2})$.
Using the Burkh\"{o}lder-Davis-Gundy inequality, the estimates for
the local volatility of continuous It\^o semimartingales  yield
\begin{align*}
\widetilde{e}'MZ_{\Delta} & =\widetilde{e}'Z_{\Delta}-\widetilde{e}'X\left(X'X\right)^{-1}X'Z_{\Delta}\\
 & =O_{p}\left(Kh^{1/2+1-\kappa}\right)-O_{p}\left(h^{1/2}\right)O_{p}\left(1\right)O_{p}\left(Kh^{1-\kappa}\right).
\end{align*}
Thus, the second term in \eqref{Eq. AAA6, Prop} is such that 
\begin{align}
K^{-1}h^{-\left(3/2-\kappa\right)} & \widetilde{e}'MZ_{\Delta}\left(Z'_{2}MZ_{2}\right)^{-1}Z_{2}Me\label{Eq. AAA7, Prop}\\
 & =B\left(K^{-1}h^{-\left(3/2-\kappa\right)}\right)O_{p}\left(Kh^{1-\kappa+1/2}\right)O_{p}\left(1\right)O_{p}\left(h^{1/2}\right)\nonumber \\
 & =BO_{p}\left(h^{1/2}\right).\nonumber 
\end{align}
 Next, let us consider \eqref{eq. A.2.4, Prop 3}. The key here is
to recognize that, on $\mathbf{D}_{K,T}$, $T_{b}$ and $T_{b}^{0}$
lies on the same window with right-hand point $N_{b}^{0}$. Thus the
difference between the two terms in \eqref{eq. A.2.4, Prop 3} is
asymptotically negligible. First, note that using \eqref{Eq. Xej spot vol eps},
\begin{align*}
\widetilde{e}'MZ_{0}\left(Z'_{0}MZ_{0}\right)^{-1}Z_{0}M\widetilde{e} & =O_{p}\left(h^{1/2}\right)O_{p}\left(1\right)O_{p}\left(h^{1/2}\right)=O_{p}\left(h\right).
\end{align*}
Applying $Z_{0}=Z{}_{2}\pm Z{}_{\Delta}$ repeatedly in \eqref{Eq. AAA6, Prop},
and noting that the cross-product terms involving $Z_{\Delta}$ are
$o_{p}\left(1\right)$ by the same reasoning as in \eqref{Eq. AAA7, Prop},
we obtain that the difference between the first and second term of
\eqref{eq. A.2.4, Prop 3} is negligible. The more intricate step
is the one arising from
\begin{align*}
e'MZ_{0} & \left(Z'_{0}MZ_{2}\pm Z'_{\Delta}MZ_{2}\right)^{-1}Z'_{0}Me-e'MZ_{0}\left(Z'_{0}MZ_{0}\right)^{-1}Z'_{0}Me\\
 & =e'MZ_{0}\left[\left(Z'_{0}MZ_{2}\pm Z'_{\Delta}MZ_{2}\right)^{-1}-\left(Z'_{0}MZ_{0}\right)^{-1}\right]Z'_{0}Me.
\end{align*}
On $\mathbf{D}_{K,T},$ $\left|N_{b}-N_{b}^{0}\right|=O_{p}(Kh^{1-\kappa})$,
and so each term involving $Z_{\Delta}$ is of higher order. By using
the continuity of probability limits, the matrix in square brackets
goes to zero at rate $h^{1-\kappa}$. Then, this expression when multiplied
by $h^{-\left(3/2-\kappa\right)}K^{-1}$, and after using the same
rearrangements as above, can be shown to satisfy {[}recall also \eqref{Eq. Xej spot vol eps}{]}
\begin{align*}
h^{-\left(3/2-\kappa\right)}K^{-1} & e'MZ_{0}\left[\left(Z'_{0}MZ_{2}\pm Z'_{\Delta}MZ_{2}\right)^{-1}-\left(Z'_{0}MZ_{0}\right)^{-1}\right]Z'_{0}Me\\
 & =h^{-\left(3/2-\kappa\right)}K^{-1}O_{p}\left(h\right)\left[\left(Z'_{0}MZ_{2}\pm Z'_{\Delta}MZ_{2}\right)^{-1}-\left(Z'_{0}MZ_{0}\right)^{-1}\right]\\
 & =h^{-\left(3/2-\kappa\right)}K^{-1}O_{p}\left(h\right)\\
 & \quad\times\left[\left(Z'_{0}MZ_{0}\pm Z'_{0}MZ'_{\Delta}\pm Z'_{\Delta}MZ_{2}\right)^{-1}-\left(Z'_{0}MZ_{0}\right)^{-1}\right]\\
 & =h^{-\left(3/2-\kappa\right)}K^{-1}O_{p}\left(h\right)o_{p}\left(h^{1-\kappa}\right)=O_{p}\left(h^{1/2}\right)o_{p}\left(1\right).
\end{align*}
 Therefore, \eqref{eq. A.2.4, Prop 3} is stochastically small uniformly
in $T_{b}\in\mathbf{D}_{K,T}$ when $T$ is large. Altogether, we
have 
\begin{align*}
h^{-1/2}\frac{g_{e}\left(T_{b}\right)}{\left|T_{b}-T_{b}^{0}\right|} & \leq2\frac{h^{-1/2}}{Kh^{1-\kappa}}\delta'_{h}\sum_{k=T_{b}+1}^{T_{b}^{0}}z_{kh}e_{kh}\\
 & \quad-h^{-1/4}O_{p}\left(1\right)O_{p}\left(h^{1/2}\right)+O_{p}\left(h^{-1/4}\right).
\end{align*}
Thus, it remains to find a bound for the first term above. By It\^o's
formula, standard estimates for the local volatility of continuous
It\^o semimartingales yield for every $T_{b},$ 
\begin{align}
\mathbb{E}\left(\left\Vert \widehat{\Sigma}_{Ze}\left(T_{2},\,T_{b}^{0}\right)-\overline{\Sigma}_{Ze}\left(T_{2},\,T_{b}^{0}\right)\right\Vert |\,\mathscr{F}_{T_{b}h}\right) & \leq Bh^{1/2},\label{eq (Spot ineq)}
\end{align}
 for some $B>0.$ Let $R_{1,h}=\sum_{k=T_{b}^{0}-\left(B+1\right)\left\lfloor T^{\kappa}\right\rfloor +1}^{T_{b}^{0}}z_{kh}\widetilde{e}_{kh}$,
$R_{2,h}\left(T_{b}\right)=\sum_{k=T_{b}+1}^{T_{2}^{0}-\left(B+1\right)\left\lfloor T^{\kappa}\right\rfloor }z_{kh}e_{kh}$
and note that $\sum_{k=T_{2}+1}^{T_{2}^{0}}z_{kh}e_{kh}=R_{1,h}+R_{2,h}\left(T_{b}\right)$.
Then, for any $C>0$, 
\begin{align}
P & \left(\sup_{T_{b}<T_{b}^{0}-KT^{\kappa}}2\frac{h^{-1/2}}{Kh^{1-\kappa}}\delta'_{h}\left\Vert \sum_{k=T_{b}+1}^{T_{b}^{0}}z_{kh}e_{kh}\right\Vert \geq C\right)\label{eq. bound on Hayek-Reini ineq}\\
 & =P\left(\sup_{T_{b}<T_{b}^{0}-KT^{\kappa}}\frac{h^{-1/2}}{Kh^{1-\kappa}}\delta'_{h}\left\Vert R_{1,h}+R_{2,h}\left(T_{b}\right)\right\Vert \geq2^{-1}C\right)\nonumber \\
 & \leq P\left(\frac{1}{Kh^{1-\kappa}}\left\Vert R_{1,h}\right\Vert >4^{-1}C\left\Vert \delta^{0}\right\Vert ^{-1}h^{1/2}\right)\nonumber \\
 & \quad+P\left(\sup_{T_{b}<T_{b}^{0}-KT^{\kappa}}\frac{K^{-1}}{h^{1-\kappa}}\left\Vert R_{2,h}\left(T_{b}\right)\right\Vert >4^{-1}C\left\Vert \delta^{0}\right\Vert ^{-1}h^{1/4}\right).\nonumber 
\end{align}
Consider first the second probability. By Markov's inequality, 
\begin{align*}
P & \left(\sup_{T_{b}<T_{b}^{0}-KT^{\kappa}}\frac{1}{Kh^{1-\kappa}}\left\Vert R_{2,h}\left(T_{b}\right)\right\Vert >4^{-1}C\left\Vert \delta^{0}\right\Vert ^{-1}h^{1/4}\right)\\
 & \leq P\left(\sup_{T_{b}<T_{b}^{0}-KT^{\kappa}}\left\Vert \frac{1}{Kh^{1-\kappa}}R_{2,h}\left(T_{b}\right)\right\Vert >4^{-1}C\left\Vert \delta^{0}\right\Vert ^{-1}h^{1/4}\right)\\
 & \leq\left(K/B\right)T^{\kappa}P\left(\left\Vert \frac{1}{Kh^{1-\kappa}}R_{2,h}\left(T_{b}\right)\right\Vert >4^{-1}C\left\Vert \delta^{0}\right\Vert ^{-1}h^{1/4}\right)\\
 & \leq\frac{\left(4\left(B+1\right)\left\Vert \delta^{0}\right\Vert \right)^{r}}{C^{r}}h^{-r/4}\frac{K}{B}T^{\kappa}\mathbb{E}\left(\left|\frac{1}{\left(B+1\right)Kh^{1-\kappa}}\left\Vert R_{2,h}\left(T_{b}\right)\right\Vert \right|^{r}\right)\\
 & \leq C_{r}\left(B+1\right)B^{-1}\left\Vert \delta^{0}\right\Vert ^{r}h^{-r/4}T^{\kappa}h^{r/2}\leq C_{r}\left\Vert \delta^{0}\right\Vert ^{r}h^{r/2-\kappa-r/4}\rightarrow0,
\end{align*}
for a sufficiently large $r>0$. We now turn to $R_{1,h}$. We have,
\begin{align*}
P & \left(\frac{1}{Kh^{1-\kappa}}\left\Vert R_{1,h}\right\Vert >2^{-1}C\left\Vert \delta^{0}\right\Vert ^{-1}h^{1/2}\right)\\
 & \leq P\biggl(\frac{\left(B+1\right)}{K}\left\Vert \left(B+1\right)^{-1}h^{-\left(1-\kappa\right)}\sum_{k=T_{b}^{0}-\left(B+1\right)\left\lfloor T^{\kappa}\right\rfloor +1}^{T_{b}^{0}}z_{kh}\widetilde{e}_{kh}\right\Vert \\
 & \quad>\frac{C}{4}\left\Vert \delta^{0}\right\Vert ^{-1}h^{1/2}\biggr)\\
 & \leq P\left(\left(B+1\right)K^{-1}O_{\mathbb{P}}\left(1\right)>4^{-1}C\left\Vert \delta^{0}\right\Vert ^{-1}\right)\rightarrow0,
\end{align*}
 by choosing $K$ large enough where we have used \eqref{eq (Spot ineq)}.
Altogether, the right-hand side of \eqref{eq. bound on Hayek-Reini ineq}
is less than $\varepsilon$, which concludes the proof. $\square$

\noindent\textit{Proof of part (iii) of Proposition \ref{Prop 3 Asym}.}
Observe that Lemma \ref{Lemma gamma OLS} applies under this setting.
Then, we have, 
\begin{align*}
\sqrt{T}\begin{bmatrix}\widehat{\beta}-\beta_{0}\\
\widehat{\delta}-\delta_{h}
\end{bmatrix} & =\begin{bmatrix}X'X & X'\widehat{Z}_{0}\\
\widehat{Z}_{0}'X & \widehat{Z}_{0}'\widehat{Z}_{0}
\end{bmatrix}^{-1}\sqrt{T}\begin{bmatrix}X'e+X'\left(Z_{0}-\widehat{Z}_{0}\right)\delta_{h}\\
\widehat{Z}_{0}'e+\widehat{Z}_{0}'\left(Z_{0}-\widehat{Z}_{0}\right)\delta_{h}
\end{bmatrix},
\end{align*}
so that we need to show that
\begin{align*}
\begin{bmatrix}X'X & X'\widehat{Z}_{0}\\
\widehat{Z}_{0}'X & \widehat{Z}_{0}'\widehat{Z}_{0}
\end{bmatrix}^{-1}\frac{1}{h^{1/2}}X'\left(Z_{0}-\widehat{Z}_{0}\right)\delta_{h} & \overset{P}{\rightarrow}0,
\end{align*}
 and the limiting distribution of $X'e/h^{1/2}$ is Gaussian. The
first claim can be proved in a manner analogous to that in the proof
of Proposition \ref{Proposition OLS Asymtptoc Distribu}. For the
second claim, we have the following decomposition from \eqref{Eq. Xej spot vol eps},
\begin{align*}
X'e & =\sum_{k=1}^{T_{b}^{0}-\left\lfloor T^{\kappa}\right\rfloor }x_{kh}\widetilde{e}_{kh}+h^{-1/4}\sum_{T_{b}^{0}-\left\lfloor T^{\kappa}\right\rfloor +1}^{T_{b}^{0}+\left\lfloor T^{\kappa}\right\rfloor }x_{kh}\widetilde{e}_{kh}+\sum_{k=T_{b}^{0}+\left\lfloor T^{\kappa}\right\rfloor +1}^{T}x_{kh}\widetilde{e}_{kh}\\
 & \triangleq R_{1,h}+R_{2,h}+R_{3,h}.
\end{align*}
 By Theorem \ref{=00005BJ=000026P,-Theorem-5.4.2=00005D}, $h^{-1/2}R_{1,h}\overset{\mathcal{L}\textrm{-s}}{\rightarrow}\mathscr{MN}\left(0,\,V_{1}\right),$
where $V_{1}\triangleq\underset{T\rightarrow\infty}{\lim}T\sum_{k=1}^{T_{b}^{0}-\left\lfloor T^{\kappa}\right\rfloor }\mathbb{E}(x_{kh}x'_{kh}\widetilde{e}_{kh}^{2}).$
Similarly, $h^{-1/2}R_{3,h}\overset{\mathcal{L}\textrm{-s}}{\rightarrow}\mathscr{MN}\left(0,\,V_{3}\right),$
where $V_{3}\triangleq\underset{T\rightarrow\infty}{\lim}T\sum_{k=T_{b}^{0}+\left\lfloor T^{\kappa}\right\rfloor +1}^{T}\mathbb{E}(x_{kh}x'_{kh}\widetilde{e}_{kh}^{2}).$
If $\kappa\in\left(0,\,1/4\right),$ $\,$ $h^{-\left(1-\kappa\right)}$
$\sum_{T_{b}^{0}-\left\lfloor T^{\kappa}\right\rfloor +1}^{T_{b}^{0}+\left\lfloor T^{\kappa}\right\rfloor }x_{kh}\widetilde{e}_{kh}\overset{P}{\rightarrow}\Sigma_{Xe,N_{b}^{0}}$
by Theorem 9.3.2 in \citeReferencesSupp{jacod/protter:12} and so
$h^{-1/2}R_{2,h}=h^{-3/4}\sum_{T_{b}^{0}-\left\lfloor T^{\kappa}\right\rfloor }^{T_{b}^{0}+\left\lfloor T^{\kappa}\right\rfloor }x_{kh}\widetilde{e}_{kh}\overset{P}{\rightarrow}0.$
If $\kappa=1/4,$ then $h^{-1/2}R_{2,h}\rightarrow\Sigma_{Xe,N_{b}^{0}}$
in probability again by Theorem 9.3.2 in \citeReferencesSupp{jacod/protter:12}.
Since by Assumption \ref{Assumption 1, CT}-(iv) $\Sigma_{Xe,t}=0$
for all $t\geq0$, whenever $\kappa\in(0,\,1/4],$ $X'e/h^{1/2}$
is asymptotically normally distributed. The rest of the proof is simple
and follows the same steps as in Proposition \ref{Proposition OLS Asymtptoc Distribu}.
$\square$

\subsubsection{Proof of Proposition \ref{Prop. Slow Time Scale}}

\noindent\textit{Proof.} By Lemma \ref{Lemma 1},
\begin{align*}
Q_{T}\left(T_{b}\right)-Q_{T}\left(T_{b}^{0}\right) & =-\delta'_{h}\left(Z_{\Delta}'Z_{\Delta}\right)\delta_{h}\pm2\delta_{h}'\left(Z'_{\Delta}e\right)+o_{p}\left(h^{3/2-\kappa}\right).
\end{align*}
 Divide both sides by $h$ to yield, 
\begin{align*}
h^{-1}\left(Q_{T}\left(T_{b}\right)-Q_{T}\left(T_{b}^{0}\right)\right) & =-h^{1/2}\left(\delta^{0}\right)'\left(\frac{Z'_{\Delta}}{\sqrt{h}}\frac{Z_{\Delta}}{\sqrt{h}}\right)\delta^{0}\\
 & \quad\pm2\left(\delta^{0}\right)'\left(\frac{Z'_{\Delta}}{\sqrt{h}}\frac{\widetilde{e}}{\sqrt{h}}\right)+o_{p}\left(h^{1/2-\kappa}\right).
\end{align*}
 Note that $z_{kh}/\sqrt{h}\sim i.n.d.\,\mathscr{N}(0,\,\Sigma_{kh})$
and $\widetilde{e}_{kh}/\sqrt{h}\sim i.n.d.\,\mathscr{N}(0,\,\sigma_{e,kh}^{2})$.
Thus,
\begin{align*}
h^{-1+\kappa/2} & \left(Q_{T}\left(T_{b}\right)-Q_{T}\left(T_{b}^{0}\right)\right)\\
 & =-\frac{h^{1/2}}{\sqrt{T^{\kappa}}}\left(\delta^{0}\right)'\left(\frac{Z'_{\Delta}}{\sqrt{h}}\frac{Z_{\Delta}}{\sqrt{h}}\right)\delta^{0}\\
 & \quad\pm\frac{2}{\sqrt{T^{\kappa}}}\left(\delta^{0}\right)'\left(\frac{Z'_{\Delta}}{\sqrt{h}}\frac{\widetilde{e}}{\sqrt{h}}\right)+o_{p}\left(h^{1/2-\kappa/2}\right)\\
 & =O_{p}\left(h^{1/2}\right)\pm\frac{2}{\sqrt{T^{\kappa}}}\left(\delta^{0}\right)'\left(\frac{Z'_{\Delta}}{\sqrt{h}}\frac{\widetilde{e}}{\sqrt{h}}\right)+o_{p}\left(h^{1/2-\kappa}\right).
\end{align*}
 Also $T_{b}=T_{b}^{0}+\left\lfloor vT^{\kappa}\right\rfloor $, and
\begin{align*}
h^{-1+\kappa/2}\left(Q_{T}\left(T_{b}\right)-Q_{T}\left(T_{b}^{0}\right)\right) & \Rightarrow2\left(\delta^{0}\right)'\mathscr{W}\left(v\right).
\end{align*}
 The continuous mapping theorem then yields the desired result. $\square$

\subsubsection{Proof of Lemma \ref{Lemma 1}}

First, we begin with the following simple identity. Throughout the
proof, $B$ is a generic constant which may change from line to line.
\begin{lem}
\label{Lemma A.3}The following identity holds
\begin{align*}
\left(\delta_{h}\right)' & \left\{ Z_{0}'MZ_{0}-\left(Z'_{0}MZ_{2}\right)\left(Z'_{2}MZ_{2}\right)^{-1}\left(Z'_{2}MZ_{0}\right)\right\} \delta_{h}\\
 & =\left(\delta_{h}\right)'\left\{ Z_{\Delta}'MZ_{\Delta}-\left(Z'_{\Delta}MZ_{2}\right)\left(Z'_{2}MZ_{2}\right)^{-1}\left(Z'_{2}MZ_{\Delta}\right)\right\} \delta_{h}.
\end{align*}
\end{lem}
\noindent\textit{Proof.} The proof follows simply from the fact that
$Z'_{0}MZ_{2}=Z'_{2}MZ_{2}\pm Z'_{\Delta}MZ_{2}$ and so
\begin{align*}
\left(\delta_{h}\right)' & \left\{ Z_{0}'MZ_{0}-\left(Z'_{2}MZ_{2}\pm Z'_{\Delta}MZ_{2}\right)\left(Z'_{2}MZ_{2}\right)^{-1}\left(Z'_{2}MZ_{0}\right)\right\} \delta_{h}\\
 & =\left(\delta_{h}\right)'\{Z_{\Delta}'MZ_{0}-\left(Z'_{\Delta}MZ_{2}\right)\left(Z'_{2}MZ_{2}\right)^{-1}\left(Z'_{2}MZ_{2}\right)\\
 & \quad-\left(Z'_{\Delta}MZ_{2}\right)\left(Z'_{2}MZ_{2}\right)^{-1}\left(Z'_{2}MZ_{\Delta}\right)\}\delta_{h}\\
 & =\left(\delta_{h}\right)'\left\{ Z_{\Delta}'MZ_{\Delta}-\left(Z'_{\Delta}MZ_{2}\right)\left(Z'_{2}MZ_{2}\right)^{-1}\left(Z'_{2}MZ_{\Delta}\right)\right\} \delta_{h}.\,\square
\end{align*}

\noindent\textit{Proof of Lemma \ref{Lemma 1}.} By the definition
of $Q_{T}\left(T_{b}\right)-Q_{T}\left(T_{0}\right)$ and Lemma \ref{Lemma A.3},
\begin{align}
 & Q_{T}\left(T_{b}\right)-Q_{T}\left(T_{0}\right)\nonumber \\
 & =-\delta_{h}'\left\{ Z_{\Delta}'MZ_{\Delta}-\left(Z'_{\Delta}MZ_{2}\right)\left(Z'_{2}MZ_{2}\right)^{-1}\left(Z'_{2}MZ_{\Delta}\right)\right\} \delta_{h}+g_{e}\left(T_{b},\,\delta_{h}\right),\label{eq. (45), Proof of  Lemma 1}
\end{align}
where
\begin{align}
g_{e}\left(T_{b},\,\delta_{h}\right) & =2\delta_{h}'\left(Z'_{0}MZ_{2}\right)\left(Z'_{2}MZ_{2}\right)^{-1}Z_{2}Me-2\delta_{h}'\left(Z'_{0}Me\right)\label{Eq ge ACRSC3}\\
 & +e'MZ_{2}\left(Z'_{2}MZ_{2}\right)^{-1}Z_{2}Me-e'MZ_{0}\left(Z'_{0}MZ_{0}\right)^{-1}Z'_{0}Me.\label{Eq ge ACRSC4}
\end{align}
Recall that $N_{b}\left(u\right)\in\mathcal{D}\left(C\right)$ implies
$T_{b}\left(u\right)=T_{b}^{0}+uT^{\kappa},\,u\in\left[-C,\,C\right]$.
We consider the case $u\leq0.$ By Theorem 9.3.2-(i) in \citeReferencesSupp{jacod/protter:12}
combined with Lemma  \ref{Lemma, Spot Uniform Approx}, we have uniformly
in $u$ as $h\downarrow0$
\begin{align}
\frac{1}{h^{1-\kappa}}\sum_{k=T_{b}^{0}+uT^{\kappa}}^{T_{b}^{0}}x_{kh}x'_{kh} & \overset{P}{\rightarrow}\Sigma_{XX,N_{b}^{0}}.\label{Eq (1), Spot Estimator}
\end{align}
Since $Z_{\Delta}'X=Z_{\Delta}'X_{\Delta}$, we will use this result
also for $Z_{\Delta}'X/h^{1-\kappa}$. With the notation of Section
\ref{Subsection Additional-Notation Spot Continuous case} {[}recall
\eqref{Eq. spot one side Ze-1 H}{]}, by the Burkh\"{o}lder-Davis-Gundy
inequality, we have that standard estimates for the local volatility,
is such that 
\begin{align}
\left\Vert \mathbb{E}\left(\widehat{\Sigma}_{ZX}\left(T_{b},\,T_{b}^{0}\right)-\Sigma_{ZX,\left(T_{b}^{0}-1\right)h}|\,\mathscr{F}_{\left(T_{b}^{0}-1\right)h}\right)\right\Vert  & \leq Bh^{1/2}.\label{eq. Eq (1), Spot Estimator}
\end{align}
 Equation \eqref{Eq (1), Spot Estimator}-\eqref{eq. Eq (1), Spot Estimator}
can be used to yield, uniformly in $T_{b}$, 
\begin{align}
\psi_{h}^{-1}Z_{\Delta}'X\left(X'X\right)^{-1}X'Z_{\Delta} & =O_{p}\left(1\right)X'Z_{\Delta},\label{eq. ACRSC2b}
\end{align}
 and
\begin{align}
Z'_{\Delta}MZ_{2} & =Z_{\Delta}'Z_{\Delta}-Z_{\Delta}'X\left(X'X\right)^{-1}X'Z_{2}=O_{p}\left(\psi_{h}\right)-O_{p}\left(\psi_{h}\right)O_{p}\left(1\right)O_{p}\left(1\right).\label{eq. ACRSC2}
\end{align}
Now, expand the first term of \eqref{eq. (45), Proof of  Lemma 1},
\begin{align}
\delta_{h}'Z_{\Delta}'MZ_{\Delta}\delta_{h} & =\delta_{h}'Z_{\Delta}'Z_{\Delta}\delta_{h}-\delta_{h}'Z_{\Delta}'X\left(X'X\right)^{-1}X'Z_{\Delta}\delta_{h}.\label{eq. first term of eq. 45}
\end{align}
By Lemma \ref{Lemma LLNs 1}, $\left(X'X\right)^{-1}=O_{p}\left(1\right)$
and recall that $\delta_{h}=h^{1/4}\delta^{0}$. Then, 
\begin{align}
\psi_{h}^{-1}\delta_{h}'Z_{\Delta}'MZ_{\Delta}\delta_{h} & =\psi_{h}^{-1}\delta_{h}'Z_{\Delta}'Z_{\Delta}\delta_{h}-\psi_{h}^{-1}\delta_{h}'Z_{\Delta}'X\left(X'X\right)^{-1}X'Z_{\Delta}\delta_{h}.\label{Eq. ZMZ-1-1-1}
\end{align}
By \eqref{eq. ACRSC2b}, the second term above is such that
\begin{align}
\left\Vert \delta^{0}\right\Vert ^{2}h^{1/2}\frac{Z_{\Delta}'X}{\psi_{h}}\left(X'X\right)^{-1}X'Z_{\Delta} & =\left\Vert \delta^{0}\right\Vert ^{2}h^{1/2}O_{p}\left(1\right)X'Z_{\Delta},\label{eq. ACRSC5}
\end{align}
uniformly in $T_{b}\left(u\right)$. Therefore, 
\begin{align}
\psi_{h}^{-1}\delta_{h}'Z_{\Delta}'MZ_{\Delta}\delta_{h} & =\psi_{h}^{-1}\delta_{h}'Z_{\Delta}'Z_{\Delta}\delta_{h}-\left\Vert \delta^{0}\right\Vert ^{2}h^{1/2}O_{p}\left(1\right)O_{p}\left(\psi_{h}\right).\label{eq. ACRSC1}
\end{align}
In the last equality the second term of $\delta'Z_{\Delta}'MZ_{\Delta}\delta$
is always of higher order. This suggests that the term involving regressors
whose parameters are allowed to shift plays a primary role in the
asymptotic analysis. The second term is a complicated function of
cross products of all regressors around the time of the change. Because
of the fast rate of convergence, these high order product estimates
around the break date will be negligible. We use this result repeatedly
in the derivations that follow. The second term of \eqref{eq. (45), Proof of  Lemma 1}
when multiplied by $\psi_{h}^{-1}$ is, uniformly in $T_{b}\left(u\right)$,
\begin{align*}
\psi_{h}^{-1}\delta{}_{h} & \left(Z'_{\Delta}MZ_{2}\right)\left(Z'_{2}MZ_{2}\right)^{-1}\left(Z'_{2}MZ_{\Delta}\right)\delta'_{h}=\left\Vert \delta^{0}\right\Vert ^{2}h^{1/2}O_{p}\left(1\right)O_{p}\left(1\right)O_{p}\left(\psi_{h}\right),
\end{align*}
where we have used the fact that $Z'_{\Delta}MZ_{2}/\psi_{h}=O_{p}\left(1\right)$
{[}cf. \eqref{eq. ACRSC2}{]}. Hence, the second term of \eqref{eq. (45), Proof of  Lemma 1},
when multiplied by $\psi_{h}^{-1}$, is $O_{p}\left(h^{3/2-\kappa}\right)$
uniformly in $T_{b}$. Finally, let us consider $g_{e}\left(T_{b},\,\delta_{h}\right).$
Recall that $\widetilde{e}_{kh}$ defined in \eqref{Eq. eps WN} is
i.n.d. with zero mean and conditional variance $\sigma_{e,k-1}^{2}h$.
Upon applying the continuity of probability limits repeatedly, one
first obtains that the difference between the two terms in \eqref{Eq ge ACRSC4}
goes to zero at a fast enough rate as in the last step of the proof
of Proposition \ref{Prop 3 Asym}-(ii). That is, for $T$ large enough,
we can find a $c_{T}$ sufficiently small such that,
\begin{align*}
\psi_{h}^{-1}\left[e'MZ_{2}\left(Z'_{2}MZ_{2}\right)^{-1}Z_{2}Me-e'MZ_{0}\left(Z'_{0}MZ_{0}\right)^{-1}Z'_{0}Me\right] & =o_{p}\left(c_{T}h\right).
\end{align*}
Next, consider the first two terms of $g_{e}\left(T_{b},\,\delta_{h}\right).$
Using $Z'_{0}MZ_{2}=Z'_{2}MZ_{2}\pm Z'_{\Delta}MZ_{2}$, it is easy
to show that
\begin{align}
2h^{1/4}\left(\delta^{0}\right)' & \left(Z'_{0}MZ_{2}\right)\left(Z'_{2}MZ_{2}\right)^{-1}Z_{2}Me-2h^{1/4}\left(\delta^{0}\right)'\left(Z'_{0}Me\right)\nonumber \\
 & =2h^{1/4}\left(\delta^{0}\right)'Z'{}_{\Delta}Me\pm2h^{1/4}\left(\delta^{0}\right)'Z'{}_{\Delta}MZ_{2}\left(Z'_{2}MZ_{2}\right)^{-1}Z'_{2}Me.\label{Eq. CRSC6}
\end{align}
 Note that, uniformly in $T_{b}\left(u\right)$,
\begin{align*}
\psi_{h}^{-1} & h^{1/4}\left(\delta^{0}\right)'Z'{}_{\Delta}MZ_{2}\\
 & =h^{1/4}\left(\delta^{0}\right)'Z'{}_{\Delta}Z{}_{\Delta}+\left(\delta^{0}\right)'h^{1/4}\frac{Z_{\Delta}'X}{\psi_{h}}\left(X'X\right)^{-1}X'Z_{2}\\
 & =h^{1/4}\left(\delta^{0}\right)'\frac{Z'_{\Delta}Z{}_{\Delta}}{\psi_{h}}+\left(\delta^{0}\right)'h^{1/4}O_{p}\left(1\right)\\
 & =h^{1/4}\left\Vert \delta^{0}\right\Vert O_{p}\left(1\right)+\left\Vert \delta^{0}\right\Vert h^{1/4}O_{p}\left(1\right),
\end{align*}
 where we have used \eqref{Eq (1), Spot Estimator} and the fact that
$\left(X'X\right)^{-1}$ and $X'Z_{2}$ are each $O_{p}\left(1\right)$.
Recall the decomposition in \eqref{Eq. Xej spot vol eps}: 
\begin{align}
X'e & =O_{p}\left(h^{1-\kappa+1/4}\right)+O_{p}\left(h^{1/2}\right).\label{eq. CRSC7}
\end{align}
Thus, the last term in \eqref{Eq. CRSC6} multiplied by $\psi_{h}^{-1}$
is such that
\begin{align*}
\psi_{h}^{-1} & 2h^{1/4}\left(\delta^{0}\right)'Z'{}_{\Delta}MZ_{2}\left(Z'_{2}MZ_{2}\right)^{-1}Z'_{2}Me\\
 & =h^{1/4}\left\Vert \delta^{0}\right\Vert O_{p}\left(1\right)O_{p}\left(1\right)\left[O_{p}\left(h^{1-\kappa+1/4}\right)+O_{p}\left(h^{1/2}\right)\right]\\
 & =\left\Vert \delta^{0}\right\Vert h^{1/4}O_{p}\left(1\right)O_{p}\left(h^{1/2}\right)=\left\Vert \delta^{0}\right\Vert O_{p}\left(h^{3/4}\right).
\end{align*}
The first term of \eqref{Eq. CRSC6} can be decomposed further as
follows
\begin{align*}
2h^{1/4}\left(\delta^{0}\right)'Z'{}_{\Delta}Me & =2h^{1/4}\left(\delta^{0}\right)'Z'{}_{\Delta}e-2h^{1/4}\left(\delta^{0}\right)Z'{}_{\Delta}X\left(X'X\right)^{-1}X'e.
\end{align*}
Then, when multiplied by $\psi_{h}^{-1}$, the second term above is
such that, uniformly in $T_{b}$, 
\begin{align*}
h^{1/4} & \left(\delta^{0}\right)'\left(Z'_{\Delta}X/\psi_{h}\right)\left(X'X\right)^{-1}X'e\\
 & =h^{1/4}\left(\delta^{0}\right)'O_{p}\left(1\right)O_{p}\left(1\right)\left[O_{p}\left(h^{1-\kappa+1/4}\right)+O_{p}\left(h^{1/2}\right)\right]=O_{p}\left(h^{3/4}\right),
\end{align*}
where we have used \eqref{Eq (1), Spot Estimator} and \eqref{eq. CRSC7}.
Combining the last results, we have uniformly in $T_{b},$
\begin{align*}
\psi_{h}^{-1}g_{e}\left(T_{b},\,\delta_{h}\right) & =2h^{1/4}\left(\delta^{0}\right)'\left(Z'_{\Delta}e/\psi_{h}\right)\\
 & \quad+O_{p}\left(h^{3/4}\right)+\left\Vert \delta^{0}\right\Vert O_{p}\left(h^{3/4}\right)+o_{p}\left(c_{T}h\right),
\end{align*}
when $T$ is large and $c_{T}$ is a sufficiently small number. Then,
\begin{align*}
\psi_{h}^{-1} & \left(Q_{T}\left(T_{b}\right)-Q_{T}\left(T_{b}^{0}\right)\right)\\
 & =-\delta_{h}\left(Z_{\Delta}'Z_{\Delta}/\psi_{h}\right)\delta_{h}\pm2\delta_{h}'\left(Z'_{\Delta}e/\psi_{h}\right)\\
 & \quad+O_{p}\left(h^{3/2-\kappa}\right)+O_{p}\left(h^{3/4}\right)+\left\Vert \delta^{0}\right\Vert O_{p}\left(h^{3/4}\right)+o_{p}\left(c_{T}h\right).
\end{align*}
Therefore, for $T$ large enough, 
\begin{align*}
\psi_{h}^{-1}\left(Q_{T}\left(T_{b}\right)-Q_{T}\left(T_{b}^{0}\right)\right) & =-\delta_{h}\left(Z_{\Delta}'Z_{\Delta}/\psi_{h}\right)\delta_{h}\pm2\delta_{h}'\left(Z'_{\Delta}e/\psi_{h}\right)+o_{p}\left(h^{1/2}\right).
\end{align*}
 This concludes the proof of Lemma \ref{Lemma 1}. $\square$

\subsubsection{Proof of Theorem \ref{Theorem 1}}

\textit{Proof.} Let us focus on the case $T_{b}\left(v\right)\leq T_{b}^{0}$
(i.e., $v\leq0$). The change of time scale is obtained by a change
in variable. On the old time scale, by Proposition \ref{Prop 3 Asym},
$N_{b}\left(v\right)$ varies on the time interval $[N_{b}^{0}-\left|v\right|h^{1-\kappa},\,N_{b}^{0}+\left|v\right|h^{1-\kappa}]$
with $v\in[-C,\,C]$. Lemma \ref{Lemma 1} shows that the conditional
first moment of $Q_{T}\left(T_{b}\left(v\right)\right)-Q_{T}\left(T_{b}^{0}\right)$
is determined by that of $-\delta_{h}'\left(Z_{\Delta}'Z_{\Delta}\right)\delta_{h}\pm2\delta_{h}'\left(Z'_{\Delta}e\right).$
Next, we rescale time with $s\mapsto t\triangleq\psi_{h}^{-1}s$ on
$\mathcal{D}\left(C\right)$. This is achieved by rescaling the criterion
function $Q_{T}\left(T_{b}\left(u\right)\right)-Q_{T}\left(T_{b}^{0}\right)$
by the factor $\psi_{h}^{-1}$. First, note that the processes $Z_{t}$
and $e_{t}^{*}$ {[}recall \eqref{Model Regressors Integral Form}
and \eqref{Eq. eps WN}{]} are rescaled as follows on $\mathcal{D}\left(C\right)$.
Let $Z_{\psi,s}\triangleq\psi_{h}^{-1/2}Z_{s}$, $W_{\psi,e,s}\triangleq\psi_{h}^{-1/2}W_{e,s}$
and note that  
\begin{align}
dZ_{\psi,s}=\psi_{h}^{-1/2}\sigma_{Z,s}dW_{Z,s}, & \qquad dW_{\psi,e,s}=\psi_{h}^{-1/2}\sigma_{e,s}dW_{e,s},\qquad\textrm{with }s\in\mathcal{D}\left(C\right).\label{eq. SCCR_App_TS1}
\end{align}
For $s\in\left[N_{b}^{0}-Ch^{1-\kappa},\,N_{b}^{0}+Ch^{1-\kappa}\right]$,
let $v=\psi_{h}^{-1}\left(N_{b}^{0}-s\right)$ and, by using the properties
of $W_{.,s}$ and the fact that $\sigma_{Z,s},\,\sigma_{e,s}$ are
$\mathscr{F}_{s}$-measurable, we have
\begin{align}
dZ_{\psi,t}=\sigma_{Z,t}dW_{Z,t}, & \qquad dW_{\psi,e,t}=\sigma_{e,t}dW_{e,t},\qquad\textrm{with }t\in\mathcal{T}\left(C\right).\label{eq. CRSC_App_TS2}
\end{align}
This can be used into the following quantities for $N_{b}\left(v\right)\in\mathcal{D}\left(C\right)$.
First, 
\begin{align*}
\psi_{h}^{-1}Z_{\Delta}'Z_{\Delta} & =\sum_{k=T_{b}\left(v\right)+1}^{T_{b}^{0}}z_{\psi,kh}z_{\psi,kh},
\end{align*}
which by \eqref{eq. SCCR_App_TS1}-\eqref{eq. CRSC_App_TS2} is such
that 
\begin{align}
\psi_{h}^{-1}Z_{\Delta}'Z_{\Delta} & =\sum_{k=T_{b}^{0}+\left\lfloor v/h\right\rfloor }^{T_{b}^{0}}z_{kh}z'_{kh},\qquad\qquad v\in\mathcal{D}^{*}\left(C\right).\label{eq. CRSC_App_TS3}
\end{align}
Using the same argument: 
\begin{align}
\psi_{h}^{-1}Z_{\Delta}'\widetilde{e} & =\sum_{k=T_{b}^{0}+\left\lfloor v/h\right\rfloor }^{T_{b}^{0}}z_{kh}\widetilde{e}_{kh},\qquad\qquad v\in\mathcal{D}^{*}\left(C\right).\label{eq. CRSC_App_TS6}
\end{align}
Now $N_{b}\left(v\right)$ varies on $\mathcal{D}^{*}\left(C\right)$.
Furthermore, for sufficiently large $T$, Lemma \ref{Lemma 1} gives
\begin{align*}
Q_{T}\left(T_{b}\right)-Q_{T}\left(T_{b}^{0}\right) & =-\delta_{h}\left(Z_{\Delta}'Z_{\Delta}\right)\delta_{h}\pm2\delta_{h}'\left(Z'_{\Delta}e\right)+o_{p}\left(h^{1/2}\right),
\end{align*}
 and thus, when multiplying by $h^{-1/2},$ we have 
\[
\overline{Q}_{T}\left(T_{b}\right)=-\left(\delta^{0}\right)'Z_{\Delta}'Z_{\Delta}\left(\delta^{0}\right)\pm2\left(\delta^{0}\right)'\left(h^{-1/2}Z'_{\Delta}\widetilde{e}\right)+o_{p}\left(1\right),
\]
 since on $\mathcal{D}^{*}\left(C\right)$, $e_{kh}\sim\textrm{i.n.d.}\,\mathscr{N}(0,\,\sigma_{h,k-1}^{2}h)$,
$\sigma_{h,k}=O(h^{-1/4})\sigma_{e,k}$ and $\widetilde{e}_{kh}$
is the normalized error {[}i.e., $\widetilde{e}{}_{kh}\sim\textrm{i.n.d.}\,\mathscr{N}(0,\,\sigma_{e,k-1}^{2}h)${]}
defined in \eqref{Eq. eps WN}. Hence, according to the re-parametrization
introduced in the main text, we examine the behavior of 
\begin{align}
\overline{Q}_{T}\left(\theta^{*}\right) & =-\left(\delta^{0}\right)'\left(\sum_{k=T_{b}+1}^{T_{b}^{0}}z_{kh}z'_{kh}\right)\delta^{0}+2\left(\delta^{0}\right)'\left(h^{-1/2}\sum_{k=T_{b}+1}^{T_{b}^{0}}z_{kh}\widetilde{e}{}_{kh}\right).\label{Eq. QT after Lemma 1}
\end{align}
For the first term, a law of large numbers will be applied which yields
convergence in probability toward some quadratic covariation process.
For the second term, we observe that the finite-dimensional convergence
follows essentially from results in \citeReferencesSupp{jacod/protter:12}
(we indicate the precise theorems below) after some adaptation to
our context. Hence, we shall then verify the asymptotic stochastic
equicontinuity of the sequence of processes $\{\overline{Q}_{T}\left(\cdot\right),\,T\geq1\}.$
Let us associate to the continuous-time index $t$ a corresponding
$\mathcal{D}^{*}\left(C\right)$-specific index $t_{v}.$ This means
that each $t_{v}$ identifies a distinct $t$ in $\mathcal{D}^{*}\left(C\right)$
through $v$ as defined above. More specifically, for each $\left(\cdot,\,v\right)\in\mathcal{D}^{*}\left(C\right)$,
define the new functions 
\begin{align*}
J_{Z,h}\left(v\right)\triangleq\sum_{k=T_{b}\left(v\right)+1}^{T_{b}^{0}}z_{kh}z'_{kh} & ,\qquad\qquad J_{e,h}\left(v\right)\triangleq\sum_{k=T_{b}\left(v\right)+1}^{T_{b}^{0}}z_{kh}\widetilde{e},
\end{align*}
for $\left(T_{b}\left(v\right)+1\right)h\leq t_{v}<\left(T_{b}\left(v\right)+2\right)h$.
For $v\leq0$, the lower limit of the summation is $T_{b}\left(v\right)+1=T_{b}^{0}+\left\lfloor v/h\right\rfloor $
and thus the number of observations in each sum increases at rate
$1/h$. The functions $\{J_{Z,h}\left(v\right)\}$ and $\{J_{e,h}\left(v\right)\}$
have discontinuous, although \textit{càdlàg}, paths and thus they
belong to $\mathbb{D}\left(\mathcal{D}^{*}\left(C\right),\,\mathbb{R}\right).$
Since $Z_{t}^{\left(j\right)}$ ($j=1,\ldots,\,p$) is a continuous
It\^o semimartingale, we have by Theorem 3.3.1 in \citeReferencesSupp{jacod/protter:12}
that $J_{Z,h}\left(v\right)\overset{\textrm{u.c.p.}}{\Rightarrow}[Z,\,Z]_{1}\left(v\right),$
where $[Z,\,Z]_{1}\left(v\right)\triangleq[Z,\,Z]_{h\left\lfloor N_{b}^{0}/h\right\rfloor }-\left[Z,\,Z\right]_{h\left\lfloor t_{v}/h\right\rfloor },$
and recall by Assumption \ref{Assumption 2} that $\left[Z,\,Z\right]_{1}\left(v\right)$
is equivalent to $\left\langle Z,\,Z\right\rangle _{1}\left(v\right)$
where $\left\langle Z,\,Z\right\rangle _{1}\left(v\right)=\left\langle Z,\,Z\right\rangle _{h\left\lfloor t_{v}/h\right\rfloor }\left(v\right)$.
Next, let $\mathscr{W}_{h}\left(v\right)=h^{-1/2}J_{e,h}\left(v\right)$
and $\mathscr{W}_{1}\left(v\right)=\int_{N_{b}^{0}+v}^{N_{b}^{0}}\sigma_{Ze,s}dW_{s}^{1*}$
where $W_{s}^{1*}$ is defined in Section \ref{subsection Description Limiting Process}.
By Theorem 5.4.2 in \citeReferencesSupp{jacod/protter:12} we have
$\mathscr{W}_{h}\left(v\right)\overset{\mathcal{L}-\mathrm{s}}{\Rightarrow}\mathscr{W}_{1}\left(v\right)$
under the Skorokhod topology. Note the that both limit processes $\left[Z,\,Z\right]_{1}\left(v\right)$
and $\mathscr{W}_{1}\left(v\right)$ are continuous. This restores
the compatibility of the Skorokhod topology with the natural linear
structure of $\mathbb{D}(\mathcal{D}^{*}\left(C\right),\,\mathbb{R}).$
For $v\leq0$, the finite-dimensional stable convergence in law for
$\overline{Q}_{T}\left(\cdot\right)$ then follows: 
\[
\overline{Q}_{T}\left(\theta^{*}\right)\overset{\mathcal{L}_{f}-\mathrm{s}}{\rightarrow}-\left(\delta^{0}\right)'\left\langle Z,\,Z\right\rangle _{1}\left(v\right)\delta^{0}+2\left(\delta^{0}\right)'\mathscr{W}_{1}\left(v\right),
\]
 where $\overset{\mathcal{L}_{f}-\mathrm{s}}{\rightarrow}$ signifies
finite-dimensional stable convergence in law. Similarly, for $v>0,$
\[
\overline{Q}_{T}\left(\theta^{*}\right)\overset{\mathcal{L}_{f}-\mathrm{s}}{\rightarrow}-\left(\delta^{0}\right)'\left\langle Z,\,Z\right\rangle _{2}\left(v\right)\delta^{0}+2\left(\delta^{0}\right)'\mathscr{W}_{2}\left(v\right).
\]
 Next, we verify the asymptotic stochastic equicontinuity of the sequence
of processes $\{\overline{Q}_{T}\left(\cdot\right),\,T\geq1\}$.\footnote{Although in this proof it is not necessary to consider a neighborhood
about $\delta^{0}$ while proving stochastic equicontinuity, this
step will be needed to justify our inference methods later. Thus,
this proof is more general and may be useful in other contexts.} For $1\leq i\leq p,$ let $\zeta_{h,k}^{\left(i\right)}\triangleq z_{kh}^{\left(i\right)}\widetilde{e}_{kh},$
$\zeta_{h,k}^{*\left(i\right)}\triangleq\mathbb{E}[z_{kh}^{\left(i\right)}\widetilde{e}_{kh}|\,\mathscr{F}_{\left(k-1\right)h}],$
and $\zeta_{h,k}^{**\left(i\right)}\triangleq\zeta_{h,k}^{\left(i\right)}-\zeta_{h,k}^{*\left(i\right)}.$
For $1\leq i,\,j\leq p,$ let $\zeta_{Z,h,k}^{\left(i,j\right)}\triangleq z_{kh}^{\left(i\right)}z_{kh}^{\left(j\right)}-\Sigma_{Z,\left(k-1\right)h}^{\left(i,j\right)}h,$
\[
\zeta_{Z,h,k}^{*\left(i,j\right)}\triangleq\mathbb{E}\left[z_{kh}^{\left(i\right)}z_{kh}^{\left(j\right)}-\Sigma_{Z,\left(k-1\right)h}^{\left(i,j\right)}h|\,\mathscr{F}_{\left(k-1\right)h}\right],
\]
 and $\zeta_{Z,h,k}^{**\left(ij\right)}\triangleq\zeta_{Z,h,k}^{\left(ij\right)}-\zeta_{Z,h,k}^{*\left(ij\right)}.$
Then, we have the following decomposition for $\overline{Q}_{T}^{\textrm{c}}\left(\theta^{*}\right)\triangleq\overline{Q}_{T}^{\textrm{}}\left(\theta^{*}\right)+\left(\delta^{0}\right)'\left\langle Z,\,Z\right\rangle _{1}\left(v\right)\delta^{0}$
(if $v\leq0$, and defined analogously for $v>0$), 
\begin{align}
\overline{Q}_{T}^{\textrm{c}}\left(\theta^{*}\right) & =\sum_{r=1}^{4}\overline{Q}_{r,T}\left(\theta^{*}\right),\label{Eq. Qc}
\end{align}
where 
\begin{align*}
\overline{Q}_{1,T}\left(\theta^{*}\right)\triangleq-\left(\delta^{0}\right)'\left(\sum_{k}\zeta^{*}{}_{Z,h,k}\right)\delta^{0}, & \qquad\overline{Q}_{2,T}\left(\theta^{*}\right)\triangleq-\left(\delta^{0}\right)'\left(\sum_{k}\zeta^{**}{}_{Z,h,k}\right)\delta^{0},\\
\overline{Q}_{3,T}\left(\theta^{*}\right)\triangleq\left(\delta^{0}\right)'\text{ }\left(h^{-1/2}\sum_{k}\zeta^{*}{}_{h,k}\right), & \quad\mathrm{and}\quad\overline{Q}_{3,T}\left(\theta^{*}\right)\triangleq\left(\delta^{0}\right)'\text{ }\left(h^{-1/2}\sum_{k}\zeta^{*}{}_{h,k}\right);
\end{align*}
where $\sum_{k}$ stands for $\sum_{k=T_{b}^{0}+\left\lfloor v/h\right\rfloor }^{T_{b}^{0}}$.
Then,
\begin{align}
\sup_{\left(\theta,\,v\right)\in\mathcal{D}^{*}\left(C\right)}\left\Vert \overline{Q}_{3,T}\left(\theta^{*}\right)\right\Vert  & \leq K\left\Vert \delta^{0}\right\Vert h^{-1/2}\sum_{k}\left\Vert \zeta{}_{h,k}^{*}\right\Vert \overset{P}{\rightarrow}0,\label{Eq. Q3 Stoc Equi}
\end{align}
which follows from \citeReferencesSupp{jacod/rosenbaum:13} given
that $\Sigma_{Ze,k}=0$ identically by Assumption \ref{Assumption 1, CT}-(iv).
As for $\overline{Q}_{1,T}\left(\theta,\,v\right)$, we prove stochastic
equicontinuity directly, using the definition in \citeReferencesSupp{andrews:37hoe}.
Choose any $\varepsilon>0$ and $\eta>0$. Consider any $\left(\theta,\,v\right),\,(\bar{\theta},\,\bar{v})$
with $v<0<\bar{v}$ (the other cases can be proven similarly) and
$\bar{\delta}=\delta+c_{p\times1},$ where $c_{p\times1}$ is a $p\times1$
vector with each entry equals to $c\in\mathbb{R},$ with $0<c\leq\tau<\infty,$
then 
\begin{align*}
| & \overline{Q}_{1,T}\left(\theta^{*}\right)-\overline{Q}_{1,T}\left(\bar{\theta}^{*}\right)|\\
 & =\left|\bar{\delta}'\left(\sum_{k=T_{b}^{0}+1}^{T_{b}\left(\bar{v}\right)}\zeta^{*}{}_{Z,h,k}\right)\bar{\delta}-\delta'\left(\sum_{k=T_{b}\left(v\right)+1}^{T_{b}^{0}}\zeta{}_{Z,h,k}^{*}\right)\delta\right|\\
 & =\left|c_{p\times1}'\left(\sum_{k=T_{b}^{0}+1}^{T_{b}^{0}+\left\lfloor \overline{v}/h\right\rfloor }\zeta{}_{Z,h,k}^{*}\right)c_{p\times1}+\delta'\left(\sum_{k=T_{b}^{0}+1}^{T_{b}\left(\bar{v}\right)}\zeta{}_{Z,h,k}^{*}-\sum_{k=T_{b}^{0}+\left\lfloor v/h\right\rfloor }^{T_{b}^{0}}\zeta{}_{Z,h,k}^{*}\right)\delta\right|\\
 & \leq K(\sum_{k=T_{b}^{0}+1}^{T_{b}\left(\bar{v}\right)}\left\Vert \zeta{}_{Z,h,k}^{*}\right\Vert \left\Vert c_{p\times1}\right\Vert ^{2}+\left\Vert \sum_{k=T_{b}^{0}+1}^{T_{b}^{0}+\left\lfloor \overline{v}/h\right\rfloor }\zeta{}_{Z,h,k}^{*}-\sum_{k=T_{b}^{0}+\left\lfloor v/h\right\rfloor }^{T_{b}^{0}}\zeta{}_{Z,h,k}^{*}\right\Vert \left\Vert \delta\right\Vert ^{2}\\
 & \leq K\left(\left(pc^{2}\right)\sum_{k=T_{b}^{0}+1}^{T_{b}^{0}+\left\lfloor \overline{v}/h\right\rfloor }\left\Vert \zeta{}_{Z,h,k}^{*}\right\Vert +\sum_{k=T_{b}^{0}+\left\lfloor v/h\right\rfloor }^{T_{b}\left(\bar{v}\right)}\left\Vert \zeta{}_{Z,h,k}^{*}\right\Vert \left\Vert \delta\right\Vert ^{2}\right).
\end{align*}
By It\^o's formula $||\zeta{}_{Z,h,k}^{*}||=O(h^{3/2})$, and so
\begin{align*}
\left|\overline{Q}_{1,T}\left(\theta^{*}\right)-\overline{Q}_{1,T}\left(\bar{\theta}^{*}\right)\right| & \leq K\left(c^{2}h^{-1}O_{p}\left(h^{3/2}\right)O\left(\tau\right)+\left\Vert \delta\right\Vert ^{2}h^{-1}O_{p}\left(h^{3/2}\right)O\left(\tau\right)\right)\\
 & \leq K\left(c^{2}O_{p}\left(h^{1/2}\right)O\left(\tau\right)+\left\Vert \delta\right\Vert ^{2}O_{p}\left(h^{1/2}\right)O\left(\tau\right)\right),
\end{align*}
 which goes to zero uniformly in $\theta^{*}\in\Theta$ as $\tau\rightarrow0$.
Next, consider $\overline{Q}_{2,T}\left(\theta^{*}\right)$ and observe
that for any $r\geq1,$ standard estimates for It\^o semimartingales
yields $\mathbb{E}(||\zeta{}_{Z,h,k}^{**}||^{r}|\,\mathscr{F}_{\left(k-1\right)h})\leq K_{r}h^{r}$.
Then, by using a maximal inequality and choosing $r>2,$ 
\begin{align}
\left(\mathbb{E}\left[\sup_{\left(\theta,\,v\right)\in\mathcal{D}^{*}\left(C\right)}\left|\overline{Q}_{2,T}\left(\theta^{*}\right)\right|\right]^{r}\right)^{1/r} & \leq K_{r}\left\Vert \delta^{0}\right\Vert ^{2}h^{-2/r}h\leq K_{r}h^{1-2/r}\rightarrow0,\label{Eq. Q2 Stoc Equi}
\end{align}
 and thus we can use Markov's inequality together with the latter
result to verify that $\overline{Q}_{2,T}\left(\theta^{*}\right)$
is stochastically equicontinuous. Turning to $\overline{Q}_{4,T}\left(\theta^{*}\right)$,
\begin{align*}
 & \left|\overline{Q}_{4,T}\left(\bar{\theta}^{*}\right)-\overline{Q}_{4,T}\left(\theta^{*}\right)\right|\\
 & =\left|\bar{\delta}'\left(h^{-1/2}\sum_{k=T_{b}^{0}+1}^{T_{b}^{0}+\left\lfloor \overline{v}/h\right\rfloor }\zeta{}_{e,h,k}^{*}\right)-\delta'\left(h^{-1/2}\sum_{k=T_{b}^{0}+\left\lfloor v/h\right\rfloor }^{T_{b}^{0}}\zeta{}_{e,h,k}^{*}\right)\right|\\
 & =\Biggl|c_{p\times1}'\left(h^{-1/2}\sum_{k=T_{b}^{0}+1}^{T_{b}^{0}+\left\lfloor \overline{v}/h\right\rfloor }\zeta{}_{e,h,k}^{*}\right)\\
 & \quad+\delta'\left(h^{-1/2}\sum_{k=T_{b}^{0}+1}^{T_{b}^{0}+\left\lfloor \overline{v}/h\right\rfloor }\zeta{}_{e,h,k}^{*}-h^{-1/2}\sum_{k=T_{b}^{0}+\left\lfloor v/h\right\rfloor }^{T_{b}^{0}}\zeta{}_{e,h,k}^{*}\right)\Biggr|\\
 & \leq K(h^{-1/2}\sum_{k=T_{b}^{0}+1}^{T_{b}^{0}+\left\lfloor \overline{v}/h\right\rfloor }\left\Vert \zeta{}_{e,h,k}^{*}\right\Vert \left\Vert c_{p\times1}\right\Vert \\
 & \quad+\left\Vert h^{-1/2}\sum_{k=T_{b}^{0}+1}^{T_{b}^{0}+\left\lfloor \overline{v}/h\right\rfloor }\zeta{}_{e,h,k}^{*}-h^{-1/2}\sum_{k=T_{b}^{0}+\left\lfloor v/h\right\rfloor }^{T_{b}^{0}}\zeta{}_{e,h,k}^{*}\right\Vert \left\Vert \delta\right\Vert )\\
 & \leq K\left(pch^{-1/2}\sum_{k=T_{b}^{0}+1}^{T_{b}^{0}+\left\lfloor \overline{v}/h\right\rfloor }\left\Vert \zeta{}_{e,h,k}^{*}\right\Vert +h^{-1/2}\sum_{k=T_{b}^{0}+\left\lfloor v/h\right\rfloor }^{T_{b}^{0}+\left\lfloor \overline{v}/h\right\rfloor }\left\Vert \zeta{}_{e,h,k}^{*}\right\Vert \left\Vert \delta\right\Vert \right).
\end{align*}
By the Burkh\"{o}lder-Davis-Gundy inequality, $||\zeta{}_{e,h,k}^{*}||\leq Kh^{3/2}$
(recall $\Sigma_{Ze,t}=0$ for all $t\geq0$), so that 
\begin{align*}
\left|\overline{Q}_{4,T}\left(\theta^{*}\right)-\overline{Q}_{4,T}\left(\bar{\theta}^{*}\right)\right| & \leq K(c^{2}h^{-1/2}h^{-1}h^{3/2}O\left(\tau\right)\\
 & \quad+\left\Vert \delta\right\Vert ^{2}h^{-1/2}h^{-1}h^{3/2}O\left(\tau\right))\\
 & \leq K\left(c^{2}O\left(\tau\right)+\left\Vert \delta\right\Vert ^{2}O\left(\tau\right)\right).
\end{align*}
 Then for every $\eta>0$, with $\mathbf{B}\left(\tau,\,\left(\theta,\,v\right)\right)$
a closed ball of radius $\tau>0$ around $\theta^{*}$, the quantity
\begin{align}
\underset{h\downarrow0}{\limsup} & P\left[\sup_{\theta^{*}\in\Theta:\,\bar{\theta}^{*}\in\mathbf{B}\left(\tau,\,\theta^{*}\right)}\left|\overline{Q}_{4,T}\left(\theta^{*}\right)-\overline{Q}_{4,T}\left(\bar{\theta}^{*}\right)\right|>\eta\right],\label{Eq. Q1 Stoc Equi}
\end{align}
 can be made arbitrary less than $\varepsilon>0$ as $h\downarrow0,$
by choosing $\tau$ small enough. Combining \eqref{Eq. Q3 Stoc Equi},
\eqref{Eq. Q2 Stoc Equi} and \eqref{Eq. Q1 Stoc Equi}, we conclude
that the process $\{\overline{Q}_{T}^{\textrm{}}\left(\theta,\,v\right),\,T\geq1\}$
is asymptotically stochastic equicontinuous. Since the finite-dimensional
convergence was demonstrated above, this suffices to guarantee the
stable convergence in law of the process $\{\overline{Q}_{T}\left(\theta,\,v\right),\,T\geq1\}$
toward a two-sided Gaussian limit process with drift $(\delta^{0})'\left[Z,\,Z\right]_{\cdot}\left(\cdot\right)\delta^{0}$,
having $P$-a.s. continuous sample paths with $\mathscr{F}$-conditional
covariance matrix given in \eqref{Eq. Cov Matrix of W process, CT}.
Because $N(\widehat{\lambda}_{b,\pi}-\lambda_{0})=O_{p}\left(1\right)$
under the new ``fast time scale'', and $\mathcal{D}^{*}\left(C\right)$
is compact, then the main assertion of the theorem follows from the
continuous mapping theorem for the argmax functional. In view of Section
\ref{Negligibility-of-the, CT}, a result which shows the negligibility
of the drift term, the proof of Theorem \ref{Theorem 1} is complete.
$\square$

\subsubsection{Proof of Theorem \ref{Theorem 2, Asymptotic Distribution immediate Stationary Regimes}}

\textit{Proof.} By Theorem \ref{Theorem 1} and using the property
of the Gaussian law of the limiting process, 
\begin{align*}
\overline{Q}_{T}\left(\theta,\,v\right) & \overset{\mathcal{L}-\mathrm{s}}{\Rightarrow}\mathscr{H}\left(v\right)=\begin{cases}
-\left(\delta^{0}\right)'\left\langle Z,\,Z\right\rangle _{1}\left(v\right)\delta^{0}+2\left(\left(\delta^{0}\right)'\Omega_{\mathscr{W},1}\left(\delta^{0}\right)\right)^{1/2}W_{1}^{*}\left(v\right), & \textrm{if }v\leq0\\
-\left(\delta^{0}\right)'\left\langle Z,\,Z\right\rangle _{2}\left(v\right)\delta^{0}+2\left(\left(\delta^{0}\right)'\Omega_{\mathscr{W},2}\left(\delta^{0}\right)\right)^{1/2}W_{2}^{*}\left(v\right), & \textrm{if }v>0.
\end{cases}
\end{align*}
By a change in variable $v=\vartheta^{-1}s$ with $\vartheta=((\delta^{0})'\left\langle Z,\,Z\right\rangle _{1}\delta^{0})^{2}/(\delta^{0})'\Omega_{\mathscr{W},1}(\delta^{0})$,
we can show that
\begin{align*}
 & \underset{v\in\mathcal{A}}{\mathrm{argmax}}\mathscr{H}\left(v\right)\\
 & \overset{d}{\equiv}\underset{s\in\mathcal{A}^{*}}{\mathrm{argmax}}\mathscr{V}\left(s\right),
\end{align*}
 where 
\begin{align*}
\mathscr{V}\left(s\right) & =\begin{cases}
-\frac{\left|s\right|}{2}+W_{1}^{*}\left(s\right), & \textrm{if }s<0\\
-\frac{\left(\delta^{0}\right)'\left\langle Z,\,Z\right\rangle _{2}\delta^{0}}{\left(\delta^{0}\right)'\left\langle Z,\,Z\right\rangle _{1}\delta^{0}}\frac{\left|s\right|}{2}+\left(\frac{\left(\delta^{0}\right)'\Omega_{\mathscr{W},2}\left(\delta^{0}\right)}{\left(\delta^{0}\right)'\Omega_{\mathscr{W},1}\left(\delta^{0}\right)}\right)^{1/2}W_{2}^{*}\left(s\right), & \textrm{if }s\geq0,
\end{cases}
\end{align*}
and we have used the facts that $W\left(s\right)\overset{d}{\equiv}W\left(-s\right),$
$W\left(cs\right)\overset{d}{\equiv}\left|c\right|^{1/2}W\left(s\right)$,
and for any $c>0$ and for any function $f\left(s\right)$, $\arg\max_{s}cf\left(s\right)=\arg\max_{s}f\left(s\right)$.
Thus, 
\begin{align*}
 & \underset{v\in\mathcal{A}}{\mathrm{argmax}}\mathscr{H}\left(v\right)\\
 & \overset{d}{\equiv}\underset{s\in\mathcal{A}^{*}}{\mathrm{argmax}}\\
 & \quad\left(\frac{\left(\left(\delta^{0}\right)'\left\langle Z,\,Z\right\rangle _{1}\delta^{0}\right)^{2}}{\left(\delta^{0}\right)'\Omega_{\mathscr{W},1}\left(\delta^{0}\right)}\right)^{-1}\mathscr{V}\left(s\right),
\end{align*}
and finally by the continuous mapping theorem for the argmax functional,
\begin{align*}
 & \frac{\left(\left(\delta^{0}\right)'\left\langle Z,\,Z\right\rangle _{1}\delta^{0}\right)^{2}}{\left(\delta^{0}\right)'\Omega_{\mathscr{W},1}\left(\delta^{0}\right)}N\left(\widehat{\lambda}_{b,\pi}-\lambda_{0}\right)\\
 & \Rightarrow\underset{s\in\mathcal{A}^{*}}{\mathrm{argmax}}\mathscr{V}\left(s\right).
\end{align*}
This concludes the proof. $\square$

\subsubsection{\label{Negligibility-of-the, CT}Negligibility of the Drift Term}

We are in the setting of Section \ref{Section Consistency-and-Rate }-\ref{Section Asymptotic Distribution: Continuous Case}.
In Proposition \ref{Proposition (1) - Consistency}-\ref{Proposition OLS Asymtptoc Distribu}
and \ref{Prop 3 Asym}, the drift processes $\mu_{\cdot,t}$ from
\eqref{Model Regressors Integral Form} are clearly of higher order
in $h$ and so they are negligible. In Theorem \ref{Theorem 1}, we
first changed the time scale and then normalized the criterion function
by the factor $h^{-1/2}$. The change of time scale now results in
\begin{align}
dZ_{\psi,s}=\psi_{h}^{-1/2}\mu_{Z,s}ds+\psi_{h}^{-1/2}\sigma_{Z,s}dW_{Z,s}, & \quad dW_{\psi,e,s}=\psi_{h}^{-1/2}\sigma_{e,s}dW_{e,s},\label{eq. SCCR_App_TS1 Drift}
\end{align}
 with $s\in\mathcal{D}\left(C\right).$ Given $s\mapsto t=\psi_{h}^{-1}s$,
we have $\psi_{h}^{-1/2}\mu_{Z,s}ds=\psi_{h}^{-1/2}\mu_{Z,s}\psi_{h}\left(ds/\psi_{h}\right)=\mu_{Z,s}\psi_{h}^{\vartheta}dt$
with $\vartheta=1/2$. Then, as in \eqref{eq. CRSC_App_TS2}, $dZ_{\psi,t}=\psi_{h}^{\vartheta}\mu_{Z,t}dt+\sigma_{Z,t}dW_{Z,t}$
and $dW_{\psi,e,t}=\sigma_{e,t}dW_{e,t}$ with $t\in\mathcal{D}^{*}\left(C\right)$.
Thus, the change of time scale effectively makes the drift $\mu_{Z,s}ds$
of even higher order. We show a stronger result in that we demonstrate
its negligibility even in the case $\vartheta=0$; hence, we show
that the limit law of \eqref{Eq. QT after Lemma 1} remains the same
when $\mu_{\cdot,t}$ are nonzero. We set, for any $1\leq i\leq p$
and $1\leq j\leq q+p$, 
\begin{align*}
\mu_{Z,k}^{*\left(i\right)} & \triangleq\int_{\left(k-1\right)h}^{kh}\mu_{Z,s}^{\left(i\right)}ds,\qquad\mu_{X,k}^{*\left(j\right)}\triangleq\int_{\left(k-1\right)h}^{kh}\mu_{X,s}^{\left(j\right)}ds,\\
z_{0,kh}^{\left(i\right)} & \triangleq\sum_{r=1}^{p}\int_{\left(k-1\right)h}^{kh}\sigma_{Z,s}^{\left(i,r\right)}dW_{Z}^{\left(r\right)},\quad\mathrm{and}\quad x_{0,kh}^{\left(j\right)}\triangleq\sum_{r=1}^{q+p}\int_{\left(k-1\right)h}^{kh}\sigma_{X,s}^{\left(j,r\right)}dW_{X}^{\left(r\right)}.
\end{align*}
Note that
\[
z_{kh}^{\left(i\right)}x_{kh}^{\left(j\right)}=\mu_{Z,k}^{*\left(i\right)}\mu_{X,k}^{*\left(j\right)}+\mu_{Z,k}^{*\left(i\right)}x_{0,kh}^{\left(j\right)}+z_{0,kh}^{\left(i\right)}\mu_{X,k}^{*\left(j\right)}+z_{0,kh}^{\left(i\right)}x_{0,kh}^{\left(j\right)}.
\]
 Recall that $\mu_{\cdot,k}^{*\left(\cdot\right)}$ is $O\left(h\right)$
uniformly in $k$, and note that $\mu_{Z,k}^{*\left(i\right)}x_{0,kh}^{\left(j\right)}+\mu_{Z,k}^{*\left(i\right)}z_{0,kh}^{\left(i\right)}$
follows a Gaussian law with zero mean and variance of order $O\left(h^{3}\right)$.
Also note that on $\mathcal{D}^{*}\left(C\right)$, $T_{b}^{0}-T_{b}-1\asymp1/h$,
where $a_{h}\asymp b_{h}$ if for some $c\geq1,$ $b_{h}/c\leq a_{h}\leq cb_{h}$.
Then, 
\begin{align*}
\sum_{k=T_{b}+1}^{T_{b}^{0}}z_{kh}^{\left(i\right)}x_{kh}^{\left(j\right)} & =\sum_{k=T_{b}+1}^{T_{b}^{0}}\mu_{Z,k}^{*\left(i\right)}\mu_{X,k}^{*\left(j\right)}+\sum_{k=T_{b}+1}^{T_{b}^{0}}\mu_{Z,k}^{*\left(i\right)}x_{0,kh}^{\left(j\right)}\\
 & \quad+\sum_{k=T_{b}+1}^{T_{b}^{0}}z_{0,kh}^{\left(i\right)}\mu_{X,k}^{*\left(j\right)}+\sum_{k=T_{b}+1}^{T_{b}^{0}}z_{0,kh}^{\left(i\right)}x_{0,kh}^{\left(j\right)}\\
 & =o\left(h^{1/2}\right)+o_{p}\left(h^{1/2}\right)+\sum_{k=T_{b}+1}^{T_{b}^{0}}z_{0,kh}^{\left(i\right)}x_{0,kh}^{\left(j\right)}.
\end{align*}
Therefore, conditionally on $\Sigma^{0}=\left\{ \mu_{\cdot,t},\,\sigma_{\cdot,t}\right\} _{t\geq0}$,
the limit law of 
\begin{align*}
\overline{Q}_{T}\left(\theta^{*}\right) & =-\left(\delta^{0}\right)'\left(\sum_{k=T_{b}+1}^{T_{b}^{0}}z_{kh}z'_{kh}\right)\delta^{0}+2\left(\delta^{0}\right)'\left(h^{-1/2}\sum_{k=T_{b}+1}^{T_{b}^{0}}z_{kh}\widetilde{e}{}_{kh}\right),
\end{align*}
is the same as the limit law of 
\begin{align*}
-\left(\delta^{0}\right)'\left(\sum_{k=T_{b}+1}^{T_{b}^{0}}z_{0,kh}z'_{0,kh}\right)\delta^{0}+2\left(\delta^{0}\right)'\left(h^{-1/2}\sum_{k=T_{b}+1}^{T_{b}^{0}}z_{0,kh}\widetilde{e}{}_{kh}\right) & ,
\end{align*}
 which completes the proof of Theorem \ref{Theorem 1}. $\square$

\subsubsection{Proof of Proposition \ref{Proposition CI}}

\noindent\textit{Proof.} Replace $\xi_{1},\,\xi_{2},\,\rho$ and
$\vartheta$ in \eqref{Equation (2) Asymptotic Distribution} by their
corresponding estimates $\xi_{1},\,\xi_{2},\,\rho$ and $\vartheta$,
respectively. Multiply both sides of \eqref{Equation (2) Asymptotic Distribution}
by $h^{\kappa}$ and apply a change in variable $v=s/h^{\kappa}$.
Consider the case $s<0$. On the ``fast time scale'', $W_{\cdot}^{*}$
is replaced by $\widehat{W}_{1,h}\left(s\right)=W_{1,h}^{*}(sh^{\kappa})$
$\left(s<0\right)$, where $W_{1,h}^{*}\left(s\right)$ is a sample-size
dependent Wiener process. It follows that 
\begin{align*}
-h^{-\kappa}\frac{\left|s\right|}{2}+h^{-\kappa}W_{1,h}^{*}\left(h^{\kappa}s\right) & =-\frac{\left|v\right|}{2}+W_{1}^{*}\left(v\right).
\end{align*}
 A similar argument can be applied when $s\geq0$. Let $\widehat{\mathscr{V}}\left(s\right)$
denote $\mathscr{V}\left(s\right)$ constructed with the proposed
estimates in place of the population parameters. Then, 
\begin{align*}
h^{-\kappa}\underset{s\in\left[\left(\pi-\widehat{\lambda}_{b}\right)\widehat{\vartheta},\,\left(1-\pi-\widehat{\lambda}_{b}^{\mathrm{}}\right)\widehat{\vartheta}\right]}{\mathrm{argmax}}\mathscr{\widehat{V}}\left(s\right) & =\underset{v\in\left[\left(\pi-\widehat{\lambda}_{b}\right)\widehat{\vartheta}/h^{2\kappa},\,\left(1-\pi-\widehat{\lambda}_{b}\right)\widehat{\vartheta}/h^{2\kappa}\right]}{\mathrm{argmax}}\mathscr{\widehat{V}}\left(v\right)\\
 & \Rightarrow\underset{v\in\left[\left(\pi-\lambda_{0}\right)\vartheta,\,\left(1-\pi-\lambda_{0}\right)\vartheta\right]}{\mathrm{argmax}}\mathscr{V}\left(v\right),
\end{align*}
which is equal to the right-hand side of \eqref{Equation (2) Asymptotic Distribution}
since
\[
\vartheta=\left\Vert \delta^{0}\right\Vert ^{2}\overline{\sigma}^{-2}\left(\left(\delta^{0}\right)'\left\langle Z,\,Z\right\rangle _{1}\delta^{0}\right)^{2}/\left(\delta^{0}\right)'\Omega_{\mathscr{W},1}\left(\delta^{0}\right).
\]
 Therefore, equation \eqref{Equation (2) Asymptotic Distribution}
holds when we use the proposed plug-in estimates. $\square$

\subsection{Proofs of the Results in Section \ref{subsection Asymptotic-Results- Pre}}

The steps are similar to those used for the case when the model does
not include predictable processes. However, we need to rely occasionally
on different asymptotic results since the latter processes have distinct
statistical properties. Recall that the dependent variable $\Delta_{h}Y_{k}$
in model \eqref{Full Model} is the increment of a discretized process
which cannot be identified as an ordinary diffusion. However, its
normalized version, $\widetilde{Y}_{\left(k-1\right)h}\triangleq h^{1/2}Y_{\left(k-1\right)h}$,
is well-defined and we exploit this property in the proof. $\Delta_{h}Y_{k}$
has first conditional moment of order $O(h^{-1/2})$, it has unbounded
variation and does not belong to the usual class of semimartingales.\footnote{For an introduction to the terminology used in this sub-section, we
refer the reader to first chapters in \citeReferencesSupp{jacod/shiryaev:03}.} The predictable process $\{Y_{\left(k-1\right)h}\}_{k=1}^{T}$ derived
from it has different properties. Its ``quadratic variation'' exists,
and thus it is finite in any fixed time interval. That is, the integrated
second moments of the regressor $Y_{\left(k-1\right)h}$ are finite,
i.e., we have
\begin{align*}
\sum_{k=1}^{T}\left(Y_{\left(k-1\right)h}h\right)^{2} & =\sum_{k=1}^{T}\left(h^{1/2}Y_{\left(k-1\right)h}h^{1/2}\right)^{2}=h\sum_{k=1}^{T}\left(\widetilde{Y}_{\left(k-1\right)h}\right)^{2}=O_{p}\left(1\right),
\end{align*}
by a standard approximation for Riemann sums and recalling that $\widetilde{Y}_{\left(k-1\right)h}$
is scaled to be $O_{p}\left(1\right).$ Then it is easy to see that
$\{\widetilde{Y}_{\left(k-1\right)h}\}_{k=1}^{T}$ has nice properties.
It is left-continuous, adapted, and of finite variation in any finite
time interval. When used as the integrand of a stochastic integral,
the integral itself makes sense. Importantly, its quadratic variation
is null and the process is orthogonal to any continuous local martingale.
These properties will be used in the sequel. In analogy to the previous
section, we use a localization procedure and thus we need the following
assumption related to Assumption \ref{Assumption Localization}.
\begin{assumption}
\label{Assumption Localization, Pre}Assumption \ref{Assumption 1, Pre}
holds, the process $\{\widetilde{Y}_{t},\,D_{t},\,Z_{t}\}_{t\geq0}$
takes value in some compact set and the processes $\{\mu_{\cdot,t},\,\sigma_{\cdot,t}\}_{t\geq0}$
(except $\{\mu_{\cdot,t}^{h}\}_{t\geq0}$) are bounded.
\end{assumption}
Recall the notation $M=I-X\left(X'X\right)^{-1}X'$, where now
\begin{align}
X & =\begin{bmatrix}h^{1/2} & Y_{0}h & \Delta_{h}D'_{1} & \Delta_{h}Z'_{1}\\
h^{1/2} & Y_{1}h & \Delta_{h}D'_{2} & \Delta_{h}Z'_{2}\\
\vdots & \vdots & \vdots & \vdots\\
h^{1/2} & Y_{Th}h & \Delta_{h}D'_{T} & \Delta_{h}Z'_{T}
\end{bmatrix}_{T\times\left(q+p+2\right)}.\label{X Expression, Pre}
\end{align}
Thus, $X'X$ is a $\left(q+p+2\right)\times\left(q+p+2\right)$ matrix
given by $\begin{bmatrix}a_{1} & a_{2} & a_{3} & a_{4}\end{bmatrix},$
where
\begin{align*}
a_{1}=\begin{bmatrix}\sum_{k=1}^{T}h\\
h^{1/2}\sum_{k=1}^{T}\left(Y_{\left(k-1\right)h}h\right)\\
\sum_{k=1}^{T}h^{1/2}\left(\Delta_{h}D{}_{k}\right)\\
\sum_{k=1}^{T}h^{1/2}\left(\Delta_{h}Z{}_{k}\right)
\end{bmatrix} & ,\qquad a_{2}=\begin{bmatrix}h^{1/2}\sum_{k=1}^{T}\left(Y_{\left(k-1\right)h}h\right)\\
\sum_{k=1}^{T}\left(Y_{\left(k-1\right)h}^{2}\cdot h^{2}\right)\\
\sum_{k=1}^{T}\left(\Delta_{h}D{}_{k}\right)\left(Y_{\left(k-1\right)h}h\right)\\
\sum_{k=1}^{T}\left(\Delta_{h}Z{}_{k}\right)\left(Y_{\left(k-1\right)h}h\right)
\end{bmatrix},
\end{align*}
\begin{align*}
a_{3}=\begin{bmatrix}\sum_{k=1}^{T}h^{1/2}\left(\Delta_{h}D'_{k}\right)\\
\sum_{k=1}^{T}\left(\Delta_{h}D'_{k}\right)\left(Y_{\left(k-1\right)h}h\right)\\
X_{D}'X_{D}\\
X'_{Z}X_{D}
\end{bmatrix} & ,\qquad a_{4}=\begin{bmatrix}\sum_{k=1}^{T}h^{1/2}\left(\Delta_{h}Z'{}_{k}\right)\\
\sum_{k=1}^{T}\left(\Delta_{h}Z'_{k}\right)\left(Y_{\left(k-1\right)h}h\right)\\
X'_{D}X_{Z}\\
X_{Z}'X_{Z}
\end{bmatrix},
\end{align*}
 where $X_{D}'X_{D}$ is a $q\times q$ matrix whose $\left(j,\,r\right)$-th
component is the approximate covariation between the $j$-th and $r$-th
element of $D$, with $X_{D}'X_{Z}$ defined similarly. In view of
the properties of $Y_{\left(k-1\right)h}$ outlined above and Assumption
\ref{Assumption Localization, Pre}, $X'X$ is $O_{p}\left(1\right)$
as $h\downarrow0.$ The limit matrix is symmetric positive definite
where the only zero elements are in the $2\times\left(q+p\right)$
upper right sub-block, and by symmetry in the $\left(q+p\right)\times2$
lower left sub-block. Furthermore, we have 
\begin{align}
X'e & =\begin{bmatrix}\sum_{k=1}^{T}h^{1/2}e_{kh}\\
\sum_{k=1}^{T}\left(Y_{\left(k-1\right)h}h\right)e_{kh}\\
\sum_{k=1}^{T}\Delta_{h}D{}_{k}e_{kh}\\
\sum_{k=1}^{T}\Delta_{h}Z{}_{k}e_{kh}
\end{bmatrix}.\label{X'e}
\end{align}
 The other statistics are omitted in order to save space. Again the
proofs are first given for the case where the drift processes $\mu_{Z,t},\,\mu_{D,t}$
of the semimartingale regressors $Z$ and $D$ are identically zero.
In the last step we extend the results to nonzero $\mu_{Z,t},\,\mu_{D,t}$.
We also start by conditioning on the processes $\mu_{Z,t},\,\mu_{D,t}$
and on all the volatility processes so that they are treated as if
they were deterministic. We begin with a preliminary lemma.
\begin{lem}
\label{Lemma Cross Product Null}For $1\leq i\leq2$, $3\leq j\leq p+2$
and $\gamma>0$, $\sum_{k=\left\lfloor s/h\right\rfloor }^{\left\lfloor t/h\right\rfloor }z_{kh}^{\left(i\right)}z_{kh}^{\left(j\right)}\overset{\mathrm{u.c.p.}}{\Rightarrow}0$,
for all $N>t>s+\gamma>s>0.$
\end{lem}
\noindent\textit{Proof.} Without loss of generality consider any
$3\leq j\leq p+2$ and $N>t>s>0.$ We have $\sum_{k=\left\lfloor s/h\right\rfloor }^{\left\lfloor t/h\right\rfloor }z_{kh}^{\left(1\right)}z_{kh}^{\left(j\right)}=\sum_{k=\left\lfloor s/h\right\rfloor }^{\left\lfloor t/h\right\rfloor }\sqrt{h}(\Delta_{h}M_{Z,k}^{\left(j\right)}),$
with further $\mathbb{E}[z_{kh}^{\left(1\right)}z_{kh}^{\left(j\right)}|\,\mathscr{F}_{\left(k-1\right)h}]=0,\,|z_{kh}^{\left(1\right)}z_{kh}^{\left(j\right)}|\leq K$
for some $K$ by Assumption \ref{Assumption Localization, Pre}. Thus
$\{z_{kh}^{\left(i\right)}z_{kh}^{\left(j\right)},\,\mathscr{F}_{kh}\}$
is a martingale difference array. Then, for any $\eta>0,$ 
\begin{align*}
P & \left(\sum_{k=\left\lfloor s/h\right\rfloor }^{\left\lfloor t/h\right\rfloor }\left|z_{kh}^{\left(1\right)}z_{kh}^{\left(j\right)}\right|^{2}>\eta\right)\\
 & \leq\frac{K}{\eta}\mathbb{E}\left(\sum_{k=\left\lfloor s/h\right\rfloor }^{\left\lfloor t/h\right\rfloor }h^{2}\left(\Delta_{h}M_{Z,k}^{\left(j\right)}\right)^{2}\right)\leq\frac{K}{\eta}hO_{p}\left(t-s\right)\rightarrow0,
\end{align*}
 where the second inequality follows from the Burkh\"{o}lder-Davis-Gundy
inequality with parameter $r=2$. This shows that the array $\{|z_{kh}^{\left(i\right)}z_{kh}^{\left(j\right)}|^{2}\}$
is asymptotically negligible. By Lemma 2.2.11 in the Appendix of \citeReferencesSupp{jacod/protter:12},
we verify the claim for $i=1.$ For the case $i=2$ note that $z_{kh}^{\left(2\right)}z_{kh}^{\left(j\right)}=(Y_{\left(k-1\right)h}h)(\Delta_{h}M_{Z,k}^{\left(j\right)}),$
and recall that $\widetilde{Y}_{\left(k-1\right)h}=h^{1/2}Y_{\left(k-1\right)h}=O_{p}\left(1\right)$.
Thus, the same proof remains valid for the case $i=2.$ The assertion
of the lemma follows. $\square$

\subsubsection{Proof of Proposition \ref{Prop 1-2 Asym - Pre}}

\noindent\textit{Proof of part (i) of Proposition \ref{Prop 1-2 Asym - Pre}.}
Following the same steps that led to \eqref{eq. A.2.5, Casini (2015a)-2},
we can write 
\begin{align}
Q_{T}\left(T_{b}\right)-Q_{T}\left(T_{0}\right) & =-\left|T_{b}-T_{b}^{0}\right|d\left(T_{b}\right)+g_{e}\left(T_{b}\right),\qquad\textrm{for all }T_{b},\label{Eq. A.2.6 Pre}
\end{align}
where
\begin{align}
d\left(T_{b}\right) & \triangleq\frac{\left(\delta_{Z}^{0}\right)'\left\{ \left(Z_{0}'MZ_{0}\right)-\left(Z'_{0}MZ_{2}\right)\left(Z'_{2}MZ_{2}\right)^{-1}\left(Z'_{2}MZ_{0}\right)\right\} \delta_{Z}^{0}}{\left|T_{b}-T_{b}^{0}\right|},\label{eq. A.2.6, Pre}
\end{align}
and we arbitrarily define $d\left(T_{b}\right)=\left(\delta_{Z}^{0}\right)'\delta_{Z}^{0}$
when $T_{b}=T_{b}^{0}$. Let $d_{T}=T\inf_{\left|T_{b}-T_{b}^{0}\right|>TK}d\left(T_{b}\right)$;
it is positive and bounded away from zero by Lemma \ref{Lemma d_CT, Pre}
below. Then 
\begin{align}
P\left(\left|\widehat{\lambda}_{b}-\lambda_{0}\right|>K\right) & =P\left(\left|\widehat{T}_{b}-T_{b}^{0}\right|>TK\right)\nonumber \\
 & \leq P\left(\sup_{\left|T_{b}-T_{b}^{0}\right|>TK}\left|g_{e}\left(T_{b}\right)\right|\geq\inf_{\left|T_{b}-T_{b}^{0}\right|>TK}\left|T_{b}-T_{b}^{0}\right|d\left(T_{b}\right)\right)\nonumber \\
 & \leq P\left(\sup_{p+2\leq T_{b}\leq T-p-2}\left|g_{e}\left(T_{b}\right)\right|\geq TK\inf_{\left|T_{b}-T_{b}^{0}\right|>TK}d\left(T_{b}\right)\right)\nonumber \\
 & =P\left(d_{T}^{-1}\sup_{p+2\leq T_{b}\leq T-p-2}\left|g_{e}\left(T_{b}\right)\right|\geq K\right).\label{eq. last line of A.2.6.b, Pre}
\end{align}
We can write the first term of $g_{e}\left(T_{b}\right)$ as
\begin{align}
2\left(\delta^{0}\right)'\left(Z'_{0}MZ_{2}\right)\left(Z'_{2}MZ_{2}\right)^{-1/2}\left(Z'_{2}MZ_{2}\right)^{-1/2}Z_{2}Me & .\label{Eq. A.2.8, Pre}
\end{align}
For the stochastic regressors, Theorem \ref{=00005BJ=000026P,-Theorem-5.4.2=00005D}
implies that for any $3\leq j\leq p+2,$ $\left(Z_{2}e\right)_{j,1}/\sqrt{h}=O_{p}\left(1\right)$
and for any $3\leq i\leq q+p+2,$ $\left(Xe\right)_{i,1}/\sqrt{h}=O_{p}\left(1\right),$
since these estimates include a positive fraction of the data. We
can use the above expression for $X'X$ to verify that $Z'_{2}MZ_{2}$
and $Z'_{0}MZ_{2}$ are $O_{p}\left(1\right)$. Then,
\begin{align*}
\sup_{T_{b}}\left(Z'_{0}MZ_{2}\right)\left(Z'_{2}MZ_{2}\right)^{-1}\left(Z'_{2}MZ_{0}\right) & \leq Z'_{0}MZ_{0}=O_{p}\left(1\right),
\end{align*}
by Lemma \ref{Lemma LLNs 1}. Next, note that the first two elements
of the vector $X'e$ and $Z'_{2}e$ are $O_{p}\left(h^{1/2}\right)$
{[}recall \eqref{X'e}{]}. By Assumption \ref{Assumption 1, CT}-(iii)
and the inequality
\begin{align*}
\sup_{T_{b}}\left\Vert \left(Z'_{2}MZ_{2}\right)^{-1/2}Z_{2}Me\right\Vert  & \leq\sup_{T_{b}}\left\Vert \left(Z'_{2}MZ_{2}\right)^{-1/2}\right\Vert \sup_{T_{b}}\left\Vert Z_{2}Me\right\Vert ,
\end{align*}
we have that $\left(Z'_{2}MZ_{2}\right)^{-1/2}Z_{2}Me$ is $O_{p}(h^{1/2})$
uniformly in $T_{b}$ since the last $q+p$ (resp., $p$) elements
of $X'e$ (resp., $Z'_{2}e$) are $o_{p}\left(1\right)$ locally uniformly
in time. Therefore, uniformly over $p+2\leq T_{b}\leq T-p-2$, the
overall expression in \eqref{Eq. A.2.8, Pre} is $O_{p}(h^{1/2})$.
As for the second term of \eqref{eq. A.2.3, Casini (2015a)-2}, $Z'_{0}Me=O_{p}(h^{1/2}).$
The first term in \eqref{eq. A.2.4, Casini (2015a)-2} is uniformly
negligible and so is the last. Therefore, combining these results
we can show that $\sup_{T_{b}}|g_{e}\left(T_{b}\right)|=O_{p}(\sqrt{h})$.
Using Lemma \ref{Lemma d_CT, Pre} below, we have $P(d_{T}^{-1}\sup_{p+2\leq T_{b}\leq T-p-2}|g_{e}\left(T_{b}\right)|\geq K)\leq\varepsilon$,
which shows that $\widehat{\lambda}_{b}\overset{P}{\rightarrow}\lambda_{0}$.
$\square$
\begin{lem}
\label{Lemma d_CT, Pre}Let $d_{B}=\inf_{\left|T_{b}-T_{b}^{0}\right|>TB}Td\left(T_{b}\right).$
There exists a $\kappa>0$ and for every $\varepsilon>0,$ there exists
a $B<\infty$ such that $P(d_{B}\geq\kappa)\leq1-\varepsilon.$
\end{lem}
\noindent\textit{Proof.} Assuming $N_{b}\leq N_{b}^{0}$ and following
the same steps as in Lemma \ref{Lemma r_CT, Prop 1-2} (but replacing
$R$ by $\overline{R}$)
\begin{align*}
Td\left(T_{b}\right) & \geq T\left(\delta_{Z}^{0}\right)'\overline{R}'\frac{X'_{\Delta}X_{\Delta}}{T_{b}^{0}-T_{b}}\left(X'_{2}X_{2}\right)^{-1}\left(X'_{0}X_{0}\right)\overline{R}\left(\delta_{Z}^{0}\right)\\
 & =\left(\delta_{Z}^{0}\right)'\overline{R}'\frac{X'_{\Delta}X_{\Delta}}{B}\left(X'_{2}X_{2}\right)^{-1}\left(X'_{0}X_{0}\right)\overline{R}\left(\delta_{Z}^{0}\right).
\end{align*}
Under Assumption \ref{Assumption 1, CT}-(iii) and in view of \eqref{X Expression, Pre},
 $X'_{\Delta}X_{\Delta}$ is positive definite: for the $p\times p$
lower-right sub-block apply Lemma \ref{Lemma LLNs 1} as in the proof
of Lemma \ref{Lemma r_CT, Prop 1-2}, whereas for the remaining elements
of $X'_{\Delta}X_{\Delta}$ the result follows from the convergence
of approximations to Riemann sums. Note that $X'_{2}X_{2}$ and $X'_{0}X_{0}$
are $O_{p}\left(1\right)$. It follows that 
\begin{align*}
Td\left(T_{b}\right) & \geq\left(\delta_{Z}^{0}\right)'\overline{R}'\frac{X'_{\Delta}X_{\Delta}}{N}\left(X'_{2}X_{2}\right)^{-1}\left(X'_{0}X_{0}\right)\overline{R}\delta_{Z}^{0}\geq\kappa>0.
\end{align*}
The result follows choosing $B>0$ such that $P\left(d_{B}\geq\kappa\right)$
is larger than $1-\varepsilon.$ $\square$

\noindent\textit{Proof of part (ii) of Proposition \ref{Prop 1-2 Asym - Pre}.}
We introduce again
\begin{align*}
\mathbf{D}_{K,T} & =\left\{ T_{b}:\,N\eta\leq N_{b}\leq N\left(1-\eta\right),\,\left|N_{b}^{0}-N_{b}\right|>KT^{-1}\right\} ,
\end{align*}
 and observe that it is enough to show that $P(\sup_{T_{b}\in\mathbf{D}_{K,T}}Q_{T}\left(T_{b}\right)\geq Q_{T}\left(T_{b}^{0}\right))<\varepsilon,$
or
\begin{align}
P & \left(\sup_{T_{b}\in\mathbf{D}_{K,T}}h^{-1}g_{e}\left(T_{b}\right)\geq\inf_{T_{b}\in\mathbf{D}_{K,T}}h^{-1}\left|T_{b}-T_{b}^{0}\right|d\left(T_{b}\right)\right)<\varepsilon.\label{A.9, Pre}
\end{align}
 By Lemma \ref{Lemma A1, Casini (2015a)},
\begin{align*}
\inf_{T_{b}\in\mathbf{D}_{K,T}}d\left(T_{b}\right) & \geq\inf_{T_{b}\in\mathbf{D}_{K,T}}\left(\delta_{Z}^{0}\right)'\overline{R}'\frac{X'_{\Delta}X_{\Delta}}{T_{b}^{0}-T_{b}}\left(X'_{2}X_{2}\right)^{-1}\left(X'_{0}X_{0}\right)\overline{R}\delta_{Z}^{0}.
\end{align*}
For the $\left(q+p\right)\times\left(q+p\right)$ lower right sub-block
of $X'_{\Delta}X_{\Delta}$ the arguments of Proposition \ref{Proposition 2, (Rate of Convergence)}
apply: $(h(T_{b}^{0}-T_{b}))^{-1}$ ~ $[X'_{\Delta}X_{\Delta}]_{\left\{ \cdot,\,\left(q+p\right)\times\left(q+p\right)\right\} }$
is bounded away from zero for all $T_{b}\in\mathbf{D}_{K,T}$ by choosing
$K$ large enough (recall $|T_{b}^{0}-T_{b}|>K$), where $\left[A\right]_{\left\{ \cdot,\,i\times j\right\} }$
is the $i\times j$ lower right sub-block of $A$. Furthermore, this
approximation is uniform in $T_{b}$ by Assumption \ref{Assumption 4 Eigenvalue}.
It remains to deal with the upper left sub-block of $X'_{\Delta}X_{\Delta}.$
Consider its $\left(1,\,1\right)$-th element. It is given by $\sum_{k=T_{b}+1}^{T_{b}^{0}}(h^{1/2})^{2}.$
Thus $(h(T_{b}^{0}-T_{b}))^{-1}\sum_{k=T_{b}+1}^{T_{b}^{0}}(h^{1/2})^{2}>0.$
The same argument applies to the $\left(2,\,2\right)$-th element
of the upper left sub-block of $X'_{\Delta}X_{\Delta}.$ The latter
results imply that $\inf_{T_{b}\in\mathbf{D}_{K,T}}Td\left(T_{b}\right)$
is bounded away from zero. It remains to show that $\sup_{T_{b}\in\mathbf{D}_{K,T}}(h|T_{b}-T_{b}^{0}|)^{-1}g_{e}\left(T_{b}\right)$
is small when $T$ is large. Recall that the terms $Z_{2}$ and $Z_{0}$
involve a positive fraction $N\eta$ of the data. We can apply Lemma
\ref{Lemma LLNs 1} to those elements which involve the stochastic
regressors only, whereas the other terms are dealt with directly using
the definition of $X'e$ in \eqref{X'e}. Consider the first term
of $g_{e}\left(T_{b}\right)$. Using the same steps which led to \eqref{eq. 21, Casini (2015a)-1},
we have
\begin{align}
 & \left|2\left(\delta_{Z}^{0}\right)'\left(Z'_{0}MZ_{2}\right)\left(Z'_{2}MZ_{2}\right)^{-1}Z_{2}Me-2\left(\delta_{Z}^{0}\right)'\left(Z'_{0}Me\right)\right|\nonumber \\
 & =\left|\left(\delta_{Z}^{0}\right)'Z'_{\Delta}Me\right|+\left|\left(\delta_{Z}^{0}\right)'\left(Z'_{\Delta}MZ_{2}\right)\left(Z'_{2}MZ_{2}\right)^{-1}\left(Z_{2}Me\right)\right|.\label{eq. 21, Pre}
\end{align}
 We can apply Lemma \ref{Lemma LLNs 1} to the terms that do not involve
$|N_{b}-N_{b}^{0}|$ but only stochastic regressors. Next consider
the first term of 
\begin{align*}
\left(h\left(T_{b}^{0}-T_{b}\right)\right)^{-1}\left(\delta_{Z}^{0}\right)'\left(Z'_{\Delta}MZ_{2}\right) & =\frac{\left(\delta_{Z}^{0}\right)'\left(Z'_{\Delta}Z_{\Delta}\right)}{h\left(T_{b}^{0}-T_{b}\right)}\\
 & \quad-\left(\delta_{Z}^{0}\right)'\left(\frac{Z'_{\Delta}X_{\Delta}}{h\left(T_{b}^{0}-T_{b}\right)}\left(X'X\right)^{-1}X'Z_{2}\right).
\end{align*}
Applying the same manipulations as those used above for the $p\times p$
lower right sub-block of $Z'_{\Delta}Z_{\Delta}$, we have $(h(T_{b}^{0}-T_{b}))^{-1}[Z'_{\Delta}Z_{\Delta}]_{\left\{ \cdot,\,p\times p\right\} }=O_{p}\left(1\right),$
since there are $T_{b}^{0}-T_{b}$ summands whose conditional first
moments are each $O\left(h\right)$. The $O_{p}\left(1\right)$ result
is uniform by Assumption \ref{Assumption 4 Eigenvalue}. The same
argument holds for the corresponding sub-block of $Z'_{\Delta}X_{\Delta}/(h(T_{b}^{0}-T_{b}))$.
Hence, as $h\downarrow0$ the second term above is $O_{p}\left(1\right).$
Next, consider the upper left $2\times2$ block of $Z'_{\Delta}Z_{\Delta}$
(the same argument holds true for $Z'_{\Delta}X_{\Delta}$). Note
that the predictable variable $Y_{\left(k-1\right)h}$ in the $\left(2,\,2\right)$-th
element can be treated as locally constant after multiplying by $h^{1/2}$
(recall $h^{1/2}Y_{\left(k-1\right)h}=\widetilde{Y}_{\left(k-1\right)h}=O_{p}\left(1\right)$
by Assumption \ref{Assumption Localization, Pre}), 
\begin{align*}
\sum_{k=T_{b}+1}^{T_{b}^{0}}\left(Y_{\left(k-1\right)h}h\right)^{2} & =\sum_{k=T_{b}+1}^{T_{b}^{0}}\left(\widetilde{Y}_{\left(k-1\right)h}h^{1/2}\right)^{2}\leq C\sum_{k=T_{b}+1}^{T_{b}^{0}}h,
\end{align*}
 where $C=\sup_{k}|\widetilde{Y}_{\left(k-1\right)h}^{2}|$ is a fixed
constant given the localization in Assumption \ref{Assumption Localization, Pre}.
Thus, when multiplied by $(h(T_{b}^{0}-T_{b}))^{-1}$, the $\left(2,\,2\right)$-th
element of $Z'_{\Delta}Z_{\Delta}$  is $O_{p}\left(1\right).$ The
same reasoning can be applied to the corresponding $\left(1,\,1\right)$-th
element. Next, let us consider the cross-products between the semimartingale
regressors and the predictable regressors. Consider any $3\leq j\leq p+2,$
\begin{align*}
\frac{1}{h\left(T_{b}^{0}-T_{b}\right)}\sum_{k=T_{b}+1}^{T_{b}^{0}}z_{kh}^{\left(2\right)}z_{kh}^{\left(j\right)} & =\frac{1}{h\left(T_{b}^{0}-T_{b}\right)}\sum_{k=T_{b}+1}^{T_{b}^{0}}\left(\widetilde{Y}_{\left(k-1\right)h}h^{1/2}\right)z_{kh}^{\left(j\right)}\\
 & =\frac{1}{T_{b}^{0}-T_{b}}\sum_{k=T_{b}+1}^{T_{b}^{0}}\widetilde{Y}_{\left(k-1\right)h}\frac{z_{kh}^{\left(j\right)}}{\sqrt{h}}.
\end{align*}
Since $z_{kh}^{\left(j\right)}/\sqrt{h}$ is i.n.d. with zero mean
and finite variance and $\widetilde{Y}_{\left(k-1\right)h}$ is $O_{p}\left(1\right)$
by Assumption \ref{Assumption Localization, Pre}, Assumption \ref{Assumption 4 Eigenvalue}
implies that we can find a $K$ large enough such that the right hand
side is $O_{p}\left(1\right)$ uniformly in $T_{b}$. The same argument
applies to $(Z'_{\Delta}Z_{\Delta})_{1,j},$ $3\leq j\leq p+2.$ This
shows that the term $(Z'_{\Delta}X_{\Delta}/(h(T_{b}^{0}-T_{b})))\left(X'X\right)^{-1}X'Z_{2}$
is bounded and so is $Z'_{\Delta}X_{\Delta}/(h(T_{b}^{0}-T_{b}))$
using the same reasoning. Thus, $(h(T_{b}^{0}-T_{b}))^{-1}\left(\delta_{Z}^{0}\right)'\left(Z'_{\Delta}MZ_{2}\right)$
is $O_{p}\left(1\right).$ By the same arguments as before, we can
use Theorem \ref{=00005BJ=000026P,-Theorem-5.4.2=00005D} to show
that the second term of \eqref{eq. 21, Pre} is $O_{p}(h^{1/2})$
when multiplied by $(h(T_{b}^{0}-T_{b}))^{-1}$ since the last term
involves a positive fraction of the data. Now, expand the $\left(p+2\right)$-dimensional
vector $Z'_{\Delta}Me$ as follows
\begin{align*}
\frac{Z'_{\Delta}Me}{h\left(T_{b}^{0}-T_{b}\right)} & =\frac{1}{h\left(T_{b}^{0}-T_{b}\right)}\sum_{k=T_{b}+1}^{T_{b}^{0}}z_{kh}e_{kh}\\
 & \quad-\frac{1}{h\left(T_{b}^{0}-T_{b}\right)}\left(\sum_{k=T_{b}+1}^{T_{b}^{0}}z_{kh}x'_{kh}\right)\left(X'X\right)^{-1}\left(X'e\right).
\end{align*}
The arguments for the last $p$ elements are the same as above and
yield {[}recall \eqref{Eq. h^-1 ZX XX Xe}{]} 
\begin{align*}
\frac{\left[Z'_{\Delta}Me\right]_{\left\{ \cdot,p\right\} }}{h\left(T_{b}^{0}-T_{b}\right)} & =o_{p}\left(K^{-1}\right)-O_{p}\left(1\right)O_{p}\left(h^{1/2}\right),
\end{align*}
where we recall that by Assumption \ref{Assumption 1, CT}-(iv) $\Sigma_{Ze,N_{b}^{0}}=0$.
Note that the convergence is uniform over $T_{b}$ by Lemma \ref{Lemma, Spot Uniform Approx}.
We now consider the first two elements of $Z'_{\Delta}e$: 
\begin{align*}
\left|\sum_{k=T_{b}+1}^{T_{b}^{0}}z_{kh}^{\left(2\right)}e_{kh}\right| & =\left|\sum_{k=T_{b}+1}^{T_{b}^{0}}h^{1/2}Y_{\left(k-1\right)h}h^{1/2}e_{kh}\right|\leq A\sum_{k=T_{b}+1}^{T_{b}^{0}}\left|\widetilde{Y}_{\left(k-1\right)h}h^{1/2}e_{kh}\right|,
\end{align*}
 for some positive $A<\infty$. Noting that $e_{kh}/\sqrt{h}\sim\textrm{i.n.d.}\mathscr{N}(0,\,\sigma_{e,k-1}^{2})$,
we have
\begin{align*}
\left(h\left(T_{b}^{0}-T_{b}\right)\right)^{-1}\sum_{k=T_{b}+1}^{T_{b}^{0}}z_{kh}^{\left(2\right)}e_{kh} & \leq C\left(\left(T_{b}^{0}-T_{b}\right)^{-1}\sum_{k=T_{b}+1}^{T_{b}^{0}}\left|e_{kh}/h^{1/2}\right|\right)
\end{align*}
 where $C=\sup_{k}|\widetilde{Y}_{\left(k-1\right)h}|$ is finite
by Assumption \ref{Assumption Localization, Pre}. Choose $K$ large
enough such that the probability that the right-hand side is larger
than $B/3N$ is less than $\varepsilon$. For the first element of
$Z'_{\Delta}e$ the argument is the same and thus 
\[
P\left(\left(h\left(T_{b}^{0}-T_{b}\right)\right)^{-1}\sum_{k=T_{b}+1}^{T_{b}^{0}}z_{kh}^{\left(1\right)}e_{kh}>\frac{B}{3N}\right)\leq\varepsilon,
\]
 when $K$ is large. For the last product in the second term of $Z'_{\Delta}Me/h$
the argument is easier. This component includes a positive fraction
of data and thus
\begin{align}
\sum_{k=1}^{T}x_{kh}^{\left(1\right)}e_{kh} & =\sum_{k=1}^{T}h^{1/2}e_{kh}=h^{1/2}O_{p}\left(1\right),\label{Eq. X1e ucp}
\end{align}
using the result $\sum_{k=1}^{\left\lfloor t/h\right\rfloor }e_{kh}\overset{\textrm{u.c.p.}}{\Rightarrow}\int_{0}^{t}\sigma_{e,s}dW_{e,s}.$
A similar argument applies to $x_{kh}^{\left(2\right)}e_{kh}$ by
using in addition the localization Assumption \ref{Assumption Localization, Pre}.
Combining the above derivations, we have
\begin{align}
\frac{1}{h\left(T_{b}^{0}-T_{b}\right)}g_{e}\left(T_{b}\right) & =\frac{1}{h\left(T_{b}^{0}-T_{b}\right)}\left(\delta_{Z}^{0}\right)'2Z'_{\Delta}e+o_{p}\left(1\right).\label{Eq. A.12, Pre}
\end{align}
In order to prove
\begin{align*}
P & \left(\sup_{T_{b}\in\mathbf{D}_{K,T}}\left(h\left(T_{b}^{0}-T_{b}\right)\right)^{-1}g_{e}\left(T_{b}\right)\geq\inf_{T_{b}\in\mathbf{D}_{K,T}}h^{-1}d\left(T_{b}\right)\right)<\varepsilon,
\end{align*}
we can use \eqref{Eq. A.12, Pre}. To this end, we shall find a $K>0$,
such that 
\begin{align}
P & \left(\sup_{T_{b}\leq T_{b}^{0}-\frac{K}{N}}\left|\mu_{\delta}^{0}\frac{2}{h}\left(T_{b}^{0}-T_{b}\right)^{-1}\sum_{k=T_{b}+1}^{T_{b}^{0}}z_{kh}^{\left(1\right)}e_{kh}\right|>\frac{B}{3N}\right)\label{A.12 +1, Pre}\\
 & \leq P\left(\sup_{T_{b}\leq T_{b}^{0}-\frac{K}{N}}\left(T_{b}^{0}-T_{b}\right)^{-1}\left|\sum_{k=T_{b}+1}^{T_{b}^{0}}\frac{e_{kh}}{\sqrt{h}}\right|>\frac{B}{6\left|\mu_{\delta}^{0}\right|N}\right)<\frac{\varepsilon}{3}.\nonumber 
\end{align}
Recalling that $e_{kh}/h^{1/2}\sim\mathscr{N}(0,\,\sigma_{e,k-1}^{2}),$
the Hájek-Réiny inequality yields 
\begin{align*}
P\left(\sup_{T_{b}\leq T_{b}^{0}-\frac{K}{N}}\left(T_{b}^{0}-T_{b}\right)^{-1}\left|\sum_{k=T_{b}+1}^{T_{b}^{0}}\frac{e_{kh}}{\sqrt{h}}\right|>\frac{B}{6\left|\mu_{\delta}^{0}\right|N}\right) & \leq A\frac{36\left(\mu_{\delta}^{0}\right)^{2}N^{2}}{B^{2}}\frac{1}{KN^{-1}}.
\end{align*}
We can choose $K$ sufficiently large such that the right-hand side
is less than $\varepsilon/3.$ The same bound holds for the second
element of $Z'_{\Delta}e$. Next, by equation \eqref{Eq. Haek Prop 2},
\begin{align*}
P\left(\sup_{T_{b}\leq T_{b}^{0}-\frac{K}{N}}\frac{1}{h\left(T_{b}^{0}-T_{b}\right)}\left\Vert 2\left(\delta_{Z}^{0}\right)'\sum_{k=T_{b}+1}^{T_{b}^{0}}\left[Z'_{\Delta}e\right]_{\left\{ \cdot,p\right\} }\right\Vert >\frac{B}{3N}\right) & <\frac{\varepsilon}{3},
\end{align*}
 since for each $j=3,\ldots,\,p$, $\{z_{kh}^{\left(j\right)}e_{kh}/h\}$
is i.n.d. with finite variance, and thus the result is implied by
the Hájek-Réiny inequality for large $K$. Using the latter results
into \eqref{Eq. A.12, Pre}, we have 
\begin{align*}
P\left(\sup_{T_{b}\leq T_{b}^{0}-\frac{K}{N}}\frac{1}{h\left(T_{b}^{0}-T_{b}\right)}\left\Vert 2\left(\delta_{Z}^{0}\right)'\sum_{k=T_{b}+1}^{T_{b}^{0}}z_{kh}e_{kh}\right\Vert >\frac{B}{N}\right) & <\varepsilon,
\end{align*}
 which verifies \eqref{A.9, Pre} and thus proves our claim. $\square$

\subsubsection{Proof of Theorem \ref{Theorem - Asy Dist - Pre}}

Part (i)-(ii) follows the same steps as in the proof of Proposition
\ref{Prop 3 Asym} part (i)-(ii) but using the results developed throughout
the proof of part (i)-(ii) of Proposition \ref{Prop 1-2 Asym - Pre}.
As for part (iii), we begin with the following lemma, where again
$\psi_{h}=h^{1-\kappa}$. Without loss of generality we set $B=1$
in Assumption \ref{Assumption 6 - Small Shifts}.
\begin{lem}
\label{Lemma 1, Pre}Under Assumption \ref{Assumption Localization, Pre},
uniformly in $T_{b}$,
\begin{align*}
\left(Q_{T}\left(T_{b}\right)-Q_{T}\left(T_{b}^{0}\right)\right)/\psi_{h} & =-\delta_{h}\left(Z_{\Delta}'Z_{\Delta}/\psi_{h}\right)\delta_{h}\pm2\delta_{h}'\left(Z'_{\Delta}\widetilde{e}/\psi_{h}\right)\\
 & \quad+O_{p}\left(h^{3/4\wedge1-\kappa/2}\right).
\end{align*}
\end{lem}
\noindent\textit{Proof.} By the definition of $Q_{T}\left(T_{b}\right)-Q_{T}\left(T_{b}^{0}\right)$
and Lemma \ref{Lemma A.3},
\begin{align}
Q_{T}\left(T_{b}\right) & -Q_{T}\left(T_{b}^{0}\right)\label{eq. (45), Pre Lemma 1}\\
 & =-\delta_{h}'\left\{ Z_{\Delta}'MZ_{\Delta}+\left(Z'_{\Delta}MZ_{2}\right)\left(Z'_{2}MZ_{2}\right)^{-1}\left(Z'_{2}MZ_{\Delta}\right)\right\} \delta_{h}\nonumber \\
 & +g_{e}\left(T_{b},\,\delta_{h}\right).\nonumber 
\end{align}
We can expand the first term of \eqref{eq. (45), Pre Lemma 1} as
follows
\begin{align}
\delta_{h}'Z_{\Delta}'MZ_{\Delta}\delta_{h} & =\delta_{h}'Z_{\Delta}'Z_{\Delta}\delta_{h}-\delta_{h}'A\delta_{h},\label{Eq. ZMZ+A, Pre}
\end{align}
where $A=Z_{\Delta}'X\left(X'X\right)^{-1}X'Z_{\Delta}$. We show
that $\delta_{h}'A\delta_{h}$ is uniformly of higher order than $\delta_{h}'Z_{\Delta}'Z_{\Delta}\delta_{h}$.
The cross-products between the semimartingale and the predictable
regressors (i.e., the $p\times2$ lower-left sub-block of $Z_{\Delta}'X$)
are $o_{p}\left(1\right)$, as can be easily verified. Lemma \ref{Lemma Cross Product Null}
provides the formal statement of the result for $Z_{\Delta}'Z_{\Delta}$.
Hence, the result carries over to $Z_{\Delta}'X$ with no changes
and by symmetry also to is the $2\times p$ upper-right block. This
allows us to treat the $2\times2$ upper-left block and the $p\times p$
lower-right block of elements, such as those in $A$ separately. By
Lemma \ref{Lemma LLNs 1}, $\left(X'X\right)^{-1}=O_{p}\left(1\right)$.
Using Proposition \ref{Prop 3 Asym}-(ii), we let $N_{b}-N_{b}^{0}=K\psi_{h}$.
By the Burkh\"{o}lder-Davis-Gundy inequality, we have standard estimates
for local volatility so that 
\begin{align*}
\left\Vert \mathbb{E}\left(\widehat{\Sigma}_{ZX}^{\left(i,j\right)}\left(T_{b},\,T_{b}^{0}\right)-\Sigma_{ZX,\left(T_{b}^{0}-1\right)h}^{\left(i,j\right)}|\,\mathscr{F}_{\left(T_{b}^{0}-1\right)h}\right)\right\Vert  & \leq Kh^{1/2},
\end{align*}
with $3\leq i\leq p+2$ and $3\leq j\leq q+p+2$ which in turn implies
$[Z_{\Delta}'X_{\Delta}]_{\left\{ \cdot,p\times p\right\} }=O_{p}(1/(h(T_{b}^{0}-T_{b})))$.
The same bound applies to the corresponding blocks of $Z_{\Delta}'Z_{\Delta}$
and $X_{\Delta}'Z_{\Delta}$. Now we focus on the $\left(2,\,2\right)$-th
element of $A$. First notice that
\begin{align*}
\left(Z_{\Delta}'X\right)_{2,2} & =\sum_{k=T_{b}+1}^{T_{b}^{0}}z_{kh}^{\left(2\right)}x_{kh}^{\left(2\right)}=\sum_{k=T_{b}+1}^{T_{b}^{0}}\left(\widetilde{Y}_{\left(k-1\right)h}\right)^{2}h.
\end{align*}
By a localization argument (cf. Assumption \ref{Assumption Localization, Pre}),
$\widetilde{Y}_{\left(k-1\right)h}$ is bounded. Then, since the number
of summands grows at a rate $T^{\kappa},$ we have $\left(Z_{\Delta}'X\right)_{2,2}=O_{p}(Kh^{1-\kappa}).$
The same proof can be used for $\left(Z_{\Delta}'X\right)_{1,1},$
which gives $\left(Z_{\Delta}'X\right)_{1,1}=O_{p}(Kh^{1-\kappa}).$
Thus, in view of \eqref{Eq. ZMZ, in Lemma 1, Pre}, we conclude that
\eqref{Eq. ZMZ+A, Pre} when divided by $\psi_{h}$ is such that 
\begin{align}
\delta_{h}'Z_{\Delta}'MZ_{\Delta}\delta_{h}/\psi_{h} & =\delta_{h}'Z_{\Delta}'Z_{\Delta}\delta_{h}/\psi_{h}-\delta_{h}'Z_{\Delta}'X\left(X'X\right)^{-1}X'Z_{\Delta}\delta_{h}/\psi_{h}\nonumber \\
 & =\psi_{h}^{-1}\left(\delta^{0}\right)'Z_{\Delta}'Z_{\Delta}\delta^{0}-\psi_{h}^{-1}h^{1/2}O_{p}\left(h^{2\left(1-\kappa\right)}\right).\label{Eq. ZMZ, in Lemma 1, Pre}
\end{align}
 For the second term of \eqref{eq. (45), Pre Lemma 1}, we have 
\begin{align}
\psi_{h}^{-1} & h^{1/2}\left(\delta^{0}\right)'\left\{ \left(Z'_{\Delta}MZ_{2}\right)\left(Z'_{2}MZ_{2}\right)^{-1}\left(Z'_{2}MZ_{\Delta}\right)\right\} \delta^{0}\label{eq. ACRSC1 Pre}\\
 & =\psi_{h}^{-1}h^{1/2}\left\Vert \delta_{0}\right\Vert ^{2}O_{p}\left(\psi_{h}\right)O_{p}\left(1\right)O_{p}\left(\psi_{h}\right)\leq K\psi_{h}^{-1}h^{1/2}O_{p}\left(h^{2\left(1-\kappa\right)}\right)\nonumber 
\end{align}
uniformly in $T_{b}$, which follows from applying the same reasoning
used for $Z_{\Delta}'\left(I-M\right)Z_{\Delta}$ above to each of
these three elements. Finally, consider the stochastic term $g_{e}\left(T_{b},\,\delta_{h}\right)$.
We have
\begin{align}
g_{e}\left(T_{b},\,\delta_{h}\right) & =2\delta_{h}'\left(Z'_{0}MZ_{2}\right)\left(Z'_{2}MZ_{2}\right)^{-1}Z_{2}Me-2\delta_{h}'\left(Z'_{0}Me\right)\label{Eq. ge in Lemma 1, Pre}\\
 & +e'MZ_{2}\left(Z'_{2}MZ_{2}\right)^{-1}Z_{2}Me-e'MZ_{0}\left(Z'_{0}MZ_{0}\right)^{-1}Z'_{0}Me.\nonumber 
\end{align}
 Recall \eqref{X'e}, and $\sum_{k=T_{b}+1}^{T_{b}^{0}}x_{kh}e_{kh}=h^{-1/4}\sum_{T_{b}+1}^{T_{b}^{0}}x_{kh}\widetilde{e}_{kh}$.
Introduce the following decomposition, 
\begin{align*}
\left(X'e\right)_{2,1} & =\sum_{k=1}^{T_{b}^{0}-\left\lfloor T^{\kappa}\right\rfloor }x_{kh}^{\left(2\right)}\widetilde{e}_{kh}+h^{-1/4}\sum_{k=T_{b}^{0}-\left\lfloor T^{\kappa}\right\rfloor +1}^{T_{b}^{0}+\left\lfloor T^{\kappa}\right\rfloor }x_{kh}^{\left(2\right)}\widetilde{e}_{kh}+\sum_{k=T_{b}^{0}+\left\lfloor T^{\kappa}\right\rfloor +1}^{T}x_{kh}^{\left(2\right)}\widetilde{e}_{kh},
\end{align*}
where $\widetilde{e}_{kh}\sim\textrm{i.n.d.\,}\mathscr{N}(0,\,\sigma_{e,k-1}^{2}h)$.
The first and third terms are $O_{p}(h^{1/2})$ in view of \eqref{Eq. X1e ucp}.
The term in the middle is $h^{3/4}\sum_{k=T_{b}^{0}-\left\lfloor T^{\kappa}\right\rfloor +1}^{T_{b}^{0}+\left\lfloor T^{\kappa}\right\rfloor }\widetilde{Y}_{\left(k-1\right)h}h^{-1/2}\widetilde{e}_{kh}$,
which involves approximately $2T^{\kappa}$ summands. Since $\widetilde{Y}_{\left(k-1\right)h}$
is bounded by the localization procedure, 
\begin{align*}
h^{3/4}\frac{T^{\kappa/2}}{T^{\kappa/2}}\sum_{k=T_{b}^{0}-\left\lfloor T^{\kappa}\right\rfloor }^{T_{b}^{0}+\left\lfloor T^{\kappa}\right\rfloor }\widetilde{Y}_{\left(k-1\right)h}\frac{\widetilde{e}_{kh}}{\sqrt{h}}=h^{3/4}T^{\kappa/2}O_{p}\left(1\right) & ,
\end{align*}
 or 
\[
h^{-1/4}\sum_{k=T_{b}^{0}-\left\lfloor T^{\kappa}\right\rfloor }^{T_{b}^{0}+\left\lfloor T^{\kappa}\right\rfloor }x_{kh}^{\left(2\right)}\widetilde{e}_{kh}=h^{3/4-\kappa/2}O_{p}\left(1\right).
\]
This implies that $\left(X'e\right)_{2,1}$ is $O_{p}(h^{1/2\wedge3/4-\kappa/2}).$
The same observation holds for $\left(X'e\right)_{1,1}.$ Therefore,
one follows the same steps as in the concluding part of the proof
of Lemma \ref{Lemma 1} {[}cf. equation \eqref{eq. ACRSC5} and the
derivations thereafter{]}. That is, for the first two terms of $g_{e}(T_{b},\,\delta_{h}),$
using $Z'_{0}MZ_{2}=Z'_{2}MZ_{2}\pm Z'_{\Delta}MZ_{2}$, we have
\begin{align}
2h^{1/4}\left(\delta^{0}\right)' & \left(Z'_{0}MZ_{2}\right)\left(Z'_{2}MZ_{2}\right)^{-1}Z_{2}Me-2h^{1/4}\left(\delta^{0}\right)'\left(Z'_{0}Me\right)\nonumber \\
 & =2h^{1/4}\left(\delta^{0}\right)'Z'{}_{\Delta}Me\pm2h^{1/4}\left(\delta^{0}\right)'Z'{}_{\Delta}MZ_{2}\left(Z'_{2}MZ_{2}\right)^{-1}Z'_{2}Me.\label{Eq. CRSC6-1}
\end{align}
The last term above when multiplied by $\psi_{h}^{-1}$ is such that
\begin{align*}
\psi_{h}^{-1}2h^{1/4}\left(\delta^{0}\right)'Z'{}_{\Delta}MZ_{2}\left(Z'_{2}MZ_{2}\right)^{-1}Z'_{2}Me & =\left\Vert \delta^{0}\right\Vert O_{p}\left(1\right)O_{p}\left(h^{1\wedge5/4-\kappa/2}\right),
\end{align*}
 where we have used the fact that $Z'_{\Delta}MZ_{2}/\psi_{h}=O_{p}\left(1\right)$.
For the first term of \eqref{Eq. CRSC6-1}, 
\begin{align*}
2h^{1/4} & \left(\delta^{0}\right)'Z'{}_{\Delta}Me/\psi_{h}\\
 & =2h^{1/4}\left(\delta^{0}\right)'Z'{}_{\Delta}e/\psi_{h}-2h^{1/4}\left(\delta^{0}\right)'Z'{}_{\Delta}X\left(X'X\right)^{-1}X'e/\psi_{h}\\
 & =2h^{1/4}\left(\delta^{0}\right)'Z'{}_{\Delta}e-2\left(\delta^{0}\right)'O_{p}\left(1\right)O_{p}\left(h^{1\wedge5/4-\kappa/2}\right).
\end{align*}
 As in the proof of Lemma \ref{Lemma 1}, we can now use part (i)
of the theorem so that the difference between the terms on the second
line of $g_{e}\left(T_{b},\,\delta_{h}\right)$ is negligible. That
is, we can find a $c_{T}$ sufficiently small such that, 
\begin{align*}
\psi_{h}^{-1}\left[e'MZ_{2}\left(Z'_{2}MZ_{2}\right)^{-1}Z_{2}Me-e'MZ_{0}\left(Z'_{0}MZ_{0}\right)^{-1}Z'_{0}Me\right] & =o_{p}\left(c_{T}h\right).
\end{align*}
This leads to
\begin{align*}
g_{e}\left(T_{b},\,\delta_{h}\right)/\psi_{h} & =2h^{1/4}\left(\delta^{0}\right)'Z'{}_{\Delta}e/\psi_{h}+O_{p}\left(h^{3/4\wedge1-\kappa/2}\right)\\
 & \quad+\left\Vert \delta^{0}\right\Vert O_{p}\left(h^{3/4\wedge1-\kappa/2}\right)+o_{p}\left(c_{T}h\right),
\end{align*}
for a sufficiently small $c_{T}$. This together with \eqref{Eq. ZMZ, in Lemma 1, Pre}
and \eqref{eq. ACRSC1 Pre} yields, 
\begin{align*}
\psi_{h}^{-1}\left(Q_{T}\left(T_{b}\right)-Q_{T}\left(T_{b}^{0}\right)\right) & =-\delta_{h}\left(Z_{\Delta}'Z_{\Delta}/\psi_{h}\right)\delta_{h}\\
 & \quad\pm2\delta_{h}'\left(Z'_{\Delta}e/\psi_{h}\right)+O_{p}\left(h^{3/4\wedge1-\kappa/2}\right)+o_{p}\left(h^{1/2}\right),
\end{align*}
when $T$ is large, where $c_{T}$ is a sufficiently small number.
This concludes the proof. $\square$

\noindent\textit{Proof of part (iii) of Theorem \ref{Theorem - Asy Dist - Pre}.}
We proceed as in the proof of Theorem \ref{Theorem 1} and, hence,
some details are omitted. We again change the time scale $s\mapsto t\triangleq\psi_{h}^{-1}s$
on $\mathcal{D}\left(C\right)$ and observe that the re-parameterization
$\theta{}_{h}$ and $\sigma_{h,t}$ does not alter the result of Lemma
\ref{Lemma 1, Pre}. In addition, we have now, 
\begin{align*}
dZ_{\psi,s}^{\left(1\right)} & =\psi_{h}^{-1/2}\left(ds\right)^{1/2}=\left(ds\right)^{1/2},\qquad\\
dZ_{\psi,s}^{\left(2\right)} & =\psi_{h}^{-1/2}Y_{s-}ds=\psi_{h}^{-1/2}\widetilde{Y}_{s-}\left(ds\right)^{1/2}=\widetilde{Y}_{s-}\left(ds\right)^{1/2},
\end{align*}
 where the first equality in the second term above follows from $\widetilde{Y}_{\left(k-1\right)h}=h^{1/2}Y_{\left(k-1\right)h}$
on the old time scale. $N_{b}^{0}\left(v\right)$ varies on the time
horizon $[N_{b}^{0}-\left|v\right|,\,N_{b}^{0}+\left|v\right|]$ as
implied by $\mathcal{D}^{*}\left(C\right)$, as defined in Section
\ref{Section Asymptotic Distribution: Continuous Case}. Again, in
order to avoid clutter, we suppress the subscript $\psi_{h}$. We
then have equation \eqref{eq. CRSC_App_TS3}-\eqref{eq. CRSC_App_TS6}.
Consider $T_{b}\leq T_{b}^{0}$ (i.e., $v\leq0$). By Lemma \ref{Lemma 1, Pre},
there exists a $\overline{T}$ such that for all $T>\overline{T},$
\begin{align*}
\overline{Q}_{T}\left(\theta^{*}\right) & =-h^{-1/2}\delta_{h}'Z_{\Delta}'Z_{\Delta}\delta_{h}+h^{-1/2}2\delta_{h}'Z_{\Delta}'e+o_{p}\left(1\right)\\
 & =-\left(\delta^{0}\right)'\left(\sum_{k=T_{b}+1}^{T_{b}^{0}}z_{kh}z'_{kh}\right)\delta^{0}\\
 & \quad+2\left(\delta^{0}\right)'\left(h^{-1/2}\sum_{k=T_{b}+1}^{T_{b}^{0}}z_{kh}\widetilde{e}_{kh}\right)+o_{p}\left(1\right),
\end{align*}
 and note that this relationship corresponds to \eqref{Eq. QT after Lemma 1}.
As in the proof of Theorem \ref{Theorem 1} it is convenient to associate
to the continuous time index $t$ in $\mathcal{D}^{*}$, a corresponding
$\mathcal{D}^{*}$-specific index $t_{v}.$ We then define the following
functions which belong to $\mathbb{D}\left(\mathcal{D}^{*},\,\mathbb{R}\right)$,
\begin{align*}
J_{Z,h}\left(v\right)\triangleq\sum_{k=T_{b}\left(v\right)+1}^{T_{b}^{0}}z_{kh}z'_{kh} & ,\qquad\qquad J_{e,h}\left(v\right)\triangleq\sum_{k=T_{b}\left(v\right)+1}^{T_{b}^{0}}z_{kh}\widetilde{e}_{kh},
\end{align*}
for $\left(T_{b}\left(v\right)+1\right)h\leq t_{v}<\left(T_{b}\left(v\right)+2\right)h$.
Recall that the lower limit of the summation is $T_{b}\left(v\right)+1=T_{b}^{0}+\left\lfloor v/h\right\rfloor $
$\left(v\leq0\right)$ and thus the number of observations in each
sum increases at rate $1/h$. We first note that the partial sums
of cross-products between the predictable and stochastic semimartingale
regressors is null because the drift processes are of higher order
(recall Lemma \ref{Lemma Cross Product Null}). Given the previous
lemma we can decompose $\overline{Q}_{T}\left(\theta,\,v\right)$
as follows,
\begin{align}
\overline{Q}_{T}\left(\theta,\,v\right) & =\left(\delta_{p}^{0}\right)'R_{1,h}\left(v\right)\delta_{p}^{0}+\left(\delta_{Z}^{0}\right)'R_{2,h}\left(v\right)\delta_{Z}^{0}\label{Eq. Decomposition of Q, Pre}\\
 & \quad+2\left(\delta^{0}\right)'\left(\frac{1}{\sqrt{h}}\sum_{k=T_{b}+1}^{T_{b}^{0}}z_{kh}\widetilde{e}{}_{kh}\right),\nonumber 
\end{align}
 where 
\begin{align*}
R_{1,h}\left(v\right) & \triangleq\sum_{k=T_{b}\left(v\right)+1}^{T_{b}^{0}}\begin{bmatrix}h & Y_{\left(k-1\right)h}h^{3/2}\\
Y_{\left(k-1\right)h}h^{3/2} & \left(Y_{\left(k-1\right)h}h\right)^{2}
\end{bmatrix},\qquad R_{2,h}\left(v\right)\triangleq\left[Z'_{\Delta}Z{}_{\Delta}\right]_{\left\{ \cdot,p\times p\right\} },
\end{align*}
 and $\delta^{0}$ has been partitioned accordingly; that is, $\delta_{p}^{0}=\left(\mu_{\delta}^{0},\,\alpha_{\delta}^{0}\right)'$
is the vector of parameters associated with the predictable regressors
whereas $\delta_{Z}^{0}$ is the vector of parameters associated with
the stochastic martingale regressors in $Z.$ By standard results
for convergence of Riemann sums, 
\begin{align}
\left(\delta_{p}^{0}\right)'R_{1,h}\left(v\right)\delta_{p}^{0} & \overset{\textrm{u.c.p.}}{\Rightarrow}\left(\delta_{p}^{0}\right)'\begin{bmatrix}N_{b}^{0}-N_{b} & \int_{N_{b}^{0}+v}^{N_{b}^{0}}\widetilde{Y}_{s}ds\\
\int_{N_{b}^{0}+v}^{N_{b}^{0}}\widetilde{Y}_{s}ds & \int_{N_{b}^{0}+v}^{N_{b}^{0}}\widetilde{Y}_{s}^{2}ds
\end{bmatrix}\delta_{p}^{0}.\label{R1 ucp}
\end{align}
 Next, since $Z_{t}^{\left(j\right)}$ ($j=3,\ldots,\,p+2$) is a
continuous It\^o semimartingale, we have by Theorem 3.3.1 in \citeReferencesSupp{jacod/protter:12},
\begin{align}
R_{2,h}\left(v\right) & \overset{\textrm{u.c.p.}}{\Rightarrow}\left\langle Z_{\Delta},\,Z_{\Delta}\right\rangle \left(v\right).\label{R2 ucp}
\end{align}
We now turn to examine the asymptotic behavior of the second term
in \eqref{Eq. Decomposition of Q, Pre} on $\mathcal{D}^{*}$. We
use the following steps. First, we present a stable central limit
theorem for each component of $Z'_{\Delta}e.$ Second, we show the
joint convergence stably in law to a continuous Gaussian process,
and finally we verify tightness of the sequence of processes, which
in turn yields the stable convergence under the uniform metric. We
begin with the second element of $Z'_{\Delta}e,$ 
\begin{align*}
\frac{1}{\sqrt{h}}\sum_{k=T_{b}\left(v\right)+1}^{T_{b}^{0}}\alpha_{\delta}^{0}z_{kh}^{\left(2\right)}\widetilde{e}{}_{kh} & =\frac{1}{\sqrt{h}}\sum_{k=T_{b}\left(v\right)+1}^{T_{b}^{0}}\alpha_{\delta}^{0}\left(Y_{\left(k-1\right)h}h\right)\widetilde{e}{}_{kh},
\end{align*}
and using $\widetilde{Y}_{\left(k-1\right)h}=h^{1/2}Y_{\left(k-1\right)h}$
{[}recall that $\widetilde{Y}_{\left(k-1\right)h}$ is bounded by
the localization Assumption \ref{Assumption Localization, Pre}{]}
we then have
\begin{align*}
\frac{1}{\sqrt{h}}\sum_{k=T_{b}\left(v\right)+1}^{T_{b}^{0}}\alpha_{\delta}^{0}\left(Y_{\left(k-1\right)h}h\right)\widetilde{e}{}_{kh} & =\sum_{k=T_{b}\left(v\right)+1}^{T_{b}^{0}}\alpha_{\delta}^{0}\left(\widetilde{Y}_{\left(k-1\right)h}\right)\widetilde{e}{}_{kh}\\
 & \overset{\textrm{u.c.p.}}{\Rightarrow}\int_{N_{b}^{0}+v}^{N_{b}^{0}}\alpha_{\delta}^{0}\widetilde{Y}_{s}dW_{e,s},
\end{align*}
which follows from the convergence of Riemann approximations for stochastic
integrals {[}cf. Proposition 2.2.8 in \citeReferencesSupp{jacod/protter:12}{]}.
For the first component, the argument is similar:
\begin{align}
\frac{1}{\sqrt{h}}\sum_{k=T_{b}\left(v\right)+1}^{T_{b}^{0}}\mu_{\delta}^{0}z_{kh}^{\left(1\right)}\widetilde{e}_{kh} & \overset{\textrm{u.c.p.}}{\Rightarrow}\int_{N_{b}^{0}+v}^{N_{b}^{0}}\mu_{\delta}^{0}dW_{e,s}.\label{ucp z1 e}
\end{align}
 Next, we consider the $p$-dimensional lower subvector of $Z'_{\Delta}e$,
which can be written as 
\begin{align}
2\left(\delta_{Z}^{0}\right)'\left(\frac{1}{\sqrt{h}}\sum_{k=T_{b}\left(v\right)+1}^{T_{b}^{0}}\widetilde{z}_{kh}\widetilde{e}{}_{kh}\right) & ,\label{Eq. Ze, Pre}
\end{align}
 where we have partitioned $z_{kh}$ as $z_{kh}=\begin{bmatrix}h^{1/2} & Y_{\left(k-1\right)h}h & \widetilde{z}'_{kh}\end{bmatrix}'.$
Then, note that the small-dispersion asymptotic re-parametrization
implies that $\widetilde{z}_{kh}\widetilde{e}{}_{kh}$ corresponds
to $z_{kh}\widetilde{e}{}_{kh}$ from Theorem \ref{Theorem 1}. Hence,
we shall apply the same arguments as in the proof of Theorem \ref{Theorem 1}
since \eqref{Eq. Ze, Pre} is simply $2\left(\delta_{Z}^{0}\right)'$
times $\mathscr{W}_{h}\left(v\right)=h^{-1/2}J_{e,h}\left(v\right)$,
where $J_{e,h}\left(v\right)\triangleq\sum_{k=T_{b}\left(v\right)+1}^{T_{b}^{0}}\widetilde{z}_{kh}\widetilde{e}$
with $\left(T_{b}\left(v\right)+1\right)h\leq t_{v}<\left(T_{b}\left(v\right)+2\right)h$.
By Theorem 5.4.2 in \citeReferencesSupp{jacod/protter:12}, $\mathscr{W}_{h}\left(v\right)\overset{\mathcal{L}-\mathrm{s}}{\Rightarrow}\mathscr{W}_{Ze}\left(v\right)$.
Since the convergence of the drift processes $R_{1,h}\left(v\right)$
and $R_{2,h}\left(v\right)$ occur in probability locally uniformly
in time while $\mathscr{W}_{h}\left(v\right)$ converges stably in
law to a continuous limit process, we have for each $\left(\theta,\,\cdot\right)$
a stable convergence in law under the uniform metric. This is a consequence
of the property of stable convergence in law {[}cf. section VIII.5c
in \citeReferencesSupp{jacod/shiryaev:03}{]}. Since the case $v>0$
is analogous, this proves the finite-dimensional convergence of the
process $\overline{Q}_{T}\left(\theta,\,\cdot\right),$ for each $\theta.$
It remains to verify stochastic equicontinuity. As for the terms in
$R_{1,h}\left(v\right)$, we can decompose
\begin{align*}
\left(\alpha_{\delta}\right)^{2} & \left(\sum_{k=T_{b}\left(v\right)+1}^{T_{b}^{0}}\left(z_{kh}^{\left(2\right)}\right)^{2}-\left(\int_{N_{b}^{0}+v}^{N_{b}^{0}}\widetilde{Y}_{s}^{2}ds\right)\right),
\end{align*}
 as 
\[
\overline{Q}_{6,T}\left(\theta,\,v\right)+\overline{Q}_{7,T}\left(\theta,\,v\right),
\]
where 
\begin{align*}
\overline{Q}_{6,T}\left(\theta,\,v\right)\triangleq\left(\alpha_{\delta}\right)^{2}\left(\sum_{k}\zeta{}_{2,h,k}^{*}\right) & ,\quad\mathrm{and}\quad\overline{Q}_{7,T}\left(\theta,\,v\right)\triangleq\left(\alpha_{\delta}\right)^{2}\left(\sum_{k}\zeta^{**}{}_{2,h,k}\right),
\end{align*}
with 
\begin{align*}
\zeta{}_{2,h,k}^{*} & \triangleq\left(z_{kh}^{\left(2\right)}\right)^{2}-\left(\int_{\left(k-1\right)h}^{kh}\widetilde{Y}_{s}^{2}ds\right)-2\widetilde{Y}_{\left(k-1\right)h}\int_{\left(k-1\right)h}^{kh}\left(\widetilde{Y}_{\left(k-1\right)h}-\widetilde{Y}_{s}\right)ds\\
 & \quad+2\mathbb{E}\left[\widetilde{Y}_{\left(k-1\right)h}\left(\widetilde{Y}_{\left(k-1\right)h}\cdot h-\int_{\left(k-1\right)h}^{kh}\widetilde{Y}_{s}ds\right)|\,\mathscr{F}_{\left(k-1\right)h}\right]\\
 & \triangleq L_{1,h,k}+L_{2,h,k},
\end{align*}
 and 
\begin{align*}
\zeta_{2,h,k}^{**} & =2\widetilde{Y}_{\left(k-1\right)h}(\widetilde{Y}_{\left(k-1\right)h}h-\int_{\left(k-1\right)h}^{kh}\widetilde{Y}_{s}ds\\
 & \quad-\mathbb{E}\left[\left(\widetilde{Y}_{\left(k-1\right)h}h-\int_{\left(k-1\right)h}^{kh}\widetilde{Y}_{s}ds\right)|\,\mathscr{F}_{\left(k-1\right)h}\right]).
\end{align*}
Then, we have the following decomposition for 
\[
\overline{Q}_{T}^{\textrm{c}}\left(\theta^{*}\right)\triangleq\overline{Q}_{T}^{\textrm{}}\left(\theta^{*}\right)+\left(\delta^{0}\right)'\Lambda\left(v\right)\delta^{0},
\]
(if $v\leq0$ and defined analogously for $v>0$), 
\[
\overline{Q}_{T}^{\textrm{c}}\left(\theta^{*}\right)=\sum_{r=1}^{9}\overline{Q}_{r,T}\left(\theta,\,v\right),
\]
where $\overline{Q}_{r,T}\left(\theta,\,v\right),$ $r=1,\ldots,\,4$,
are defined in \eqref{Eq. Qc} and $\overline{Q}_{5,T}\left(\theta,\,v\right)\triangleq\left(\mu_{\delta}\right)^{2}(\sum_{k}\zeta{}_{1,h,k})$,
$\overline{Q}_{8,T}\left(\theta,\,v\right)\triangleq\left(\mu_{\delta}\right)^{2}(h^{-1/2}\sum_{k}\xi_{1,h,k})$,
$\overline{Q}_{9,T}\left(\theta,\,v\right)\triangleq(\alpha_{\delta})^{2}(h^{-1/2}\sum_{k}\xi_{2,h,k})$
where $\zeta_{1,h,k}\triangleq(z_{kh}^{\left(1\right)})^{2}-h,$ $\xi_{1,h,k}\triangleq h^{1/2}\widetilde{e}_{kh}$
and $\xi_{2,h,k}\triangleq(\widetilde{Y}_{\left(k-1\right)h}h^{1/2})\widetilde{e}_{kh}$.
Moreover, recall that $\sum_{k}$ stands for $\sum_{T_{b\left(v\right)+1}}^{T_{b}^{0}}$
for $N_{b}\left(v\right)\in\mathcal{D}^{*}\left(C\right)$. Let us
consider $\overline{Q}_{6,T}\left(\theta,\,v\right)$ first. For $s\in\left[\left(k-1\right)h,\,kh\right]$,
by the Burkh\"{o}lder-Davis-Gundy inequality
\begin{align*}
\left|\mathbb{E}\left[\widetilde{Y}_{\left(k-1\right)h}\left(\widetilde{Y}_{\left(k-1\right)h}-\widetilde{Y}_{s}\right)|\,\mathscr{F}_{\left(k-1\right)h}\right]\right| & \leq Kh,
\end{align*}
from which we can deduce that, using a maximal inequality for any
$r>1,$ 
\begin{align}
\left[\mathbb{E}\left(\sup_{\left(\theta,\,v\right)}\left|\left(\alpha_{\delta}\right)^{2}\sum_{k}L_{2,h,k}\right|\right)^{r}\right]^{1/r} & \leq K_{r}\left(\sup_{\left(\theta,\,v\right)}\left(\alpha_{\delta}\right)^{2r}\sum_{k}h^{r}\right)^{1/r}=K_{r}h^{\frac{r-1}{r}}.\label{A.24}
\end{align}
 By a Taylor series expansion for the mapping $f:\,y\rightarrow y^{2}$,
and $s\in\left[\left(k-1\right)h,\,kh\right]$, 
\begin{align*}
\mathbb{E}\left|\widetilde{Y}_{\left(k-1\right)h}^{2}-\widetilde{Y}_{s}^{2}-2\widetilde{Y}_{\left(k-1\right)h}\left(\widetilde{Y}_{\left(k-1\right)h}-\widetilde{Y}_{s}\right)\right| & \leq K\mathbb{E}\left[\left(\widetilde{Y}_{\left(k-1\right)h}-\widetilde{Y}_{s}\right)^{2}\right]\leq Kh,
\end{align*}
 where the second inequality follows from the Burkh\"{o}lder-Davis-Gundy
inequality. Thus, using a maximal inequality as in \eqref{A.24},
we have for $r>1$
\begin{align}
\left[\mathbb{E}\left(\sup_{\left(\theta,\,v\right)}\left|\left(\alpha_{\delta}\right)^{2}\sum_{k}L_{1,h,k}\right|\right)^{r}\right]^{1/r} & =K_{r}h^{\frac{r-1}{r}}.\label{Eq A23}
\end{align}
 \eqref{A.24} and \eqref{Eq A23} imply that $\overline{Q}_{6,T}\left(\cdot,\,\cdot\right)$
is stochastically equicontinuous. Next, note that $\overline{Q}_{7,T}\left(\theta,\,v\right)$
is a sum of martingale differences times $h^{1/2}$ (recall the definition
of $\Delta_{h}\widetilde{V}_{k}=h^{1/2}\Delta_{h}V_{k}(\nu,\,\delta_{Z,1},\,\delta_{Z,2})$).
Therefore by Assumption \ref{Assumption Localization, Pre}, for any
$0\leq s<t\leq N$, $V_{t}-V_{s}=O_{p}\left(1\right)$ uniformly and
therefore, 
\begin{align}
\sup_{\left(\theta,\,v\right)}\left|\overline{Q}_{7,T}\left(\theta,\,v\right)\right| & \leq KO_{p}\left(h^{1/2}\right).\label{Eq. A25}
\end{align}
Given \eqref{R1 ucp} and \eqref{A.24}-\eqref{Eq. A25}, we deduce
that
\[
\sup_{\left(\theta,\,v\right)}\left\{ \left|\overline{Q}_{6,T}\left(\theta,\,v\right)\right|+\left|\overline{Q}_{7,T}\left(\theta,\,v\right)\right|\right\} =o_{p}\left(1\right).
\]
As for the term involving $R_{1,h}\left(v\right)$, it is easy to
see that $\sup_{\left(\theta,\,v\right)}|\overline{Q}\left(\theta,\,v\right)|\overset{}{\rightarrow}0.$
Next, we can use some of the results in the proof of Theorem \ref{Theorem 1}.
In particular, the asymptotic stochastic equicontinuity of the sequence
of processes $\{2\left(\delta_{Z}\right)'\mathscr{W}_{h}\left(v\right)\}$
follows from the same property as those applied to $\{\overline{Q}_{3,T}\left(\theta,\,v\right)\}$
and $\{\overline{Q}_{4,T}\left(\theta,\,v\right)\}$. The stochastic
equicontinuity of 
\[
\left(\delta_{Z}\right)'\left(R_{2,h}\left(\theta,\,v\right)-\left\langle Z_{\Delta},\,Z_{\Delta}\right\rangle \left(v\right)\right)\delta_{Z},
\]
also follows from the same proof. Recall $\overline{Q}_{1,T}\left(\theta,\,v\right)+\overline{Q}_{2,T}\left(\theta,\,v\right)$
as defined in \eqref{Eq. Qc}. Thus, stochastic equicontinuity follows
from \eqref{Eq. Q2 Stoc Equi} and the equation right before that.
Next, let us consider $\overline{Q}_{9,T}\left(\theta,\,v\right).$
We use the alternative definition (ii) of stochastic equicontinuity
in \citeReferencesSupp{andrews:37hoe}. Consider any sequence $\{\left(\theta,\,v\right)\}$
and $\{(\bar{\theta},\,\bar{v})\}$ (we omit the dependence on $h$
for simplicity). Assume $N_{b}\leq N_{b}^{0}\leq\bar{N}_{b}$ (the
other cases can be proven similarly) and let $Nd_{h}\triangleq\bar{N}_{b}-N_{b}$.
Then,
\begin{align}
\left|\overline{Q}_{9,T}\left(\theta,\,v\right)-\overline{Q}_{9,T}\left(\bar{\theta},\,\bar{v}\right)\right| & =\left|\alpha_{\delta}\sum_{k=T_{b}\left(v\right)+1}^{T_{b}^{0}}\widetilde{Y}_{\left(k-1\right)h}\widetilde{e}_{kh}-\bar{\alpha}_{\delta}\sum_{k=T_{b}^{0}}^{T_{b}\left(\bar{v}\right)}\widetilde{Y}_{\left(k-1\right)h}\widetilde{e}_{kh}\right|\nonumber \\
 & \leq\left|\alpha_{\delta}\right|\left|\sum_{k=T_{b}\left(v\right)+1}^{T_{b}^{0}}\widetilde{Y}_{\left(k-1\right)h}\widetilde{e}{}_{kh}\right|\label{Eq. AAA6}\\
 & \quad+\left|\bar{\alpha}_{\delta}\right|\left|\sum_{k=T_{b}^{0}}^{T_{b}\left(\bar{v}\right)}\widetilde{Y}_{\left(k-1\right)h}\widetilde{e}{}_{kh}\right|.\nonumber 
\end{align}
The second term is such that, by the Burkh\"{o}lder-Davis-Gundy inequality
for any $r\geq1,$ 
\begin{align*}
\mathbb{E} & \left[\sup_{0\leq u\leq d_{h}}\left|\sum_{k=T_{b}^{0}}^{T_{b}^{0}+\left\lfloor Nu/h\right\rfloor }\widetilde{Y}_{\left(k-1\right)h}\widetilde{e}{}_{kh}\right|^{r}|\,\mathscr{F}_{N_{b}^{0}}\right]\\
 & \leq K_{r}\left(Nd_{h}\right)^{r/2}\mathbb{E}\left[\frac{1}{Nd_{h}}\left(\sum_{k=T_{b}^{0}}^{T_{b}^{0}+\left\lfloor Nd_{h}/h\right\rfloor }\int_{\left(k-1\right)h}^{kh}\left(\widetilde{Y}_{s}\right)^{2}ds\right)^{r/2}|\,\mathscr{F}_{N_{b}^{0}}\right]\leq K_{r}d_{h}^{r/2}.
\end{align*}
By the law of iterated expectations, and using the property that $d_{h}\downarrow0$
in probability, we can find a $T$ large enough such that for any
$B>0$ 
\begin{align*}
\left(\mathbb{E}\left[\sup_{0\leq u\leq d_{h}}\left|\sum_{k=T_{b}^{0}}^{T_{b}^{0}+\left\lfloor Nu/h\right\rfloor }\widetilde{Y}_{\left(k-1\right)h}\widetilde{e}{}_{kh}\right|^{r}|\,\mathscr{F}_{N_{b}^{0}}\right]\right)^{1/r} & \leq K_{r}d_{h}^{1/2}P\left(Nd_{h}>B\right)\rightarrow0.
\end{align*}
 The argument for the first term in \eqref{Eq. AAA6} is analogous.
By Markov's inequality and combining the above steps, we have that
for any $\varepsilon>0$ and $\eta>0$ there exists some $\overline{T}$
such that for $T>\overline{T}$, 
\begin{align*}
P\left(\left|\overline{Q}_{9,T}\left(\theta,\,v\right)-\overline{Q}_{9,T}\left(\bar{\theta},\,\bar{v}\right)\right|>\eta\right) & <\varepsilon.
\end{align*}
 Thus, the sequence $\{\overline{Q}_{9,T}\left(\cdot,\,\cdot\right)\}$
is stochastically equicontinuous. Noting that the same proof can be
repeated for $\overline{Q}_{8,T}\left(\cdot,\,\cdot\right)$, we conclude
that the sequence of processes $\{\overline{Q}_{T}^{\textrm{c}}\left(\theta^{*}\right),\,T\geq1\}$
in \eqref{Eq. Decomposition of Q, Pre} is stochastically equicontinuous.
Furthermore, by \eqref{R1 ucp} and \eqref{R2 ucp} we obtain, 
\begin{align*}
\left(\delta_{p}^{0}\right)'R_{1,h}\left(\theta,\,v\right)\delta_{p}^{0}+\left(\delta_{Z}^{0}\right)'\left(R_{2,h}\left(\left(\theta,\,v\right)\right)\right)\delta_{Z}^{0} & \overset{\textrm{u.c.p.}}{\Rightarrow}\left(\delta^{0}\right)'\Lambda\left(v\right)\delta_{\textrm{}}^{0}.
\end{align*}
 This suffices to guarantee the $\mathscr{G}$-stable convergence
in law of the process $\{\overline{Q}_{T}\left(\cdot,\,\cdot\right),\,T\geq1\}$
towards a process $\mathscr{W}\left(\cdot\right)$ with drift $\Lambda\left(\cdot\right)$
which, conditional on $\mathscr{G}$, is a two-sided Gaussian martingale
process with covariance matrix given in \eqref{Eq. Covariance Matrix of Limit process, Pre}.
By definition, $\mathcal{D}^{*}\left(C\right)$ is compact and $Th(\widehat{\lambda}_{b,\pi}-\lambda_{0})=O_{p}\left(1\right),$
which together with the fact that the limit process is a continuous
Gaussian process enable one to deduce the main assertion from the
continuous mapping theorem for the argmax functional. $\square$

\subsubsection{Proof of Proposition \ref{Proposition Consistency QLS, Pre}}

\begin{singlespace}
We begin with a few lemmas. Let $\widetilde{Y}_{t}^{*}\triangleq\widetilde{Y}_{\left\lfloor t/h\right\rfloor h}$.
The first result states that the observed process $\{\widetilde{Y}_{t}^{*}\}$
converges to the non-stochastic process $\{\widetilde{Y}_{t}^{0}\}$
defined in \eqref{Eq. Dynamical System} as $h\downarrow0$. Assumption
\ref{Assumption Localization, Pre} is maintained throughout and the
constant $K>0$ may vary from line to line.
\end{singlespace}
\begin{lem}
As $h\downarrow0$, $\sup_{0\leq t\leq N}|\widetilde{Y}_{t}^{*}-\widetilde{Y}_{t}^{0}|=o_{p}\left(1\right)$.\label{Lemma QLS 1} 
\end{lem}
\noindent\textit{Proof.} Let us introduce a parameter $\gamma_{h}$
with the property $\gamma_{h}\downarrow0$ and $h^{1/2}/\gamma_{h}\rightarrow B$
where $B<\infty.$ By construction, for $t<N_{b}^{0}$,
\begin{align*}
\widetilde{Y}_{t}-\widetilde{Y}_{t}^{0} & =\int_{0}^{t}\alpha_{1}^{0}\left(\widetilde{Y}_{s}-\widetilde{Y}_{s}^{0}\right)ds+B\gamma_{h}\left(\nu^{0}\right)'D{}_{t}\\
 & \quad+B\gamma_{h}\left(\delta_{Z,1}^{0}\right)'\int_{0}^{t}dZ_{s}+B\gamma_{h}\int_{0}^{t}\sigma_{e,s}dW_{e,s}.
\end{align*}
 We can use Cauchy-Schwarz's inequality, so that
\begin{align*}
\left|\widetilde{Y}_{t}-\widetilde{Y}_{t}^{0}\right|^{2} & \leq2K\Bigl[\left|\int_{0}^{t}\alpha_{1}\left(\widetilde{Y}_{s}-\widetilde{Y}_{s}^{0}\right)ds\right|^{2}\\
 & +\Bigl(\left|\nu^{0\prime}D{}_{t}\right|^{2}+\left|\delta_{Z,1}^{0\prime}\int_{0}^{t}dZ_{s}\right|^{2}+\left|\int_{0}^{t}\sigma_{e,s}dW_{e,s}\right|^{2}\Bigl)\left(B\gamma_{h}\right)^{2}\\
 & \leq2Kt\Bigl[\left|\alpha_{1}^{0}\right|^{2}\int_{0}^{t}\left|\widetilde{Y}_{s}-\widetilde{Y}_{s}^{0}\right|^{2}ds+\Bigl(\sup_{0\leq s\leq t}\left|\nu^{0\prime}D{}_{s}\right|^{2}\\
 & \quad+\sup_{0\leq s\leq t}\left|\delta_{Z,1}^{0\prime}\int_{0}^{t}dZ_{s}\right|^{2}+\sup_{0\leq s\leq t}\left|\int_{0}^{s}\sigma_{e,u}dW_{e,u}\right|^{2}\Bigl)\left(B\gamma_{h}\right)^{2}\Bigl].
\end{align*}
 By Gronwall\textquoteright s inequality,
\begin{align*}
\left|\widetilde{Y}_{t}-\widetilde{Y}_{t}^{0}\right|^{2} & \leq2\left(B\gamma_{h}\right)^{2}C\exp\left(\int_{0}^{t}2K^{2}tds\right)\\
 & \leq2\left(B\gamma_{h}\right)^{2}C\exp\left(2K^{2}t^{2}\right),
\end{align*}
 where $C<\infty$ is a bound on the sum of the supremum terms in
the last equation above. The bound follows from Assumption \ref{Assumption Localization, Pre}.
Then, 
\[
\sup_{0\leq t\leq N}\left|\widetilde{Y}_{t}-\widetilde{Y}_{t}^{0}\right|\leq K\sqrt{2}B\gamma_{h}\exp\left(K^{2}N^{2}\right)\rightarrow0,
\]
 as $h\downarrow0$ (and so $\gamma_{h}\downarrow0$). The assertion
then follows from $\left\lfloor t/h\right\rfloor h\rightarrow t$
as $h\downarrow0$. The case with $t\geq N_{b}^{0}$ is proved in
a similar fashion. $\square$
\begin{lem}
\label{Lemma QLS Riemann}As $h\downarrow0$, uniformly in $\left(\mu_{1},\,\alpha_{1}\right)$,
$\left(N/T\right)\sum_{k=1}^{T_{b}^{0}}(\mu_{1}+\alpha_{1}\widetilde{Y}_{\left(k-1\right)h})\overset{P}{\rightarrow}\int_{0}^{N_{b}^{0}}(\mu_{1}+\alpha_{1}\widetilde{Y}_{s}^{0})ds.$
\end{lem}
\noindent\textit{Proof.} Note that
\begin{align*}
\sup_{\mu_{1},\alpha_{1}} & \left|\frac{N}{T}\sum_{k=1}^{T_{b}^{0}}\left(\mu_{1}+\alpha_{1}\widetilde{Y}_{\left(k-1\right)h}\right)-\int_{0}^{N\lambda_{0}}\left(\mu_{1}+\alpha_{1}\widetilde{Y}_{s}^{0}\right)\right|\\
 & =\sup_{\mu_{1},\alpha_{1}}\left|\int_{0}^{N_{b}^{0}}\left(\mu_{1}+\alpha_{1}\widetilde{Y}_{s}^{*}\right)ds-\int_{0}^{N_{b}^{0}}\left(\mu_{1}+\alpha_{1}\widetilde{Y}_{s}^{0}\right)ds\right|\\
 & \leq\sup_{\alpha_{1}}\int_{0}^{N_{b}^{0}}\left|\alpha_{1}\right|\left|\widetilde{Y}_{s}^{*}-\widetilde{Y}_{s}^{0}\right|ds\leq KO_{p}\left(\gamma_{h}\right)\sup_{\alpha_{1}}\left|\alpha_{1}\right|,
\end{align*}
which goes to zero as $h\downarrow0$ by Lemma \ref{Lemma QLS 1}
(recall $h^{1/2}/\gamma_{h}\rightarrow B$) and by Assumption \ref{Assumption Localization, Pre}.
$\square$
\begin{lem}
\label{Lemma QLS 3-1}For each $3\leq j\leq p+2$ and each $\theta$,
as $h\downarrow0$, 
\begin{align*}
\sum_{k=1}^{\left\lfloor N_{b}^{0}/h\right\rfloor }\left(\mu_{1}+\alpha_{1}\widetilde{Y}_{\left(k-1\right)h}\right)\delta_{Z,1}^{\left(j\right)}\Delta_{h}Z_{k}^{\left(j\right)} & \overset{P}{\rightarrow}\int_{0}^{N\lambda_{0}}\left(\mu_{1}+\alpha_{1}\widetilde{Y}_{\left(k-1\right)h}^{0}\right)dZ_{s}^{\left(j\right)}.
\end{align*}
\end{lem}
\noindent\textit{Proof.} Note that 
\begin{align*}
\sum_{k=1}^{\left\lfloor N_{b}^{0}/h\right\rfloor }\left(\mu_{1}+\alpha_{1}\widetilde{Y}_{\left(k-1\right)h}\right)\delta_{Z,1}^{\left(j\right)}\Delta_{h}Z_{k}^{\left(j\right)} & =\int_{0}^{N_{b}^{0}}\left(\mu_{1}+\alpha_{1}\widetilde{Y}_{s}^{*}\right)dZ_{s}^{\left(j\right)}.
\end{align*}
 By Markov's inequality and the dominated convergence theorem, for
every $\varepsilon>0$ and every $\eta>0$: 
\begin{align*}
P\biggl( & \left|\int_{0}^{N_{b}^{0}}\alpha_{1}\left(\widetilde{Y}_{s}^{*}-\widetilde{Y}_{s}^{0}\right)\delta_{Z,1}^{\left(j\right)}dZ_{s}^{\left(j\right)}\right|>\eta\biggr)\\
 & \leq\frac{\left(\sup_{0\leq s\leq N}\sum_{r=1}^{p}\left(\sigma_{Z,s}^{\left(j,r\right)}\right)^{2}\right)^{1/2}}{\eta}\left|\alpha_{1}\right|\left|\delta_{Z,1}^{\left(j\right)}\right|\left(\int_{0}^{N_{b}^{0}}\mathbb{E}\left[\left(\widetilde{Y}_{s}^{*}-\widetilde{Y}_{s}^{0}\right)^{2}\right]ds\right)^{1/2},
\end{align*}
which goes to zero as $h\downarrow0$ in view of Lemma \ref{Lemma QLS 1}
and Assumption \ref{Assumption Localization, Pre}. $\square$
\begin{lem}
\label{Uniform miu * (dX- miu0), Pre, cf 3.5}As $h\downarrow0$,
uniformly in $\mu_{1},\,\alpha_{1}$, 
\begin{align*}
\sum_{k=1}^{T_{b}^{0}}\left(\mu_{1}+\alpha_{1}\widetilde{Y}_{\left(k-1\right)h}\right)\left(\widetilde{Y}_{kh}-\widetilde{Y}_{\left(k-1\right)h}-\left(\mu_{1}^{0}+\alpha_{1}^{0}\widetilde{Y}_{\left(k-1\right)h}\right)h\right) & \overset{P}{\rightarrow}0.
\end{align*}
\end{lem}
\noindent\textit{Proof.} By definition {[}recall the notation in
\eqref{Small-Dispersion Model}{]}, 
\begin{align*}
\widetilde{Y}_{kh}-\widetilde{Y}_{\left(k-1\right)h} & =\int_{\left(k-1\right)h}^{kh}\left(\mu_{1}^{0}+\alpha_{1}^{0}\widetilde{Y}_{s}\right)ds+\Delta_{h}\widetilde{V}_{k}\left(\nu^{0},\,\delta_{Z,1}^{0},\,\delta_{Z,2}^{0}\right).
\end{align*}
 Then,
\begin{align*}
\sum_{k=1}^{T_{b}^{0}} & \left(\mu_{1}+\alpha_{1}\widetilde{Y}_{\left(k-1\right)h}\right)\left(\widetilde{Y}_{kh}-\widetilde{Y}_{\left(k-1\right)h}-\left(\mu_{1}^{0}+\alpha_{1}^{0}\widetilde{Y}_{\left(k-1\right)h}\right)h\right)\\
 & =\sum_{k=1}^{T_{b}^{0}}\int_{\left(k-1\right)h}^{kh}\left(\mu_{1}+\alpha_{1}\widetilde{Y}_{\left(k-1\right)h}\right)\left(\mu_{1}^{0}+\alpha_{1}^{0}\widetilde{Y}_{s}-\left(\mu_{1}^{0}+\alpha_{1}^{0}\widetilde{Y}_{\left(k-1\right)h}\right)\right)\\
 & \quad+\sum_{k=1}^{T_{b}^{0}}\int_{\left(k-1\right)h}^{kh}\left(\mu_{1}+\alpha_{1}\widetilde{Y}_{\left(k-1\right)h}\right)\Delta_{h}\widetilde{V}_{k}\left(\nu^{0},\,\delta_{Z,1}^{0},\,\delta_{Z,2}^{0}\right)\\
 & =\int_{0}^{N_{b}^{0}}\left(\mu_{1}+\alpha_{1}\widetilde{Y}_{\left(k-1\right)h}^{*}\right)\left(\alpha_{1}^{0}\left(\widetilde{Y}_{s}-\widetilde{Y}_{\left(k-1\right)h}^{*}\right)\right)ds\\
 & \quad+B\gamma_{h}\int_{0}^{N_{b}^{0}}\left(\mu_{1}+\alpha_{1}\widetilde{Y}_{s}^{*}\right)dV_{s}.
\end{align*}
For the first term on the right-hand side of the last equation 
\begin{align*}
\sup_{\mu_{1},\alpha_{1}} & \left|\int_{0}^{N_{b}^{0}}\left(\mu_{1}+\alpha_{1}\widetilde{Y}_{s}^{*}\right)\left(\alpha_{1}^{0}\left(\widetilde{Y}_{s}-\widetilde{Y}_{s}^{*}\right)\right)ds\right|\\
 & \leq\left|\alpha_{1}^{0}\right|\left|\int_{0}^{N_{b}^{0}}\sup_{\mu_{1},\alpha_{1}}\left(\mu_{1}+\alpha_{1}\widetilde{Y}_{s}^{*}\right)\left(\widetilde{Y}_{s}-\widetilde{Y}_{s}^{0}+\widetilde{Y}_{s}^{0}-\widetilde{Y}_{s}^{*}\right)ds\right|\\
 & \leq\left|\alpha_{1}^{0}\right|K\left(\int_{0}^{N_{b}^{0}}\sup_{0\leq s\leq N_{b}^{0}}\left|\widetilde{Y}_{s}-\widetilde{Y}_{s}^{0}\right|+\sup_{0\leq s\leq N_{b}^{0}}\left|\widetilde{Y}_{s}^{0}-\widetilde{Y}_{s}^{*}\right|ds\right),
\end{align*}
 which is $o_{p}\left(1\right)$ as $h\downarrow0$ from Lemma \ref{Lemma QLS 1}
and Assumption \ref{Assumption Localization, Pre}. Next, consider
the vector of regressors $Z$, and note that for any $3\leq j\leq p+2$,
\begin{align*}
B\gamma_{h} & \sup_{\mu_{1},\alpha_{1}}\left|\int_{0}^{N_{b}^{0}}\left(\mu_{1}+\alpha_{1}\widetilde{Y}_{s}^{*}\right)dZ_{s}^{\left(j\right)}\right|\\
 & \leq B\gamma_{h}\sup_{\mu_{1},\alpha_{1}}\left|\int_{0}^{N_{b}^{0}}\left(\mu_{1}+\alpha_{1}\widetilde{Y}_{s}^{*}\right)\sum_{r=1}^{p}\sigma_{Z,s}^{\left(j,r\right)}dW_{Z}^{\left(r\right)}\right|.
\end{align*}
Let
\[
R_{j,h}=R_{j,h}\left(\mu_{1},\,\alpha_{1}\right)\triangleq\int_{0}^{N_{b}^{0}}B\gamma_{h}\left(\mu_{1}+\alpha_{1}\widetilde{Y}_{s}^{*}\right)\sum_{r=1}^{p}\sigma_{Z,s}^{\left(j,r\right)}dW_{Z}^{\left(r\right)},
\]
 (we index $R_{j}$ by $h$ because $\widetilde{Y}_{s}^{*}$ depends
on $h$). Then, we want to show that, for every $\varepsilon>0$ and
$K>0,$ 
\begin{align}
P\left(\sup_{\mu_{1},\alpha_{1}}\left|R_{j,h}\left(\mu_{1},\,\alpha_{1}\right)\right|>K\right) & \leq\varepsilon.\label{QLS Uniform miuxWz}
\end{align}
In view of Chebyshev's inequality and the It\^o's isometry,
\begin{align*}
P\left(\left|R_{j,h}\right|>K\right) & \leq\left(\frac{B\gamma_{h}}{K}\right)^{2}\mathbb{E}\left[\left|\int_{0}^{N_{b}^{0}}\left(R_{j,h}/\left(B\gamma_{h}\right)\right)\right|^{2}\right],\\
 & \leq\left[\sup_{0\leq s\leq N}\sum_{r=1}^{p}\left(\sigma_{Z,s}^{\left(j,r\right)}\right)^{2}\right]\left(\frac{B\gamma_{h}}{K}\right)^{2}\int_{0}^{N_{b}^{0}}\mathbb{E}\left[\left|\mu_{1}+\alpha_{1}\widetilde{Y}_{s}^{*}\right|^{2}ds\right],
\end{align*}
 so that by the boundness of the processes (cf. Assumption \ref{Assumption Localization, Pre})
and the compactness of $\Theta_{0},$ we have for some $A<\infty$,
\begin{align}
P\left(\left|R_{j,h}\right|>K\right) & \leq A\left[\sup_{0\leq s\leq T}\sum_{r=1}^{p}\left(\sigma_{Z,s}^{\left(j,r\right)}\right)^{2}\right]\left(\frac{B\gamma_{h}}{K}\right)^{2}\rightarrow0,\label{QLS miuxZ Pointwise Convergence}
\end{align}
 since $\gamma_{h}\downarrow0.$ This demonstrates pointwise convergence.
It remains to show the stochastic equicontinuity of the sequence of
processes $\{R_{j,h}\left(\cdot\right)\}.$ Choose $2m>p$ and note
that standard estimates for continuous It\^o semimartingales result
in $\mathbb{E}[|R_{j,h}|^{2m}]\leq K$ which follows using the same
steps that led to \eqref{QLS miuxZ Pointwise Convergence} with the
Burkh\"{o}lder-Davis-Gundy inequality  in place of the It\^o's isometry.
Let $g(\widetilde{Y}_{s}^{*},\,\widetilde{\theta})\triangleq\mu_{1,1}+\alpha_{1,1}\widetilde{Y}_{s}^{*}$,
$\widetilde{\theta}_{1}\triangleq(\mu_{1,1},\alpha_{1,1})'$ and $\widetilde{\theta}_{1}\triangleq(\mu_{2,1},\alpha_{2,1})'$.
For any $\widetilde{\theta}_{1},\,\widetilde{\theta}_{2},$ first
use the Burkh\"{o}lder-Davis-Gundy inequality to yield, 

\begin{align*}
\mathbb{E} & \left[\left|R_{j,h}\left(\widetilde{\theta}_{2}\right)-R_{j,h}\left(\widetilde{\theta}_{1}\right)\right|^{2m}\right]\\
 & \leq\left(B\gamma_{h}\right)^{2m}K_{m}\left[\sup_{0\leq s\leq N}\sum_{r=1}^{p}\left(\sigma_{Z,s}^{\left(j,r\right)}\right)^{2}\right]^{m}\\
 & \quad\times\mathbb{E}\left[\left(\int_{0}^{N_{b}^{0}}\left(g\left(\widetilde{Y}_{s}^{*},\,\widetilde{\theta}_{2}\right)-g\left(\widetilde{Y}_{s}^{*},\,\widetilde{\theta}_{1}\right)\right)^{2}ds\right)^{m}\right]\\
 & \leq\left(B\gamma_{h}\right)^{2m}K_{m}\left[\sup_{0\leq s\leq N}\sum_{r=1}^{p}\left(\sigma_{Z,s}^{\left(j,r\right)}\right)^{2}\right]^{m}\\
 & \quad\times\mathbb{E}\left[\left(\int_{0}^{N_{b}^{0}}\left(\left(\mu_{1,2}-\mu_{1,1}\right)+\left(\alpha_{1,2}-\alpha_{1,1}\right)\widetilde{Y}_{s}^{*}\right)^{2}ds\right)^{m}\right]
\end{align*}
\begin{align}
 & \leq\left(B\gamma_{h}\right)^{2m}K_{m}\left[\sup_{0\leq s\leq N}\sum_{r=1}^{p}\left(\sigma_{Z,s}^{\left(j,r\right)}\right)^{2}\right]^{m}\nonumber \\
 & \quad\times\mathbb{E}\left[\left(\int_{0}^{N_{b}^{0}}\left(\left(\mu_{1,2}-\mu_{1,1}\right)+\left(\alpha_{1,2}-\alpha_{1,1}\right)C\right)^{2}ds\right)^{m}\right]\nonumber \\
 & \leq\left(B\gamma_{h}\right)^{2m}K_{m}\mathbb{E}\left[\left(\int_{0}^{N_{b}^{0}}\left(2\left(\mu_{1,2}-\mu_{1,1}\right)^{2}+2C\left(\alpha_{1,2}-\alpha_{1,1}\right)^{2}\right)ds\right)^{m}\right]\nonumber \\
 & \leq2^{m}\left(B\gamma_{h}\right)^{2m}K_{m}\left\Vert 2\left(\widetilde{\theta}_{2}-\widetilde{\theta}_{1}\right)\right\Vert ^{2m}\left(\int_{0}^{N_{b}^{0}}ds\right)^{m}\label{QLS miuxZ Stochastic equic}\\
 & \quad+2^{m}\left(B\gamma_{h}\right)^{2m}K\left(\widetilde{\theta}_{1},\,\widetilde{\theta}_{2},\,m,\,C\right),\nonumber 
\end{align}
where $C=\sup_{s\geq0}|\widetilde{Y}_{s}^{*}|,$ $K(\widetilde{\theta}_{1},\,\widetilde{\theta}_{2},\,m,\,C)$
is some constant that depends on its arguments and we have used the
fact that $\left(a+b\right)^{2}\leq2a^{2}+2b^{2}$. Thus, since $\gamma_{h}\downarrow0$,
the mapping $R_{j,h}\left(\cdot\right)$ satisfies a Lipschitz-type
condition {[}cf. Section 2 in \citeReferencesSupp{andrews:92}{]}.
This is sufficient for the asymptotic stochastic equicontiuity of
$\{R_{j,h}\left(\cdot\right)\}$. Therefore, using Theorem 20 in Appendix
I of \citeReferencesSupp{ibragimov/has:81}, \eqref{QLS miuxZ Pointwise Convergence}
and \eqref{QLS miuxZ Stochastic equic} yield \eqref{QLS Uniform miuxWz}.
Since the same result can be shown to remain valid for each term in
the stochastic element $\Delta_{h}V_{k}(\nu,\,\delta_{Z,1},\,\delta_{Z,2}),$
this establishes the claim. $\square$

\noindent\textit{Proof of Proposition \ref{Proposition Consistency QLS, Pre}.}
To avoid clutter, we prove the case for which the true parameters
are $\left(\mu_{1}^{0},\,\alpha_{1}^{0}\right)'$. The extension to
parameters being local-to-zero is straightforward. The least-squares
estimates of $\left(\mu_{1}^{0},\,\alpha_{1}^{0}\right)'$ are given
by,
\begin{align}
\widehat{\mu}_{1}\widehat{N}_{b} & =\widetilde{Y}_{\widehat{N}_{b}}-\widetilde{Y}_{0}-\widehat{\alpha}_{1}h\sum_{k=1}^{\widehat{T}_{b}}\widetilde{Y}_{\left(k-1\right)h}\label{Eq Expression miu1}\\
\widehat{\alpha}_{1} & =\frac{\sum_{k=1}^{\widehat{T}_{b}}\left(\widetilde{Y}_{kh}-\widetilde{Y}_{\left(k-1\right)h}\right)\widetilde{Y}_{\left(k-1\right)h}}{h\sum_{k=1}^{\widehat{T}_{b}}\widetilde{Y}_{\left(k-1\right)h}^{2}-\widehat{N}_{b}^{-1}\left(h\sum_{k=1}^{\widehat{T}_{b}}\widetilde{Y}_{\left(k-1\right)h}\right)^{2}}\label{Eq. Expression alpha1}\\
 & \quad-\frac{\widehat{N}_{b}^{-1}\left(\widetilde{Y}_{\widehat{N}_{b}}-\widetilde{Y}_{0}\right)h\sum_{k=1}^{\widehat{T}_{b}}\widetilde{Y}_{\left(k-1\right)h}}{h\sum_{k=1}^{\widehat{T}_{b}}\widetilde{Y}_{\left(k-1\right)h}^{2}-\widehat{N}_{b}^{-1}\left(h\sum_{k=1}^{\widehat{T}_{b}}\widetilde{Y}_{\left(k-1\right)h}\right)^{2}}.\nonumber 
\end{align}
 Then, assuming $\widehat{T}_{b}<T_{b}^{0}$,
\begin{align*}
\widehat{\alpha}_{1} & =\frac{\sum_{k=1}^{\widehat{T}_{b}}\left(\mu_{1}^{0}h+\alpha_{1}^{0}\widetilde{Y}_{\left(k-1\right)h}h+\Delta_{h}\widetilde{V}_{h,k}\right)\widetilde{Y}_{\left(k-1\right)h}}{h\sum_{k=1}^{\widehat{T}_{b}}\widetilde{Y}_{\left(k-1\right)h}^{2}-\widehat{N}_{b}^{-1}\left(h\sum_{k=1}^{\widehat{T}_{b}}\widetilde{Y}_{\left(k-1\right)h}\right)^{2}}\\
 & \quad-\frac{\left(\mu_{1}^{0}+\alpha_{1}^{0}\widehat{N}_{b}^{-1}\sum_{k=1}^{\widehat{T}_{b}}\widetilde{Y}_{\left(k-1\right)h}h+\widehat{N}_{b}^{-1}B\gamma_{h}\left(V_{\widehat{N}_{b}}-V_{0}\right)\right)h\sum_{k=1}^{\widehat{T}_{b}}\widetilde{Y}_{\left(k-1\right)h}}{h\sum_{k=1}^{\widehat{T}_{b}}\widetilde{Y}_{\left(k-1\right)h}^{2}-\widehat{N}_{b}^{-1}\left(h\sum_{k=1}^{\widehat{T}_{b}}\widetilde{Y}_{\left(k-1\right)h}\right)^{2}}\\
 & +o_{p}\left(1\right),
\end{align*}
and thus
\begin{align*}
\widehat{\alpha}_{1} & =\frac{\sum_{k=1}^{T_{b}^{0}}\left(\mu_{1}^{0}h+\alpha_{1}^{0}\widetilde{Y}_{\left(k-1\right)h}h+\Delta_{h}\widetilde{V}_{k}\right)\widetilde{Y}_{\left(k-1\right)h}}{h\sum_{k=1}^{\widehat{T}_{b}}\widetilde{Y}_{\left(k-1\right)h}^{2}-\widehat{N}_{b}^{-1}\left(h\sum_{k=1}^{\widehat{T}_{b}}\widetilde{Y}_{\left(k-1\right)h}\right)^{2}}\\
 & \quad-\frac{\left(\mu_{1}^{0}+\alpha_{1}^{0}\widehat{N}_{b}^{-1}\sum_{k=1}^{T_{b}^{0}}\widetilde{Y}_{\left(k-1\right)h}h+\widehat{N}_{b}^{-1}B\gamma_{h}\left(V_{N_{b}^{0}}-V_{0}\right)\right)h\sum_{k=1}^{T_{b}^{0}}\widetilde{Y}_{\left(k-1\right)h}}{h\sum_{k=1}^{\widehat{T}_{b}}\widetilde{Y}_{\left(k-1\right)h}^{2}-\widehat{N}_{b}^{-1}\left(h\sum_{k=1}^{\widehat{T}_{b}}\widetilde{Y}_{\left(k-1\right)h}\right)^{2}}\\
 & \quad\times h\sum_{k=1}^{T_{b}^{0}}\widetilde{Y}_{\left(k-1\right)h}\\
 & \quad-\frac{\sum_{k=\widehat{T}_{b}+1}^{T_{b}^{0}}\left(\mu_{1}^{0}h+\alpha_{1}^{0}\widetilde{Y}_{\left(k-1\right)h}h+\Delta_{h}\widetilde{V}_{k}\right)\widetilde{Y}_{\left(k-1\right)h}}{h\sum_{k=1}^{\widehat{T}_{b}}\widetilde{Y}_{\left(k-1\right)h}^{2}-\widehat{N}_{b}^{-1}\left(h\sum_{k=1}^{\widehat{T}_{b}}\widetilde{Y}_{\left(k-1\right)h}\right)^{2}}\\
 & \quad+\frac{\widehat{N}_{b}^{-1}\left(\sum_{k=\widehat{T}_{b}+1}^{T_{b}^{0}}\mu_{1}^{0}h+\alpha_{1}^{0}\sum_{k=\widehat{T}_{b}+1}^{T_{b}^{0}}\widetilde{Y}_{\left(k-1\right)h}h+B\gamma_{h}\left(V_{N_{b}^{0}}-V_{\widehat{N}_{b}}\right)\right)}{h\sum_{k=1}^{\widehat{T}_{b}}\widetilde{Y}_{\left(k-1\right)h}^{2}-\widehat{N}_{b}^{-1}\left(h\sum_{k=1}^{\widehat{T}_{b}}\widetilde{Y}_{\left(k-1\right)h}\right)^{2}}\\
 & \quad\times h\sum_{k=\widehat{T}_{b}+1}^{T_{b}^{0}}\widetilde{Y}_{\left(k-1\right)h}.
\end{align*}
 By part (ii) of Theorem \ref{Theorem - Asy Dist - Pre}, $N_{b}^{0}-\widehat{N}_{b}=O_{p}(h^{1-\kappa})$,
and thus it is easy to see that the third and fourth terms converge
to zero in probability at a rate slower than $h^{1-\kappa}$. For
the first and second terms, recalling that $\Delta_{h}\widetilde{V}_{h,k}=h^{1/2}\Delta V_{h,k}$
from \eqref{Small-Dispersion Model}, we have by ordinary convergence
of approximations to Riemann sums, Lemma \ref{Lemma QLS Riemann}
and the continuity of probability limits, 
\begin{align*}
\alpha_{1}^{0}\sum_{k=1}^{T_{b}^{0}}\widetilde{Y}_{\left(k-1\right)h}h\overset{P}{\rightarrow}\alpha_{1}^{0}\int_{0}^{N_{b}^{0}}\widetilde{Y}_{s}ds, & \qquad\sum_{k=1}^{T_{b}^{0}}\mu_{1}^{0}h\overset{P}{\rightarrow}\mu_{1}^{0}\int_{0}^{N_{b}^{0}}ds,
\end{align*}
and by Lemma \ref{Lemma QLS 3-1}, $\sum_{k=1}^{T_{b}^{0}}\widetilde{Y}_{\left(k-1\right)h}\Delta_{h}\widetilde{V}_{k}\overset{P}{\rightarrow}0.$
Thus, we deduce that
\begin{align}
\widehat{\alpha}_{1} & =\alpha_{1}^{0}+O_{p}\left(B\gamma_{h}\right).\label{Eq alpha1 consistent}
\end{align}
 Using \eqref{Eq alpha1 consistent} into \eqref{Eq Expression miu1},
\begin{align*}
\widehat{\mu}_{1}\widehat{N}_{b} & =\widetilde{Y}_{\widehat{N}_{b}}-\widetilde{Y}_{0}-\alpha_{1}^{0}h\sum_{k=1}^{\widehat{T}_{b}}\widetilde{Y}_{\left(k-1\right)h}-O_{p}\left(B\gamma_{h}\right),\\
 & =\widetilde{Y}_{\widehat{N}_{b}}-\widetilde{Y}_{0}-\alpha_{1}^{0}h\sum_{k=1}^{T_{b}^{0}}\widetilde{Y}_{\left(k-1\right)h}-\alpha_{1}^{0}h\sum_{k=\widehat{T}_{b}+1}^{T_{b}^{0}}\widetilde{Y}_{\left(k-1\right)h}-o_{p}\left(1\right).
\end{align*}
 By part (ii) of Theorem \ref{Theorem - Asy Dist - Pre}, the number
of terms in the second sum above increases at rate $T^{\kappa}$ and
thus, $\alpha_{1}^{0}h\sum_{k=\widehat{T}_{b}+1}^{T_{b}^{0}}\widetilde{Y}_{\left(k-1\right)h}=KO_{p}(h^{1-\kappa}),$
where we have also used standard estimates for the drift arising from
the Burkh\"{o}lder-Davis-Gundy inequality. This gives 
\begin{align*}
\widehat{\mu}_{1}\widehat{N}_{b} & =\widetilde{Y}_{N_{b}^{0}}-\widetilde{Y}_{0}-\alpha_{1}^{0}\int_{0}^{N_{b}^{0}}\widetilde{Y}_{s}ds-\alpha_{1}^{0}O_{p}\left(h^{1-\kappa}\right)-o_{p}\left(1\right).
\end{align*}
Noting that
\begin{align*}
\widetilde{Y}_{N_{b}^{0}}-\widetilde{Y}_{0} & =\mu_{1}^{0}N_{b}^{0}+\alpha_{1}^{0}\int_{0}^{N_{b}^{0}}\widetilde{Y}_{s}ds+O_{p}\left(B\gamma_{h}\right)\left(V_{N_{b}^{0}}-V_{0}\right),
\end{align*}
we have $\widehat{\mu}_{1}N_{b}^{0}=\mu_{1}^{0}N_{b}^{0}+O_{p}\left(B\gamma_{h}\right)(V_{N_{b}^{0}}-V_{0}),$
which yields
\begin{align}
\widehat{\mu}_{1} & =\mu_{1}^{0}+O_{p}\left(B\gamma_{h}\right).\label{Eq. miu1 consistent}
\end{align}
Thus, as $h\downarrow0$, $\widehat{\mu}_{1}$ is consistent for $\mu_{1}^{0}$.
The case where $\widehat{T}_{b}>T_{b}^{0}$ can be treated in the
same fashion and is omitted. The consistency proof for $(\widehat{\mu}_{2},\,\widehat{\alpha}_{2})'$
is analogous and also omitted. The second step is to construct the
least-squares residuals and scaling them up. The residuals are constructed
as follows,
\begin{align*}
\widehat{u}_{kh} & =\begin{cases}
h^{-1/2}\left(\Delta_{h}\widetilde{Y}_{k}-\widehat{\mu}_{1}\widetilde{x}_{kh}^{\left(1\right)}-\widehat{\alpha}_{1}\widetilde{x}_{kh}^{\left(2\right)}\right), & k\leq\widehat{T}_{b}\\
h^{-1/2}\left(\Delta_{h}\widetilde{Y}_{k}-\widehat{\mu}_{2}\widetilde{x}_{kh}^{\left(1\right)}-\widehat{\alpha}_{2}\widetilde{x}_{kh}^{\left(2\right)}\right), & k>\widehat{T}_{b},
\end{cases}
\end{align*}
where $\widetilde{x}_{kh}^{\left(1\right)}=h$ and $\widetilde{x}_{kh}^{\left(2\right)}=\widetilde{Y}_{\left(k-1\right)h}h$.
This yields, for $k\leq T_{b}^{0}\leq\widehat{T}_{b},$ 
\begin{align*}
\widehat{u}_{kh} & =h^{-1/2}\left(\mu_{1}^{0}h+\alpha_{1}^{0}\widetilde{Y}_{\left(k-1\right)h}h+B\gamma_{h}\Delta_{h}V_{k}-\widehat{\mu}_{1}h-\widehat{\alpha}_{1}\widetilde{Y}_{\left(k-1\right)h}h\right),
\end{align*}
 and using \eqref{Eq alpha1 consistent} and \eqref{Eq. miu1 consistent},
\begin{align}
\widehat{u}_{kh} & =h^{-1/2}(\mu_{1}^{0}h+\alpha_{1}^{0}\widetilde{Y}_{\left(k-1\right)h}h+B\gamma_{h}\Delta_{h}V_{k}-\mu_{1}^{0}h\nonumber \\
 & \quad-O_{p}\left(h^{3/2}\right)-\alpha_{1}^{0}\widetilde{Y}_{\left(k-1\right)h}h-O_{p}\left(h^{3/2}\right))\nonumber \\
 & =h^{-1/2}B\gamma\Delta_{h}V_{k}-O_{p}\left(h\right).\label{Eq. Utkh 1}
\end{align}
Similarly, for $T_{b}^{0}\leq\widehat{T}_{b}\leq k,$ 
\begin{align}
\widehat{u}_{kh} & =h^{-1/2}B\gamma_{h}\Delta_{h}V_{k}-O_{p}\left(h\right),\label{Eq. Utkh 2}
\end{align}
whereas for $\widehat{T}_{b}<k\leq T_{b}^{0},$
\begin{align}
\widehat{u}_{kh} & =h^{-1/2}(\mu_{1}^{0}h+\alpha_{1}^{0}\widetilde{Y}_{\left(k-1\right)h}h+B\gamma_{h}\Delta_{h}V_{k}-\mu_{2}^{0}h\nonumber \\
 & \quad-O_{p}\left(h^{3/2}\right)-\alpha_{2}^{0}\widetilde{Y}_{\left(k-1\right)h}h-O_{p}\left(h^{3/2}\right))\nonumber \\
 & =h^{-1/2}\left(-\mu_{\delta}^{0}h-\alpha_{\delta}^{0}\widetilde{Y}_{\left(k-1\right)h}h+B\gamma_{h}\Delta_{h}V_{k}-O_{p}\left(h^{3/2}\right)\right)\nonumber \\
 & =-\mu_{\delta}^{0}h^{1/2}-\alpha_{\delta}^{0}\widetilde{Y}_{\left(k-1\right)h}h^{1/2}+h^{-1/2}B\gamma_{h}\Delta_{h}V_{k}-O_{p}\left(h\right).\label{Eq. Utkh}
\end{align}
 Next, note that $\sum_{k=\widehat{T}_{b}+1}^{T_{b}^{0}}\mu_{\delta}^{0}h^{1/2}\leq Kh^{1/2-\kappa}$
and $\sum_{k=\widehat{T}_{b}+1}^{T_{b}^{0}}\alpha_{\delta}^{0}\widetilde{Y}_{\left(k-1\right)h}h^{1/2}\leq Kh^{1/2-\kappa}$
since by Theorem \ref{Theorem - Asy Dist - Pre}-(ii) there are $T^{\kappa}$
terms in each sum. Moreover, recall that $e_{kh}=\Delta_{h}e_{k}^{*}\sim\mathscr{N}(0,\,\sigma_{e,k-1}^{2}h)$
and thus\footnote{The same bound holds for the corresponding sum involving the other
terms in $\Delta_{h}V_{k}$.}
\[
\sum_{k=\widehat{T}_{b}+1}^{T_{b}^{0}}e_{kh}=\sqrt{h}\sum_{k=\widehat{T}_{b}+1}^{T_{b}^{0}}h^{-1/2}e_{kh}=h^{1/2-\kappa}o_{p}\left(1\right).
\]
 Therefore, $\sum_{k=\widehat{T}_{b}+1}^{T_{b}^{0}}\widehat{u}_{kh}=Ko_{p}(h^{1/2-\kappa}).$
Since $\kappa\in\left(0,\,1/2\right),$ this shows that the residuals
$\widehat{u}_{kh}$ from equation \eqref{Eq. Utkh} are asymptotically
negligible. That is, asymptotically the estimator of $((\beta_{S}^{0})',\,(\delta_{Z,1}^{0})',\,(\delta_{Z,2}^{0})')'$
minimizes (assuming $\widehat{T}_{b}\leq T_{b}^{0}$), 
\begin{align*}
\sum_{k=1}^{\widehat{T}_{b}}\left(\widehat{u}_{kh}-\widetilde{x}'_{kh}\beta_{S}\right)^{2} & +\sum_{k=T_{b}^{0}+1}^{T}\left(\widehat{u}_{kh}-\widetilde{x}'_{kh}\beta_{S}-\widetilde{z}'_{0,kh}\delta_{S}\right)^{2}+o_{p}\left(1\right),
\end{align*}
 where $X=[\widetilde{X}^{\left(1\right)}$, $\beta^{0}=\begin{bmatrix}\mu_{1}^{0} & \alpha_{1}^{0} & \left(\beta_{S}^{0}\right)'\end{bmatrix}',$
and $Z_{0}$ and $\delta_{S}^{0}$ are partitioned accordingly. The
subscript $S$ indicates that these are the parameters of the stochastic
semimartingale regressors. This is exactly the same regression model
as in Proposition \ref{Proposition OLS Asymtptoc Distribu}. Hence,
the consistency result for the slope coefficients of the semimartingale
regressors follows from the same proof. The following regression model
estimated by least-squares provides consistent estimates for $\beta_{S}^{0}$
and $\delta_{S}^{0}$: $\widehat{U}=\widetilde{X}\widehat{\beta}_{S}+\widehat{Z}_{0}\widehat{\delta}_{S}+\textrm{residuals},$
where 
\begin{align*}
\widehat{Z}_{0} & =\begin{bmatrix}\widetilde{z}_{1}^{\left(1\right)} & \cdots & \widetilde{z}_{1}^{\left(p\right)}\\
\vdots & \ddots & \vdots\\
\widetilde{z}_{\widehat{T}_{b}h}^{\left(1\right)} & \cdots & \widetilde{z}_{\widehat{T}_{b}h}^{\left(p\right)}\\
\widetilde{z}_{\left(T_{b}^{0}+1\right)h}^{\left(1\right)} & \cdots & \widetilde{z}_{\left(T_{b}^{0}+1\right)h}^{\left(p\right)}\\
\vdots & \ddots & \vdots\\
\widetilde{z}_{N}^{\left(1\right)} & \cdots & \widetilde{z}_{N}^{\left(p\right)}
\end{bmatrix},
\end{align*}
 and $\widehat{U}=(\widehat{u}_{kh};\,k=1,\ldots,\,\widehat{T}_{b},\,T_{b}^{0}+1,\ldots,\,N)$.
Therefore, using \eqref{Eq. Utkh 1} and \eqref{Eq. Utkh 2}, we have
\begin{align*}
h^{-1/2}\begin{bmatrix}\widehat{\beta}_{S}-\beta^{0}\\
\widehat{\delta}_{S}-\delta^{0}
\end{bmatrix} & =\begin{bmatrix}\widetilde{X}'\widetilde{X} & \widetilde{X}'\widehat{Z}_{0}\\
\widehat{Z}_{0}'\widetilde{X} & \widehat{Z}_{0}'\widehat{Z}_{0}
\end{bmatrix}^{-1}\\
 & \quad\times h^{-1/2}\begin{bmatrix}\widetilde{X}'e & \widetilde{X}'\left(Z_{0}-\widehat{Z}_{0}\right)\delta^{0}+\widetilde{X}'AO_{p}\left(h\right)\\
\widehat{Z}_{0}'e & \widehat{Z}_{0}'\left(Z_{0}-\widehat{Z}_{0}\right)\delta^{0}+\widehat{Z}_{0}'AO_{p}\left(h\right)
\end{bmatrix},
\end{align*}
 for some matrix $A=O_{p}\left(1\right)$. It then follows by the
same proof as in Proposition \ref{Proposition OLS Asymtptoc Distribu}
that 
\begin{align}
\begin{bmatrix}\widetilde{X}'\widetilde{X} & \widetilde{X}'\widehat{Z}_{0}\\
\widehat{Z}_{0}'\widetilde{X} & \widehat{Z}_{0}'\widehat{Z}_{0}
\end{bmatrix}^{-1}\widetilde{X}'AO_{p}\left(h^{1/2}\right) & =o_{p}\left(1\right),\label{AA4}
\end{align}
 and 
\begin{align}
\begin{bmatrix}\widetilde{X}'\widetilde{X} & \widetilde{X}'\widehat{Z}_{0}\\
\widehat{Z}_{0}'\widetilde{X} & \widehat{Z}_{0}'\widehat{Z}_{0}
\end{bmatrix}^{-1}\frac{1}{h^{1/2}}\widetilde{X}'\left(Z_{0}-\widehat{Z}_{0}\right)\delta^{0} & =O_{p}\left(1\right)o_{p}\left(1\right)=o_{p}\left(1\right).\label{AA5}
\end{align}
The same arguments can be used for $\widehat{Z}_{0}'(Z_{0}-\widehat{Z}_{0})\delta^{0}$
and $\widehat{Z}_{0}'AO_{p}\left(h\right).$ Therefore, in view of
\eqref{Eq alpha1 consistent} and \eqref{Eq. miu1 consistent}, we
obtain $\widehat{\mu}_{1}=\mu_{1}^{0}+o_{p}\left(1\right)$ and $\widehat{\alpha}_{1}=\alpha_{1}^{0}+o_{p}\left(1\right)$,
respectively, whereas \eqref{AA4} and \eqref{AA5} imply $\widehat{\beta}_{S}=\beta_{S}^{0}+o_{p}\left(1\right)$
and $\widehat{\delta}_{S}=\delta_{S}^{0}+o_{p}\left(1\right)$, respectively.
Under the setting where the magnitude of the shifts is local to zero,
we observe that by Proposition \ref{Prop 3 Asym}, $\widehat{N}_{b}-\widehat{N}_{b}^{0}=O_{p}(h^{1-\kappa})$
and one can follow the same steps that led to \eqref{Eq alpha1 consistent}
and \eqref{Eq. miu1 consistent} and proceed as above. The final result
is $\widehat{\theta}=\theta^{0}+o_{p}\left(1\right)$, which completes
the proof. $\square$ 

\subsubsection{\label{Negligibility-of-the, Pre}Negligibility of the Drift Term}

Recall Lemma \ref{Lemma Cross Product Null} and apply the same proof
as in Section \ref{Negligibility-of-the, CT}. Of course, the negligibility
only applies to the drift processes $\mu_{\cdot,t}$ from \eqref{Model Regressors Integral Form}
(i.e., only the drift processes of the semimartingale regressors)
and not to $\mu_{1}^{0},\,\mu_{2}^{0},\,\alpha_{1}^{0}$ or $\alpha_{2}^{0}$.
The steps are omitted since they are the same.

\section{\label{sec:Additional-Simulations-Results}Additional Simulations
Results about HDR Confidence Sets}

We continue with the analysis of finite-sample from Section \ref{Section Small-Sample-Effectiveness-of}.
We consider discrete-time DGPs of the form
\begin{align}
y_{t}=D'_{t}\nu^{0}+Z'_{t}\beta^{0}+Z'_{t}\delta_{Z}^{0}\boldsymbol{1}_{\left\{ t>T_{b}^{0}\right\} }+e_{t}, & \qquad\qquad t=1,\ldots,\,T,\label{Eq. DGP Simulation Study}
\end{align}
with $T=100$ and, without loss of generality, $\nu^{0}=0$ (except
for M5-M6, M8-M9). We consider eight versions of \eqref{Eq. DGP Simulation Study}:
M3 involves a break in the simultaneous mean and variance of an $i.i.d.$
series with $Z_{t}=1$ for all $t$, $D_{t}$ absent, and $e_{t}=(1+\mathbf{1}_{\left\{ t>T_{b}^{0}\right\} })u_{t}$
with $u_{t}\sim i.i.d.\,\mathscr{N}\left(0,\,1\right)$; M4 is the
same as M1 but with stationary Gaussian AR(1) disturbances $e_{t}=0.3e_{t-1}+u_{t}$,
$u_{t}\sim i.i.d.\,\mathscr{N}\left(0,\,0.49\right)$; M5 is a partial
structural change model with $D_{t}=1$ for all $t$, $\nu^{0}=1$
and $Z_{t}=0.5Z_{t}+u_{t}$ with $u_{t}\sim i.i.d.\,\mathscr{N}\left(0,\,0.75\right)$
independent of $e_{t}\sim i.i.d.\,\mathscr{N}\left(0,\,1\right)$;
M6 is similar to M5 but with $u_{t}\sim i.i.d.\,\mathscr{N}\left(0,\,1\right)$
and heteroskedastic disturbances given by $e_{t}=v_{t}\left|Z_{t}\right|$
where $v_{t}$ is a sequence of $i.i.d.$ $\mathscr{N}\left(0,\,1\right)$
random variables independent of $\left\{ Z_{t}\right\} $; M7 is the
same as M4 but with $u_{t}$ drawn from a $\textrm{t}_{\upsilon}$
distribution with $\upsilon=5$ degrees of freedom; M8 is a model
with a lagged dependent variable with $D_{t}=y_{t-1}$, $Z_{t}=1$,
$e_{t}\sim i.i.d.\,\mathscr{N}\left(0,\,0.49\right)$, $\nu^{0}=0.3$
and $Z'_{t}\delta_{Z}^{0}\boldsymbol{1}_{\left\{ t>T_{b}^{0}\right\} }$
is replaced by $Z'_{t}\left(1-\nu^{0}\right)\delta_{Z}^{0}\boldsymbol{1}_{\left\{ t>T_{b}^{0}\right\} }$;
M9 has FIGARCH(1,d,1) errors given by $e_{t}=\sigma_{t}u_{t},\,u_{t}\sim i.i.d.\,\mathscr{N}\left(0,\,1\right)$
and $\sigma_{t}=0.1+(1-0.2L\left(1-L\right)^{d})e_{t}^{2}$ where
$d=0.6$ is the order of differencing and $L$ the lag operator, $D_{t}=1$,
$\nu^{0}=1$ and $Z_{t}\sim i.i.d.\,\mathscr{N}\left(1,\,1.44\right)$
independent of $e_{t}$. M10 is similar to M6 but with an ARFIMA$\mathrm{(0.3,\,\mathit{d},\,0)}$
regressor $Z_{t}$ with order of differencing $d=0.5$, $\mathrm{Var}\left(Z_{t}\right)=1$
and $e_{t}\sim\mathscr{N}\left(0,\,1\right)$ independent of $\left\{ Z_{t}\right\} $.
We set $\beta^{0}=1$ in all models, except in M8 where $\beta^{0}=0$.
The Results are reported in Table \ref{Table M2}-\ref{Table M10 LongMemory}.

\pagebreak{}

\begin{table}[H]
\caption{\label{Table M2}Small-sample coverage rate and length of the confidence
set for model M3}

\begin{centering}
{\scriptsize{}}%
\begin{tabular}{cccccccccc}
\hline 
 &  & \multicolumn{2}{c}{{\scriptsize{}$\delta_{Z}^{0}=0.3$}} & \multicolumn{2}{c}{{\scriptsize{}$\delta_{Z}^{0}=0.6$}} & \multicolumn{2}{c}{{\scriptsize{}$\delta_{Z}^{0}=1$}} & \multicolumn{2}{c}{{\scriptsize{}$\delta_{Z}^{0}=2$}}\tabularnewline
 &  & {\scriptsize{}$\textrm{Cov.}$} & {\scriptsize{}$\textrm{Lgth.}$} & {\scriptsize{}$\textrm{Cov.}$} & {\scriptsize{}$\textrm{Lgth.}$} & {\scriptsize{}$\textrm{Cov.}$} & {\scriptsize{}$\textrm{Lgth.}$} & {\scriptsize{}$\textrm{Cov.}$} & {\scriptsize{}$\textrm{Lgth.}$}\tabularnewline
\hline 
{\scriptsize{}$\lambda_{0}=0.5$} & {\scriptsize{}HDR} & {\scriptsize{}0.970} & {\scriptsize{}86.65} & {\scriptsize{}0.937} & {\scriptsize{}76.29} & {\scriptsize{}0.901} & {\scriptsize{}55.59} & {\scriptsize{}0.934} & {\scriptsize{}26.11}\tabularnewline
 & {\scriptsize{}Bai (1997)} & {\scriptsize{}0.854} & {\scriptsize{}70.60} & {\scriptsize{}0.843} & {\scriptsize{}58.27} & {\scriptsize{}0.857} & {\scriptsize{}40.70} & {\scriptsize{}0.923} & {\scriptsize{}14.24}\tabularnewline
 & {\scriptsize{}$\widehat{U}_{T}.\textrm{neq}$} & {\scriptsize{}0.961} & {\scriptsize{}88.95} & {\scriptsize{}0.961} & {\scriptsize{}80.33} & {\scriptsize{}0.961} & {\scriptsize{}61.15} & {\scriptsize{}0.964} & {\scriptsize{}32.16}\tabularnewline
 & {\scriptsize{}ILR} & {\scriptsize{}0.989} & {\scriptsize{}92.53} & {\scriptsize{}0.985} & {\scriptsize{}84.06} & {\scriptsize{}0.977} & {\scriptsize{}58.05} & {\scriptsize{}0.958} & {\scriptsize{}12.31}\tabularnewline
{\scriptsize{}$\lambda_{0}=0.35$} & {\scriptsize{}HDR} & {\scriptsize{}0.976} & {\scriptsize{}89.81} & {\scriptsize{}0.961} & {\scriptsize{}83.26} & {\scriptsize{}0.935} & {\scriptsize{}64.87} & {\scriptsize{}0.934} & {\scriptsize{}26.11}\tabularnewline
 & {\scriptsize{}Bai (1997)} & {\scriptsize{}0.823} & {\scriptsize{}69.86} & {\scriptsize{}0.822} & {\scriptsize{}55.87} & {\scriptsize{}0.844} & {\scriptsize{}38.91} & {\scriptsize{}0.932} & {\scriptsize{}14.24}\tabularnewline
 & {\scriptsize{}$\widehat{U}_{T}.\textrm{neq}$} & {\scriptsize{}0.963} & {\scriptsize{}89.84} & {\scriptsize{}0.963} & {\scriptsize{}82.26} & {\scriptsize{}0.961} & {\scriptsize{}65.87} & \multicolumn{1}{c}{{\scriptsize{}0.964}} & {\scriptsize{}32.16}\tabularnewline
 & {\scriptsize{}ILR} & {\scriptsize{}0.990} & {\scriptsize{}93.48} & {\scriptsize{}0.985} & {\scriptsize{}88.69} & {\scriptsize{}0.982} & {\scriptsize{}68.23} & {\scriptsize{}0.977} & {\scriptsize{}15.45}\tabularnewline
{\scriptsize{}$\lambda_{0}=0.2$} & {\scriptsize{}HDR} & {\scriptsize{}0.978} & {\scriptsize{}90.39} & {\scriptsize{}0.975} & {\scriptsize{}85.89} & {\scriptsize{}0.934} & {\scriptsize{}70.05} & {\scriptsize{}0.957} & {\scriptsize{}29.63}\tabularnewline
 & {\scriptsize{}Bai (1997)} & {\scriptsize{}0.782} & {\scriptsize{}70.24} & {\scriptsize{}0.805} & {\scriptsize{}56.37} & {\scriptsize{}0.831} & {\scriptsize{}37.66} & {\scriptsize{}0.928} & {\scriptsize{}14.80}\tabularnewline
 & {\scriptsize{}$\widehat{U}_{T}.\textrm{neq}$} & {\scriptsize{}0.968} & {\scriptsize{}91.11} & {\scriptsize{}0.968} & {\scriptsize{}87.62} & {\scriptsize{}0.972} & {\scriptsize{}78.17} & {\scriptsize{}0.967} & {\scriptsize{}46.24}\tabularnewline
 & {\scriptsize{}ILR} & {\scriptsize{}0.980} & {\scriptsize{}93.32} & {\scriptsize{}0.981} & {\scriptsize{}91.60} & {\scriptsize{}0.978} & {\scriptsize{}81.60} & {\scriptsize{}0.981} & {\scriptsize{}22.60}\tabularnewline
\hline 
\end{tabular}{\scriptsize\par}
\par\end{centering}
\noindent\begin{minipage}[t]{1\columnwidth}%
{\scriptsize{}The model is $y_{t}=\beta^{0}+\delta_{Z}^{0}\mathbf{1}_{\left\{ t>\left\lfloor T\lambda_{0}\right\rfloor \right\} }+e_{t},\,e_{t}=\left(1+\mathbf{1}_{\left\{ t>\left\lfloor T\lambda_{0}\right\rfloor \right\} }\right)u_{t},\,u_{t}\sim i.i.d.\,\mathscr{N}\left(0,\,1\right),\,T=100$.
The notes of Table \ref{Table M1} apply.}%
\end{minipage}
\end{table}

\begin{table}[H]
\caption{\label{Table M3}Small-sample coverage rate and length of the confidence
set for model M4}

\begin{centering}
{\footnotesize{}}%
\begin{tabular}{cccccccccc}
\hline 
 &  & \multicolumn{2}{c}{{\scriptsize{}$\delta_{Z}^{0}=0.3$}} & \multicolumn{2}{c}{{\scriptsize{}$\delta_{Z}^{0}=0.6$}} & \multicolumn{2}{c}{{\scriptsize{}$\delta_{Z}^{0}=1$}} & \multicolumn{2}{c}{{\scriptsize{}$\delta_{Z}^{0}=2$}}\tabularnewline
 &  & {\scriptsize{}$\textrm{Cov.}$} & {\scriptsize{}$\textrm{Lgth.}$} & {\scriptsize{}$\textrm{Cov.}$} & {\scriptsize{}$\textrm{Lgth.}$} & {\scriptsize{}$\textrm{Cov.}$} & {\scriptsize{}$\textrm{Lgth.}$} & {\scriptsize{}$\textrm{Cov.}$} & {\footnotesize{}$\textrm{Lgth.}$}\tabularnewline
\hline 
{\scriptsize{}$\lambda_{0}=0.5$} & {\scriptsize{}HDR} & {\scriptsize{}0.904} & {\scriptsize{}72.44} & {\scriptsize{}0.901} & {\scriptsize{}57.37} & {\scriptsize{}0.919} & {\scriptsize{}29.70} & {\scriptsize{}0.971} & {\scriptsize{}5.85}\tabularnewline
 & {\scriptsize{}Bai (1997)} & {\scriptsize{}0.833} & {\scriptsize{}66.34} & {\scriptsize{}0.834} & {\scriptsize{}41.32} & {\scriptsize{}0.895} & {\scriptsize{}18.63} & {\scriptsize{}0.969} & {\scriptsize{}5.49}\tabularnewline
 & {\scriptsize{}$\widehat{U}_{T}.\textrm{eq}$} & {\scriptsize{}0.958} & {\scriptsize{}87.16} & {\scriptsize{}0.968} & {\scriptsize{}71.47} & {\scriptsize{}0.958} & {\scriptsize{}45.82} & {\scriptsize{}0.957} & {\scriptsize{}28.01}\tabularnewline
 & {\scriptsize{}ILR} & {\scriptsize{}0.932} & {\scriptsize{}79.38} & {\scriptsize{}0.944} & {\scriptsize{}53.48} & {\scriptsize{}0.966} & {\scriptsize{}21.98} & {\scriptsize{}0.993} & {\scriptsize{}4.87}\tabularnewline
{\scriptsize{}$\lambda_{0}=0.35$} & {\scriptsize{}HDR} & {\scriptsize{}0.910} & {\scriptsize{}70.98} & {\scriptsize{}0.902} & {\scriptsize{}53.88} & {\scriptsize{}0.917} & {\scriptsize{}28.07} & {\scriptsize{}0.973} & {\scriptsize{}5.99}\tabularnewline
 & {\scriptsize{}Bai (1997)} & {\scriptsize{}0.849} & {\scriptsize{}65.13} & {\scriptsize{}0.840} & {\scriptsize{}40.43} & {\scriptsize{}0.900} & {\scriptsize{}18.69} & {\scriptsize{}0.974} & {\scriptsize{}5.49}\tabularnewline
 & {\scriptsize{}$\widehat{U}_{T}.\textrm{eq}$} & {\scriptsize{}0.960} & {\scriptsize{}87.46} & {\scriptsize{}0.961} & {\scriptsize{}72.79} & {\scriptsize{}0.962} & {\scriptsize{}46.44} & \multicolumn{1}{c}{{\scriptsize{}0.961}} & {\scriptsize{}28.03}\tabularnewline
 & {\scriptsize{}ILR} & {\scriptsize{}0.942} & {\scriptsize{}80.94} & {\scriptsize{}0.946} & {\scriptsize{}55.20} & {\scriptsize{}0.965} & {\scriptsize{}23.55} & {\scriptsize{}0.993} & {\scriptsize{}4.93}\tabularnewline
{\scriptsize{}$\lambda_{0}=0.2$} & {\scriptsize{}HDR} & {\scriptsize{}0.905} & {\scriptsize{}72.26} & {\scriptsize{}0.913} & {\scriptsize{}50.61} & {\scriptsize{}0.933} & {\scriptsize{}25.07} & {\scriptsize{}0.973} & {\scriptsize{}6.35}\tabularnewline
 & {\scriptsize{}Bai (1997)} & {\scriptsize{}0.829} & {\scriptsize{}65.56} & {\scriptsize{}0.899} & {\scriptsize{}41.42} & {\scriptsize{}0.932} & {\scriptsize{}19.62} & {\scriptsize{}0.966} & {\scriptsize{}5.55}\tabularnewline
 & {\scriptsize{}$\widehat{U}_{T}.\textrm{eq}$} & {\scriptsize{}0.962} & {\scriptsize{}88.77} & {\scriptsize{}0.968} & {\scriptsize{}78.61} & {\scriptsize{}0.963} & {\scriptsize{}57.87} & {\scriptsize{}0.965} & {\scriptsize{}29.88}\tabularnewline
 & {\scriptsize{}ILR} & {\scriptsize{}0.938} & {\scriptsize{}83.24} & {\scriptsize{}0.951} & {\scriptsize{}63.66} & {\scriptsize{}0.972} & {\scriptsize{}28.94} & {\scriptsize{}0.994} & {\scriptsize{}5.16}\tabularnewline
\hline 
\end{tabular}{\footnotesize\par}
\par\end{centering}
\noindent\begin{minipage}[t]{1\columnwidth}%
{\scriptsize{}The model is $y_{t}=\beta^{0}+\delta_{Z}^{0}\mathbf{1}_{\left\{ t>\left\lfloor T\lambda_{0}\right\rfloor \right\} }+e_{t},\,e_{t}=0.3e_{t-1}+u_{t},\,u_{t}\sim i.i.d.\,\mathscr{N}\left(0,\,0.49\right),\,T=100$.
The notes of Table \ref{Table M1} apply.}%
\end{minipage}
\end{table}

\begin{table}[H]
\caption{\label{Table M4}Small-sample coverage rate and length of the confidence
set for model M5}

\begin{centering}
{\scriptsize{}}%
\begin{tabular}{cccccccccc}
\hline 
 &  & \multicolumn{2}{c}{{\scriptsize{}$\delta_{Z}^{0}=0.3$}} & \multicolumn{2}{c}{{\scriptsize{}$\delta_{Z}^{0}=0.6$}} & \multicolumn{2}{c}{{\scriptsize{}$\delta_{Z}^{0}=1$}} & \multicolumn{2}{c}{{\scriptsize{}$\delta_{Z}^{0}=2$}}\tabularnewline
 &  & {\scriptsize{}$\textrm{Cov.}$} & {\scriptsize{}$\textrm{Lgth.}$} & {\scriptsize{}$\textrm{Cov.}$} & {\scriptsize{}$\textrm{Lgth.}$} & {\scriptsize{}$\textrm{Cov.}$} & {\scriptsize{}$\textrm{Lgth.}$} & {\scriptsize{}$\textrm{Cov.}$} & {\scriptsize{}$\textrm{Lgth.}$}\tabularnewline
\hline 
{\scriptsize{}$\lambda_{0}=0.5$} & {\scriptsize{}HDR} & {\scriptsize{}0.915} & {\scriptsize{}77.14} & {\scriptsize{}0.912} & {\scriptsize{}61.71} & {\scriptsize{}0.910} & {\scriptsize{}30.64} & {\scriptsize{}0.912} & {\scriptsize{}7.15}\tabularnewline
 & {\scriptsize{}Bai (1997)} & {\scriptsize{}0.805} & {\scriptsize{}65.94} & {\scriptsize{}0.821} & {\scriptsize{}44.07} & {\scriptsize{}0.850} & {\scriptsize{}20.71} & {\scriptsize{}0.887} & {\scriptsize{}5.96}\tabularnewline
 & {\scriptsize{}$\widehat{U}_{T}.\textrm{eq}$} & {\scriptsize{}0.950} & {\scriptsize{}85.23} & {\scriptsize{}0.951} & {\scriptsize{}67.40} & {\scriptsize{}0.951} & {\scriptsize{}39.87} & {\scriptsize{}0.955} & {\scriptsize{}17.46}\tabularnewline
 & {\scriptsize{}ILR} & {\scriptsize{}0.961} & {\scriptsize{}84.37} & {\scriptsize{}0.966} & {\scriptsize{}59.94} & {\scriptsize{}0.977} & {\scriptsize{}26.09} & {\scriptsize{}0.986} & {\scriptsize{}7.14}\tabularnewline
{\scriptsize{}$\lambda_{0}=0.35$} & {\scriptsize{}HDR} & {\scriptsize{}0.915} & {\scriptsize{}75.53} & {\scriptsize{}0.911} & {\scriptsize{}58.88} & {\scriptsize{}0.905} & {\scriptsize{}29.77} & {\scriptsize{}0.912} & {\scriptsize{}7.27}\tabularnewline
 & {\scriptsize{}Bai (1997)} & {\scriptsize{}0.821} & {\scriptsize{}64.69} & {\scriptsize{}0.826} & {\scriptsize{}42.93} & {\scriptsize{}0.849} & {\scriptsize{}20.77} & {\scriptsize{}0.888} & {\scriptsize{}5.99}\tabularnewline
 & {\scriptsize{}$\widehat{U}_{T}.\textrm{eq}$} & {\scriptsize{}0.948} & {\scriptsize{}85.48} & {\scriptsize{}0.948} & {\scriptsize{}68.95} & {\scriptsize{}0.948} & {\scriptsize{}41.40} & \multicolumn{1}{c}{{\scriptsize{}0.954}} & {\scriptsize{}17.57}\tabularnewline
 & {\scriptsize{}ILR} & {\scriptsize{}0.959} & {\scriptsize{}84.67} & {\scriptsize{}0.964} & {\scriptsize{}61.55} & {\scriptsize{}0.973} & {\scriptsize{}27.70} & {\scriptsize{}0.987} & {\scriptsize{}7.13}\tabularnewline
{\scriptsize{}$\lambda_{0}=0.2$} & {\scriptsize{}HDR} & {\scriptsize{}0.911} & {\scriptsize{}74.46} & {\scriptsize{}0.931} & {\scriptsize{}56.22} & {\scriptsize{}0.935} & {\scriptsize{}29.22} & {\scriptsize{}0.929} & {\scriptsize{}7.85}\tabularnewline
 & {\scriptsize{}Bai (1997)} & {\scriptsize{}0.820} & {\scriptsize{}64.06} & {\scriptsize{}0.870} & {\scriptsize{}42.86} & {\scriptsize{}0.896} & {\scriptsize{}22.11} & {\scriptsize{}0.887} & {\scriptsize{}6.16}\tabularnewline
 & {\scriptsize{}$\widehat{U}_{T}.\textrm{eq}$} & {\scriptsize{}0.952} & {\scriptsize{}86.80} & {\scriptsize{}0.956} & {\scriptsize{}75.20} & {\scriptsize{}0.952} & {\scriptsize{}51.99} & {\scriptsize{}0.952} & {\scriptsize{}19.92}\tabularnewline
 & {\scriptsize{}ILR} & {\scriptsize{}0.961} & {\scriptsize{}86.03} & {\scriptsize{}0.964} & {\scriptsize{}68.69} & {\scriptsize{}0.978} & {\scriptsize{}36.34} & {\scriptsize{}0.985} & {\scriptsize{}7.51}\tabularnewline
\hline 
\end{tabular}{\scriptsize\par}
\par\end{centering}
\noindent\begin{minipage}[t]{1\columnwidth}%
{\scriptsize{}The model is $y_{t}=\nu^{0}+Z_{t}\beta^{0}+Z_{t}\delta_{Z}^{0}\mathbf{1}_{\left\{ t>\left\lfloor T\lambda_{0}\right\rfloor \right\} }+e_{t},\,X_{t}=0.5X_{t-1}+u_{t},\,u_{t}\sim i.i.d.\,\mathscr{N}\left(0,\,0.75\right),\,e_{t}\sim\textrm{i.i.d.}\mathscr{N}\left(0,\,1\right),\,T=100$.
The notes of Table \ref{Table M1} apply.}%
\end{minipage}
\end{table}

\begin{table}[H]
\caption{\label{Table M10}Small-sample coverage rate and length of the confidence
set for model M6}

\begin{centering}
{\scriptsize{}}%
\begin{tabular}{cccccccccc}
\hline 
 &  & \multicolumn{2}{c}{{\scriptsize{}$\delta_{Z}^{0}=0.3$}} & \multicolumn{2}{c}{{\scriptsize{}$\delta_{Z}^{0}=0.6$}} & \multicolumn{2}{c}{{\scriptsize{}$\delta_{Z}^{0}=1$}} & \multicolumn{2}{c}{{\scriptsize{}$\delta_{Z}^{0}=2$}}\tabularnewline
 &  & {\scriptsize{}$\textrm{Cov.}$} & {\scriptsize{}$\textrm{Lgth.}$} & {\scriptsize{}$\textrm{Cov.}$} & {\scriptsize{}$\textrm{Lgth.}$} & {\scriptsize{}$\textrm{Cov.}$} & {\scriptsize{}$\textrm{Lgth.}$} & {\scriptsize{}$\textrm{Cov.}$} & {\scriptsize{}$\textrm{Lgth.}$}\tabularnewline
\hline 
{\scriptsize{}$\lambda_{0}=0.5$} & {\scriptsize{}HDR} & {\scriptsize{}0.920} & {\scriptsize{}77.03} & {\scriptsize{}0.923} & {\scriptsize{}70.69} & {\scriptsize{}0.930} & {\scriptsize{}60.02} & {\scriptsize{}0.969} & {\scriptsize{}35.03}\tabularnewline
 & {\scriptsize{}Bai (1997)} & {\scriptsize{}0.690} & {\scriptsize{}56.73} & {\scriptsize{}0.716} & {\scriptsize{}41.63} & {\scriptsize{}0.783} & {\scriptsize{}27.53} & {\scriptsize{}0.885} & {\scriptsize{}12.70}\tabularnewline
 & {\scriptsize{}$\widehat{U}_{T}.\textrm{eq}$} & {\scriptsize{}0.962} & {\scriptsize{}87.76} & {\scriptsize{}0.962} & {\scriptsize{}78.32} & {\scriptsize{}0.962} & {\scriptsize{}63.80} & {\scriptsize{}0.962} & {\scriptsize{}40.82}\tabularnewline
 & {\scriptsize{}ILR} & {\scriptsize{}0.790} & {\scriptsize{}71.07} & {\scriptsize{}0.805} & {\scriptsize{}59.66} & {\scriptsize{}0.824} & {\scriptsize{}40.78} & {\scriptsize{}0.909} & {\scriptsize{}11.63}\tabularnewline
{\scriptsize{}$\lambda_{0}=0.35$} & {\scriptsize{}HDR} & {\scriptsize{}0.928} & {\scriptsize{}76.41} & {\scriptsize{}0.925} & {\scriptsize{}68.21} & {\scriptsize{}0.933 } & {\scriptsize{}56.17} & {\scriptsize{}0.964 } & {\scriptsize{}31.73}\tabularnewline
 & {\scriptsize{}Bai (1997)} & {\scriptsize{}0.691} & {\scriptsize{}55.18} & {\scriptsize{}0.720} & {\scriptsize{}40.25} & {\scriptsize{}0.757 } & {\scriptsize{}26.90} & {\scriptsize{}0.883} & {\scriptsize{}12.62}\tabularnewline
 & {\scriptsize{}$\widehat{U}_{T}.\textrm{eq}$} & {\scriptsize{}0.953} & {\scriptsize{}87.76} & {\scriptsize{}0.953} & {\scriptsize{}78.55} & {\scriptsize{}0.953} & {\scriptsize{}64.81} & \multicolumn{1}{c}{{\scriptsize{}0.953}} & {\scriptsize{}41.98}\tabularnewline
 & {\scriptsize{}ILR} & {\scriptsize{}0.795} & {\scriptsize{}71.34} & {\scriptsize{}0.804} & {\scriptsize{}60.48} & {\scriptsize{}0.832} & {\scriptsize{}30.42} & {\scriptsize{}0.903} & {\scriptsize{}10.78}\tabularnewline
{\scriptsize{}$\lambda_{0}=0.2$} & {\scriptsize{}HDR} & {\scriptsize{}0.915} & {\scriptsize{}75.86} & {\scriptsize{}0.919} & {\scriptsize{}66.79} & {\scriptsize{}0.926} & {\scriptsize{}52.50} & {\scriptsize{}0.957} & {\scriptsize{}27.46}\tabularnewline
 & {\scriptsize{}Bai (1997)} & {\scriptsize{}0.707} & {\scriptsize{}55.03} & {\scriptsize{}0.770} & {\scriptsize{}39.77} & {\scriptsize{}0.828} & {\scriptsize{}26.82} & {\scriptsize{}0.901} & {\scriptsize{}12.68}\tabularnewline
 & {\scriptsize{}$\widehat{U}_{T}.\textrm{eq}$} & {\scriptsize{}0.951} & {\scriptsize{}88.48} & {\scriptsize{}0.952} & {\scriptsize{}82.09} & {\scriptsize{}0.954} & {\scriptsize{}71.84} & {\scriptsize{}0.950} & {\scriptsize{}50.72}\tabularnewline
 & {\scriptsize{}ILR} & {\scriptsize{}0.795} & {\scriptsize{}72.01} & {\scriptsize{}0.809} & {\scriptsize{}62.75} & {\scriptsize{}0.829} & {\scriptsize{}45.18} & {\scriptsize{}0.913} & {\scriptsize{}12.62}\tabularnewline
\hline 
\end{tabular}{\scriptsize\par}
\par\end{centering}
\noindent\begin{minipage}[t]{1\columnwidth}%
{\scriptsize{}The model is $y_{t}=\nu^{0}+Z_{t}\beta^{0}+Z_{t}\delta_{Z}^{0}\mathbf{1}_{\left\{ t>\left\lfloor T\lambda_{0}\right\rfloor \right\} }+e_{t},\,e_{t}=v_{t}\left|Z_{t}\right|,\,v_{t}\sim i.i.d.\,\mathscr{N}\left(0,\,1\right),\,Z_{t}=0.5Z_{t-1}+u_{t},\,u_{t}\sim i.i.d.\,\mathscr{N}\left(0,\,1\right)\,T=100$.
The notes of Table \ref{Table M1} apply.}%
\end{minipage}
\end{table}

\begin{table}[H]
\caption{\label{Table M5}Small-sample coverage rate and length of the confidence
set for model M7}

\begin{centering}
{\scriptsize{}}%
\begin{tabular}{cccccccccc}
\hline 
 &  & \multicolumn{2}{c}{{\scriptsize{}$\delta_{Z}^{0}=0.3$}} & \multicolumn{2}{c}{{\scriptsize{}$\delta_{Z}^{0}=0.6$}} & \multicolumn{2}{c}{{\scriptsize{}$\delta_{Z}^{0}=1$}} & \multicolumn{2}{c}{{\scriptsize{}$\delta_{Z}^{0}=2$}}\tabularnewline
 &  & {\scriptsize{}$\textrm{Cov.}$} & {\scriptsize{}$\textrm{Lgth.}$} & {\scriptsize{}$\textrm{Cov.}$} & {\scriptsize{}$\textrm{Lgth.}$} & {\scriptsize{}$\textrm{Cov.}$} & {\scriptsize{}$\textrm{Lgth.}$} & {\scriptsize{}$\textrm{Cov.}$} & {\scriptsize{}$\textrm{Lgth.}$}\tabularnewline
\hline 
{\scriptsize{}$\lambda_{0}=0.5$} & {\scriptsize{}HDR} & {\scriptsize{}0.918} & {\scriptsize{}75.64} & {\scriptsize{}0.910} & {\scriptsize{}67.46} & {\scriptsize{}0.931} & {\scriptsize{}48.54} & {\scriptsize{}0.957} & {\scriptsize{}12.50}\tabularnewline
 & {\scriptsize{}Bai (1997)} & {\scriptsize{}0.834} & {\scriptsize{}70.13} & {\scriptsize{}0.824} & {\scriptsize{}52.16} & {\scriptsize{}0.861} & {\scriptsize{}28.69} & {\scriptsize{}0.948} & {\scriptsize{}8.45}\tabularnewline
 & {\scriptsize{}$\widehat{U}_{T}.\textrm{eq}$} & {\scriptsize{}0.959} & {\scriptsize{}88.62} & {\scriptsize{}0.959} & {\scriptsize{}78.87} & {\scriptsize{}0.959} & {\scriptsize{}58.60} & {\scriptsize{}0.952} & {\scriptsize{}30.15}\tabularnewline
 & {\scriptsize{}ILR} & {\scriptsize{}0.969} & {\scriptsize{}86.75} & {\scriptsize{}0.959} & {\scriptsize{}67.91} & {\scriptsize{}0.967} & {\scriptsize{}34.13} & {\scriptsize{}0.995} & {\scriptsize{}9.17}\tabularnewline
{\scriptsize{}$\lambda_{0}=0.35$} & {\scriptsize{}HDR} & {\scriptsize{}0.926} & {\scriptsize{}74.78} & {\scriptsize{}0.914} & {\scriptsize{}64.86} & {\scriptsize{}0.924} & {\scriptsize{}45.69} & {\scriptsize{}0.956} & {\scriptsize{}12.25}\tabularnewline
 & {\scriptsize{}Bai (1997)} & {\scriptsize{}0.851} & {\scriptsize{}69.35} & {\scriptsize{}0.847} & {\scriptsize{}51.17} & {\scriptsize{}0.878} & {\scriptsize{}28.59} & {\scriptsize{}0.944} & {\scriptsize{}8.47}\tabularnewline
 & {\scriptsize{}$\widehat{U}_{T}.\textrm{eq}$} & {\scriptsize{}0.964} & {\scriptsize{}88.82} & {\scriptsize{}0.960} & {\scriptsize{}79.74} & {\scriptsize{}0.964} & {\scriptsize{}60.26} & \multicolumn{1}{c}{{\scriptsize{}0.964}} & {\scriptsize{}30.64}\tabularnewline
 & {\scriptsize{}ILR} & {\scriptsize{}0.972} & {\scriptsize{}88.69} & {\scriptsize{}0.975} & {\scriptsize{}73.95} & {\scriptsize{}0.981} & {\scriptsize{}39.08} & {\scriptsize{}0.992} & {\scriptsize{}9.08}\tabularnewline
{\scriptsize{}$\lambda_{0}=0.2$} & {\scriptsize{}HDR} & {\scriptsize{}0.909} & {\scriptsize{}78.12} & {\scriptsize{}0.921} & {\scriptsize{}61.87} & {\scriptsize{}0.933} & {\scriptsize{}40.66} & {\scriptsize{}0.961} & {\scriptsize{}11.70}\tabularnewline
 & {\scriptsize{}Bai (1997)} & {\scriptsize{}0.824} & {\scriptsize{}65.23} & {\scriptsize{}0.867} & {\scriptsize{}51.35} & {\scriptsize{}0.915} & {\scriptsize{}29.83} & {\scriptsize{}0.955} & {\scriptsize{}8.70}\tabularnewline
 & {\scriptsize{}$\widehat{U}_{T}.\textrm{eq}$} & {\scriptsize{}0.961} & {\scriptsize{}89.71} & {\scriptsize{}0.960} & {\scriptsize{}83.68} & {\scriptsize{}0.961} & {\scriptsize{}69.25} & {\scriptsize{}0.960} & {\scriptsize{}35.78}\tabularnewline
 & {\scriptsize{}ILR} & {\scriptsize{}0.966} & {\scriptsize{}91.48} & {\scriptsize{}0.971} & {\scriptsize{}82.78} & {\scriptsize{}0.984} & {\scriptsize{}51.93} & {\scriptsize{}0.995} & {\scriptsize{}10.87}\tabularnewline
\hline 
\end{tabular}{\scriptsize\par}
\par\end{centering}
\noindent\begin{minipage}[t]{1\columnwidth}%
{\scriptsize{}The model is $y_{t}=\beta^{0}+\delta_{Z}^{0}\mathbf{1}_{\left\{ t>\left\lfloor T\lambda_{0}\right\rfloor \right\} }+e_{t},\,e_{t}=0.3e_{t-1}+u_{t},\,u_{t}\sim i.i.d.\,\textrm{t}_{\upsilon},\,\upsilon=5,\,T=100$.
The notes of Table \ref{Table M1} apply.}%
\end{minipage}
\end{table}

\begin{table}[H]
\caption{\label{Table M6}Small-sample coverage rate and length of the confidence
set for model M8}

\begin{centering}
{\scriptsize{}}%
\begin{tabular}{cccccccccc}
\hline 
 &  & \multicolumn{2}{c}{{\scriptsize{}$\delta_{Z}^{0}=0.3$}} & \multicolumn{2}{c}{{\scriptsize{}$\delta_{Z}^{0}=0.6$}} & \multicolumn{2}{c}{{\scriptsize{}$\delta_{Z}^{0}=1$}} & \multicolumn{2}{c}{{\scriptsize{}$\delta_{Z}^{0}=2$}}\tabularnewline
 &  & {\scriptsize{}$\textrm{Cov.}$} & {\scriptsize{}$\textrm{Lgth.}$} & {\scriptsize{}$\textrm{Cov.}$} & {\scriptsize{}$\textrm{Lgth.}$} & {\scriptsize{}$\textrm{Cov.}$} & {\scriptsize{}$\textrm{Lgth.}$} & {\scriptsize{}$\textrm{Cov.}$} & {\scriptsize{}$\textrm{Lgth.}$}\tabularnewline
\hline 
{\scriptsize{}$\lambda_{0}=0.5$} & {\scriptsize{}HDR} & {\scriptsize{}0.918} & {\scriptsize{}75.08} & {\scriptsize{}0.913} & {\scriptsize{}60.44} & {\scriptsize{}0.931} & {\scriptsize{}32.30} & {\scriptsize{}0.965} & {\scriptsize{}6.34}\tabularnewline
 & {\scriptsize{}Bai (1997)} & {\scriptsize{}0.778} & {\scriptsize{}60.94} & {\scriptsize{}0.815} & {\scriptsize{}38.14} & {\scriptsize{}0.885} & {\scriptsize{}17.29} & {\scriptsize{}0.949} & {\scriptsize{}5.34}\tabularnewline
 & {\scriptsize{}$\widehat{U}_{T}.\textrm{eq}$} & {\scriptsize{}0.949} & {\scriptsize{}84.56} & {\scriptsize{}0.950} & {\scriptsize{}67.64} & {\scriptsize{}0.953} & {\scriptsize{}42.95} & {\scriptsize{}0.950} & {\scriptsize{}30.25}\tabularnewline
 & {\scriptsize{}ILR} & {\scriptsize{}0.943} & {\scriptsize{}83.69} & {\scriptsize{}0.946} & {\scriptsize{}63.24} & {\scriptsize{}0.956} & {\scriptsize{}32.85} & {\scriptsize{}0.982} & {\scriptsize{}10.49}\tabularnewline
{\scriptsize{}$\lambda_{0}=0.35$} & {\scriptsize{}HDR} & {\scriptsize{}0.919} & {\scriptsize{}74.16} & {\scriptsize{}0.916} & {\scriptsize{}58.53} & {\scriptsize{}0.931} & {\scriptsize{}32.10} & {\scriptsize{}0.965} & {\scriptsize{}6.48}\tabularnewline
 & {\scriptsize{}Bai (1997)} & {\scriptsize{}0.799} & {\scriptsize{}60.25} & {\scriptsize{}0.814} & {\scriptsize{}37.94} & {\scriptsize{}0.872} & {\scriptsize{}17.49} & {\scriptsize{}0.952} & {\scriptsize{}5.35}\tabularnewline
 & {\scriptsize{}$\widehat{U}_{T}.\textrm{eq}$} & {\scriptsize{}0.951} & {\scriptsize{}85.01} & {\scriptsize{}0.948} & {\scriptsize{}69.14} & {\scriptsize{}0.957} & {\scriptsize{}48.40} & \multicolumn{1}{c}{{\scriptsize{}0.949}} & {\scriptsize{}30.31}\tabularnewline
 & {\scriptsize{}ILR} & {\scriptsize{}0.946} & {\scriptsize{}84.12} & {\scriptsize{}0.944} & {\scriptsize{}63.99} & {\scriptsize{}0.960} & {\scriptsize{}33.45} & {\scriptsize{}0.977} & {\scriptsize{}8.71}\tabularnewline
{\scriptsize{}$\lambda_{0}=0.2$} & {\scriptsize{}HDR} & {\scriptsize{}0.912} & {\scriptsize{}73.43} & {\scriptsize{}0.929} & {\scriptsize{}56.18} & {\scriptsize{}0.949} & {\scriptsize{}31.23} & {\scriptsize{}0.965} & {\scriptsize{}6.96}\tabularnewline
 & {\scriptsize{}Bai (1997)} & {\scriptsize{}0.795} & {\scriptsize{}59.43} & {\scriptsize{}0.864} & {\scriptsize{}38.17} & {\scriptsize{}0.910} & {\scriptsize{}18.52} & {\scriptsize{}0.954} & {\scriptsize{}5.34}\tabularnewline
 & {\scriptsize{}$\widehat{U}_{T}.\textrm{eq}$} & {\scriptsize{}0.950} & {\scriptsize{}86.94} & {\scriptsize{}0.951} & {\scriptsize{}76.52} & {\scriptsize{}0.946} & {\scriptsize{}55.72} & {\scriptsize{}0.947} & {\scriptsize{}38.80}\tabularnewline
 & {\scriptsize{}ILR} & {\scriptsize{}0.945} & {\scriptsize{}83.94} & {\scriptsize{}0.953} & {\scriptsize{}63.55} & {\scriptsize{}0.963} & {\scriptsize{}32.41} & {\scriptsize{}0.982} & {\scriptsize{}15.01}\tabularnewline
\hline 
\end{tabular}{\scriptsize\par}
\par\end{centering}
\noindent\begin{minipage}[t]{1\columnwidth}%
{\scriptsize{}The model is $y_{t}=\delta_{Z}^{0}\left(1-\nu^{0}\right)\mathbf{1}_{\left\{ t>\left\lfloor T\lambda_{0}\right\rfloor \right\} }+\nu^{0}y_{t-1}+e_{t},\,e_{t}\sim i.i.d.\,\mathscr{N}\left(0,\,0.49\right),\,\nu^{0}=0.3,\,T=100$.
The notes of Table \ref{Table M1} apply.}%
\end{minipage}
\end{table}

\begin{table}[H]
\caption{\label{Table M8}Small-sample coverage rate and length of the confidence
sets for model M9}

\begin{centering}
{\scriptsize{}}%
\begin{tabular}{cccccccccc}
\hline 
 &  & \multicolumn{2}{c}{{\scriptsize{}$\delta_{Z}^{0}=0.3$}} & \multicolumn{2}{c}{{\scriptsize{}$\delta_{Z}^{0}=0.6$}} & \multicolumn{2}{c}{{\scriptsize{}$\delta_{Z}^{0}=1$}} & \multicolumn{2}{c}{{\scriptsize{}$\delta_{Z}^{0}=2$}}\tabularnewline
 &  & {\scriptsize{}$\textrm{Cov.}$} & {\scriptsize{}$\textrm{Lgth.}$} & {\scriptsize{}$\textrm{Cov.}$} & {\scriptsize{}$\textrm{Lgth.}$} & {\scriptsize{}$\textrm{Cov.}$} & {\scriptsize{}$\textrm{Lgth.}$} & {\scriptsize{}$\textrm{Cov.}$} & {\scriptsize{}$\textrm{Lgth.}$}\tabularnewline
\hline 
{\scriptsize{}$\lambda_{0}=0.5$} & {\scriptsize{}HDR} & {\scriptsize{}0.903} & {\scriptsize{}61.09} & {\scriptsize{}0.927} & {\scriptsize{}31.14} & {\scriptsize{}0930} & {\scriptsize{}18.33} & {\scriptsize{}0.930} & {\scriptsize{}9.10}\tabularnewline
 & {\scriptsize{}Bai (1997)} & {\scriptsize{}0.791} & {\scriptsize{}37.86} & {\scriptsize{}0.831} & {\scriptsize{}17.73} & {\scriptsize{}0.855} & {\scriptsize{}10.43} & {\scriptsize{}0.868} & {\scriptsize{}5.30}\tabularnewline
 & {\scriptsize{}$\widehat{U}_{T}.\textrm{eq}$} & {\scriptsize{}0.947} & {\scriptsize{}65.23} & {\scriptsize{}0.947} & {\scriptsize{}39.76} & {\scriptsize{}0.947} & {\scriptsize{}28.82} & {\scriptsize{}0.947} & {\scriptsize{}20.36}\tabularnewline
 & {\scriptsize{}ILR} & {\scriptsize{}0.909} & {\scriptsize{}72.62} & {\scriptsize{}0.946} & {\scriptsize{}45.06} & {\scriptsize{}0.962} & {\scriptsize{}23.97} & {\scriptsize{}0.978} & {\scriptsize{}9.34}\tabularnewline
{\scriptsize{}$\lambda_{0}=0.35$} & {\scriptsize{}HDR} & {\scriptsize{}0.904} & {\scriptsize{}60.58} & {\scriptsize{}0.918} & {\scriptsize{}30.96} & {\scriptsize{}0.904} & {\scriptsize{}18.16} & {\scriptsize{}0.928} & {\scriptsize{}0.34}\tabularnewline
 & {\scriptsize{}Bai (1997)} & {\scriptsize{}0.791} & {\scriptsize{}37.70} & {\scriptsize{}0.829} & {\scriptsize{}18.04} & {\scriptsize{}0.852} & {\scriptsize{}10.61} & {\scriptsize{}0.870} & {\scriptsize{}5.34}\tabularnewline
 & {\scriptsize{}$\widehat{U}_{T}.\textrm{eq}$} & {\scriptsize{}0.942} & {\scriptsize{}66.27} & {\scriptsize{}0.942} & {\scriptsize{}40.63} & \multicolumn{1}{c}{{\scriptsize{}0.942}} & {\scriptsize{}29.39} & {\scriptsize{}0.942} & {\scriptsize{}20.67}\tabularnewline
 & {\scriptsize{}ILR} & {\scriptsize{}0.922} & {\scriptsize{}72.20} & {\scriptsize{}0.947} & {\scriptsize{}45.27} & {\scriptsize{}0.959} & {\scriptsize{}24.93} & {\scriptsize{}0.973} & {\scriptsize{}8.55}\tabularnewline
{\scriptsize{}$\lambda_{0}=0.2$} & {\scriptsize{}HDR} & {\scriptsize{}0.920} & {\scriptsize{}61.37} & {\scriptsize{}0.946} & {\scriptsize{}31.00} & {\scriptsize{}0.942} & {\scriptsize{}20.44} & {\scriptsize{}0.944} & {\scriptsize{}9.04}\tabularnewline
 & {\scriptsize{}Bai (1997)} & {\scriptsize{}0.791} & {\scriptsize{}39.23} & {\scriptsize{}0.841} & {\scriptsize{}19.28} & {\scriptsize{}0.876} & {\scriptsize{}11.99} & {\scriptsize{}0.886} & {\scriptsize{}6.16}\tabularnewline
 & {\scriptsize{}$\widehat{U}_{T}.\textrm{eq}$} & {\scriptsize{}0.934} & {\scriptsize{}71.42} & {\scriptsize{}0.931} & {\scriptsize{}47.53} & {\scriptsize{}0.934} & {\scriptsize{}34.12} & {\scriptsize{}0934} & {\scriptsize{}24.06}\tabularnewline
 & {\scriptsize{}ILR} & {\scriptsize{}0.920} & {\scriptsize{}72.68} & {\scriptsize{}0.935} & {\scriptsize{}49.61} & {\scriptsize{}0.959} & {\scriptsize{}27.90} & {\scriptsize{}0.972} & {\scriptsize{}10.01}\tabularnewline
\hline 
\end{tabular}{\scriptsize\par}
\par\end{centering}
\noindent\begin{minipage}[t]{1\columnwidth}%
{\footnotesize{}The model is}{\scriptsize{} $y_{t}=\nu^{0}+Z_{t}\beta^{0}+Z_{t}\delta_{Z}^{0}\mathbf{1}_{\left\{ t>\left\lfloor T\lambda_{0}\right\rfloor \right\} }+e_{t},\,Z_{t}\sim i.i.d.\,\mathscr{N}\left(1,\,1.44\right),\,\left\{ e_{t}\right\} $
follows a FIGARCH(1,0.6,1) process and $T=100$. The notes of Table
\ref{Table M1} apply.}%
\end{minipage}
\end{table}

\begin{table}[H]
\caption{\label{Table M10 LongMemory}Small-sample coverage rate and length
of the confidence set for model M10}

\begin{centering}
{\scriptsize{}}%
\begin{tabular}{cccccccccc}
\hline 
 &  & \multicolumn{2}{c}{{\scriptsize{}$\delta_{Z}^{0}=0.3$}} & \multicolumn{2}{c}{{\scriptsize{}$\delta_{Z}^{0}=0.6$}} & \multicolumn{2}{c}{{\scriptsize{}$\delta_{Z}^{0}=1$}} & \multicolumn{2}{c}{{\scriptsize{}$\delta_{Z}^{0}=2$}}\tabularnewline
 &  & {\scriptsize{}$\textrm{Cov.}$} & {\scriptsize{}$\textrm{Lgth.}$} & {\scriptsize{}$\textrm{Cov.}$} & {\scriptsize{}$\textrm{Lgth.}$} & {\scriptsize{}$\textrm{Cov.}$} & {\scriptsize{}$\textrm{Lgth.}$} & {\scriptsize{}$\textrm{Cov.}$} & {\scriptsize{}$\textrm{Lgth.}$}\tabularnewline
\hline 
{\scriptsize{}$\lambda_{0}=0.5$} & {\scriptsize{}HDR} & {\scriptsize{}0.952} & {\scriptsize{}74.84} & {\scriptsize{}0.930} & {\scriptsize{}36.02} & {\scriptsize{}0.921} & {\scriptsize{}13.11} & {\scriptsize{}0.916} & {\scriptsize{}4.34}\tabularnewline
 & {\scriptsize{}Bai (1997)} & {\scriptsize{}0.809} & {\scriptsize{}45.33} & {\scriptsize{}0.844} & {\scriptsize{}17.11} & {\scriptsize{}0.864} & {\scriptsize{}8.27} & {\scriptsize{}0.883} & {\scriptsize{}3.61}\tabularnewline
 & {\scriptsize{}$\widehat{U}_{T}.\textrm{eq}$} & {\scriptsize{}0.959} & {\scriptsize{}72.69} & {\scriptsize{}0.959} & {\scriptsize{}39.81} & {\scriptsize{}0.959} & {\scriptsize{}24.25} & {\scriptsize{}0.959} & {\scriptsize{}14.79}\tabularnewline
 & {\scriptsize{}ILR} & {\scriptsize{}0.929} & {\scriptsize{}83.23} & {\scriptsize{}0.951} & {\scriptsize{}69.67} & {\scriptsize{}0.971} & {\scriptsize{}44.40} & {\scriptsize{}0.987} & {\scriptsize{}10.44}\tabularnewline
{\scriptsize{}$\lambda_{0}=0.35$} & {\scriptsize{}HDR} & {\scriptsize{}0.934} & {\scriptsize{}73.08} & {\scriptsize{}0.937} & {\scriptsize{}35.37} & {\scriptsize{}0.923} & {\scriptsize{}13.68} & {\scriptsize{}0.920} & {\scriptsize{}4.55}\tabularnewline
 & {\scriptsize{}Bai (1997)} & {\scriptsize{}0.821} & {\scriptsize{}45.70} & {\scriptsize{}0.838} & {\scriptsize{}17.78} & {\scriptsize{}0.867} & {\scriptsize{}8.53} & {\scriptsize{}0.889} & {\scriptsize{}3.71}\tabularnewline
 & {\scriptsize{}$\widehat{U}_{T}.\textrm{eq}$} & {\scriptsize{}0.964} & {\scriptsize{}76.14} & {\scriptsize{}0.964} & {\scriptsize{}44.61} & {\scriptsize{}0.965} & {\scriptsize{}27.33} & \multicolumn{1}{c}{{\scriptsize{}0.964}} & {\scriptsize{}15.84}\tabularnewline
 & {\scriptsize{}ILR} & {\scriptsize{}0.934} & {\scriptsize{}81.32} & {\scriptsize{}0.959} & {\scriptsize{}62.98} & {\scriptsize{}0.977} & {\scriptsize{}34.38} & {\scriptsize{}0.984} & {\scriptsize{}9.12}\tabularnewline
{\scriptsize{}$\lambda_{0}=0.2$} & {\scriptsize{}HDR} & {\scriptsize{}0.941} & {\scriptsize{}71.46} & {\scriptsize{}0.959} & {\scriptsize{}59.03} & {\scriptsize{}0.950} & {\scriptsize{}15.39} & {\scriptsize{}0.919} & {\scriptsize{}5.03}\tabularnewline
 & {\scriptsize{}Bai (1997)} & {\scriptsize{}0.818} & {\scriptsize{}47.82} & {\scriptsize{}0.872} & {\scriptsize{}20.44} & {\scriptsize{}0.878} & {\scriptsize{}9.60} & {\scriptsize{}0.873} & {\scriptsize{}3.92}\tabularnewline
 & {\scriptsize{}$\widehat{U}_{T}.\textrm{eq}$} & {\scriptsize{}0.971} & {\scriptsize{}82.40} & {\scriptsize{}0.971} & {\scriptsize{}59.03} & {\scriptsize{}0.971} & {\scriptsize{}39.02} & {\scriptsize{}0.972} & {\scriptsize{}20.42}\tabularnewline
 & {\scriptsize{}ILR} & {\scriptsize{}0.928} & {\scriptsize{}83.26} & {\scriptsize{}0.952} & {\scriptsize{}70.03} & {\scriptsize{}0.964} & {\scriptsize{}42.65} & {\scriptsize{}0.982} & {\scriptsize{}10.30}\tabularnewline
\hline 
\end{tabular}{\scriptsize\par}
\par\end{centering}
\noindent\begin{minipage}[t]{1\columnwidth}%
{\scriptsize{}The model is $y_{t}=\nu^{0}+Z_{t}\beta^{0}+Z_{t}\delta_{Z}^{0}\mathbf{1}_{\left\{ t>\left\lfloor T\lambda_{0}\right\rfloor \right\} }+e_{t},\,e_{t}\sim i.i.d.\,\mathscr{N}\left(0,\,1\right),\,Z_{t}\sim\mathrm{ARFIMA}\left(0.3,\,0.6,\,0\right),\,T=100$.
The notes of Table \ref{Table M1} apply.}%
\end{minipage}
\end{table}


\end{singlespace}

\normalsize
\indent


\setcounter{figure}{0}
\renewcommand{\thefigure}{S-\arabic{figure}}
\renewcommand{\thepage}{S-\arabic{page}}   
\renewcommand{\thesection}{S.\Alph{section}}   
\renewcommand{\theequation}{S.\arabic{equation}}

\nociteReferencesSupp{casini/perron_CR_Single_Break} \nociteReferencesSupp{casini/perron_CR_Single_Break_Extended}

\newpage{}

\bibliographystyleReferencesSupp{elsarticle-harv}  
\bibliographyReferencesSupp{References_Supp}
\end{document}